\documentclass[reqno,11pt]{amsart} 
\usepackage{amsfonts, amsmath, amssymb, amsthm}
\usepackage[left=1.61cm,right=1.61cm,top=1.8cm,bottom=0.9cm,heightrounded]{geometry}
\usepackage{fancybox}
\usepackage{hhline, float}
\usepackage{mathrsfs}
\usepackage{nccmath}
\usepackage[dvipsnames]{xcolor}
\usepackage[colorlinks=true]{hyperref}
\hypersetup{
  citecolor=red,
  linkcolor=blue,
  urlcolor=blue,
  pdftitle={On the sharp quantitative stability of critical points of the Hardy--Littlewood--Sobolev inequality in Euclidean n-space with n at least 3},
  pdfauthor={Wei Dai; Yichen Hu; Shaolong Peng},
  pdfsubject={Sharp quantitative stability for the energy-critical Hartree equation},
  pdfkeywords={Hardy--Littlewood--Sobolev inequality; Hartree equation; sharp stability; Struwe decomposition; nondegeneracy}
}

\numberwithin{equation}{section}
\allowdisplaybreaks[4]
\AtBeginDocument{%
  \setlength{\jot}{0.5pt}%
  \setlength{\abovedisplayskip}{2.2pt plus .8pt minus .8pt}%
  \setlength{\belowdisplayskip}{2.2pt plus .8pt minus .8pt}%
  \setlength{\abovedisplayshortskip}{1pt plus .5pt minus .5pt}%
  \setlength{\belowdisplayshortskip}{1.4pt plus .5pt minus .5pt}%
}

\newtheorem{theorem}{Theorem}[section]
\newtheorem{proposition}[theorem]{Proposition}
\newtheorem{corollary}[theorem]{Corollary}
\newtheorem{lemma}[theorem]{Lemma}
\theoremstyle{definition}

\newtheorem{rmk}{Remark}[section]
\newtheorem{defn}{Definition}[section]

\newcommand{\N}{\mathbb{N}}
\newcommand{\R}{\mathbb{R}}

\newcommand{\hphi}{\hat{\phi}}

\newcommand{\dx}{\textup{d}x}

\newcommand{\dy}{\textup{d}y}

\usepackage{tcolorbox}

\begin{document}
\title[Sharp quantitative stability of critical points of HLS inequality]{On the sharp quantitative stability of critical points of the Hardy--Littlewood--Sobolev inequality in $\mathbb{R}^{n}$ with $n\geq3$}

\author{Wei Dai}
\address{School of Mathematical Sciences, Beihang University (BUAA), Beijing 100083, People's Republic of China, and Key Laboratory of Mathematics, Informatics and Behavioral Semantics, Ministry of Education, Beijing 100191, People's Republic of China}
\email{weidai@buaa.edu.cn}

\author{Yichen Hu}
\address{School of Mathematical Sciences, Dalian University of Technology, Dalian 116024, Liaoning, People's Republic of China}
\email{huyc24@dlut.edu.cn}

\author{Shaolong Peng}
\address{School of Mathematical Sciences, Beihang University (BUAA), Beijing 100083, People's Republic of China, and Key Laboratory of Mathematics, Informatics and Behavioral Semantics, Ministry of Education, Beijing 100191, People's Republic of China}
\email{slpeng@amss.ac.cn}

\thanks{\noindent Wei Dai is supported by the NNSF of China (Nos. 12222102 and 12571113), the National Science and Technology Major Project (2022ZD0116401), and the Fundamental Research Funds for the Central Universities. Yichen Hu is supported by NNSF of China (No. 12501135), the Fundamental Research Funds for the Central Universities (DUT24RC(3)110), and the Open Research Fund of
	Hubei Key Laboratory of Mathematical Sciences (Central China
	Normal University), Wuhan 430079, P. R. China. Shaolong
	Peng was supported by the National Natural Science Foundation of China (No. 12401148), and Beijing Natural
	Science Foundation (No. 1262014).}


\begin{abstract}
Assume $n\geq3$ and $u\in\dot H^1(\mathbb R^n)$. Recently, Piccione,
Yang and Zhao \cite{Piccione-Yang-Zhao} established a nonlocal version
of Struwe's decomposition in \cite{Struwe-1984}, i.e., if
\(\Gamma(u):=\left\|\Delta u+D_{n,\alpha}
\int_{\mathbb R^n}\frac{|u|^{p_\alpha}(y)}{|x-y|^\alpha}\,dy\,
|u|^{p_\alpha-2}u
\right\|_{\dot H^{-1}}\to0\) and $u\geq0$, then
$\operatorname{dist}(u,\mathcal T)\to0$, where
$\operatorname{dist}(u,\mathcal T)$ denotes the
$\dot H^1(\mathbb R^n)$-distance of $u$ from the manifold of sums of
Talenti bubbles. In this paper, we establish the nonlocal version of
the quantitative estimates of Struwe's decomposition in
Ciraolo, Figalli and Maggi \cite{CFM} for one bubble and $n\geq3$,
Figalli and Glaudo \cite{Figalli-Glaudo2020} for $3\leq n\leq5$ and
Deng, Sun and Wei \cite{DSW} for $n\geq6$ and two or more bubbles.
We prove that, for $n\geq3$, $\nu \in \mathbb{N}^{+}$ and $0<\alpha<n$, there are constants
$\delta=\delta(n,\nu,\alpha)>0$ and $C=C(n,\nu,\alpha)>0$ such that whenever $u\in\dot H^1(\mathbb R^n)$ satisfies 
\[
	\left\|\nabla u-\sum_{i=1}^\nu \nabla \tilde{U}_i\right\|_{L^2} \leq \delta,
\]
where $\tilde U_i=U[\tilde z_i,\tilde\lambda_i]$,
$\tilde z_i\in\mathbb R^n$, $\tilde\lambda_i>0$, and
$\{\tilde U_i:1\leq i\leq\nu\}$ is $\delta$-interacting, it holds that
\[
\operatorname{dist}(u,\mathcal T)
\leq C
\begin{cases}
\Gamma(u)\bigl(1+|\log\Gamma(u)|\bigr)^{1/2},
&\nu\geq2,\ (n,\alpha)\in\mathcal C_{\log},\\
\Gamma(u)^{\theta_\alpha},
&\nu\geq2,\ n\geq7,\ \dfrac{n+2}{2}<\alpha<n,\\
\Gamma(u),&\text{otherwise},
\end{cases}
\]
where $\nu$ denotes the number of bubbles,
\[
a_\alpha:=\min\{\alpha,n-2\},\qquad
b_\alpha:=\min\left\{\alpha,n-2,\frac{n+2}{2}\right\},\qquad
\theta_\alpha:=\frac{b_\alpha}{a_\alpha},
\]
and
\[
\mathcal C_{\log}
:=\{(6,\alpha):4\leq\alpha<6\}
\cup
\left\{\left(n,\frac{n+2}{2}\right):n\geq7\right\}.
\]
Furthermore, the logarithmic rate is sharp when
$(n,\alpha)\in\mathcal C_{\log}$, and the power rate
$\Gamma(u)^{\theta_\alpha}$ is sharp when
$n\geq7$ and $(n+2)/2<\alpha<n$.
It should be emphasized that, in our paper, we have developed new
techniques to deal with the strong singular case $4<\alpha<n$, which
can not be handled by reduction methods in previous works. We believe
that our method can also be applied to other problems related to the
physically interesting Hartree equation.
\end{abstract}

\maketitle

\section{Introduction}
\subsection{Motivation and main results}

We consider the following $\dot{H}^{1}$-energy-critical nonlocal Hartree equation on $\mathbb{R}^n$:
\begin{equation}\label{eq1}
   -\Delta u(x)-D_{n,\alpha} \int_{\mathbb{R}^n} \frac{|u|^{p_{\alpha}}(y)}{|x-y|^\alpha} \mathrm{d} y\, g_{p_\alpha}(u(x))=0, \quad x \in \mathbb{R}^n,
\end{equation}
where $n\geq3$, $0<\alpha<n$,
$
2^*:=\frac{2n}{n-2}$,
$p_{\alpha}:=\frac{2n-\alpha}{n-2}$,
and, for every $p>1$, we use the continuous convention
\begin{equation*}
 g_p(t):=
 \begin{cases}
 |t|^{p-2}t,&t\ne0,\\
 0,&t=0.
 \end{cases}
\end{equation*}
Thus every occurrence of $|t|^{p-2}t$ below is understood as
$g_p(t)$ at $t=0$.  The normalization constant is
$
D_{n,\alpha}=\frac{\Gamma\left(n-\frac{\alpha}{2}\right)}{(n(n-2))^{\frac{n-\alpha}{2}} \pi^{\frac{n}{2}} \Gamma\left(\frac{n-\alpha}{2}\right)}.
$
We also write
\[
 \kappa_n:=\frac{1}{(n-2)|\mathbb S^{n-1}|},\qquad
 \mathfrak c_{n,\alpha}:=
 D_{n,\alpha}^{\frac{1}{2(p_\alpha-1)}}
\]
for the Newtonian-kernel and equation-normalization constants,
respectively.
The Hartree equation originates from the Hartree--Fock theory of quantum mechanics.  A finite-energy weak solution of \eqref{eq1} is a function $u\in\dot H^1(\mathbb R^n)$ satisfying
\[
 \int_{\mathbb R^n}\nabla u\cdot\nabla\varphi\,\dx
 =D_{n,\alpha}\int_{\mathbb R^n}
 \bigl(|x|^{-\alpha}*|u|^{p_\alpha}\bigr)
 g_{p_\alpha}(u)\,\varphi\,\dx
 \qquad(\varphi\in C_c^\infty(\mathbb R^n)).
\]
 The right-hand side defines a continuous functional on $\dot H^1(\mathbb R^n)$ by the Sobolev and Hardy--Littlewood--Sobolev inequalities.  Every nontrivial nonnegative finite-energy weak solution of \eqref{eq1} is an Aubin--Talenti bubble; see \cite{CDL,DFHQW,DLQ,DLL,DQ} and also \cite{Aubin,Talenti}.  Put
\[
 A_n:=[n(n-2)]^{\frac{n-2}{4}}.
\]
With the normalization in \eqref{eq1}, these bubbles are
\[
U[z, \lambda](x):=A_n\left(\frac{\lambda}{1+\lambda^2|x-z|^2}\right)^{\frac{n-2}{2}}.
\]
With $v(x)=\mathfrak c_{n,\alpha}u(x)$, equation \eqref{eq1} becomes
\begin{equation}\label{eq2}
     -\Delta v(x)-\int_{\mathbb{R}^n} \frac{|v|^{p_{\alpha}}(y)}{|x-y|^\alpha} \mathrm{d} y\, g_{p_\alpha}(v(x))=0, \quad x \in \mathbb{R}^n.
\end{equation}
Thus every nontrivial nonnegative finite-energy weak solution of \eqref{eq2} has the form
\[
\mathcal{W}[z, \lambda](x)=\mathfrak c_{n,\alpha}U[z,\lambda](x)
=\mathfrak c_{n,\alpha}[n(n-2)]^{\frac{n-2}{4}}\left(\frac{\lambda}{1+\lambda^2|x-z|^2}\right)^{\frac{n-2}{2}}, \qquad \lambda>0, \,\,\, z \in \mathbb{R}^n.
\]
The $\dot H^1$-energy-critical nonlocal Hartree equation \eqref{eq1} (equivalently, \eqref{eq2}) is the Euler--Lagrange equation associated with the Hardy--Littlewood--Sobolev inequality; see \cite{DFHQW,DLQ,DLL,DQ,FL,Lieb}.  We fix the normalization
\begin{equation}\label{eq:def-SHL}
S_{\mathrm{HL}}:=
\inf_{0\ne v\in\dot H^1(\mathbb R^n)}
\frac{\displaystyle\int_{\mathbb R^n}|\nabla v|^2\,\dx}
{\displaystyle
\left(\int_{\mathbb R^n}\int_{\mathbb R^n}
\frac{|v(x)|^{p_\alpha}|v(y)|^{p_\alpha}}{|x-y|^\alpha}
\,\dx\,\dy\right)^{1/p_\alpha}}.
\end{equation}
Equivalently,
\[
\left(\int_{\mathbb R^n}\int_{\mathbb R^n}
\frac{|v(x)|^{p_\alpha}|v(y)|^{p_\alpha}}{|x-y|^\alpha}
\,\dx\,\dy\right)^{1/(2p_\alpha)}
\leq S_{\mathrm{HL}}^{-1/2}\|\nabla v\|_{L^2}.
\]
Since $\mathcal W[z,\lambda]$ is an extremizer in \eqref{eq:def-SHL}, testing the coefficient-one equation \eqref{eq2} by $\mathcal W[z,\lambda]$ gives
\begin{equation*}
\|\mathcal W[z,\lambda]\|_{\dot H^1}^2
=S_{\mathrm{HL}}^{p_\alpha/(p_\alpha-1)}
=S_{\mathrm{HL}}^{(2n-\alpha)/(n+2-\alpha)}.
\end{equation*}
To distinguish the two equation normalizations, define
\begin{equation*}
\begin{split}
\mathcal T_\nu^{(D)}&:=
\left\{\sum_{i=1}^{\nu}U[z_i,\lambda_i]:
z_i\in\mathbb R^n,\ \lambda_i>0\right\},\\
\mathcal T_\nu^{(1)}&:=
\left\{\sum_{i=1}^{\nu}\mathcal W[z_i,\lambda_i]:
z_i\in\mathbb R^n,\ \lambda_i>0\right\}
=\mathfrak c_{n,\alpha}\mathcal T_\nu^{(D)}.
\end{split}
\end{equation*}
Thus $\mathcal T_\nu^{(D)}$ is used for the equation with coefficient $D_{n,\alpha}$, while $\mathcal T_\nu^{(1)}$ is used for the coefficient-one equation.  We write
\[
\operatorname{dist}_{\dot H^1}(v,\mathcal A)
:=\inf_{w\in\mathcal A}\|v-w\|_{\dot H^1},
\qquad \|v\|_{\dot H^1}:=\|\nabla v\|_{L^2}.
\]

\medskip

Struwe \cite{Struwe-1984} proved the corresponding sequential compactness statement for the Yamabe equation: under a fixed energy window, if a nonnegative sequence $(u_m)\subset\dot H^1(\mathbb R^n)$ satisfies $\|\Delta u_m+u_m^{(n+2)/(n-2)}\|_{\dot H^{-1}}\to0$, then, after extraction, $u_m$ converges in $\dot H^1$ to a sum of asymptotically orthogonal bubbles with the number of bubbles determined by that energy window. In the same spirit, Piccione, Yang and Zhao \cite{Piccione-Yang-Zhao} proved the following nonlocal profile decomposition for nonnegative functions.
\begin{theorem}[Theorem 1.4 in \cite{Piccione-Yang-Zhao}]
\label{thm:qualitative-nonlocal-profile}
    Let $n\geq3$, $0<\alpha<n$, and let $\nu\geq1$ be an integer. Let $\left(u_m\right)_{m \in \mathbb{N}} \subseteq \dot{H}^{1}\left(\mathbb{R}^n\right)$ be a sequence of nonnegative functions such that
\[
\left(\nu-\frac{1}{2}\right) S_{\mathrm{HL}}^{\frac{2 n-\alpha}{n+2-\alpha}} \leq \left\| u_m\right\|_{\dot{H}^{1}\left(\mathbb{R}^n\right)} ^2 \leq\left(\nu+\frac{1}{2}\right) S_{\mathrm{HL}}^{\frac{2n-\alpha}{n+2-\alpha}}
\]
where $S_{\mathrm{HL}}$ is defined by \eqref{eq:def-SHL}, and assume that
\[
\left\|-\Delta u_{m}(x)-\int_{\mathbb{R}^n} \frac{u_{m}^{p_{\alpha}}(y)}{|x-y|^\alpha} \mathrm{d} y u_{m}^{p_{\alpha}-1}(x)\right\|_{\dot{H}^{-1}} \rightarrow 0, \qquad \text { as } \quad m \rightarrow \infty.
\]
Then there exist $\nu$-tuples of points $\bigl(z_1^{(m)},\allowbreak \ldots,\allowbreak z_\nu^{(m)}\bigr)_{m \in \mathbb{N}}$ in $\mathbb{R}^n$ and $\nu$-tuples of positive real numbers $\bigl(\lambda_1^{(m)},\allowbreak \ldots,\allowbreak \lambda_\nu^{(m)}\bigr)_{m \in \mathbb{N}}$ such that
\[
\left\|\nabla\left(u_m-\sum_{i=1}^\nu \mathcal{W}\left[z_i^{(m)}, \lambda_i^{(m)}\right]\right)\right\|_{L^2} \rightarrow 0, \qquad \text { as } \quad m \rightarrow \infty.
\]
Moreover, for every $i\ne j$,
\[
 \frac{\lambda_i^{(m)}}{\lambda_j^{(m)}}
 +\frac{\lambda_j^{(m)}}{\lambda_i^{(m)}}
 +\lambda_i^{(m)}\lambda_j^{(m)}
   \lvert z_i^{(m)}-z_j^{(m)}\rvert^2
 \longrightarrow\infty.
\]
\end{theorem}

In \cite{CFM}, Ciraolo, Figalli and Maggi obtained the first quantitative version of Struwe's one-bubble decomposition in every dimension. For two or more bubbles, Figalli and Glaudo \cite{Figalli-Glaudo2020} proved the same linear estimate when $3\leq n\leq5$ and showed that it fails for $n\geq6$. Deng, Sun and Wei \cite{DSW} then proved, for $n\geq6$, the sharp estimate for the Yamabe bubble manifold. More specifically, they proved that, for $n\geq6$ and $\nu \in \mathbb{N}^{+}$, there are constants
$\delta=\delta(n,\nu,\alpha)>0$ and $C=C(n,\nu,\alpha)>0$ such that whenever $u\in\dot H^1(\mathbb R^n)$ satisfies 
\[
\left\|\nabla u-\sum_{i=1}^\nu \nabla \tilde{U}_i\right\|_{L^2} \leq \delta,
\]
where $\tilde U_i=U[\tilde z_i,\tilde\lambda_i]$,
$\tilde z_i\in\mathbb R^n$, $\tilde\lambda_i>0$, and
$\{\tilde U_i:1\leq i\leq\nu\}$ is $\delta$-interacting(see Definition 1.1 and the paragraph immediately
following it below), it holds that
    \[\operatorname{dist}_{\dot H^1}(u,\mathcal T)\leq C\begin{cases} \Gamma_{\mathrm{Yam}}(u)\left|\log \Gamma_{\mathrm{Yam}}(u)\right|^{\frac{1}{2}},\quad&\text{if }n=6\text{ and }\nu\geq 2,\\
    \Gamma_{\mathrm{Yam}}(u)^{\frac{n+2}{2(n-2)}},\quad&\text{if }n\geq 7\text{ and }\nu\geq 2,
    \\\Gamma_{\mathrm{Yam}}(u),\quad&\text{otherwise, } \end{cases}\]
where $\Gamma_{\mathrm{Yam}}(u):=\|\Delta u+|u|^{4/(n-2)}u\|_{\dot H^{-1}}$ and  $\mathcal T$ is  the manifold of sums of the positive Yamabe bubbles.
Dimension is central to the analysis of the Yamabe equation $\Delta u+u^{\frac{n+2}{n-2}}=0$, with dimension $6$ acting as a threshold in several estimates; see, for example, Aubin \cite{Aubin} and Schoen \cite{Schoen}. It also plays a decisive role in prescribed scalar-curvature problems; see \cite{Li,Druet,CY,BahriCoron,ACCH,MM} and the references therein.

\medskip

Inspired by \cite{CFM,DSW,Figalli-Glaudo2020}, we establish the sharp quantitative form of the nonlocal Struwe decomposition proved by Piccione, Yang and Zhao \cite{Piccione-Yang-Zhao}.

\smallskip

To this end, we need the following definition.
\begin{defn}[Interaction of Aubin--Talenti bubbles] Let $U\left[z_i, \lambda_i\right]$ and $U\left[z_j, \lambda_j\right]$ be two bubbles. Their interaction parameter is
\begin{equation*}
    q\left(z_i, z_j, \lambda_i, \lambda_j\right)=\left(\frac{\lambda_i}{\lambda_j}+\frac{\lambda_j}{\lambda_i}+\lambda_i \lambda_j\left|z_i-z_j\right|^2\right)^{-\frac{n-2}{2}}.
\end{equation*}
\end{defn}

We shall denote $q_{i j}=q_{j i}=q\left(z_i, z_j, \lambda_i, \lambda_j\right)$. Let $\left\{U_i: 1 \leq i \leq \nu\right\}$ be a family of Aubin--Talenti bubbles. Set
\begin{equation*}
Q:=
\begin{cases}
0,&\nu=1,\\
\displaystyle\max \left\{q_{i j}: 1 \leq i \neq j \leq \nu\right\},&\nu\geq2.
\end{cases}
\end{equation*}
We say that the family is $\delta$-interacting if $Q\leq\delta$.
We shall use the following interaction exponents
\begin{equation}\label{eq:corrected-exponents}
    a_{\alpha}:=\min\{\alpha,n-2\},\qquad
    b_{\alpha}:=\min\left\{\alpha,n-2,\frac{n+2}{2}\right\},\qquad
    \theta_{\alpha}:=\frac{b_{\alpha}}{a_{\alpha}},\qquad
    c_n:=\min\left\{n-2,\frac{n+2}{2}\right\}.
\end{equation}
\begin{equation*}
\mathcal C_{\log}:=
 \{(6,\alpha):4\leq\alpha<6\}
 \cup
 \left\{\left(n,\frac{n+2}{2}\right):n\geq7\right\}.
\end{equation*}
Whenever it occurs below, the function
$t\left(1+|\log t|\right)^{1/2}$ is extended continuously by the value
$0$ at $t=0$.
For every integer $\nu\geq2$ and every $R>e$, set
\begin{equation}\label{eq:H-alpha}
    \mathfrak{H}_{\nu,n,\alpha}(R):=
    \begin{cases}
        R^{-b_\alpha}(\log R)^{\frac12},&(n,\alpha)\in\mathcal C_{\log},\\
        R^{-b_\alpha},&n\geq7,\ \dfrac{n+2}{2}<\alpha<n,\\
        R^{-a_\alpha},&\text{for the remaining cases}.
    \end{cases}
\end{equation}

\medskip

Our main result on sharp quantitative stability is the following theorem.
\begin{theorem}\label{mainthm1}
Let $\nu\geq1$, $n\geq3$, and $0<\alpha<n$. There are constants
$\delta=\delta(n,\nu,\alpha)>0$ and $C=C(n,\nu,\alpha)>0$ such that
whenever $u\in\dot H^1(\mathbb R^n)$ satisfies
\begin{equation}\label{1.1}
    \left\|\nabla u-\sum_{i=1}^\nu \nabla \tilde{U}_i\right\|_{L^2} \leq \delta,
\end{equation}
where $\tilde U_i=U[\tilde z_i,\tilde\lambda_i]$,
$\tilde z_i\in\mathbb R^n$, $\tilde\lambda_i>0$, and
$\{\tilde U_i:1\leq i\leq\nu\}$ is $\delta$-interacting, there exist
$z_i\in\mathbb R^n$ and $\lambda_i>0$, $1\leq i\leq\nu$, such that, with
\begin{equation}\label{introdueq1.1}
\left\|\nabla u-\sum_{i=1}^\nu\nabla U_i\right\|_{L^2}
\leq C
\begin{cases}
\Gamma(u)\bigl(1+|\log\Gamma(u)|\bigr)^{1/2},
&\nu\geq2,\ (n,\alpha)\in\mathcal C_{\log},\\
\Gamma(u)^{\theta_\alpha},
&\nu\geq2,\ n\geq7,\ \dfrac{n+2}{2}<\alpha<n,\\
\Gamma(u),&\text{otherwise},
\end{cases}
\end{equation}
where $
    U_i:=U[z_i,\lambda_i]$,
\(\Gamma(u):=\left\|\Delta u
+D_{n,\alpha}\bigl(|x|^{-\alpha}*|u|^{p_\alpha}\bigr)
g_{p_\alpha}(u)\right\|_{\dot H^{-1}}\).
Moreover, if $\nu\geq2$, then, for every $i\ne j$,
\begin{equation}\label{eq3}
\int_{\mathbb R^n}
\bigl(|x|^{-\alpha}*U[z_i,\lambda_i]^{p_\alpha}\bigr)
U[z_i,\lambda_i]^{p_\alpha-1}U[z_j,\lambda_j] \,\dx
\leq C\,\Gamma(u)^{\frac{n-2}{a_\alpha}}.
\end{equation}

\end{theorem}

\begin{rmk}\label{rmk:one-bubble-case}
In the notation introduced in the proof below, for $\nu=1$ one has
$Q=0$ and $h=0$. Proposition
\ref{pro4.2.1} and the nonlinear estimate give
\[
\|\rho\|_{\dot H^1}
\leq C\left(\|f\|_{\dot H^{-1}}
+\|\rho\|_{\dot H^1}^{1+\eta}\right).
\]
Since $\|\rho\|_{\dot H^1}\leq\delta$, choosing the initial distance
$\delta$ so that $C\delta^\eta\leq1/2$ yields
$\|\rho\|_{\dot H^1}\leq2C\|f\|_{\dot H^{-1}}$, which is the
$\nu=1$ assertion.
\end{rmk}

\begin{rmk}
Our results are nonlocal counterparts of the quantitative stability results in \cite{CFM,DSW,Figalli-Glaudo2020}. The construction in Section \ref{Sec6} proves that both nonlinear branches of \eqref{introdueq1.1} are sharp. More precisely, for $(n,\alpha)\in\mathcal C_{\log}$, it proves the sharpness of the logarithmic correction, while, for $n\geq7$ and $(n+2)/2<\alpha<n$, it proves the sharpness of the power estimate $\Gamma(u)^{\theta_\alpha}$; equivalently, the exponent $\theta_{\alpha}$ in \eqref{eq:corrected-exponents} is optimal.
\end{rmk}

\begin{rmk}
Liu, Zhang and Zou \cite{LZZ} proved quantitative stability for $n=3$ and $5/2<\alpha<3$. In the multi-bubble problem, Piccione, Yang and Zhao \cite{Piccione-Yang-Zhao} treated $3\leq n<6-\alpha$ with $0<\alpha\leq4$, while Yang and Zhao \cite{YZ} obtained estimates for $n\geq6-\alpha$ and $0<\alpha<\min\{4,n\}$. Theorem \ref{mainthm1} covers all $n\geq3$ and $0<\alpha<n$ and gives the optimal dependence on $\Gamma(u)$. The new ingredient is the treatment of the strongly singular range $4<\alpha<n$, where the earlier reduction arguments do not apply directly.
\end{rmk}

\begin{rmk}
Equation \eqref{eq3} gives
\[
\lim_{\alpha\downarrow0}
\frac{n-2}{\min\{\alpha,n-2\}}=+\infty,
\]
and, for $i\ne j$,
\[
\int_{\mathbb{R}^n}\left(\frac{1}{|x|^\alpha}
*U[z_i,\lambda_i]^{p_\alpha}\right)
U[z_i,\lambda_i]^{p_\alpha-1}U[z_j,\lambda_j]\dx
=O\left(\Gamma(u)^{\frac{n-2}{\min\{\alpha,n-2\}}}\right).
\]
\end{rmk}

\begin{theorem}\label{mainthm3}
Assume
 $
 n\geq6 $, $\frac{n+2}{2}\leq\alpha<n $.
There exist $R_0=R_0(n,\alpha)>e$ and constants
$0<c=c(n,\alpha)\leq C=C(n,\alpha)<\infty$ such that,
for every $R\geq R_0$, there is $\rho_R\in\dot H^1(\mathbb R^n)$.
Define
\[
 u_R:=U[-Re_1,1]+U[Re_1,1]+\rho_R,\qquad
 e_1=(1,0,\ldots,0).
\]
Then the following quantitative small-defect properties hold:
\[
\begin{aligned}
 c\mathfrak{H}_{2,n,\alpha}(R)
 \leq
 \operatorname{dist}_{\dot H^1}\bigl(u_R,\mathcal T_2^{(D)}\bigr)
 \leq \|\rho_R\|_{\dot H^1}
 \leq C\mathfrak{H}_{2,n,\alpha}(R)\longrightarrow0,\qquad
 0\leq\Gamma(u_R)&\leq CR^{-a_\alpha}\longrightarrow0.
\end{aligned}
\]
Moreover,
\begin{equation*}
\inf_{\substack{z_1,z_2\in\mathbb R^n\\ \lambda_1,\lambda_2>0}}
\left\|\nabla u_R-\sum_{i=1}^2\nabla U[z_i,\lambda_i]\right\|_{L^2}
\geq c
\begin{cases}
\Gamma(u_R)\left(1+|\log\Gamma(u_R)|\right)^{1/2},
&(n,\alpha)\in\mathcal C_{\log},\\[2mm]
\Gamma(u_R)^{\theta_\alpha},
&n\geq7,\ \dfrac{n+2}{2}<\alpha<n,
\end{cases}
\end{equation*}
where
\[
 \Gamma(u_R):=\left\|\Delta u_R
 +D_{n,\alpha}\left(|x|^{-\alpha}*|u_R|^{p_\alpha}\right)
 g_{p_\alpha}(u_R)\right\|_{\dot H^{-1}},
 \qquad \theta_\alpha=\frac{b_\alpha}{a_\alpha},
\]
and $\mathfrak{H}_{2,n,\alpha}$ is defined in \eqref{eq:H-alpha}.
\end{theorem}

As a direct consequence of Theorem \ref{mainthm1}, we can obtain the following corollary.
\begin{corollary}\label{mainthm2}
Suppose $n \geq 3$, $0<\alpha<n$, and let $\nu\geq1$ be an integer. There exists a constant $C=C(n, \nu,\alpha)>0$ such that the following statement holds. For any nonnegative function $u \in \dot{H}^{1}\left(\mathbb{R}^n\right)$ such that
\[
\left(\nu-\frac{1}{2}\right) S_{\mathrm{HL}}^{\frac{2n-\alpha}{n+2-\alpha}} \leq\|u\|_{\dot{H}^{1}\left(\mathbb{R}^n\right)}^2 \leq\left(\nu+\frac{1}{2}\right) S_{\mathrm{HL}}^{\frac{2n-\alpha}{n+2-\alpha}},
\]
there exist $z_i\in\mathbb R^n$ and $\lambda_i>0$, $1\leq i\leq\nu$. Writing $\mathcal W_i:=\mathcal W[z_i,\lambda_i]$, one has
\begin{equation*}
    \begin{split}
        \left\|\nabla u-\sum_{i=1}^\nu \nabla \mathcal{W}_i\right\|_{L^2} \leq C \begin{cases}
    \hat{\Gamma}(u)\left(1+\left|\log\hat{\Gamma}(u)\right|\right)^{\frac{1}{2}},
    &\nu\geq2,\ (n,\alpha)\in\mathcal C_{\log},\\
     \hat{\Gamma}(u)^{\theta_{\alpha}},
    &\nu\geq2,\ n\geq7,\ \dfrac{n+2}{2}<\alpha<n,\\
     \hat{\Gamma}(u),&\text{otherwise},
\end{cases}
    \end{split}
\end{equation*}
where $ \hat{\Gamma}(u):=\left\|\Delta u(x)+\int_{\mathbb{R}^{n}}\frac{|u|^{p_{\alpha}}(y)}{|x-y|^{\alpha}}dy\,g_{p_\alpha}(u(x))\right\|_{\dot{H}^{-1}}$.
Furthermore, if $\nu \geq 2$, then for any $i \neq j$, the interaction between the bubbles can be estimated as
\begin{equation*}
    \begin{split}
        \int_{\mathbb{R}^n}\left(\frac{1}{|x|^\alpha} * \mathcal{W}\left[z_i, \lambda_i\right]^{p_{\alpha}}\right) \mathcal{W}\left[z_i, \lambda_i\right]^{p_{\alpha}-1} \mathcal{W}\left[z_j, \lambda_j\right]\dx
        \leq C\,\hat\Gamma(u)^{\frac{n-2}{a_\alpha}}.
    \end{split}
\end{equation*}
\end{corollary}

\smallskip

For quantitative stability of Sobolev, isoperimetric, and related geometric and functional inequalities, as well as applications to the calculus of variations and PDEs, see, e.g., \cite{BE,BF,BPV,CMM,DMMN,DS,DE,Figalli-Glaudo2020,FMP0,FN,FZ,Frank,Frank1,FMP,Maggi} and the references therein.


\medskip

\subsection{Sketch of the proof}
Write \(A\lesssim B\) if \(A\leq C(n,\nu,\alpha)B\), and
\(A\asymp B\) if \(A\lesssim B\lesssim A\).  For \(i\ne j\), set
\[
R_{ij}:=\max\left\{\sqrt{\frac{\lambda_i}{\lambda_j}},
\sqrt{\frac{\lambda_j}{\lambda_i}},
\sqrt{\lambda_i\lambda_j}|z_i-z_j|\right\},\qquad
R:=\frac12\min_{i\ne j}R_{ij},
\]
and
\begin{equation*}
\mathfrak d_{ij}:=
\frac{\lambda_i/\lambda_j-\lambda_j/\lambda_i+
\lambda_i\lambda_j|z_i-z_j|^2}
{\lambda_i/\lambda_j+\lambda_j/\lambda_i+
\lambda_i\lambda_j|z_i-z_j|^2}.
\end{equation*}

\paragraph{\textbf{Modulation and bubble geometry.}}
We first use the standard best-approximation construction for weakly
interacting Aubin--Talenti bubbles; see
\cite[Appendix A]{BahriCoron} and
\cite[(1.12)--(1.14)]{DSW}.  After decreasing \(\delta\), the infimum
defining
\(\operatorname{dist}_{\dot H^1}(u,\mathcal T_\nu^{(D)})\) is attained
at a sum \(\sigma=\sum_iU_i\).  Thus, after relabeling, there exists a
modulus \(\omega_{\rm mod}(t)\to0\) as \(t\downarrow0\), depending only
on \((n,\nu)\), such that, with \(\rho=u-\sigma\),
\begin{equation}\label{introeq2}
\max_i\left(
\left|\log\frac{\lambda_i}{\widetilde\lambda_i}\right|
+\widetilde\lambda_i|z_i-\widetilde z_i|\right)
\leq\omega_{\rm mod}(\delta),\quad
\|\rho\|_{\dot H^1}\leq\delta, \quad
Q\leq\omega_{\rm mod}(\delta),\quad
\int_{\mathbb R^n}\nabla\rho\cdot\nabla Z_i^a\dx=0.
\end{equation}
The orthogonality relations are the Euler--Lagrange equations of the
finite-dimensional minimization, and the estimate for \(Q\) follows
from parameter matching and the interaction assumption on the original
family.

A weakly interacting family may contain towers, clusters, and mixtures
of the two.  Section \ref{Sec2} therefore organizes every sequence with
\(R\to\infty\) into bubble trees, using the order \(i\succ j\), its
reverse, and the incomparable case.  These three alternatives are used
throughout the interaction estimates.  Since the desired estimates are
uniform in the bubble parameters, a contradiction argument converts the
subsequential tree decomposition into constants depending only on
\((n,\nu,\alpha)\).

\paragraph{\textbf{The reduced equation.}}
Set
\[
\mathcal I_\alpha(g)(x):=
D_{n,\alpha}\int_{\mathbb R^n}\frac{g(y)}{|x-y|^\alpha}\,\dy .
\]
For
$
\mathscr F(v):=\mathcal I_\alpha(|v|^{p_\alpha})
g_{p_\alpha}(v),
$
we define
\begin{align}
h&:=\mathscr F(\sigma)-\sum_i\mathscr F(U_i), \label{h}\\
N_1(\eta,\xi)&:=D\mathscr F(\xi)[\eta],\qquad
N(\eta,\xi):=\mathscr F(\xi+\eta)-\mathscr F(\xi)
-D\mathscr F(\xi)[\eta]. \label{N}
\end{align}
Then
\begin{equation*}
\Delta\rho+N_1(\rho,\sigma)+h+N(\rho,\sigma)+f=0.
\end{equation*}
The principal analytic problem is to estimate \(h\) in
\(\dot H^{-1}\), uniformly over all bubble trees, and then to show that
the two terms depending on \(\rho\) are lower order in both the energy
estimate and the finite-dimensional projections.

Put
\[
 \beta:=p_\alpha-1=\frac{n+2-\alpha}{n-2},\qquad
 \gamma:=2^*-p_\alpha=\frac{\alpha}{n-2}.
\]
Applying the single-bubble identity
\(D_{n,\alpha}(|x|^{-\alpha}*U_i^{p_\alpha})=U_i^\gamma\)
to each component of the multi-bubble sum
\(\sigma=\sum_iU_i\) gives the exact interaction decomposition
\[
h=\sum_iU_i^\gamma\bigl(\sigma^\beta-U_i^\beta\bigr)
+D_{n,\alpha}
\left(|x|^{-\alpha}*
\left(\sigma^{p_\alpha}-\sum_iU_i^{p_\alpha}\right)\right)
\sigma^\beta
=:A+B.
\]
Here \(A\) contains a local mixed-bubble interaction, whereas \(B\)
contains interactions simultaneously in the source and observation
variables of the Riesz potential.

\paragraph{\textbf{Projected equation and conclusion.}}
Lemmas \ref{le7.1} and \ref{le3.1} give
\[
 \|h\|_{\dot H^{-1}}
 \lesssim\mathfrak H_{\nu,n,\alpha}(R).
\]
After removing the translation and scale modes by
\eqref{introeq2}, Proposition \ref{pro4.1.1} gives a uniform a priori
estimate for the projected linearized operator.  Proposition
\ref{pro4.2.1} consequently provides a uniformly bounded projected
inverse.  The nonlinear estimates in Section \ref{Sec5} imply, for
some \(\eta=\eta(n,\alpha)>0\),
\[
 \|N(\rho,\sigma)\|_{\dot H^{-1}}
 \lesssim\|\rho\|_{\dot H^1}^{1+\eta}.
\]
Absorbing this superlinear term gives
\[
 \|\rho\|_{\dot H^1}\lesssim
 \mathfrak H_{\nu,n,\alpha}(R)+\|f\|_{\dot H^{-1}}.
\]

It remains to control the bubble parameters.  Lemma \ref{le2.1} in the
weakly singular range and Lemma \ref{le2.1strong} in the strongly
singular range give the leading scale projection of \(h\), of order
\(Q^{a_\alpha/(n-2)}\).  Proposition \ref{pro5.5} shows that the
projections containing \(\rho\) are lower order.  A scale direction
selected from the bubble tree therefore satisfies
\[
 Q^{a_\alpha/(n-2)}
 \lesssim\|f\|_{\dot H^{-1}}
 +o\!\left(Q^{a_\alpha/(n-2)}\right).
\]
After absorption, \(Q\) is controlled by the residual.  Substitution
into the estimate for \(\rho\), together with
\eqref{eq:corrected-exponents}--\eqref{eq:H-alpha}, proves Theorem
\ref{mainthm1}.

\paragraph{\textbf{Sharpness.}}
For \(n\geq6\), \((n+2)/2\leq\alpha<n\), and \(\nu=2\), we let
\[
 \sigma=U[-Re_1,1]+U[Re_1,1]
\]
and construct a projected correction \(\rho\) satisfying
\[
 \Gamma(\sigma+\rho)\lesssim R^{-a_\alpha},\qquad
 \|\rho\|_{\dot H^1}\asymp\mathfrak H_{2,n,\alpha}(R),
\]
in Section
\ref{Sec6}.
The modulation argument identifies \(\|\rho\|_{\dot H^1}\) with the
distance to the two-bubble manifold, up to uniform constants.  Since
\(R^{-b_\alpha}=(R^{-a_\alpha})^{b_\alpha/a_\alpha}\), this proves the
optimality of the logarithmic rate and of
\(\theta_\alpha=b_\alpha/a_\alpha\).

\subsection{Main difficulties and new ingredients}

The preceding reduction has the same formal structure as a
finite-dimensional multi-bubble argument for a local critical
equation.  Its implementation is different for the Hartree equation,
because the mixed-bubble source
\(\sigma^{p_\alpha}(y)-\sum_iU_i^{p_\alpha}(y)\) is coupled through
the kernel \(|x-y|^{-\alpha}\) with the coefficient
\(\sigma^\beta(x)\), and because \(p_\alpha-1<1\) when
\(\alpha>4\).  We record the principal difficulties and the estimates
developed to overcome them.

\paragraph{\textbf{Dominance regions and the endpoint \((5,4)\).}}
For \(0<\alpha\leq4\), one has \(\beta\geq1\).  On the dominance region
\(\Omega_\ell=\{U_\ell=\max_iU_i\}\), the mean-value theorem yields
\[
 0\leq A\leq C\left(
 U_\ell^{2^*-2}\sum_{j\ne\ell}U_j+
 \sum_{j\ne\ell}U_j^\gamma U_\ell^\beta\right).
\]
For every case covered by Lemma \ref{le7.1} except
\((n,\alpha)=(5,4)\), Lemma \ref{A.4} and the
Hardy--Littlewood--Sobolev inequality yield
\[
 \|A\|_{L^{\frac{2n}{n+2}}}
 +\|B\|_{L^{\frac{2n}{n+2}}}
 \leq C R^{-a_\alpha}.
\]
At the endpoint \((n,\alpha)=(5,4)\), the same argument gives only
\[
 \|B\|_{L^{10/7}}
 \leq C R^{-3}(1+\log R)^{3/5},
\]
whereas Lemma \ref{le7.1} requires the sharp bound \(CR^{-3}\).
The logarithmic factor arises when the critical \(L^{5/3}\)-estimate
of \(U_iU_j\) is inserted directly into the global
Hardy--Littlewood--Sobolev inequality, before the full convolution
product \(\mathcal I_4(U_iU_j)U_k\) is estimated.  Since
\(p_\alpha=2\) and \(\beta=1\), one instead uses the exact identity
\[
 B=2\sum_{i<j}\sum_k\mathcal I_4(U_iU_j)U_k.
\]
For each \(i\ne j\), the \(y\)-domain in the convolution of
\(U_i(y)U_j(y)\) is split into \(\{U_i\geq U_j\}\) and
\(\{U_j>U_i\}\).  After rescaling around the dominant bubble, a
truncated convolution bound is applied separately in the core,
intermediate, and far regions.  The ordered, reverse ordered, and
incomparable configurations then all give
\[
 \|\mathcal I_4(U_iU_j)U_k\|_{L^{10/7}}\lesssim R_{ij}^{-3}
 \quad (k=i,j),
\]
and the three-bubble terms have the same or smaller order.  This proves
\(\|h\|_{L^{10/7}}\lesssim R^{-3}\), without a logarithmic loss.

\paragraph{\textbf{The concave coefficient in \(h\).}}
The strongly singular range \(4<\alpha<n\) is qualitatively different:
\(0<\beta<1\), so \(t\mapsto t^\beta\) is concave and its derivative is
singular at the origin.  In particular, a termwise linear expansion of
\(\sigma^\beta-U_i^\beta\) is not uniform.  The local interactions must
be kept in the combined form \(A\).  On \(\Omega_\ell\),
\[
 0\leq A\leq C\left(
 U_\ell^{2^*-2}\sum_{j\ne\ell}U_j+
 \sum_{j\ne\ell}U_j^\gamma U_\ell^\beta\right).
\]
The first term has the Yamabe interaction order, while the second has a
different order in a bubble core and in a neck.  Lemma \ref{A.4} gives
the three competing exponents
\[
 b_\alpha=\min\left\{\alpha,n-2,\frac{n+2}{2}\right\};
\]
the balanced cases are estimated directly in \(\dot H^{-1}\) by
Lemma \ref{A.6}, which gives the factor \((\log R)^{1/2}\).

\paragraph{\textbf{The range \(\alpha>n-2\) in Lemma \ref{le3.1}.}}
The convolution term \(B\) requires an additional argument when
\(\alpha>n-2\).  A direct Hardy--Littlewood--Sobolev estimate only gives
the exponent
\[
 \eta_\alpha:=\frac{2n-\alpha}{2},
\]
which is smaller than the exponent required in Lemma \ref{le3.1} in
part of this range.  For
\[
 K_{\ell j}:=\mathcal I_\alpha
 \bigl(U_\ell^\beta U_j\mathbf1_{\Omega_\ell}\bigr),
\]
we decompose the \(y\)-domain in the integral defining
\(K_{\ell j}(x)\) and the \(x\)-domain in its
\(\dot H^{-1}\)-pairing according to the bubble tree, and then into
core, neck, and far regions.  In an ordered pair the potential is
rescaled at the parent and child scales; in an incomparable pair it is
split at the separation scale.  Lemma \ref{A.4} controls the resulting
critical products, and Lemma \ref{A.6} treats the balanced Newtonian
products.  This restores
\[
 \|h\|_{\dot H^{-1}}
 \lesssim\mathfrak H_{\nu,n,\alpha}(R)
\]
in all cases of Lemma \ref{le3.1}.  Thus the exponent \(b_\alpha\) and
the logarithmic set \(\mathcal C_{\log}\) are consequences of the
actual core--neck competition, rather than of a non-sharp global HLS
bound.

\paragraph{\textbf{The critical dominance-region estimate.}}
The basic device behind these computations is Lemma \ref{A.4}.  If
\(a+b=n\), set
\[
 I_{12}:=\int_{\{U_1\geq U_2\}}
 \frac{\lambda_1^a}{(1+\lambda_1^2|x-z_1|^2)^a}
 \frac{\lambda_2^b}{(1+\lambda_2^2|x-z_2|^2)^b}\,\dx .
\]
Lemma \ref{A.4} proves
\[
\begin{aligned}
 a>b&\quad\Longrightarrow\quad I_{12}\asymp R_{12}^{-2b},\\
 a=b=\frac n2&\quad\Longrightarrow\quad
 I_{12}=O\!\left(R_{12}^{-n}(1+\log R_{12})\right),\\
 a<b&\quad\Longrightarrow\quad I_{12}=O(R_{12}^{-n}).
\end{aligned}
\]
Moreover, it identifies the rescaled leading term in each of the three
tree configurations.  Standard full-space two-bubble estimates do not
retain this dominance information.  Lemma \ref{A.4} is therefore used
not only in Lemmas \ref{le7.1} and \ref{le3.1}, but also in the scale
projections and the estimates of the nonlinear projected terms.

\paragraph{\textbf{The sublinear exponent for \(\alpha>4\).}}
For \(\alpha\leq4\), the nonlinear remainder is controlled by the
standard superlinear power inequalities available when
\(p_\alpha\geq2\).  For \(\alpha>4\), however,
\(1<p_\alpha<2\); consequently \(g_{p_\alpha}\) is only
\((p_\alpha-1)\)-H\"older continuous at the origin and its derivative is
singular there.  In particular, a quadratic Taylor estimate is false.
Lemma
\ref{le:strong-nonlinear} uses
\[
 E_1:=\{\sigma>2|\rho|\},\qquad
 E_2:=\{\sigma\leq2|\rho|\}.
\]
On \(E_1\) one expands around the positive quantity \(\sigma\), while
on \(E_2\) every negative power of \(\sigma\) is replaced by the
corresponding power of \(|\rho|\).  Applying this decomposition in both
variables of the Riesz potential and then using the bilinear
Hardy--Littlewood--Sobolev inequality gives
\[
 \|N(\rho,\sigma)\|_{\dot H^{-1}}
 \lesssim\|\rho\|_{\dot H^1}^{1+\eta},
 \qquad
 \eta=\min\{p_\alpha-1,2^*-2\}>0.
\]
Thus the remainder is superlinear in \(\rho\), although the power in
the equation is subquadratic.  Proposition \ref{pro5.5} refines this
estimate after projection onto every mode \(Z_r^b\).  In particular,
if \(\varepsilon=\|\rho\|_{\dot H^1}\), the strongly singular scale
projection is bounded by terms of the form
\[
 \varepsilon^2+R^{-\eta_\alpha}\varepsilon^{p_\alpha},
 \qquad
 2b_\alpha>a_\alpha,\qquad
 \eta_\alpha+p_\alpha b_\alpha>a_\alpha,
\]
and is therefore of smaller order than the leading scale interaction.

\paragraph{\textbf{The fixed-point argument.}}
For \(n\geq6\), \((n+2)/2\leq\alpha<n\), and \(\nu=2\), Section
\ref{Sec6} constructs
\[
u=U[-Re_1,1]+U[Re_1,1]+\rho
\]
satisfying
\[
\Gamma(u)\lesssim R^{-a_\alpha},\qquad
\operatorname{dist}_{\dot H^1}(u,\mathcal T_2^{(D)})
\gtrsim\mathfrak H_{2,n,\alpha}(R).
\]
Let \(X_R\) be the orthogonal complement of all two-bubble kernel
modes and let \(\mathcal L_R\) be the projected inverse.  The correction
is a fixed point of
\[
 \mathcal A_R(\rho)
 :=-\mathcal L_R\bigl(h+N(\rho,\sigma)\bigr)
\]
on the weakly compact convex ball
\[
 \mathcal E_R
 :=\{\rho\in X_R:\|\rho\|_{\dot H^1}
 \leq M\mathfrak H_{2,n,\alpha}(R)\}.
\]
In the range \(\alpha>4\), the preceding nonlinear estimate is H\"older
rather than Lipschitz in \(\rho\); the estimates yield invariance of
\(\mathcal E_R\), but not the contraction estimate required by the
Banach fixed-point theorem.  We instead prove that \(\mathcal A_R\) is
weakly continuous.  The required local compactness of the Riesz
potential is obtained by splitting its kernel into a far tail, a
near-diagonal part, and a truncated kernel.  The first two parts are
uniformly small, whereas, on each \(B_L\times K\), the truncated
integral operator is compact from \(L^s(B_L)\) into \(L^\ell(K)\) for
the exponents chosen in Proposition \ref{prop:sharp-correction}.  The
weak topology of \(\mathcal E_R\) is metrizable, and the
Schauder--Tychonoff theorem \cite{Tychonoff1935} yields a fixed point.

The finite-dimensional multipliers are then estimated by the same
projection system as in Proposition \ref{pro4.1.1}.  The scale
projection of \(h\) has order \(R^{-a_\alpha}\), while all projected
nonlinear terms are \(o(R^{-a_\alpha})\).  Finally, after projecting
onto \(X_R\), the leading part of \(h\) gives
\(\|\rho\|_{\dot H^1}\gtrsim\mathfrak H_{2,n,\alpha}(R)\); the nonlinear remainder
is negligible by Lemma \ref{le:strong-nonlinear}.  Since
\(R^{-b_\alpha}=(R^{-a_\alpha})^{b_\alpha/a_\alpha}\), this proves the
optimality of both the logarithmic rate and the exponent
\(\theta_\alpha=b_\alpha/a_\alpha\).

\paragraph{\textbf{Comparison with previous work.}}
For the local Yamabe equation, quantitative multi-bubble stability was
developed in \cite{CFM,Figalli-Glaudo2020,DSW}.  The bubble-tree
geometry used here is based on \cite{DSW}, but the Hartree interaction
is nonlocal: the mixed-bubble function
\(\sigma^{p_\alpha}(y)-\sum_iU_i^{p_\alpha}(y)\) lies inside the
Riesz convolution and is coupled with the exterior factor
\(\sigma^\beta(x)\).  Consequently, pointwise local interaction
estimates alone do not determine either the exponent \(b_\alpha\) or
the logarithmic thresholds.  Lemmas \ref{A.4} and \ref{A.6}, the
endpoint \((5,4)\) argument, and the joint decomposition of the
\(y\)-integration domain and the \(x\)-estimate domain for
\(\alpha>n-2\) provide the required nonlocal substitutes.
In the balanced cases, estimating the interaction directly in
\(\dot H^{-1}\) through the Newtonian representation also avoids the
maximum-principle and \(L^\infty\) estimates that occur in pointwise
local reductions.

Liu, Zhang and Zou \cite{LZZ} considered \(n=3\) and
\(5/2<\alpha<3\).  Piccione, Yang and Zhao
\cite{Piccione-Yang-Zhao} treated \(3\leq n<6-\alpha\) with
\(\alpha\leq4\), and Yang and Zhao \cite{YZ} treated
\(n\geq6-\alpha\) with \(0<\alpha<\min\{4,n\}\).  These arguments do
not cover the strongly singular range \(4<\alpha<n\).  Theorem
\ref{mainthm1} covers every \(n\geq3\) and \(0<\alpha<n\).  The new
range is obtained by combining the dominance-region estimates for the
concave coefficient \(\sigma^{p_\alpha-1}\), the joint
\(y\)- and \(x\)-region decomposition of the Riesz interaction, and
the two-region estimate for the nonlinear remainder.  The construction
in Section \ref{Sec6} further proves that the logarithmic rate and, in
the superbalanced range, the power \(\theta_\alpha\) are optimal.

\subsection{Organization of the paper.}
Section \ref{Sec2} establishes the bubble-tree decomposition.  Section
\ref{Sec3} derives the projected equations and proves the main theorem
from the analytic estimates.  Section \ref{Sec4} estimates the
inhomogeneous interaction \(h\).  Section \ref{Sec5} develops the
uniform linear theory, the nonlinear estimates, and the projected
error bounds.  Section \ref{Sec6} constructs the sharp two-bubble
examples.  Appendix \ref{appA} contains the one-bubble nondegeneracy
and the critical interaction lemmas, while Appendix \ref{appB}
contains the projection expansions in both the weakly and strongly
singular ranges.

\section{Configuration of Bubble Trees}\label{Sec2}

In this section we describe a useful partition of $\{U_i:=U[z_i,\lambda_i]:i\in I\}$ and record its basic structure. This decomposition was introduced by Deng, Sun and Wei in \cite{DSW}. Throughout this section, we assume $\nu\geq2$.

\smallskip

Let $I=\{1,\ldots,\nu\}$. Suppose that $\{U_{i,(k)}:=U[z_{i,(k)},\lambda_{i,(k)}]:i\in I\}_{k=1}^{\infty}$ is a sequence of $\nu$-bubble configurations such that
\[
 Q_{(k)}:=\max_{\substack{i,j\in I\\i\ne j}}q_{ij,(k)}\longrightarrow0
 \qquad\text{as }k\to\infty,
\]
or, equivalently,
\begin{equation}\label{confi1.1}
   R_{(k)}=\frac{1}{2} \min \left\{R_{i j}^{(k)}: \forall i, j \in I, i \neq j\right\} \rightarrow \infty, \quad \text { as } \quad k \rightarrow \infty .
\end{equation}
After relabeling and passing to subsequences finitely many times, we may assume
\begin{equation}\label{confi2}
    \lambda_{1,(k)} \leq \cdots \leq \lambda_{\nu,(k)},
\end{equation}
\begin{equation}\label{confi3}
   \text { either } \quad \lim _{k \rightarrow \infty} z_{i j,(k)} \,\, \text { exists} \qquad  \text{or } \quad \lim _{k \rightarrow \infty}\left|z_{i j,(k)}\right|=\infty,
\end{equation}
where $z_{i j,(k)}:=\lambda_{i,(k)}\left(z_{j,(k)}-z_{i,(k)}\right)$ for $j \in I \backslash\{i\}$.
There is a geometric interpretation of $z_{ij,(k)}$. In the rescaled
$z_i$-centered coordinates $x_i=\lambda_i(x-z_i)$, one has
$U_i(x)=\lambda_i^{(n-2)/2}U(x_i)$, where $U=U[0,1]$. In this paragraph
we write \(U_i,z_i,\lambda_i\) for
\(U_{i,(k)},z_{i,(k)},\lambda_{i,(k)}\). In the $x_i$-coordinates, the
other bubbles become
\[
U_j(x)=\left(\frac{(n(n-2))^{1 / 2} \lambda_j}{1+\left(\lambda_j / \lambda_i\right)^2\left|x_i-z_{i j}\right|^2}\right)^{\frac{n-2}{2}}=\lambda_i^{\frac{n-2}{2}} U\left[z_{i j}, \lambda_j / \lambda_i\right]\left(x_i\right).
\]
Then $z_{i j}$ is the new center of $U_j$, $j \in I \backslash\{i\}$.
Under the $x_i$-coordinates, put
\[
    \lambda'_p:=\frac{\lambda_p}{\lambda_i},\qquad
    z'_p:=\lambda_i(z_p-z_i).
\]
For every $p\ne q$,
\[
 \frac{\lambda'_p}{\lambda'_q}=\frac{\lambda_p}{\lambda_q},
 \qquad
 \sqrt{\lambda'_p\lambda'_q}\,|z'_p-z'_q|
 =\sqrt{\lambda_p\lambda_q}\,|z_p-z_q|.
\]
Thus each of the three quantities in the definition of $R_{pq}$ is
unchanged, and hence $R'_{pq}=R_{pq}$.

We define a partial order $\prec$ on $I=\{1, \ldots, \nu\}$ by
\[
\begin{aligned}
& i \prec j \quad \Longleftrightarrow \quad i<j \text { and } \lim _{k \rightarrow \infty} z_{i j,(k)} \text { exists, } \\
& i \preceq j \quad \Longleftrightarrow \quad i \prec j \text { or } i=j .
\end{aligned}
\]
\begin{lemma}\cite[Lemma 4.4]{DSW}
    $\prec$ is a strict partial order.
\end{lemma}
\begin{lemma}\cite[Lemma 4.5]{DSW}
Suppose that $\mathscr U=\{U_{i,(k)}:i=1,\ldots,\nu,\ k\geq1\}$ is a
sequence of $\nu$-bubble configurations satisfying
\eqref{confi1.1}--\eqref{confi3}. After discarding finitely many terms,
\begin{equation}\label{confi4}
 C_{\mathscr U}:=1+
 \max\left(\{0\}\cup
 \left\{\sup_{k\geq1}|z_{ij,(k)}|:
 i,j\in I,\ i\prec j\right\}\right)<\infty .
\end{equation}
Moreover,
\begin{equation}\label{confi7}
 i\prec j\quad\Longleftrightarrow\quad
 \sqrt{\lambda_{j,(k)}/\lambda_{i,(k)}}
 \leq R_{ij}^{(k)}
 \leq C_{\mathscr U}\sqrt{\lambda_{j,(k)}/\lambda_{i,(k)}},
 \qquad k\geq1,
\end{equation}
\begin{equation}\label{confi8}
i \nprec j \text { and } j \nprec i \quad \Leftrightarrow \quad \sqrt{\lambda_{j,(k)} / \lambda_{i,(k)}}+\sqrt{\lambda_{i,(k)} / \lambda_{j,(k)}}=o(1) R_{i j}^{(k)} .
\end{equation}
Here $o(1)$ denotes a quantity tending to zero as $k\to\infty$. The
constant $C_{\mathscr U}$ is attached to the extracted sequence and is not
asserted to be uniform over all configurations. This causes no loss in the
uniform estimates below: those estimates are proved by contradiction along a
single sequence with $R_{(k)}\to\infty$. We call $i$ and $j$ incomparable
when \eqref{confi8} holds.
\end{lemma}

Recall that a tree is a partially ordered set, say $(\mathrm{T}, \prec)$, such that for any $t \in \mathrm{~T}$, the set $\{s \in \mathrm{~T}: s \prec t\}$ is well-ordered by the relation $\prec$. The following lemma shows that $I$ can be decomposed into several trees.

\begin{lemma}\cite[Lemma 4.6]{DSW}
    For every sequence $\left\{U_{i,(k)}:i=1,\ldots,\nu\right\}$ satisfying \eqref{confi1.1}--\eqref{confi3}, there exists an integer $\beta^*$, depending on the sequence, such that $I=\{1,\ldots,\nu\}$ can be partitioned into trees $\mathrm T_\beta$, $\beta=1,\ldots,\beta^*$.
\end{lemma}
For each $i\in I$, define the set of successors of $i$ by
\begin{equation*}
    S(i)=\left\{j \in I: i<j, \lim _{k \rightarrow \infty} z_{i j,(k)} \text { exists }\right\}.
\end{equation*}

\section{Proof of the Main Theorem}\label{Sec3}

In this section, we prove Theorem \ref{mainthm1} using the interaction and error estimates established in the subsequent sections and appendices.

\smallskip

For $i=1, \ldots, \nu$, let us define
\begin{equation}\label{eq2.1}
    \begin{split}
        Z_i^a(x) & :=\left.\frac{1}{\lambda_i} \frac{\partial U\left[z, \lambda_i\right]}{\partial z^a}\right|_{z=z_i}=(n-2) U\left[z_i, \lambda_i\right](x) \frac{\lambda_i\left(x^a-z_i^a\right)}{1+\lambda_i^2\left|x-z_i\right|^2}, \\
Z_i^{n+1}(x) & :=\left.\lambda_i \frac{\partial U\left[z_i, \lambda\right]}{\partial \lambda}\right|_{\lambda=\lambda_i}=\frac{n-2}{2} U\left[z_i, \lambda_i\right](x) \frac{1-\lambda_i^2\left|x-z_i\right|^2}{1+\lambda_i^2\left|x-z_i\right|^2},
    \end{split}
\end{equation}
where $z^a$ is the $a$-th component of $z$ for $a=1, \ldots, n$.

\begin{lemma}
\label{lem:uniform-remainders}
Fix $n$, $\nu\geq2$, and $\alpha$.  Let $E$ be a scalar quantity
defined for every $\nu$-bubble configuration with $0<Q<1$, and let
$s>0$.  Assume that for every sequence of configurations with
$Q_k\to0$ one has
\[
 \frac{|E_k|}{Q_k^s}\longrightarrow0.
\]
Then there is a function $\omega_1:(0,1)\to[0,\infty)$ such that
\[
 \lim_{t\downarrow0}\omega_1(t)=0,
 \qquad
 |E|\leq\omega_1(Q)Q^s.
\]
The same conclusion holds for a finite family of remainders with one
common modulus.
\end{lemma}

\begin{proof}
For $0<t<1$, define
\[
 \omega_1(t):=\sup\left\{
 \frac{|E|}{Q^s}:0<Q\leq t
 \right\}.
\]
If $\omega_1(t)$ failed to be finite for every sufficiently small
$t$, there would be configurations with $Q_k\leq1/k$ and
$|E_k|/Q_k^s\geq k$, contradicting the hypothesis.  If
$\omega_1(t)\not\to0$, there would be $c>0$, numbers $t_k\downarrow0$,
and configurations with $Q_k\leq t_k$ and
$|E_k|/Q_k^s\geq c/2$, which is the same contradiction.  For finitely
many remainders, take the maximum of their moduli.
\end{proof}

Throughout the remainder of the paper, \(n\), \(\alpha\), and \(\nu\)
are fixed.  Unless a little-\(o\) term is explicitly associated with a
prescribed or extracted sequence, every remainder in an interaction
estimate is understood uniformly over all admissible bubble
configurations and all mode indices.  In particular, remainders denoted
by \(o(q_{ij}^s)\), \(o(Q^s)\), and \(o(R^{-c})\) mean, respectively,
\[
\begin{aligned}
 \lim_{t\downarrow0}
 \sup_{0<q_{ij}\leq t}\frac{|E|}{q_{ij}^s}=0,\qquad
 \lim_{t\downarrow0}
 \sup_{0<Q\leq t}\frac{|E|}{Q^s}=0,\qquad
 \lim_{T\to\infty}
 \sup_{R\geq T}R^c|E|&=0,
\end{aligned}
\]
where each supremum is taken over all admissible configurations and
mode indices.  For every finite collection of remainders, one common
modulus is understood.  When \(\nu\geq2\), Lemma
\ref{lem:uniform-remainders} and \(Q\asymp R^{-(n-2)}\) show that the
formulations in terms of \(Q\) and \(R\) are equivalent.
In particular, the remainders furnished by Lemma \ref{A.2}, Lemma
\ref{le:weak-projection-expansion}, and Lemma
\ref{le:strong-projection-proof} satisfy this convention.  Little-\(o\)
terms attached to a fixed sequence in a bubble-tree, subsequence, or
contradiction argument retain their usual sequential meaning.

Suppose $u=\sigma+\rho$, where
$\sigma=\sum_{i=1}^{\nu}U_i$ is the minimizing approximation supplied
by the standard modulation construction fixed above.  Then
\begin{equation}\label{2.1}
    \begin{split}
            &\Delta u + D_{n,\alpha}\int_{\mathbb{R}^{n}}\frac{|u|^{p_{\alpha}}(y)}{|x-y|^{\alpha}}\dy\,g_{p_\alpha}(u(x))
            \\=&\Delta (\sigma+\rho) + D_{n,\alpha}\int_{\mathbb{R}^{n}}\frac{|\sigma+\rho|^{p_{\alpha}}(y)   }{|x-y|^{\alpha}}\dy|\sigma+\rho|^{p_{\alpha}-2}(x)(\sigma+\rho)(x)
    \\=&\Delta \rho +N_1(\rho,\sigma) +h+N(\rho,\sigma),
    \end{split}
\end{equation}
    where
    \begin{align*}
        h =D_{n,\alpha}\bigg(\int_{\mathbb{R}^{n}}\frac{\sigma^{p_{\alpha}}(y)   }{|x-y|^{\alpha}}\dy \sigma^{p_{\alpha}-1}(x)-\sum_{i=1}^{\nu}\int_{\mathbb{R}^{n}}\frac{U_{i}^{p_{\alpha}}(y)   }{|x-y|^{\alpha}}\dy U^{p_{\alpha}-1}_{i}(x)\bigg),
    \end{align*}
    \begin{align*}
         N_{1}(\eta,\xi) := D_{n,\alpha} \bigg(\int_{\mathbb{R}^{n}}\frac{p_{\alpha}|\xi|^{p_{\alpha}-2}(y)\xi(y)\eta(y)   }{|x-y|^{\alpha}}\dy |\xi|^{p_{\alpha}-2}(x)\xi(x) +\int_{\mathbb{R}^{n}}\frac{|\xi|^{p_{\alpha}}(y)  }{|x-y|^{\alpha}}\dy (p_{\alpha}-1)|\xi|^{p_{\alpha}-2}(x) \eta(x)   \bigg),
     \end{align*}
    and
     \begin{align*}
        N(\eta,\xi):= &D_{n,\alpha}\bigg(\int_{\mathbb{R}^{n}}\frac{|\eta+\xi|^{p_{\alpha}}(y)   }{|x-y|^{\alpha}}\dy |\eta+\xi|^{p_{\alpha}-2}(x)(\eta+\xi)(x) - \int_{\mathbb{R}^{n}}\frac{|\xi|^{p_{\alpha}}(y)   }{|x-y|^{\alpha}}\dy |\xi|^{p_{\alpha}-2}(x)\xi(x)
        \\&-\int_{\mathbb{R}^{n}}\frac{p_{\alpha}|\xi|^{p_{\alpha}-2}(y)\xi(y)\eta(y)   }{|x-y|^{\alpha}}\dy |\xi|^{p_{\alpha}-2}(x)\xi(x) -\int_{\mathbb{R}^{n}}\frac{|\xi|^{p_{\alpha}}(y)  }{|x-y|^{\alpha}}\dy (p_{\alpha}-1)|\xi|^{p_{\alpha}-2}(x) \eta(x) \bigg).
    \end{align*}
    Let $ f =-\Delta u - D_{n,\alpha}\int_{\mathbb{R}^{n}}\frac{|u|^{p_{\alpha}}(y)}{|x-y|^{\alpha}}\dy\,g_{p_\alpha}(u(x))$. Then \eqref{2.1} can be rewritten as
    \begin{equation}\label{2.2}
        \begin{split}
            &\Delta \rho +N_1(\rho,\sigma) +h+N(\rho,\sigma)+f=0.
        \end{split}
    \end{equation}

\begin{lemma}\label{lem:residual-smallness}
Let
\[
 \mathcal K(v):=
 \bigl(|x|^{-\alpha}*|v|^{p_\alpha}\bigr)
 |v|^{p_\alpha-2}v,
 \qquad
 \mu_\alpha:=\min\{1,p_\alpha-1\}.
\]
For every $M>0$ there is $C_M>0$ such that, whenever
$v,w\in\dot H^1(\mathbb R^n)$ and
$\|v\|_{\dot H^1}+\|w\|_{\dot H^1}\leq M$, one has
\begin{equation}\label{eq:hartree-residual-continuity}
 \|\mathcal K(v)-\mathcal K(w)\|_{\dot H^{-1}}
 \leq C_M\left(
 \|v-w\|_{\dot H^1}+
 \|v-w\|_{\dot H^1}^{\mu_\alpha}
 \right).
\end{equation}
Consequently, after the constant $\delta$ in Theorem
\ref{mainthm1} is decreased depending only on $(n,\nu,\alpha)$,
\begin{equation}\label{eq:residual-at-most-one}
 \|f\|_{\dot H^{-1}}=\Gamma(u)\leq1.
\end{equation}
\end{lemma}

\begin{proof}
Put
$
 s_\alpha:=\frac{2n}{2n-\alpha}$,
and $p_\alpha s_\alpha=2^*$. The bilinear Hardy--Littlewood--Sobolev inequality states that
\begin{equation}\label{eq:bilinear-HLS-residual}
 \left|\int_{\mathbb R^n}
 (|x|^{-\alpha}*F)G\,\dx\right|
 \leq C\|F\|_{L^{s_\alpha}}\|G\|_{L^{s_\alpha}}.
\end{equation}
For $p:=p_\alpha>1$, a direct computation, using
$s_\alpha^{-1}=p/(2^*)$, and H\"older inequality give
$$
 \||v|^p-|w|^p\|_{L^{s_\alpha}}
 \leq C M^{p-1}\|v-w\|_{L^{2^*}},\qquad \\
 \||v|^{p-1}\varphi\|_{L^{s_\alpha}}
 \leq M^{p-1}\|\varphi\|_{L^{2^*}},
$$
and
\begin{equation}\label{eq:residual-concave-power}
 \|(|v|^{p-2}v-|w|^{p-2}w)\varphi\|_{L^{s_\alpha}}
 \leq C
 \begin{cases}
 M^{p-2}\|v-w\|_{L^{2^*}}\|\varphi\|_{L^{2^*}},&p\geq2,\\
 \|v-w\|_{L^{2^*}}^{p-1}\|\varphi\|_{L^{2^*}},&1<p<2.
 \end{cases}
\end{equation}
Decompose
\[
 \begin{split}
 \mathcal K(v)-\mathcal K(w)
 =\bigl(|x|^{-\alpha}*(|v|^p-|w|^p)\bigr)
 |v|^{p-2}v+\bigl(|x|^{-\alpha}*|w|^p\bigr)
 \bigl(|v|^{p-2}v-|w|^{p-2}w\bigr).
 \end{split}
\]
Pairing this identity with $\varphi\in\dot H^1$, applying
\eqref{eq:bilinear-HLS-residual}--\eqref{eq:residual-concave-power},
and using Sobolev inequality proves
\eqref{eq:hartree-residual-continuity}.

For the modulated sum $\sigma$, the single-bubble equation and the
definition of $h$ give the exact identity
\[
 -\Delta\sigma-D_{n,\alpha}\mathcal K(\sigma)=-h.
\]
If $\nu=1$, then $h=0$.  If $\nu\geq2$, the estimates in Lemmas
\ref{le7.1} and \ref{le3.1}, together with
$L^{2n/(n+2)}\hookrightarrow\dot H^{-1}$, imply
\begin{equation*}
 \|h\|_{\dot H^{-1}}\longrightarrow0
 \qquad\text{uniformly as }Q\longrightarrow0.
\end{equation*}
Furthermore,
$\|\sigma\|_{\dot H^1}\leq\nu\|U[0,1]\|_{\dot H^1}$ and
$\rho=u-\sigma$ satisfies $\|\rho\|_{\dot H^1}\leq\delta$.  From the
definitions of $f$ and $\rho$, we have the following identity in
$\dot H^{-1}(\mathbb R^n)$:
\begin{equation}\label{eq:residual-exact-decomposition}
\begin{aligned}
 f
 =-\Delta u-D_{n,\alpha}\mathcal K(u)=-h-\Delta\rho
   -D_{n,\alpha}\bigl(\mathcal K(u)-\mathcal K(\sigma)\bigr).
\end{aligned}
\end{equation}
After decreasing $\delta$ so that $\delta\leq1$, we have
\[
 \|u\|_{\dot H^1}+\|\sigma\|_{\dot H^1}
 \leq2\nu\|U[0,1]\|_{\dot H^1}+1=:M_{n,\nu}.
\]
Applying \eqref{eq:hartree-residual-continuity} with
$(v,w)=(u,\sigma)$ therefore yields
\[
 \|\mathcal K(u)-\mathcal K(\sigma)\|_{\dot H^{-1}}
 \leq C_{M_{n,\nu}}
 \left(\|\rho\|_{\dot H^1}
       +\|\rho\|_{\dot H^1}^{\mu_\alpha}\right).
\]
Taking the $\dot H^{-1}$-norm in
\eqref{eq:residual-exact-decomposition}, and using
$\|\rho\|_{\dot H^1}\leq\delta$, gives
\begin{align*}
 \|f\|_{\dot H^{-1}}
 \leq\|h\|_{\dot H^{-1}}
 +\|\rho\|_{\dot H^1}
 +D_{n,\alpha}\|\mathcal K(u)-\mathcal K(\sigma)\|_{\dot H^{-1}}\leq\|h\|_{\dot H^{-1}}
 +C(\delta+\delta^{\mu_\alpha}).
\end{align*}
Here \(C=C(n,\nu,\alpha)\), since \(D_{n,\alpha}\) is fixed and
\(M_{n,\nu}\) is independent of all bubble parameters.
Since \(Q\leq\omega_{\rm mod}(\delta)\) by \eqref{introeq2}, the
right-hand side is at most $1$ after decreasing $\delta$, proving
\eqref{eq:residual-at-most-one}.
\end{proof}

     Testing \eqref{2.2} against $Z_{r}^{n+1}$ gives
     \begin{align*}
         \left|\int h Z_{r}^{n+1}\dx\right|
         \leq \left|\int N(\rho,\sigma) Z_{r}^{n+1}\dx\right|
         +\left|\left\langle f,Z_r^{n+1}\right\rangle_{\dot H^{-1},\dot H^1}\right|
         + \left|\int N_{1}(\rho,\sigma) Z_{r}^{n+1}\dx\right|.
     \end{align*}
     Here we have used the orthogonal conditions \eqref{introeq2}. It follows from \eqref{eq2.1} that $\left\|\nabla Z_r^{n+1}\right\|_{L^2}^2 \lesssim\left\|U_r\right\|_{L^{2^*}}^{2^*} $, up to a dimensional constant independent of $z_{r}$ and $\lambda_{r}$. Thus $\left|\langle f,Z_r^{n+1}\rangle\right| \leq C \left\|f\right\|_{\dot{H}^{-1}}$. Hence
       \begin{equation}\label{2.4}
              \left|\int h Z_{r}^{n+1}\dx\right| \leq \left|\int N(\rho,\sigma) Z_{r}^{n+1}\dx\right| +C\left\|f\right\|_{\dot{H}^{-1}} + \left|\int N_{1}(\rho,\sigma) Z_{r}^{n+1}\dx\right|.
       \end{equation}

     \begin{lemma}\label{le2.1}
         Suppose that $\nu\geq2$, $n\geq3$, $0<\alpha\leq4$, and $u$
         satisfies \eqref{1.1} with $\delta$ small enough. The following
         expansions hold for every $1\leq r\leq\nu$. If $n\geq6$, then
         $$\int h Z_{r}^{n+1}\dx
                 = \sum_{j=1,j\ne r}^{\nu}\int U_{j}^{2^*-p_{\alpha}}  U_{r}^{p_{\alpha}-1}Z_{r}^{n+1}(x) \dx+o\left(R^{-a_\alpha}\right)
                 = -c_{1,\alpha}\sum_{j\ne r}\mathfrak d_{rj}q_{rj}^{\frac{\alpha}{n-2}}
                 +o\left(Q^{\frac{\alpha}{n-2}}\right),
                 \quad \alpha<n-2,$$
         \begin{equation*}
             \int h Z_{r}^{n+1}\dx
                 = 3\sum_{j=1,j\ne r}^{\nu}\int U_jU_rZ_r^{n+1}(x)\dx+o(R^{-4})
                 = -c_3\sum_{j\ne r}\mathfrak d_{rj}q_{rj}+o(Q),
                 \quad n=6,\ \alpha=4,
         \end{equation*}
         where $c_{1,\alpha},c_3>0$. If $3\leq n\leq 5$, then
         \begin{equation*}
             \begin{split}
                 \int h Z_{r}^{n+1}\dx
                = &\sum_{j=1,j\ne r}^{\nu}\int U_{j}^{2^*-p_{\alpha}}  U_{r}^{p_{\alpha}-1}Z_{r}^{n+1}(x) \dx +\int_{\mathbb{R}^{n} }(2^*-1)U_{r}^{2^*-2}\sum_{j\ne r}^{\nu}U_{j}Z_{r}^{n+1}(x)\dx
                 +o\left(R^{-a_\alpha}\right)\\
                = &-\sum_{j\ne r}\mathfrak d_{rj}
                \left(c_{1,\alpha}q_{rj}^{\frac{\alpha}{n-2}}
                +c_2q_{rj}\right)+o\left(Q^{\frac{a_\alpha}{n-2}}\right),
             \end{split}
         \end{equation*}
         Here
         \[
         c_{1,\alpha},c_2\geq0,\qquad
         \begin{cases}
         c_{1,\alpha}>0,&\alpha<n-2,\\
         c_2>0,&\alpha>n-2,\\
         c_{1,\alpha}+c_2>0,&\alpha=n-2.
         \end{cases}
         \]
         In all the preceding cases, there exist a constant
         $c=c(n,\nu,\alpha)>0$ and a nondecreasing function
         $\omega_2:[0,1)\to[0,\infty)$, depending only on
         $(n,\nu,\alpha)$, independent of the bubble parameters, and satisfying
         $\omega_2(0)=0$ and $\lim_{t\downarrow0}\omega_2(t)=0$, such that
         there exists an index $r$ for which
         \begin{equation}\label{eq:weak-projection-lower}
         \left|\int hZ_r^{n+1}\dx\right|
         \geq\bigl(c-\omega_2(Q)\bigr)
         Q^{\frac{a_\alpha}{n-2}}.
         \end{equation}
         After decreasing $\delta$ if necessary,
         $\omega_2(Q)\leq c/2$.
      \end{lemma}
     \begin{proof}
Put
\[
 \beta:=p_\alpha-1,
 \qquad
 \gamma:=2^*-p_\alpha,
 \qquad
 \kappa:=\frac{a_\alpha}{n-2}
 =\min\left\{\frac{\alpha}{n-2},1\right\}.
\]
For $k\ne j$, set
\[
 \mathcal J^{(\alpha)}_{kj}
 :=\int_{\mathbb R^n}U_j^\gamma U_k^\beta Z_k^{n+1}\dx,
 \qquad
 \mathcal J^{(0)}_{kj}
 :=\int_{\mathbb R^n}U_k^{2^*-2}U_jZ_k^{n+1}\dx.
\]
Lemma \ref{le:weak-projection-expansion}, namely
\eqref{eq:weak-range-master-expansion}, gives
\begin{equation}\label{eq:le21-master-expansion}
 \int_{\mathbb R^n}hZ_k^{n+1}\dx
 =\sum_{j\ne k}\mathcal J^{(\alpha)}_{kj}
 +(2^*-1)\sum_{j\ne k}\mathcal J^{(0)}_{kj}
 +\varepsilon_k,
 \qquad
 \max_{1\leq k\leq\nu}|\varepsilon_k|=o(R^{-a_\alpha}).
\end{equation}
Lemma \ref{A.2} gives, uniformly for $k\ne j$,
$$
 \mathcal J^{(\alpha)}_{kj}
 =-C_{2,n,\alpha}\mathfrak d_{kj}
 q_{kj}^{\frac{\alpha}{n-2}}
 +o\left(q_{kj}^{\frac{\alpha}{n-2}}\right),
\qquad
 \mathcal J^{(0)}_{kj}
 =-C_{1,n}\mathfrak d_{kj}q_{kj}+o(q_{kj}).
$$
For
\[
 S_{ij}:=\frac{\lambda_i}{\lambda_j}
 +\frac{\lambda_j}{\lambda_i}
 +\lambda_i\lambda_j|z_i-z_j|^2
 =q_{ij}^{-\frac2{n-2}},
\]
one has
\[
 R_{ij}^2\leq S_{ij}\leq3R_{ij}^2,
 \qquad
 3^{-\frac{n-2}{2}}(2R)^{-(n-2)}
 \leq Q\leq(2R)^{-(n-2)},
\]
and hence
\begin{equation}\label{eq:le21-RQ}
 R^{-a_\alpha}\asymp Q^\kappa.
\end{equation}
If $n\geq6$ and $\alpha<n-2$, then
$
 Q=o(Q^\kappa)$.
 \eqref{eq:le21-master-expansion}--\eqref{eq:le21-RQ} imply
\begin{align*}
 \int hZ_k^{n+1}\dx
 =\sum_{j\ne k}\mathcal J^{(\alpha)}_{kj}
 +o(R^{-a_\alpha})=-C_{2,n,\alpha}\sum_{j\ne k}\mathfrak d_{kj}
 q_{kj}^{\frac\alpha{n-2}}
 +o\left(Q^{\frac\alpha{n-2}}\right).
\end{align*}
If $(n,\alpha)=(6,4)$, then
\[
 \beta=\gamma=1,\qquad 2^*=3,\qquad
 \mathcal J^{(\alpha)}_{kj}=\mathcal J^{(0)}_{kj},
\]
and
\[
 \int hZ_k^{n+1}\dx
 =3\sum_{j\ne k}\int U_jU_kZ_k^{n+1}\dx+o(R^{-4})
 =-3C_{1,6}\sum_{j\ne k}\mathfrak d_{kj}q_{kj}+o(Q).
\]
If $3\leq n\leq5$, then
\begin{align*}
 \int hZ_k^{n+1}\dx
 &=\sum_{j\ne k}\mathcal J^{(\alpha)}_{kj}
 +(2^*-1)\sum_{j\ne k}\mathcal J^{(0)}_{kj}
 +o(R^{-a_\alpha})\\
 &=-\sum_{j\ne k}\mathfrak d_{kj}
 \left(C_{2,n,\alpha}q_{kj}^{\frac\alpha{n-2}}
 +(2^*-1)C_{1,n}q_{kj}\right)
 +o(Q^\kappa).
\end{align*}
Therefore
\[
 c_{1,\alpha}:=C_{2,n,\alpha},
 \qquad
 c_2:=(2^*-1)C_{1,n},
 \qquad
 c_3:=3C_{1,6}.
\]
For $i\ne j$, define
\[
 w_{ij}=w_{ji}:=
 \begin{cases}
 C_{2,n,\alpha}q_{ij}^{\frac\alpha{n-2}},
 &n\geq6,\ \alpha<n-2,\\[1mm]
 3C_{1,6}q_{ij},
 &(n,\alpha)=(6,4),\\[1mm]
 C_{2,n,\alpha}q_{ij}^{\frac\alpha{n-2}}
 +(2^*-1)C_{1,n}q_{ij},
 &3\leq n\leq5.
 \end{cases}
\]
There exist $c_-,c_+>0$ such that
\begin{equation*}
 c_-q_{ij}^\kappa\leq w_{ij}\leq c_+q_{ij}^\kappa.
\end{equation*}
Thus
\begin{equation}\label{eq:le21-unified-projection}
 \int hZ_k^{n+1}\dx
 =-\sum_{j\ne k}\mathfrak d_{kj}w_{kj}+r_k,
 \qquad
 \max_{1\leq k\leq\nu}|r_k|=o(Q^\kappa).
\end{equation}
By the uniformity convention above, the function $\omega_2$ in the
statement may be chosen so that
\[
 \max_{1\leq k\leq\nu}|r_k|
 \leq\omega_2(Q)Q^\kappa.
\]
After a permutation,
$
 \lambda_1\leq\lambda_2\leq\cdots\leq\lambda_\nu.
$
For $i<j$,
\[
 \mathfrak d_{ji}
 =\frac{\lambda_j/\lambda_i-\lambda_i/\lambda_j
 +\lambda_i\lambda_j|z_i-z_j|^2}{S_{ij}}
 =1-\frac{2\lambda_i/\lambda_j}{S_{ij}},
\]
and hence
\[
 0\leq1-\mathfrak d_{ji}
 \leq2S_{ij}^{-1}
 =2q_{ij}^{\frac2{n-2}}
 \leq2Q^{\frac2{n-2}}
 \leq2\omega_{\rm mod}(\delta)^{\frac2{n-2}}.
\]
For
$
 2\omega_{\rm mod}(\delta)^{\frac2{n-2}}\leq\frac12,
$
one has
\begin{equation}\label{eq:le21-positive-direction}
 \mathfrak d_{ji}\geq\frac12,
 \qquad 1\leq i<j\leq\nu.
\end{equation}
Set
\[
 M:=\max_{1\leq k\leq\nu}
 \left|\sum_{j\ne k}\mathfrak d_{kj}w_{kj}\right|,
 \qquad
 A_l:=\sum_{i=1}^{l-1}w_{il},
 \qquad
 T_l:=\sum_{j=l}^{\nu}A_j,
 \qquad T_{\nu+1}:=0.
\]
For $2\leq l\leq\nu$, \eqref{eq:le21-positive-direction} and
$|\mathfrak d_{lj}|\leq1$ give
\begin{align*}
 \frac12A_l
 \leq\sum_{i=1}^{l-1}\mathfrak d_{li}w_{il}=\sum_{j\ne l}\mathfrak d_{lj}w_{lj}
 -\sum_{j=l+1}^{\nu}\mathfrak d_{lj}w_{lj}\leq M+\sum_{j=l+1}^{\nu}w_{lj}
 \leq M+T_{l+1}.
\end{align*}
Therefore
\[
 A_l\leq2M+2T_{l+1},
 \qquad
 T_l=A_l+T_{l+1}\leq2M+3T_{l+1}.
\]
The descending induction gives
\[
 T_l\leq\left(3^{\nu-l+1}-1\right)M,
 \qquad 2\leq l\leq\nu;
\qquad
 \sum_{1\leq i<j\leq\nu}w_{ij}
 =T_2\leq(3^{\nu-1}-1)M.
\]
Choose $(i_0,j_0)$ such that $q_{i_0j_0}=Q$.  Then
\[
 c_-Q^\kappa\leq w_{i_0j_0}
 \leq T_2\leq(3^{\nu-1}-1)M,
\qquad
 M\geq c_*Q^\kappa,
 \qquad
 c_*:=\frac{c_-}{3^{\nu-1}-1}>0.
 \]
Choose $r$ such that
$
 M=\left|\sum_{j\ne r}\mathfrak d_{rj}w_{rj}\right|
$.
Equation \eqref{eq:le21-unified-projection} yields
\[
 \left|\int hZ_r^{n+1}\dx\right|
 \geq M-|r_r|
 \geq\bigl(c_*-\omega_2(Q)\bigr)Q^\kappa.
\]
By \eqref{introeq2}, \(Q\leq\omega_{\rm mod}(\delta)\).  Choose
\(\delta\) so that
\(\omega_2(\omega_{\rm mod}(\delta))\leq c_*/2\).  Then the right-hand side is
at least $(c_*/2)Q^\kappa$.  Taking $c=c_*$ and using
$\kappa=a_\alpha/(n-2)$ gives
\eqref{eq:weak-projection-lower}.
    \end{proof}

     \begin{lemma}\label{le2.1strong}
         Suppose that $\nu\geq2$, $4<\alpha<n$ and $u$ satisfies \eqref{1.1} with $\delta$ small enough. Then, for every scale direction $Z_r^{n+1}$,
         \begin{equation}\label{eq:strong-projection-expansion}
             \begin{split}
             \int h Z_r^{n+1}\dx
             =&\sum_{j\ne r}\int U_j^{2^*-p_\alpha}U_r^{p_\alpha-1}Z_r^{n+1}\dx
             +(2^*-1)\sum_{j\ne r}\int U_r^{2^*-2}U_jZ_r^{n+1}\dx+o(R^{-a_\alpha})\\
             =&-\sum_{j\ne r}\mathfrak d_{rj}
             \left(c_{1}q_{rj}^{\frac{\alpha}{n-2}}+c_{2}q_{rj}\right)
             +o\left(Q^{\frac{a_\alpha}{n-2}}\right).
             \end{split}
         \end{equation}
         Here
         \[
         c_1,c_2\geq0,\qquad
         \begin{cases}
         c_1>0,&\alpha<n-2,\\
         c_2>0,&\alpha>n-2,\\
         c_1+c_2>0,&\alpha=n-2.
         \end{cases}
         \]
         Consequently, there exist a constant
         $c=c(n,\nu,\alpha)>0$ and a nondecreasing function
         $\omega_3:[0,1)\to[0,\infty)$, depending only on
         $(n,\nu,\alpha)$, independent of the bubble parameters, and satisfying
         $\omega_3(0)=0$ and $\lim_{t\downarrow0}\omega_3(t)=0$, such that
         there exists a scale direction $Z_r^{n+1}$ for which
         \begin{equation}\label{eq:strong-projection-lower}
             \left|\int hZ_r^{n+1}\dx\right|
             \geq\bigl(c-\omega_3(Q)\bigr)
             Q^{\frac{a_\alpha}{n-2}}.
         \end{equation}
         After decreasing $\delta$ if necessary,
         $\omega_3(Q)\leq c/2$.
     \end{lemma}
     \begin{proof}
Set \(\kappa=a_\alpha/(n-2)=\min\{\alpha/(n-2),1\}\). Lemma
\ref{le:strong-projection-proof} gives
\[
 \int_{\mathbb R^n}hZ_k^{n+1}\dx
 =-\sum_{j\ne k}\mathfrak d_{kj}w_{kj}+\varepsilon_k,\qquad
 w_{ij}:=c_1q_{ij}^{\alpha/(n-2)}+c_2q_{ij},
 \qquad \max_k|\varepsilon_k|=o(Q^\kappa).
\]
By the uniformity convention above, the function $\omega_3$ in the
statement may be chosen so that
$
 \max_k|\varepsilon_k|\leq\omega_3(Q)Q^\kappa.
$
The alternatives for \(c_1,c_2\) imply
\[
 c_-q_{ij}^{\kappa}\leq w_{ij}=w_{ji}\leq c_+q_{ij}^{\kappa}
 \qquad(i\ne j)
\]
for some \(c_\pm>0\).  Relabel so that
\(\lambda_1\leq\cdots\leq\lambda_\nu\).  For \(i<j\), the computation
in \eqref{eq:le21-positive-direction} applies verbatim and yields
\[
 \mathfrak d_{ji}
 =1-\frac{2\lambda_i/\lambda_j}
 {\lambda_i/\lambda_j+\lambda_j/\lambda_i+
   \lambda_i\lambda_j|z_i-z_j|^2}
 \geq\frac12
\]
after decreasing \(\delta\).  Hence the finite-dimensional descending
induction used in the proof of \eqref{eq:weak-projection-lower}, with the
present \(w_{ij}\), gives an index \(r\) such that
$
 \left|\sum_{j\ne r}\mathfrak d_{rj}w_{rj}\right|
 \geq c_*Q^\kappa$.
Consequently,
\[
 \left|\int hZ_r^{n+1}\dx\right|
 \geq\bigl(c_*-\omega_3(Q)\bigr)Q^\kappa.
\]
Since \(Q\leq\omega_{\rm mod}(\delta)\) by \eqref{introeq2}, choosing
\(\omega_3(\omega_{\rm mod}(\delta))\leq c_*/2\), and taking
\(c=c_*\), proves
\eqref{eq:strong-projection-lower}.
     \end{proof}

     We now return to \eqref{2.4}. Since $R\approx Q^{-\frac{1}{n-2}}$, Proposition \ref{pro5.5} provides the two estimates needed to obtain
     \begin{equation*}
         \left|\int N(\rho,\sigma) Z_{r}^{n+1}\dx\right|+ \left|\int N_{1}(\rho,\sigma) Z_{r}^{n+1}\dx\right| = o\left(Q^{\frac{a_\alpha}{n-2}}\right)+ O\left(\left\|  f \right\|_{\dot{H}^{-1}}\right).
     \end{equation*}
     \begin{lemma}\label{le2.2}
         Suppose that $\nu\geq2$, $n\geq3$, $0<\alpha<n$, and $u$
         satisfies \eqref{1.1} with $\delta$ small enough. Then
         \begin{equation*}
             \begin{split}
            Q^{\frac{a_\alpha}{n-2}} \lesssim \left\|  f \right\|_{\dot{H}^{-1}}.
             \end{split}
         \end{equation*}
     \end{lemma}
     \begin{proof}
Let $\kappa:=a_\alpha/(n-2)$.  If $\nu=1$, then $Q=0$ and the
conclusion is immediate.  Assume $\nu\geq2$, set
\[
 \iota_\alpha:=
 \begin{cases}
 2,&0<\alpha\leq4,\\
 3,&4<\alpha<n.
 \end{cases}
\]
For $0<\alpha\leq4$, Lemma \ref{le2.1} gives an index $r$ and a
constant $c_0>0$ with modulus $\omega_2$.  For $4<\alpha<n$, Lemma
\ref{le2.1strong} gives the same conclusion, possibly with a different
$c_0$, with modulus $\omega_3$.  Thus, in either case,
\[
 \left|\int hZ_r^{n+1}\,\dx\right|
 \geq\bigl(c_0-\omega_{\iota_\alpha}(Q)\bigr)Q^\kappa.
\]
At $(n,\alpha)=(6,4)$ the leading coefficient
$3C_{1,6}$ is strictly positive; hence the constant $c_0$ in the
conclusion of Lemma \ref{le2.1} is positive.

Testing \eqref{2.2} by $Z_r^{n+1}$, using
\eqref{introeq2}, and using the parameter-independent bound
$\|Z_r^{n+1}\|_{\dot H^1}\leq C_n$, we obtain
\begin{equation}\label{eq:interaction-absorption-first}
 \bigl(c_0-\omega_{\iota_\alpha}(Q)\bigr)Q^\kappa
 \leq C\|f\|_{\dot H^{-1}}
 +\left|\int N(\rho,\sigma)Z_r^{n+1}\,\dx\right|
 +\left|\int N_1(\rho,\sigma)Z_r^{n+1}\,\dx\right|.
\end{equation}
Proposition \ref{pro5.5} gives
\begin{equation}\label{eq:interaction-absorption-second}
 \left|\int N(\rho,\sigma)Z_r^{n+1}\,\dx\right|
 +\left|\int N_1(\rho,\sigma)Z_r^{n+1}\,\dx\right|
 \leq\omega_5(R)Q^\kappa+C\|f\|_{\dot H^{-1}}.
\end{equation}
By \eqref{introeq2},
\(Q\leq\omega_{\rm mod}(\delta)\), and
\(Q\asymp R^{-(n-2)}\).  Since
\(\omega_{\iota_\alpha}(t)\to0\) as \(t\downarrow0\) and
\(\omega_5(T)\to0\) as \(T\to\infty\), decrease \(\delta\) so that
\[
 \omega_{\iota_\alpha}(Q)+\omega_5(R)\leq\frac{c_0}{2}.
\]
Substituting \eqref{eq:interaction-absorption-second} into
\eqref{eq:interaction-absorption-first} gives
$
 \frac{c_0}{2}Q^\kappa\leq C\|f\|_{\dot H^{-1}}$,
which proves the lemma.
     \end{proof}

     We are now ready to prove the main result.

\begin{proof}[Proof of Theorem \ref{mainthm1}]
The case $\nu=1$ is proved in Remark
\ref{rmk:one-bubble-case}: since $Q=h=0$, Proposition
\ref{pro4.2.1} and the nonlinear remainder estimate give the desired
linear bound after absorption.  Hence assume $\nu\geq2$.  By
Proposition \ref{prop:error-estimate}, we have
\[
\left\|\nabla \rho\right\|_{L^2} \lesssim \mathfrak{H}_{\nu,n,\alpha}(R)+\|f\|_{\dot{H}^{-1}},
\]
where $\mathfrak{H}_{\nu,n,\alpha}(R)$ is defined in
\eqref{eq:H-alpha}.  Put
\[
 t:=\|f\|_{\dot H^{-1}}=\Gamma(u).
\]
Lemma \ref{lem:residual-smallness} gives $0\leq t\leq1$, and Lemma
\ref{le2.2}, together with $R^{-(n-2)}\asymp Q$, gives
\[
\|\nabla \rho\|_{L^2}  \lesssim  \begin{cases}
t\left(1+\left|\log t\right|\right)^{\frac12},
&\nu\geq2,\ (n,\alpha)\in\mathcal C_{\log},\\
t^{\theta_\alpha},
&\nu\geq2,\ n\geq7,\ \dfrac{n+2}{2}<\alpha<n,\\
t,&\text{for the remaining cases}.
\end{cases}
\]
This is \eqref{introdueq1.1}.  Finally, the single-bubble equation gives
\[
\begin{split}
\int_{\mathbb{R}^n}
\left(|x|^{-\alpha} * U_i^{p_{\alpha}}\right)
U_i^{p_{\alpha}-1} U_j\,\dx=D_{n,\alpha}^{-1}
\int_{\mathbb R^n}U_i^{2^*-1}U_j\,\dx
=D_{n,\alpha}^{-1}\langle U_i,U_j\rangle_{\dot H^1}
\leq Cq_{ij}\leq CQ.
\end{split}
\]
The last inequality follows from the standard two-bubble overlap
estimate in Lemma \ref{A.3}.  Lemma
\ref{le2.2} gives
$Q\leq Ct^{(n-2)/a_\alpha}$, and hence proves \eqref{eq3}.
This proves Theorem \ref{mainthm1}.
\end{proof}

\medskip

\begin{proof}[Proof of Corollary \ref{mainthm2}]
The proof is identical to that of Corollary 3.4 in \cite{Figalli-Glaudo2020}.
\end{proof}

\medskip

 It remains to establish Proposition \ref{prop:error-estimate} and Proposition \ref{pro5.5}. These depend crucially on the estimate of $\rho$, which will be done in Section \ref{Sec5}.

\section{Estimates for the Inhomogeneous Term \texorpdfstring{$h$}{h}}\label{Sec4}

In \eqref{2.1}, the linearized part is $\Delta\rho+N_1(\rho,\sigma)$ and the inhomogeneous part is $h+N(\rho,\sigma)+f$.  The term $h$ encodes the leading bubble interactions, whereas $N(\rho,\sigma)$ is of higher order in $\rho$.  We begin with two estimates for $h$.

Let $I=\{1,2, \ldots, \nu\}$. Recall that $x_i=\lambda_i\left(x-z_i\right)$,
\[
R_{ij}:=\max\left\{\sqrt{\lambda_i / \lambda_j},
\sqrt{\lambda_j / \lambda_i},
\sqrt{\lambda_i \lambda_j}\,|z_i-z_j|\right\},\qquad i\ne j,
\]
and
\[
R:=\frac{1}{2}\min_{i \neq j \in I}\left\{R_{ij}\right\}.
\]

\paragraph{\textbf{Uniformity convention.}}
Let $w:(1,\infty)\to(0,\infty)$ and let $\mathcal E$ be a
nonnegative, scale-invariant functional of a $\nu$-bubble configuration.
The following two assertions are equivalent:
\begin{align}
 &\exists\,C,R_0>0\quad
   \mathcal E(\mathcal U)\leq Cw(R(\mathcal U))
   \quad\hbox{whenever }R(\mathcal U)\geq R_0;       \label{eq:uniform-principle-a}\\
 &\text{for every sequence }\mathcal U_{(k)}\text{ with }R_{(k)}\to\infty,
   \text{ some subsequence satisfies }
   \sup_k\frac{\mathcal E(\mathcal U_{(k)})}{w(R_{(k)})}<\infty .
                                                               \label{eq:uniform-principle-b}
\end{align}
Indeed, \eqref{eq:uniform-principle-a} immediately implies
\eqref{eq:uniform-principle-b}.  Conversely, if
\eqref{eq:uniform-principle-a} is false, then for every integer $k\geq1$
there is a configuration $\mathcal U_{(k)}$ such that
\[
 R_{(k)}\geq k,
 \qquad
 \mathcal E(\mathcal U_{(k)})>k\,w(R_{(k)}).
\]
No subsequence of this sequence can satisfy
\eqref{eq:uniform-principle-b}, which is a contradiction.  Therefore, in
the configuration arguments below we may pass to the subsequences of
Section~\ref{Sec2}.  Constants such as $C_{\mathscr U}$ may depend on the
extracted sequence, but the constant obtained from
\eqref{eq:uniform-principle-a} depends only on $(n,\nu,\alpha)$.  We apply
this principle with $w(R)=R^{-s}$ or
$w(R)=R^{-s}(\log R)^{1/2}$, always taking $R_0>e$ in the latter case.

\begin{lemma}\label{le7.1}
Assume that $\nu\geq2$, $n\geq3$, $0<\alpha\leq4$, and either
$3\leq n\leq5$ or $n\geq6$ with
$\alpha\ne(n+2)/2$. Then there exists a constant
$C=C(n,\nu,\alpha)>0$ such that
\begin{equation}\label{eq:510}
\begin{split}
\|h\|_{L^{\frac{2n}{n+2}}(\mathbb{R}^{n})}
\leq C R^{-\min\{\alpha,n-2\}}
\end{split}
\end{equation}
\end{lemma}
\begin{proof}
By \eqref{eq:uniform-principle-a}--\eqref{eq:uniform-principle-b}, it is
enough to consider an arbitrary sequence of configurations satisfying
$R_{(k)}\to\infty$ and to prove the asserted estimate, with a constant
independent of $k$, along a subsequence.  After relabeling and passing to
a subsequence as in Section~\ref{Sec2}, relations
\eqref{confi2}--\eqref{confi8} hold.  We suppress the index $(k)$ below.
Thus, whenever the proof uses $i\succ j$, $i\prec j$, or the
incomparability of $i$ and $j$, these are the three exhaustive
subsequential alternatives defined in Section~\ref{Sec2}.

Put
\[
    \beta:=p_\alpha-1=\frac{n+2-\alpha}{n-2}\geq1,\qquad
    \gamma:=2^*-p_\alpha=\frac{\alpha}{n-2},\qquad
    r:=\frac{2n}{n+2}.
\]
Using the single bubble identity
\[
    D_{n,\alpha}\int_{\mathbb R^n}\frac{U_i^{p_\alpha}(y)}{|x-y|^\alpha}\dy=U_i^\gamma,
\]
we write the expression inside the norm as
\begin{align*}
    h=\sum_{i=1}^{\nu}U_i^\gamma(\sigma^\beta-U_i^\beta)
    +D_{n,\alpha}\int_{\mathbb R^n}
    \frac{\sigma^{p_\alpha}(y)-\sum_iU_i^{p_\alpha}(y)}{|x-y|^\alpha}\dy\,\sigma^\beta
    :=A+B .
\end{align*}
Choose a measurable partition \((\Omega_\ell)_{\ell=1}^{\nu}\) of
\(\mathbb R^n\) such that
\[
 \Omega_\ell\subset\{x:U_\ell(x)=\max_iU_i(x)\};
\]
ties may be resolved by the smallest index.  Since \(\beta\geq1\) and
\(\sigma\leq\nu U_\ell\) on \(\Omega_\ell\), the \(i=\ell\) term is
estimated by the mean value theorem, and the terms \(i\ne\ell\) are
dominated by \(U_\ell\):
\[
    0\leq A\leq C\left(U_\ell^{2^*-2}\sum_{j\ne\ell}U_j
    +\sum_{j\ne\ell}U_j^\gamma U_\ell^\beta\right)
    \qquad\hbox{on }\Omega_\ell .
\]
Lemma \ref{A.4} gives, for $j\ne\ell$,
\begin{equation}\label{eq:lemma31-A-yamabe}
    \|U_\ell^{2^*-2}U_j\|_{L^r(\Omega_\ell)}
    \leq C R_{\ell j}^{-c_n}
    (1+\log R_{\ell j})^{\frac1r\mathbf1_{\{n=6\}}},
\end{equation}
and
\begin{equation}\label{eq:lemma31-A-singular}
    \|U_j^\gamma U_\ell^\beta\|_{L^r(\Omega_\ell)}
    \leq C\begin{cases}
        R_{\ell j}^{-\alpha},&\alpha<\frac{n+2}{2},\\[1mm]
        R_{\ell j}^{-\frac{n+2}{2}}(\log R_{\ell j})^{\frac{n+2}{2n}},&\alpha=\frac{n+2}{2},\\[1mm]
        R_{\ell j}^{-\frac{n+2}{2}},&\alpha>\frac{n+2}{2}.
    \end{cases}
\end{equation}
In the full parameter range of the lemma, the powers in
\eqref{eq:lemma31-A-yamabe}--\eqref{eq:lemma31-A-singular} are at
least \(a_\alpha=\min\{\alpha,n-2\}\), and every logarithmic factor is
absorbed by a strictly larger power; hence, since
\(R_{\ell j}\geq2R\) and the number of pairs is finite,
\begin{equation}\label{eq:lemma31-A-final}
 \|A\|_{L^r}\leq C R^{-a_\alpha}.
\end{equation}

We next estimate \(B\).  Since \(p_\alpha=\beta+1\geq2\), the mean value
theorem gives, on \(\Omega_\ell\),
\[
 0\leq \sigma^{p_\alpha}-\sum_iU_i^{p_\alpha}
 \leq C U_\ell^\beta\sum_{j\ne\ell}U_j.
\]
Moreover, \(\sigma^\beta\leq C\sum_kU_k^\beta\).  Splitting both
variables according to the partition \((\Omega_\ell)\), and applying
the bilinear Hardy--Littlewood--Sobolev inequality with equal exponent
\(s=2n/(2n-\alpha)\), we obtain, for every
\(\varphi\in L^{r'}(\mathbb R^n)\) with
\(\|\varphi\|_{L^{r'}}=1\),
\begin{align*}
\left|\int_{\mathbb R^n}
\left(D_{n,\alpha}\int_{\mathbb R^n}
\frac{\sigma^{p_\alpha}(y)-\sum_iU_i^{p_\alpha}(y)}{|x-y|^\alpha}\dy\right)
\sigma^\beta(x)\varphi(x)\dx\right|\leq C\sum_{\ell}\sum_{j\ne\ell}\sum_k
\left(\int_{\Omega_\ell}(U_\ell^\beta U_j)^s\dy\right)^{1/s}
\left(\int_{\Omega_k}(U_k^\beta|\varphi|)^s\dx\right)^{1/s},
\end{align*}
where
$s=\frac{2n}{2n-\alpha}$,
    $\frac1s=\frac1t+\frac1{r'}$, and
    $\beta t=2^*$ .
H\"older inequality gives
\[
    \left(\int_{\Omega_k}(U_k^\beta|\varphi|)^s\dx\right)^{1/s}
    \leq \|U_k^\beta\|_{L^t}\|\varphi\|_{L^{r'}}\leq C.
\]
For the first factor, Lemma \ref{A.4} gives
\[
\left(\int_{\Omega_\ell}(U_\ell^\beta U_j)^s\dy\right)^{1/s}
\leq C
\begin{cases}
R_{\ell j}^{-(n-2)},&0<\alpha<4,\\[1mm]
R_{\ell j}^{-(n-2)}(\log R_{\ell j})^{1/s},&\alpha=4.
\end{cases}
\]
If \(\alpha<4\), then
\(R^{-(n-2)}=O(R^{-a_\alpha})\).  If \(\alpha=4\) and \(n\geq7\), then
\[
 R^{-(n-2)}(1+\log R)^{1/s}=O(R^{-4})
 =O(R^{-a_\alpha}).
\]
Consequently, in every case covered by the lemma except
\((n,\alpha)=(5,4)\),
\begin{equation}\label{eq:lemma31-B-final}
 \|B\|_{L^r}\leq C R^{-a_\alpha}.
\end{equation}

It remains only to estimate \(B\) at the endpoint
\((n,\alpha)=(5,4)\).  Put
\[
 \mathcal I_4f(x):=D_{5,4}\int_{\mathbb R^5}\frac{f(y)}{|x-y|^4}\dy,
 \qquad r=\frac{10}{7}.
\]
Since \(p_\alpha=2\), \(\beta=1\), and
\(\mathcal I_4(U_i^2)=U_i^{4/3}\), the second term in the decomposition
of \(h\) is exactly
\begin{equation}\label{eq:54-B-decomposition}
 B=2\sum_{i<j}\sum_k\mathcal I_4(U_iU_j)U_k.
\end{equation}
For \(i\ne j\), define
\[
 \Theta_{ij}:=\{U_i\geq U_j\},\qquad
 K_{ij}:=\mathcal I_4(U_iU_j\mathbf1_{\Theta_{ij}}),
 \qquad D:=R_{ij}.
\]
We shall also use the following pair-product consequence of Lemma
\ref{A.4}:
\begin{align}
 \|U_i^{4/3}U_j\|_{L^{10/7}}^{10/7}
 =\int_{\Theta_{ij}}U_i^{40/21}U_j^{10/7}\dx
   +\int_{\Theta_{ji}}U_i^{40/21}U_j^{10/7}\dx\leq C\left(D^{-30/7}+D^{-5}\right)
 \leq C D^{-30/7},                                      \label{eq:54-pair-product}
\end{align}
and therefore
\begin{equation*}
 \|U_i^{4/3}U_j\|_{L^{10/7}}\leq CD^{-3}.
\end{equation*}
We shall prove
\begin{equation}\label{eq:54-conv-one}
 \|K_{ij}U_i\|_{L^r}+\|K_{ij}U_j\|_{L^r}\leq CD^{-3}.
\end{equation}
The elementary convolution estimate used in each configuration is
\begin{equation}\label{eq:54-elementary-convolution}
 \int_{|Y|\leq L}|X-Y|^{-4}\langle Y\rangle^{-3}\,dY
 \leq C\begin{cases}
 \langle X\rangle^{-2},&|X|\leq2L,\\
 L^2|X|^{-4},&|X|>2L.
 \end{cases}
\end{equation}
For $|X|\leq2L$, this follows by splitting into
$|Y|\leq|X|/2$, $|X-Y|\leq|X|/2$, and the remaining region.  If
$|X|>2L$, then $|X-Y|\geq|X|/2$ and the second line follows from
$\int_{|Y|\leq L}\langle Y\rangle^{-3}dY\leq CL^2$.

\noindent{\bf Case 1: $i\succ j$.}
Set
\[
 \vartheta=\left(\frac{\lambda_i}{\lambda_j}\right)^{1/2},\qquad
 \mu=\vartheta^{-2},\qquad b_{ji}=\lambda_j(z_i-z_j).
\]
The ordered relations imply
\[
 |b_{ji}|\leq C_{\mathscr U},\qquad
 \vartheta\leq D\leq C_{\mathscr U}\vartheta,
 \qquad y_j=\mu y_i+b_{ji}.
\]
For $y\in\Theta_{ij}$,
$
 1+|\mu y_i+b_{ji}|^2\geq\mu(1+|y_i|^2)$,
and hence $|y_i|\leq CD$.  On this region
\[
 U_j(y)\leq C\lambda_i^{3/2}D^{-3},\qquad
 U_i(y)\leq C\lambda_i^{3/2}\langle y_i\rangle^{-3}.
\]
After the change $Y=y_i$, formula
\eqref{eq:54-elementary-convolution} gives
\begin{equation*}
 K_{ij}(x)\leq C\lambda_i^2
 \begin{cases}
 D^{-3}\langle x_i\rangle^{-2},&|x_i|\leq CD,\\
 D^{-1}|x_i|^{-4},&|x_i|>CD.
 \end{cases}
\end{equation*}
Multiplying by $U_i$ and using polar coordinates yields
\begin{align*}
 \|K_{ij}U_i\|_{L^r}
 \leq C D^{-3}
 \left(\int_0^{CD}(1+s)^{4-5r}ds\right)^{1/r}+C D^{-1}
 \left(\int_{CD}^{\infty}s^{4-7r}ds\right)^{1/r}
 \leq C(D^{-3}+D^{-9/2}).
\end{align*}
Using instead the bound for $U_j$ on the first region and its exact
$x_i$-coordinate decay after $|x_i|\simeq D^2$, one obtains
\begin{align*}
 \|K_{ij}U_j\|_{L^r}
 \leq C D^{-6}
 \left(\int_0^{CD}(1+s)^{4-2r}ds\right)^{1/r}
 +C D^{-4}
 \left(\int_{CD}^{\infty}s^{4-4r}ds\right)^{1/r}\leq CD^{-9/2}.
\end{align*}
Thus \eqref{eq:54-conv-one} holds in Case 1.

\noindent{\bf Case 2: $i\prec j$.}
Now set
\[
 \vartheta=\left(\frac{\lambda_j}{\lambda_i}\right)^{1/2},\qquad
 \mu=\vartheta^{-2},\qquad b_{ij}=\lambda_i(z_j-z_i).
\]
Then
\[
 |b_{ij}|\leq C_{\mathscr U},\qquad
 \vartheta\leq D\leq C_{\mathscr U}\vartheta,
 \qquad y_i=\mu y_j+b_{ij}.
\]
The inequality $U_i\geq U_j$ implies
$
 \mu(1+|y_j|^2)\geq1+|\mu y_j+b_{ij}|^2$,
and therefore
\begin{equation}\label{eq:54-case2-theta-region}
 \Theta_{ij}\subset\{|y_j|\geq cD\}.
\end{equation}
Let \(X=x_j\) and \(Y=y_j\).  Since
\(\mu\asymp D^{-2}\) and \(|b_{ij}|\leq C_{\mathscr U}\), there are
constants \(0<c_0<1<C_0\), depending only on \(\mathscr U\), such that
the definition of \(K_{ij}\),
\eqref{eq:54-case2-theta-region}, and
\[
 U_i(y)U_j(y)
 =C\lambda_j^3\mu^{3/2}
 (1+|\mu Y+b_{ij}|^2)^{-3/2}(1+|Y|^2)^{-3/2}
\]
give the following two pointwise bounds.  On
\(c_0D\leq|Y|\leq C_0D^2\),
\[
 \mu^{3/2}(1+|\mu Y+b_{ij}|^2)^{-3/2}
 (1+|Y|^2)^{-3/2}
 \leq CD^{-3}|Y|^{-3}.
\]
Choose \(C_0\) so large that
\(\mu|Y|\geq2(1+|b_{ij}|)\) whenever
\(|Y|\geq(C_0/2)D^2\).  On this latter region,
\[
 1+|\mu Y+b_{ij}|^2\geq c\mu^2|Y|^2,
\]
and hence
\[
 \mu^{3/2}(1+|\mu Y+b_{ij}|^2)^{-3/2}
 (1+|Y|^2)^{-3/2}
 \leq C\mu^{-3/2}|Y|^{-6}
 \leq CD^3|Y|^{-6}.
\]
These inequalities and
\eqref{eq:54-case2-theta-region}
give
\[
 \lambda_j^{-2}K_{ij}(x)\leq C\bigl(\mathcal K_1(X)+\mathcal K_2(X)\bigr),
\]
where
\[
\begin{split}
 \mathcal K_1(X)
 &:=D^{-3}\int_{c_0D\leq|Y|\leq C_0D^2}
 |X-Y|^{-4}|Y|^{-3}\,dY,\\
 \mathcal K_2(X)
 &:=D^3\int_{|Y|\geq (C_0/2)D^2}
 |X-Y|^{-4}|Y|^{-6}\,dY.
\end{split}
\]
For \(|X|\leq CD\), the changes of variables \(X=D\xi,Y=D\zeta\)
in \(\mathcal K_1\), and \(X=D^2\xi,Y=D^2\zeta\) in
\(\mathcal K_2\), yield
\[
 \mathcal K_1(X)\leq CD^{-5},\qquad
 \mathcal K_2(X)\leq CD^{-7}\leq CD^{-5}.
\]
For \(CD<|X|\leq CD^2\), decompose the domain of
\(\mathcal K_1\) into
\[
 A_1:=\{|Y|\leq|X|/2\},\qquad
 A_2:=\{|X-Y|\leq|X|/2\},\qquad
 A_3:=\{|Y|>|X|/2,\ |X-Y|>|X|/2\}.
\]
Direct radial integration gives
\[
\begin{split}
 D^{-3}\int_{A_1}|X-Y|^{-4}|Y|^{-3}\,dY
 &\leq CD^{-3}|X|^{-4}\int_0^{|X|/2}r\,dr
 \leq CD^{-3}|X|^{-2},\\
 D^{-3}\int_{A_2}|X-Y|^{-4}|Y|^{-3}\,dY
 &\leq CD^{-3}|X|^{-3}\int_0^{|X|/2}dr
 \leq CD^{-3}|X|^{-2},\\
 D^{-3}\int_{A_3}|X-Y|^{-4}|Y|^{-3}\,dY
 &\leq CD^{-3}\left(
 |X|^{-7}|B_{2|X|}|+
 \int_{2|X|}^{\infty}r^{-3}\,dr\right)\leq CD^{-3}|X|^{-2}.
\end{split}
\]
The \(D^2\)-rescaling gives
\[
 \mathcal K_2(X)\leq CD^{-7}
 \leq CD^{-3}|X|^{-2}.
\]
Finally, if \(|X|>CD^2\), with the constants chosen so that
\(|X-Y|\geq|X|/2\) on the support of \(\mathcal K_1\), then
\[
 \mathcal K_1(X)
 \leq CD^{-3}|X|^{-4}\int_D^{CD^2}r\,dr
 \leq CD|X|^{-4}.
\]
For \(\mathcal K_2\), split into
\[
 \{|Y|\leq|X|/2\},\quad
 \{|X-Y|\leq|X|/2\},\quad
 \{|Y|>|X|/2,\ |X-Y|>|X|/2\}.
\]
The corresponding bounds are
\[
 CD^3(D^2)^{-1}|X|^{-4},\qquad
 CD^3|X|^{-5},\qquad
 CD^3|X|^{-5},
\]
each of which is at most \(CD|X|^{-4}\).  We have therefore proved
\begin{equation*}
 K_{ij}(x)\leq C\lambda_j^2
 \begin{cases}
 D^{-5},&|x_j|\leq CD,\\
 D^{-3}|x_j|^{-2},&CD<|x_j|\leq CD^2,\\
 D|x_j|^{-4},&|x_j|>CD^2.
 \end{cases}
\end{equation*}
In these coordinates
\[
 U_j\leq C\lambda_j^{3/2}\langle x_j\rangle^{-3},\qquad
 U_i\leq C\lambda_j^{3/2}
 \begin{cases}
 D^{-3},&|x_j|\leq CD^2,\\
 D^3|x_j|^{-3},&|x_j|>CD^2.
 \end{cases}
\]
Consequently,
\begin{align*}
 \|K_{ij}U_j\|_{L^r}
 &\leq C\bigg[D^{-5}
 \left(\int_0^{CD}(1+s)^{4-3r}ds\right)^{1/r}
 +D^{-3}\left(\int_D^{CD^2}s^{4-5r}ds\right)^{1/r}+D
 \left(\int_{D^2}^{\infty}s^{4-7r}ds\right)^{1/r}\bigg]
 \leq CD^{-9/2},\\
 \|K_{ij}U_i\|_{L^r}
 &\leq C\bigg[D^{-8}
 \left(\int_0^{CD}s^4ds\right)^{1/r}
 +D^{-6}\left(\int_D^{CD^2}s^{4-2r}ds\right)^{1/r}+D^4
 \left(\int_{D^2}^{\infty}s^{4-7r}ds\right)^{1/r}\bigg]
 \leq CD^{-3}.
\end{align*}
This proves \eqref{eq:54-conv-one} in Case 2.

\noindent{\bf Case 3: $i$ and $j$ are incomparable.}
Put
\[
 d=|z_i-z_j|,\qquad D_0=\sqrt{\lambda_i\lambda_j}\,d,
 \qquad\theta=\sqrt{\frac{\lambda_i}{\lambda_j}},
 \qquad
 L_i=\frac{D_0\theta}{2},\qquad L_j=\frac{D_0}{2\theta}.
\]
Then
\begin{equation}\label{eq:54-case3-incomp-scales}
 D_0\asymp D,\qquad
 \max\{\theta,\theta^{-1}\}=o(D_0),\qquad L_i,L_j\to\infty.
\end{equation}
Let
\[
 \mathcal B_i=B_{d/2}(z_i),\qquad
 \mathcal B_j=B_{d/2}(z_j),\qquad
 \mathcal B_\infty=(\mathcal B_i\cup\mathcal B_j)^c,
\]
and set $\mathcal E_m=\Theta_{ij}\cap\mathcal B_m$ for
$m=i,j,\infty$.  This is a disjoint partition of $\Theta_{ij}$.  Define
\[
 T_m=\mathcal I_4(U_iU_j\mathbf1_{\mathcal E_m}),\qquad
 m=i,j,\infty,
\]
so that $K_{ij}=T_i+T_j+T_\infty$.

For $m\in\{i,j\}$, let $m'$ denote the other index. Points in
$\mathcal E_m$ satisfy $|y-z_{m'}|\geq d/2$ and $|y_m|<L_m$.
Hence
\[
 U_iU_j\mathbf1_{\mathcal E_m}
 \leq C\lambda_m^3D_0^{-3}\langle y_m\rangle^{-3}
 \mathbf1_{\{|y_m|<L_m\}}.
\]
The change to $y_m$ coordinates and
\eqref{eq:54-elementary-convolution} give
\begin{equation*}
 T_m(x)\leq C\lambda_m^2D_0^{-3}\langle x_m\rangle^{-2}.
\end{equation*}
Thus
\begin{align*}
 \|T_mU_m\|_{L^r}^r
 &\leq CD_0^{-3r}\lambda_m^{\frac72r-5}
 \int_{\mathbb R^5}\langle X\rangle^{-5r}dX
 \leq CD_0^{-3r},
\end{align*}
and therefore
\begin{equation}\label{eq:54-case3-local-same}
 \|T_iU_i\|_{L^r}+\|T_jU_j\|_{L^r}\leq CD_0^{-3}.
\end{equation}

We next estimate a crossed local term, say $T_iU_j$.  On
$\mathcal B_i$, $U_j\leq C\lambda_i^{3/2}D_0^{-3}$; hence
\[
 \|T_iU_j\|_{L^r(\mathcal B_i)}
 \leq CD_0^{-6}
 \left(\int_0^{L_i}\rho^4\langle\rho\rangle^{-2r}d\rho\right)^{1/r}
 \leq CD_0^{-6}L_i^{3/2}.
\]
On $\mathcal B_j$, one has $|x_i|\geq L_i$ and
$dx=\lambda_j^{-5}dx_j$, so
\[
 \|T_iU_j\|_{L^r(\mathcal B_j)}
 \leq CD_0^{-3}L_i^{-2}\theta^4
 \left(\int_0^{L_j}\rho^4\langle\rho\rangle^{-3r}d\rho\right)^{1/r}
 \leq CD_0^{-3}L_i^{-2}\theta^4L_j^{1/2}.
\]
Finally, on $\mathcal B_\infty$, both $|x_i|$ and
$|x_i-c_{ij}|$ are at least $L_i$, and rescaling by $L_i$ gives
\begin{align*}
 \|T_iU_j\|_{L^r(\mathcal B_\infty)}
 \leq CD_0^{-3}\theta^3
 \left(\int_{\substack{|X|\geq L_i\\|X-c_{ij}|\geq L_i}}
 |X|^{-2r}|X-c_{ij}|^{-3r}dX\right)^{1/r}\leq CD_0^{-3}\theta^3L_i^{-3/2}.
\end{align*}
Substituting $L_i=D_0\theta/2$ and $L_j=D_0/(2\theta)$ into these
three bounds yields
$$
 \|T_iU_j\|_{L^r}\leq CD_0^{-9/2}\theta^{3/2}\leq CD_0^{-3}.
$$
Interchanging $i,j$ gives
\begin{equation}\label{eq:54-case3-local-cross-sym}
 \|T_jU_i\|_{L^r}\leq CD_0^{-9/2}\theta^{-3/2}\leq CD_0^{-3}.
\end{equation}

It remains to control $T_\infty$.  On $\mathcal B_\infty$ the two
distances to the centers are comparable; in $i$- and $j$-coordinates,
respectively,
\[
 U_iU_j\mathbf1_{\mathcal B_\infty}
 \leq C\lambda_i^3\theta^3|y_i|^{-6}\mathbf1_{\{|y_i|\geq L_i\}},
\qquad
 U_iU_j\mathbf1_{\mathcal B_\infty}
 \leq C\lambda_j^3\theta^{-3}|y_j|^{-6}\mathbf1_{\{|y_j|\geq L_j\}}.
\]
For $L\geq1$ set
\[
 G_L(X)=\int_{|Y|\geq L}|X-Y|^{-4}|Y|^{-6}dY.
\]
If \(|X|\leq2L\), the substitution \(X=L\xi\), \(Y=L\eta\)
gives \(|\xi|\leq2\).  On \(1\leq|\eta|\leq4\), the factor
\(|\eta|^{-6}\) is bounded and
\(|\xi-\eta|^{-4}\in L^1_{\mathrm{loc}}(\mathbb R^5)\); on
\(|\eta|>4\), one has
\[
 |\xi-\eta|^{-4}|\eta|^{-6}\leq16|\eta|^{-10}.
\]
Consequently, the following bound is uniform for \(|\xi|\leq2\):
\[
 G_L(X)
 =L^{-5}\int_{|\eta|\geq1}
 |\xi-\eta|^{-4}|\eta|^{-6}\,d\eta
 \leq CL^{-5}.
\]
If \(|X|>2L\), decompose \(\{|Y|\geq L\}\) into
\[
 A_1:=\{|Y|\leq|X|/2\},\quad
 A_2:=\{|X-Y|\leq|X|/2\},\quad
 A_3:=\{|Y|>|X|/2,\ |X-Y|>|X|/2\}.
\]
Then
\[
\begin{split}
 \int_{A_1}|X-Y|^{-4}|Y|^{-6}\,dY
 &\leq C|X|^{-4}\int_L^{|X|/2}r^{-2}\,dr
 \leq CL^{-1}|X|^{-4},\\
 \int_{A_2}|X-Y|^{-4}|Y|^{-6}\,dY
 &\leq C|X|^{-6}\int_0^{|X|/2}dr
 \leq C|X|^{-5},\\
 \int_{A_3}|X-Y|^{-4}|Y|^{-6}\,dY
 &\leq C|X|^{-4}\int_{|X|/2}^{\infty}r^{-2}\,dr
 \leq C|X|^{-5}.
\end{split}
\]
Since \(|X|^{-5}\leq L^{-1}|X|^{-4}\), these inequalities prove
\begin{equation*}
 G_L(X)\leq C\begin{cases}
 L^{-5},&|X|\leq2L,\\
 L^{-1}|X|^{-4},&|X|>2L.
 \end{cases}
\end{equation*}
Consequently,
\begin{align*}
 \left(\int_{\mathbb R^5}G_L(X)^r\langle X\rangle^{-3r}dX\right)^{1/r}
 \leq C\bigg[L^{-5}
 \left(\int_0^{2L}\rho^4\langle\rho\rangle^{-3r}d\rho\right)^{1/r}+L^{-1}
 \left(\int_{2L}^{\infty}\rho^{4-4r}\langle\rho\rangle^{-3r}d\rho\right)^{1/r}
 \bigg]\leq CL^{-9/2}.
\end{align*}
It follows that
\[
 \|T_\infty U_i\|_{L^r}\leq CD_0^{-9/2}\theta^{-3/2},
 \qquad
 \|T_\infty U_j\|_{L^r}\leq CD_0^{-9/2}\theta^{3/2}.
\]
Together with \eqref{eq:54-case3-incomp-scales} and
\eqref{eq:54-case3-local-same}--\eqref{eq:54-case3-local-cross-sym},
these inequalities prove \eqref{eq:54-conv-one} in Case 3.

It remains to consider $k\notin\{i,j\}$.  H\"older inequality with
respect to the Riesz kernel and the single-bubble identity give
\begin{equation*}
 K_{ij}\leq\mathcal I_4(U_i^2)^{1/2}
 \mathcal I_4(U_j^2)^{1/2}=U_i^{2/3}U_j^{2/3}.
\end{equation*}
Partition according to the largest of $U_i,U_j,U_k$.  On the
$U_k$-region,
\[
 U_i^{2/3}U_j^{2/3}U_k
 \leq\frac12(U_i^{4/3}U_k+U_j^{4/3}U_k),
\]
while on the $U_i$-region,
\[
 U_i^{2/3}U_j^{2/3}U_k
 \leq C U_i^{2/3}(U_j^{5/3}+U_k^{5/3}),
\]
with the symmetric inequality on the $U_j$-region.  Lemma \ref{A.4}
and \eqref{eq:54-pair-product} therefore yield
\begin{equation}\label{eq:54-conv-three}
 \|K_{ij}U_k\|_{L^r}\leq CR^{-3},\qquad k\notin\{i,j\}.
\end{equation}
Since \(\mathcal I_4(U_iU_j)=K_{ij}+K_{ji}\), estimate
\eqref{eq:54-conv-one} controls the terms with
\(k\in\{i,j\}\), and \eqref{eq:54-conv-three} controls those with
\(k\notin\{i,j\}\).  Summing the exact decomposition
\eqref{eq:54-B-decomposition}, using the finiteness of the index set
and \(R_{ij}\geq2R\), gives
\(\|B\|_{L^r}\leq CR^{-3}\).  At this endpoint
\(a_\alpha=\min\{4,3\}=3\), so this estimate and
\eqref{eq:lemma31-A-final} prove \eqref{eq:510}.  In every other case,
\eqref{eq:lemma31-A-final} and \eqref{eq:lemma31-B-final} prove the same
conclusion.
\end{proof}

\begin{lemma}\label{le3.1}\label{le:strong-h-estimate}
Assume
\[
 \nu\geq2,\qquad
 \bigl[(n,\alpha)=(6,4)
 \quad\hbox{or}\quad
 4<\alpha<n\bigr].
\]
There exist $R_0=R_0(n,\nu,\alpha)>e$ and
$C=C(n,\nu,\alpha)>0$ with the following properties for every
$\nu$-bubble configuration satisfying $R\geq R_0$. If
$4<\alpha<n$ and $(n,\alpha)\notin\mathcal C_{\log}$, then
\begin{equation*}
\|h\|_{L^{\frac{2n}{n+2}}(\mathbb{R}^{n})}
\leq C R^{-b_\alpha}.
\end{equation*}
If $(n,\alpha)\in\mathcal C_{\log}$, then
\begin{equation*}
\|h\|_{\dot H^{-1}(\mathbb{R}^{n})}
\leq C R^{-b_\alpha}(\log R)^{\frac12}.
\end{equation*}
\end{lemma}

\begin{proof}
By \eqref{eq:uniform-principle-a}--\eqref{eq:uniform-principle-b}, it is
enough to consider an arbitrary sequence of configurations with
$R_{(k)}\to\infty$ and to prove the stated estimate, with a constant
independent of $k$, on a subsequence. Apply Section~\ref{Sec2}, relabel the
bubbles, and pass to a subsequence on which \eqref{confi2}--\eqref{confi8}
hold. In the rest of the proof the index $(k)$ is suppressed. A constant
may depend on the fixed number $C_{\mathscr U}$ from \eqref{confi4}, but not
on $k$. The uniformity convention then produces constants depending only
on $(n,\nu,\alpha)$.

Put
\[
 \beta=p_\alpha-1=\frac{n+2-\alpha}{n-2}\in(0,1],\qquad
 \gamma=2^*-p_\alpha=\frac{\alpha}{n-2},
 \qquad r=\frac{2n}{n+2},
\]
and let $\Omega_\ell=\{U_\ell=\max_iU_i\}$.  Using the
single-bubble identity, decompose
\begin{align*}
 h&=\sum_iU_i^\gamma(\sigma^\beta-U_i^\beta)
 +D_{n,\alpha}\int_{\mathbb R^n}
 \frac{\sigma^{p_\alpha}(y)-\sum_iU_i^{p_\alpha}(y)}{|x-y|^\alpha}
 \dy\,\sigma^\beta:=A+B.
\end{align*}
On $\Omega_\ell$, the mean value theorem for the $i=\ell$ term and
$U_i\leq U_\ell$ for $i\ne\ell$ give
\begin{equation}\label{eq:strong-A-point}
 0\leq A\leq C\left(
 U_\ell^{2^*-2}\sum_{j\ne\ell}U_j
 +\sum_{j\ne\ell}U_j^\gamma U_\ell^\beta\right).
\end{equation}
Lemma \ref{A.4} yields
\begin{equation}\label{eq:strong-A-yamabe-detail}
 \|U_\ell^{2^*-2}U_j\|_{L^r(\Omega_\ell)}
 \leq CR_{\ell j}^{-c_n}
 (1+\log R_{\ell j})^{\frac1r\mathbf1_{\{n=6\}}},
\end{equation}
and
\begin{equation*}
 \|U_j^\gamma U_\ell^\beta\|_{L^r(\Omega_\ell)}
 \leq C\begin{cases}
 R_{\ell j}^{-\alpha},&\alpha<(n+2)/2,\\
 R_{\ell j}^{-(n+2)/2}(1+\log R_{\ell j})^{(n+2)/(2n)},
      &\alpha=(n+2)/2,\\
 R_{\ell j}^{-(n+2)/2},&\alpha>(n+2)/2.
 \end{cases}
\end{equation*}
After summing over the finitely many pairs, these estimates give
\begin{equation}\label{eq:strong-A-short}
 \|A\|_{L^r}\leq CR^{-b_\alpha},\qquad
 (n,\alpha)\notin\mathcal C_{\log}.
\end{equation}
For the logarithmic pairs, Lemma \ref{A.6} is applied to the balanced
product in \eqref{eq:strong-A-point}.  Together with Sobolev duality
and \eqref{eq:strong-A-yamabe-detail}, it gives
\begin{equation}\label{eq:strong-A-short-log}
 \|A\|_{\dot H^{-1}}
 \leq CR^{-b_\alpha}(1+\log R)^{1/2},\qquad
 (n,\alpha)\in\mathcal C_{\log}.
\end{equation}

It remains to estimate $B$.  Put
\[
 F=\sigma^{p_\alpha}-\sum_iU_i^{p_\alpha},\qquad
 \mathcal I_\alpha f(x)=D_{n,\alpha}\int_{\mathbb R^n}
 \frac{f(y)}{|x-y|^\alpha}\dy.
\]
On $\Omega_\ell$, the mean value theorem gives
\begin{equation}\label{eq:strong-B-point}
 0\leq F\mathbf1_{\Omega_\ell}
 \leq C\sum_{j\ne\ell}U_\ell^\beta U_j\mathbf1_{\Omega_\ell},
 \qquad \sigma^\beta\leq\sum_kU_k^\beta.
\end{equation}
We first record the direct HLS estimate.  Let
\[
 s=\frac{2n}{2n-\alpha},\qquad
 t=\frac{2n}{n+2-\alpha},\qquad
 \eta_\alpha=\frac{2n-\alpha}{2}.
\]
For $\varphi\in L^{r'}$ with $\|\varphi\|_{L^{r'}}=1$, HLS inequality and
H\"older inequality implies
\begin{align*}
 \left|\int_{\mathbb R^n}B\varphi\dx\right|
 \leq C\sum_{\ell}\sum_{j\ne\ell}\sum_k
 \|U_\ell^\beta U_j\|_{L^s(\Omega_\ell)}
 \|U_k^\beta\varphi\|_{L^s}
 \leq C\sum_{\ell}\sum_{j\ne\ell}
 \|U_\ell^\beta U_j\|_{L^s(\Omega_\ell)}.
\end{align*}
Here $\|U_k^\beta\varphi\|_{L^s}\leq C$ by H\"older inequality.  Lemma
\ref{A.4}, applied on the region $U_\ell\geq U_j$, gives
$\|U_\ell^\beta U_j\|_{L^s(\Omega_\ell)}
 \leq CR_{\ell j}^{-\eta_\alpha}$.
Consequently,
\begin{equation*}
 \|B\|_{L^r}\leq CR^{-\eta_\alpha}.
\end{equation*}
This proves the desired estimate when $4<\alpha\leq n-2$.

The remaining parameter set is
$ \{(n,\alpha):4<\alpha<n,\ \alpha>n-2\}\cup\{(6,4)\}$.
For $\ell\ne j$ define
\[
 K_{\ell j}=\mathcal I_\alpha
 (U_\ell^\beta U_j\mathbf1_{\Omega_\ell}).
\]
By \eqref{eq:strong-B-point}, it is enough first to prove
\begin{equation}\label{eq:strong-B-same-index}
 \|K_{\ell j}U_\ell^\beta\|_{L^r}
 +\|K_{\ell j}U_j^\beta\|_{L^r}
 \leq CR_{\ell j}^{-b_\alpha}.
\end{equation}
Put
\[
 a=n+2-\alpha,\qquad b_0=\frac{n+2}{2},\qquad
 e_\alpha=n-\alpha+b_0=n+a-b_0,\qquad
 \kappa=a-b_0+2,\qquad D=R_{\ell j}.
\]
In the parameter range now under consideration,
\begin{equation}\label{eq:strong-parameter-identities}
 2<a\leq4,\qquad a\leq b_0,\qquad
 b_\alpha=\min\{n-2,b_0\},\qquad
 e_\alpha-\kappa=n-2.
\end{equation}
The equality $a=b_0$ occurs only for $(n,\alpha)=(6,4)$.

We shall repeatedly use the following consequences of direct integration.
For $T\geq2$ and $q>0$, set
\[
 \Phi_q(T):=
 \left(\int_0^T(1+s)^{n-1-qr}\,ds\right)^{1/r}.
\]
Then
\begin{equation}\label{eq:strong-power-integral}
 \Phi_q(T)\leq C
 \begin{cases}
  T^{b_0-q},&q<b_0,\\
  (1+\log T)^{1/r},&q=b_0,\\
  1,&q>b_0.
 \end{cases}
\end{equation}
Moreover, if $2\leq A<B$, then
\begin{equation}\label{eq:strong-annular-integral}
 \left(\int_A^B s^{n-1-qr}\,ds\right)^{1/r}
 \leq C
 \begin{cases}
  B^{b_0-q},&q<b_0,\\
  (\log(B/A))^{1/r},&q=b_0,\\
  A^{b_0-q},&q>b_0,
 \end{cases}
\end{equation}
and, for $q>b_0$,
\begin{equation}\label{eq:strong-exterior-integral}
 \left(\int_T^\infty s^{n-1-qr}\,ds\right)^{1/r}
 =C_qT^{b_0-q}.
\end{equation}
Indeed, $n/r=b_0$, and the three assertions follow by evaluating the
antiderivative of $s^{n-1-qr}$; the case $qr=n$ gives the logarithm.
We repeatedly use the following convolution estimate:
\begin{equation}\label{eq:strong-B-elementary-convolution}
 \int_{|Y|\leq L}|X-Y|^{-\alpha}\langle Y\rangle^{-a}dY
 \leq C\begin{cases}
 \langle X\rangle^{-2},&|X|\leq2L,\\
 L^{\alpha-2}|X|^{-\alpha},&|X|>2L.
 \end{cases}
\end{equation}
For the first line, split into $|Y|\leq|X|/2$,
$|X-Y|\leq|X|/2$, and the complement.  For the second line use
$|X-Y|\geq|X|/2$ and
$\int_{|Y|\leq L}\langle Y\rangle^{-a}dY\leq CL^{\alpha-2}$.
We now keep the three alternatives in the partial order separate.

\noindent{\bf Case 1: $\ell\succ j$.}
Set
\[
 \vartheta=\left(\frac{\lambda_\ell}{\lambda_j}\right)^{1/2},
 \qquad\mu=\vartheta^{-2},\qquad
 b_{j\ell}=\lambda_j(z_\ell-z_j),\qquad
 \Theta_{\ell j}=\{U_\ell\geq U_j\}.
\]
Then
\[
 |b_{j\ell}|\leq C_{\mathscr U},\qquad
 \vartheta\leq D\leq C_{\mathscr U}\vartheta,\qquad
 y_j=\mu y_\ell+b_{j\ell}.
\]
For $D$ large, $\mu\leq1/4$.  If $y\in\Theta_{\ell j}$, then
\begin{align*}
 1+|\mu y_\ell+b_{j\ell}|^2
 \geq\mu(1+|y_\ell|^2),\qquad
 \mu(1-2\mu)|y_\ell|^2
 \leq1+2|b_{j\ell}|^2.
\end{align*}
Since $\mu^{-1/2}=\vartheta\asymp D$, this proves
\begin{equation*}
 \Theta_{\ell j}\subset\{|y_\ell|\leq CD\}.
\end{equation*}
On this set,
\[
 U_j(y)\leq C\lambda_\ell^{(n-2)/2}D^{-(n-2)},\qquad
 U_\ell(y)^\beta
 \leq C\lambda_\ell^{a/2}\langle y_\ell\rangle^{-a}.
\]
Changing to $y_\ell$ coordinates in $K_{\ell j}$ and applying
\eqref{eq:strong-B-elementary-convolution}, we obtain
\begin{equation}\label{eq:strong-B-order-greater-point}
 K_{\ell j}(x)\leq C\lambda_\ell^{\alpha/2}
 \begin{cases}
 D^{-(n-2)}\langle x_\ell\rangle^{-2},&|x_\ell|\leq CD,\\
 D^{-(n-\alpha)}|x_\ell|^{-\alpha},&|x_\ell|>CD.
 \end{cases}
\end{equation}
Multiplication by $U_\ell^\beta$ gives
\begin{align}
 \|K_{\ell j}U_\ell^\beta\|_{L^r}
 &\leq C D^{-(n-2)}
 \left(\int_0^{CD}(1+s)^{n-1-(a+2)r}ds\right)^{1/r}+C D^{-(n-\alpha)}
 \left(\int_{CD}^{\infty}s^{n-1-(n+2)r}ds\right)^{1/r}\notag\\
 &\leq C\left(D^{-(n-2)}+D^{-e_\alpha}
 +D^{-(n-2)}(1+\log D)^{1/r}
 \mathbf1_{\{\alpha=(n+6)/2\}}\right).          \label{eq:strong-B-order-greater-l}
\end{align}
For the output at the smaller scale, use
$U_j^\beta\leq C\lambda_\ell^{a/2}D^{-a}$ until
$|x_\ell|\simeq D^2$, and its exact decay afterwards.  Splitting at
$D$ and $D^2$ gives
\begin{align}
 \|K_{\ell j}U_j^\beta\|_{L^r}
 &\leq C D^{-(n-2+a)}
 \left(\int_0^{CD}s^{n-1}(1+s)^{-2r}ds\right)^{1/r}\notag\\
 &\quad+C D^{-(n-\alpha+a)}
 \left(\int_{CD}^{CD^2}s^{n-1-\alpha r}ds\right)^{1/r}
 +C D^{-e_\alpha}\leq CD^{-e_\alpha}.                              \label{eq:strong-B-order-greater-j}
\end{align}
At $(n,\alpha)=(6,4)$ one can see the endpoint contribution directly.
In the $x_\ell$ coordinates,
\[
 U_j(x)=C\lambda_\ell^2\mu^2
 \langle\mu x_\ell+b_{j\ell}\rangle^{-4}.
\]
For a fixed large $C_0$ this gives
\[
 U_j(x)\leq C\lambda_\ell^2
 \begin{cases}
 D^{-4},&|x_\ell|\leq C_0D^2,\\
 D^4|x_\ell|^{-4},&|x_\ell|>C_0D^2.
 \end{cases}
\]
Combining this with \eqref{eq:strong-B-order-greater-point},
\begin{align*}
 \|K_{\ell j}U_j\|_{L^{3/2}}
 &\leq CD^{-8}
 \left(\int_0^{CD}s^5\langle s\rangle^{-3}ds\right)^{2/3}+CD^{-6}
 \left(\int_{CD}^{C_0D^2}s^{-1}ds\right)^{2/3}
 +CD^2\left(\int_{C_0D^2}^{\infty}s^{-7}ds\right)^{2/3}\\
 &\leq CD^{-6}(1+\log D)^{2/3}\leq CD^{-4}.
\end{align*}
Thus
\eqref{eq:strong-B-order-greater-l} and
\eqref{eq:strong-B-order-greater-j} prove
\eqref{eq:strong-B-same-index} in Case 1.

\noindent{\bf Case 2: $\ell\prec j$.}
Set
\[
 \vartheta=\left(\frac{\lambda_j}{\lambda_\ell}\right)^{1/2},
 \qquad\mu=\vartheta^{-2},\qquad
 b_{\ell j}=\lambda_\ell(z_j-z_\ell).
\]
As in Case 1,
\[
 |b_{\ell j}|\leq C_{\mathscr U},\qquad
 \vartheta\leq D\leq C_{\mathscr U}\vartheta,\qquad
 y_\ell=\mu y_j+b_{\ell j}.
\]
This time the condition $U_\ell\geq U_j$ implies
\begin{equation*}
 \Theta_{\ell j}\subset\{|y_j|\geq cD\}.
\end{equation*}
In $y_j$ coordinates the source has two ranges:
\begin{align*}
 U_\ell^\beta U_j\mathbf1_{\Theta_{\ell j}}
 &\leq C\lambda_j^{(2n-\alpha)/2}
 \bigg[D^{-a}|y_j|^{-(n-2)}
 \mathbf1_{\{cD\leq|y_j|\leq CD^2\}}+D^a|y_j|^{-(2n-\alpha)}
 \mathbf1_{\{|y_j|>CD^2\}}\bigg].
\end{align*}
Accordingly write
\[
 K_{\ell j}(x)\leq C\lambda_j^{\alpha/2}(I_1(x_j)+I_2(x_j)),
\]
where, for fixed constants $c,C>0$,
\begin{align*}
 I_1(X)=D^{-a}\int_{cD\leq|Y|\leq CD^2}
 |X-Y|^{-\alpha}|Y|^{-(n-2)}dY,\qquad
 I_2(X)=D^a\int_{|Y|>CD^2}
 |X-Y|^{-\alpha}|Y|^{-(2n-\alpha)}dY.
\end{align*}
We estimate these two integrals in the three output regions.  First,
if $|X|\leq C_1D$, the change of variables $X=D\xi$, $Y=DZ$
gives
\begin{align*}
 I_1(X)+I_2(X)
 \leq CD^{-n}\int_{|Z|\geq c}
 |\xi-Z|^{-\alpha}|Z|^{-(n-2)}dZ\leq CD^{-n}.
\end{align*}
The last integral is uniform for bounded $\xi$: near $Z=\xi$ use
$\alpha<n$, and for large $Z$ use $\alpha>2$.

Next suppose $C_1D<|X|\leq C_2D^2$.  Split the $Y$-space into
\[
 \mathcal A_1=\{|Y|\leq|X|/2\},\qquad
 \mathcal A_2=\{|Y|>|X|/2,\ |X-Y|\leq|X|/2\},
\]
and the remaining set $\mathcal A_3$.  Direct radial integration gives
\begin{align*}
 \int_{\mathbb R^n}|X-Y|^{-\alpha}|Y|^{-(n-2)}dY
 \leq C|X|^{-\alpha}\int_0^{|X|/2}s\,ds+C|X|^{-(n-2)}
 \int_0^{|X|/2}s^{n-1-\alpha}ds+C\int_{|X|/2}^{\infty}s^{1-\alpha}ds
 \leq C|X|^{2-\alpha}.
\end{align*}
For the exterior source, rescaling at $D^2$ gives
\[
 \int_{|Y|>CD^2}|X-Y|^{-\alpha}|Y|^{-(2n-\alpha)}dY
 \leq CD^{-2n},
\]
and hence
\begin{equation*}
 I_1(X)+I_2(X)\leq CD^{-a}|X|^{2-\alpha},\qquad
 C_1D<|X|\leq C_2D^2.
\end{equation*}

Finally assume $|X|>C_2D^2$.  Since
$|X-Y|\geq|X|/2$ on the support of $I_1$,
\[
 I_1(X)\leq CD^{-a}|X|^{-\alpha}
 \int_{cD}^{CD^2}s\,ds
 \leq CD^{\alpha-(n-2)}|X|^{-\alpha}.
\]
For $I_2$, split once more into
$|Y|\leq|X|/2$, $|X-Y|\leq|X|/2$, and their complement.  Then
\begin{align*}
 &\quad\int_{|Y|>CD^2}|X-Y|^{-\alpha}|Y|^{-(2n-\alpha)}dY\\
 &\leq C|X|^{-\alpha}
 \int_{CD^2}^{|X|/2}s^{\alpha-n-1}ds
 +C|X|^{-(2n-\alpha)}
 \int_0^{|X|/2}s^{n-1-\alpha}ds+C\int_{2|X|}^{\infty}s^{-n-1}ds
 \leq CD^{-2(n-\alpha)}|X|^{-\alpha}.
\end{align*}
Therefore
\begin{equation*}
 I_1(X)+I_2(X)
 \leq CD^{\alpha-(n-2)}|X|^{-\alpha},\qquad |X|>C_2D^2.
\end{equation*}
Thus
\begin{equation}\label{eq:strong-B-order-less-point}
 K_{\ell j}(x)\leq C\lambda_j^{\alpha/2}
 \begin{cases}
 D^{-n},&|x_j|\leq CD,\\
 D^{-a}|x_j|^{-(\alpha-2)},&CD<|x_j|\leq CD^2,\\
 D^{\alpha-(n-2)}|x_j|^{-\alpha},&|x_j|>CD^2.
 \end{cases}
\end{equation}
Since
$
 U_j^\beta\leq C\lambda_j^{a/2}\langle x_j\rangle^{-a}$,
polar integration in the three regions gives
\begin{align}
 \|K_{\ell j}U_j^\beta\|_{L^r}
 &\leq C\bigg[D^{-n}
 \left(\int_0^{CD}(1+s)^{n-1-ar}ds\right)^{1/r}\notag\\
 &\quad+D^{-a}\left(\int_{CD}^{CD^2}s^{n-1-nr}ds\right)^{1/r}
 +D^{\alpha-n+2}
 \left(\int_{CD^2}^{\infty}s^{n-1-(n+2)r}ds\right)^{1/r}\bigg]\notag\\
 &\leq C\left(D^{-e_\alpha}+D^{-(2n-\alpha)}
 +D^{-n}(1+\log D)^{1/r}\mathbf1_{\{(n,\alpha)=(6,4)\}}\right).
                                                        \label{eq:strong-B-order-less-j}
\end{align}
Moreover,
\[
 U_\ell^\beta\leq C\lambda_j^{a/2}
 \begin{cases}
 D^{-a},&|x_j|\leq CD^2,\\
 D^a|x_j|^{-a},&|x_j|>CD^2.
 \end{cases}
\]
Using this in \eqref{eq:strong-B-order-less-point},
\begin{align}
 \|K_{\ell j}U_\ell^\beta\|_{L^r}
 &\leq C\bigg[D^{-(n+a)}(CD)^{b_0}
 +D^{-2a}\left(\int_{CD}^{CD^2}
 s^{n-1-(\alpha-2)r}ds\right)^{1/r}+D^4
 \left(\int_{CD^2}^{\infty}s^{n-1-(n+2)r}ds\right)^{1/r}\bigg]\notag\\
 &\leq C\left(D^{-e_\alpha}+D^{-(n-2)}
 +D^{-(n-2)}(1+\log D)^{1/r}
 \mathbf1_{\{\alpha=(n+6)/2\}}\right).          \label{eq:strong-B-order-less-l}
\end{align}
The logarithmic terms in the last two displays have faster algebraic
decay than $D^{-b_\alpha}$.  Hence
\eqref{eq:strong-B-order-less-j}--
\eqref{eq:strong-B-order-less-l} prove
\eqref{eq:strong-B-same-index} in Case 2.

\noindent{\bf Case 3: $\ell$ and $j$ are incomparable.}
Put
\[
 d_{\ell j}=|z_\ell-z_j|,\qquad
 D_0=\sqrt{\lambda_\ell\lambda_j}\,d_{\ell j},\qquad
 \theta=\sqrt{\frac{\lambda_\ell}{\lambda_j}},
\qquad
 L_\ell=\frac{D_0\theta}{2},\qquad
 L_j=\frac{D_0}{2\theta}.
\]
In this case
\begin{equation}\label{eq:strong-case3-scales}
 D_0\asymp D,\qquad
 \max\{\theta,\theta^{-1}\}=o(D_0),\qquad
 L_\ell,L_j\to\infty.
\end{equation}
Let
\[
 \mathcal B_\ell=B_{d_{\ell j}/2}(z_\ell),\qquad
 \mathcal B_j=B_{d_{\ell j}/2}(z_j),\qquad
 \mathcal B_\infty=(\mathcal B_\ell\cup\mathcal B_j)^c,
\]
and set $\mathcal E_m=\Theta_{\ell j}\cap\mathcal B_m$ for
$m=\ell,j,\infty$.  We justify the geometry that will be used below.
The dominance condition is equivalent to
\[
 \frac{\lambda_\ell}{1+|y_\ell|^2}
 \geq\frac{\lambda_j}{1+|y_j|^2}.
\]
If $y\in\mathcal B_\ell$, then $|y_j|<3L_j$ and hence
\[
 \frac{\lambda_j}{1+|y_j|^2}
 \geq c\frac{\lambda_\ell}{D_0^2}.
\]
Thus $y\in\mathcal E_\ell$ implies $|y_\ell|\leq CD_0$; the ball
condition also gives $|y_\ell|<L_\ell$.  If
$y\in\mathcal B_j$, then $|y_\ell|>L_\ell$, so
\[
 \frac{\lambda_\ell}{1+|y_\ell|^2}
 \leq C\frac{\lambda_j}{D_0^2}.
\]
Therefore $y\in\mathcal E_j$ implies $|y_j|\geq cD_0$, while
$|y_j|<L_j$ follows from $y\in\mathcal B_j$.  Finally, on
$\mathcal B_\infty$ the triangle inequality gives
\[
 \frac13|y-z_j|\leq|y-z_\ell|\leq3|y-z_j|.
\]
We have proved
\begin{align}
 \mathcal E_\ell\subset\{|y_\ell|\leq\min(L_\ell,CD_0)\},
\quad
 \mathcal E_j\subset\{cD_0\leq|y_j|<L_j\},
 \quad
 y\in\mathcal B_\infty&\Longrightarrow
 |y-z_\ell|\asymp|y-z_j|.                         \label{eq:strong-case3-Binf}
\end{align}
We retain the two possible scale directions inside Case 3.

\smallskip
\noindent{\it Case 3(a): $\lambda_\ell\geq\lambda_j$.}
Then $\mathcal E_j=\varnothing$.  Indeed, put
$t:=\lambda_\ell/\lambda_j\geq1$ and
$A:=(\lambda_jd_{\ell j}/2)^2=L_j^2$.  By
\eqref{eq:strong-case3-scales}, $A>1$ for all sufficiently large
indices.  If $y\in\mathcal B_j$, then
$|y-z_j|<d_{\ell j}/2<|y-z_\ell|$, and hence
\[
 \frac{\lambda_\ell/(1+\lambda_\ell^2|y-z_\ell|^2)}
 {\lambda_j/(1+\lambda_j^2|y-z_j|^2)}
 <\frac{t(1+A)}{1+t^2A}\leq1.
\]
The last inequality is equivalent to
$(t-1)(tA-1)\geq0$.  Thus $U_\ell(y)<U_j(y)$ on
$\mathcal B_j$, which is incompatible with
$y\in\Theta_{\ell j}$.  Consequently
$\mathcal E_j=\Theta_{\ell j}\cap\mathcal B_j=\varnothing$.
Split
\[
 U_\ell^\beta U_j\mathbf1_{\Theta_{\ell j}}=f_1+f_2,
\]
where $f_1$ is supported in $|y_\ell|\leq2D_0$ and $f_2$ in its
complement.  The dominance relation gives
\begin{align*}
 f_1\leq C\lambda_\ell^{(2n-\alpha)/2}D_0^{-(n-2)}
 \langle y_\ell\rangle^{-a}\mathbf1_{\{|y_\ell|\leq2D_0\}},\quad
 f_2\leq U_\ell^{p_\alpha}\mathbf1_{\{|y_\ell|>2D_0\}}.
\end{align*}
Writing $K^{(q)}=\mathcal I_\alpha f_q$, the elementary convolution
estimate and its exterior analogue give
\begin{align}
 K^{(1)}(x)&\leq C\lambda_\ell^{\alpha/2}
 \begin{cases}
 D_0^{-(n-2)}\langle x_\ell\rangle^{-2},&|x_\ell|\leq CD_0,\\
 D_0^{-(n-\alpha)}|x_\ell|^{-\alpha},&|x_\ell|>CD_0,
 \end{cases}                                               \label{eq:strong-B-incomp-f1}\\
 K^{(2)}(x)&\leq C\lambda_\ell^{\alpha/2}
 \begin{cases}
 D_0^{-n},&|x_\ell|\leq D_0,\\
 D_0^{-(n-\alpha)}|x_\ell|^{-\alpha},&|x_\ell|>D_0.
 \end{cases}                                               \label{eq:strong-B-incomp-f2}
\end{align}
Multiplying by $U_\ell^\beta$ and integrating at the cutoff $D_0$
gives
\begin{align}
 \|(K^{(1)}+K^{(2)})U_\ell^\beta\|_{L^r}
 &\leq C D_0^{-(n-2)}
 \left(\int_0^{CD_0}(1+s)^{n-1-(a+2)r}ds\right)^{1/r}+C D_0^{-(n-\alpha)}
 \left(\int_{D_0}^{\infty}s^{n-1-(n+2)r}ds\right)^{1/r}\notag\\
 &\quad+C D_0^{-n}
 \left(\int_0^{D_0}(1+s)^{n-1-ar}ds\right)^{1/r}\leq C\left(D_0^{-(n-2)}+D_0^{-e_\alpha}
 +\mathcal L_{\rm loc}(D_0)\right),                 \label{eq:strong-B-incomp-l}
\end{align}
where
\[
 \mathcal L_{\rm loc}(D_0)=
 D_0^{-(n-2)}(1+\log D_0)^{1/r}\mathbf1_{\{\alpha=(n+6)/2\}}
 +D_0^{-n}(1+\log D_0)^{1/r}\mathbf1_{\{(n,\alpha)=(6,4)\}}.
\]
For the crossed output $U_j^\beta$, write
\[
 c_{\ell j}=\lambda_\ell(z_j-z_\ell),\qquad
 |c_{\ell j}|=D_0\theta,
\]
and use
\[
 U_j^\beta=C\lambda_\ell^{a/2}\theta^{-a}
 \left\langle\theta^{-2}(x_\ell-c_{\ell j})\right\rangle^{-a}.
\]
Let
\[
 E_\ell=\{|x_\ell|\leq D_0/4\},\qquad
 E_j=\{|x_j|\leq D_0/(4\theta)\}.
\]
On $E_\ell$, $U_j^\beta\leq C\lambda_\ell^{a/2}D_0^{-a}$.
Since $b_0>2$, \eqref{eq:strong-power-integral} gives
\begin{align}
 \|(K^{(1)}+K^{(2)})U_j^\beta\|_{L^r(E_\ell)}
 &\leq C D_0^{-(n-2+a)}\Phi_2(D_0)
      +C D_0^{-(n+a-b_0)}\notag\\
 &\leq C D_0^{-(n-2+a)+(b_0-2)}+CD_0^{-e_\alpha}
 =CD_0^{-e_\alpha}.                         \label{eq:strong-case3a-Eell}
\end{align}
On $E_j$, one has $|x_\ell|\asymp D_0\theta$ and
$dx_\ell=\theta^{2n}dx_j$.  If $a<b_0$, then
\begin{equation*}
 \|(K^{(1)}+K^{(2)})U_j^\beta\|_{L^r(E_j)}
 \leq CD_0^{-n}\Phi_a(D_0/\theta)
 \leq CD_0^{-e_\alpha}\theta^{a-b_0}
 \leq CD_0^{-e_\alpha}.
\end{equation*}
If $a=b_0$, then $(n,\alpha)=(6,4)$ and
\begin{equation*}
 \|(K^{(1)}+K^{(2)})U_j^\beta\|_{L^r(E_j)}
 \leq CD_0^{-n}\bigl(1+\log(D_0/\theta)\bigr)^{1/r}
 \leq CD_0^{-n}(1+\log D_0)^{1/r}.
\end{equation*}

For completeness, we evaluate every contribution on
$(E_\ell\cup E_j)^c$.  The two balls have centers at distance
$D_0\theta$ in the $x_\ell$ coordinates.  Splitting at one half of that
distance and applying \eqref{eq:strong-B-incomp-f1}--
\eqref{eq:strong-B-incomp-f2} gives four nonnegative quantities
\begin{align*}
 J_1=CD_0^{-(n-\alpha+a)}
 \left(\int_{D_0/4}^{D_0\theta/2}
 s^{n-1-\alpha r}\,ds\right)^{1/r},\quad
 J_2=CD_0^{-(n-\alpha)}\theta^{-a}(D_0\theta)^{-\alpha}
 \left(\int_{D_0\theta/4}^{D_0\theta/2}
 s^{n-1}(s/\theta^2)^{-ar}\,ds\right)^{1/r},
 \end{align*}
 \begin{align*}
 J_3^{(1)}=CD_0^{-(n-\alpha)}\theta^{-a}
 (D_0\theta)^{b_0-\alpha}(D_0/\theta)^{-a},\qquad
 J_3^{(2)}=CD_0^{-(n-\alpha)}\theta^{-a}\theta^{2a}
 (D_0\theta)^{b_0-(n+2)}.
\end{align*}
Here $\alpha-b_0=b_0-a\geq0$.  If $\alpha>b_0$, then
\eqref{eq:strong-annular-integral} and
$\alpha+a=n+2=2b_0$ imply
\begin{equation*}
 J_1\leq
 CD_0^{-(n-\alpha+a)}D_0^{b_0-\alpha}
 =CD_0^{-(n+a-b_0)}=CD_0^{-e_\alpha}.
\end{equation*}
If $\alpha=b_0$, then $(n,\alpha)=(6,4)$ and instead
\begin{equation*}
 J_1\leq CD_0^{-n}(1+\log\theta)^{1/r}
 \leq CD_0^{-n}(1+\log D_0)^{1/r}.
\end{equation*}
The other three quantities contain no unevaluated integral:
\begin{align}
 J_2
 &=CD_0^{-(n-\alpha)-\alpha+b_0-a}
   \theta^{-a-\alpha+b_0-a+2a}
 =CD_0^{-e_\alpha}\theta^{b_0-\alpha},
                                                        \label{eq:strong-case3a-J2}\\
 J_3^{(1)}
 &=CD_0^{-(n-\alpha)+b_0-\alpha-a}
   \theta^{-a+b_0-\alpha+a}
 =CD_0^{-e_\alpha}\theta^{b_0-\alpha},
                                                        \label{eq:strong-case3a-J31}\\
 J_3^{(2)}
 &=CD_0^{-(n-\alpha)+b_0-(n+2)}
   \theta^{-a+2a+b_0-(n+2)}
 =CD_0^{-e_\alpha}\theta^{a-b_0}.       \label{eq:strong-case3a-J32}
\end{align}
Because $\theta\geq1$, $b_0-\alpha\leq0$, and $a-b_0\leq0$,
\eqref{eq:strong-case3a-J2}--\eqref{eq:strong-case3a-J32} are bounded by
$CD_0^{-e_\alpha}$.  Equations
\eqref{eq:strong-case3a-Eell}--\eqref{eq:strong-case3a-J32}  give
\begin{equation*}
 \|(K^{(1)}+K^{(2)})U_j^\beta\|_{L^r}
 \leq C\left(D_0^{-e_\alpha}+\mathcal L_{\rm loc}(D_0)\right).
\end{equation*}
Combining this estimate with \eqref{eq:strong-B-incomp-l} yields
\begin{equation}\label{eq:strong-case3a-result}
 \|K_{\ell j}U_\ell^\beta\|_{L^r}
 +\|K_{\ell j}U_j^\beta\|_{L^r}
 \leq C\left(D_0^{-(n-2)}+D_0^{-e_\alpha}
 +\mathcal L_{\rm loc}(D_0)\right).
\end{equation}

\smallskip
\noindent{\it Case 3(b): $\lambda_\ell<\lambda_j$.}
Set $\theta_*=\sqrt{\lambda_j/\lambda_\ell}>1$.  Then
\[
 L_j=\frac{D_0\theta_*}{2},\qquad
 L_\ell=\frac{D_0}{2\theta_*},\qquad \theta_*=o(D_0).
\]
Define
\[
 g_m=U_\ell^\beta U_j\mathbf1_{\mathcal E_m},\qquad
 P_m=\mathcal I_\alpha(g_m),\qquad m=\ell,j,\infty.
\]
The partition is exact: $K_{\ell j}=P_\ell+P_j+P_\infty$.  From
\eqref{eq:strong-case3-Binf},
\begin{align*}
 g_\ell&\leq C\lambda_\ell^{(2n-\alpha)/2}D_0^{-(n-2)}
 \langle y_\ell\rangle^{-a}\mathbf1_{\{|y_\ell|<L_\ell\}},
 \\
 g_j&\leq C\lambda_j^{(2n-\alpha)/2}D_0^{-a}|y_j|^{-(n-2)}
 \mathbf1_{\{cD_0\leq|y_j|<L_j\}}, \\
 g_\infty&\leq C\lambda_j^{(2n-\alpha)/2}\theta_*^a
 |y_j|^{-(2n-\alpha)}\mathbf1_{\{|y_j|\geq L_j\}}.
\end{align*}
For the exterior source we use
\begin{equation}\label{eq:strong-case3-tail-convolution}
 \int_{|Y|\geq L}|X-Y|^{-\alpha}|Y|^{-(2n-\alpha)}dY
 \leq C\begin{cases}
 L^{-n},&|X|\leq2L,\\
 L^{-(n-\alpha)}|X|^{-\alpha},&|X|>2L.
 \end{cases}
\end{equation}
It follows from \eqref{eq:strong-B-elementary-convolution} and
\eqref{eq:strong-case3-tail-convolution} that
\begin{align}
 P_\ell(x)&\leq C\lambda_\ell^{\alpha/2}D_0^{-(n-2)}
 \begin{cases}
 \langle x_\ell\rangle^{-2},&|x_\ell|\leq2L_\ell,\\
 L_\ell^{\alpha-2}|x_\ell|^{-\alpha},&|x_\ell|>2L_\ell,
 \end{cases}                                               \notag\\
 P_j(x)&\leq C\lambda_j^{\alpha/2}
 \begin{cases}
 D_0^{-n},&|x_j|\leq cD_0,\\
 D_0^{-a}|x_j|^{-(\alpha-2)},&cD_0<|x_j|\leq4L_j,\\
 D_0^{-a}L_j^2|x_j|^{-\alpha},&|x_j|>4L_j,
 \end{cases}                                               \label{eq:strong-case3b-Pj}\\
 P_\infty(x)&\leq C\lambda_j^{\alpha/2}\theta_*^a
 \begin{cases}
 L_j^{-n},&|x_j|\leq2L_j,\\
 L_j^{-(n-\alpha)}|x_j|^{-\alpha},&|x_j|>2L_j.
 \end{cases}                                               \label{eq:strong-case3b-Pinf}
\end{align}
We first estimate the two outputs carrying the same index as the
potential.  Since $L_\ell=D_0/(2\theta_*)$, equations
\eqref{eq:strong-power-integral} and
\eqref{eq:strong-exterior-integral} give
\begin{align}
 \|P_\ell U_\ell^\beta\|_{L^r}
 &\leq CD_0^{-(n-2)}\Phi_{a+2}(2L_\ell)
 +CD_0^{-(n-2)}L_\ell^{\alpha-2-b_0}.       \label{eq:strong-case3b-PlUl-raw}
\end{align}
The second term on the right equals
\begin{equation*}
 D_0^{-(n-2)}L_\ell^{\alpha-2-b_0}
 =CD_0^{-e_\alpha}\theta_*^\kappa .
\end{equation*}
For the first term there are exactly three possibilities.  If
$\kappa<0$, then $a+2<b_0$ and
\begin{equation*}
 D_0^{-(n-2)}\Phi_{a+2}(2L_\ell)
 \leq CD_0^{-(n-2)}(D_0/\theta_*)^{-\kappa}
 =CD_0^{-e_\alpha}\theta_*^\kappa.
\end{equation*}
If $\kappa=0$, equivalently
$\alpha=(n+6)/2$, then
\begin{equation*}
 D_0^{-(n-2)}\Phi_{a+2}(2L_\ell)
 \leq CD_0^{-(n-2)}
 (1+\log(D_0/\theta_*))^{1/r}.
\end{equation*}
If $\kappa>0$, then $a+2>b_0$ and
\begin{equation*}
 D_0^{-(n-2)}\Phi_{a+2}(2L_\ell)
 \leq CD_0^{-(n-2)}.
\end{equation*}

The exponent in the middle integral for $P_jU_j^\beta$ is $n$, because
$(\alpha-2)+a=n$.  Thus the three pointwise ranges in
\eqref{eq:strong-case3b-Pj} yield
\begin{align}
 \|P_jU_j^\beta\|_{L^r}
 \leq CD_0^{-n}\Phi_a(cD_0)
 +CD_0^{-a}
 \left(\int_{cD_0}^{4L_j}s^{n-1-nr}\,ds\right)^{1/r}+CD_0^{-a}L_j^2
 \left(\int_{4L_j}^{\infty}s^{n-1-(n+2)r}\,ds\right)^{1/r}.
                                                        \label{eq:strong-case3b-PjUj-raw}
\end{align}
If $a<b_0$, the first term is bounded by $CD_0^{-e_\alpha}$;
if $a=b_0$, it is bounded by
$CD_0^{-n}(1+\log D_0)^{1/r}$.  Since $n>b_0$, the second term in
\eqref{eq:strong-case3b-PjUj-raw} satisfies
\begin{equation*}
 D_0^{-a}
 \left(\int_{cD_0}^{4L_j}s^{n-1-nr}\,ds\right)^{1/r}
 \leq CD_0^{-a+b_0-n}=CD_0^{-e_\alpha}.
\end{equation*}
The last term satisfies
\begin{equation*}
 D_0^{-a}L_j^2(4L_j)^{-b_0}
 =CD_0^{-e_\alpha}\theta_*^\kappa.
\end{equation*}

We next calculate the crossed output $P_\ell U_j^\beta$.  In the
$x_\ell$ coordinates put
\[
 L=L_\ell,\qquad c=\lambda_\ell(z_j-z_\ell),\qquad |c|=2L,
\qquad
 V_j(X)=\theta_*^a\langle\theta_*^2(X-c)\rangle^{-a}.
\]
Split the output into
\[
 \mathcal C_\ell=\{|X|\leq L/2\},\qquad
 \mathcal C_j=\{|X-c|\leq L/2\},\qquad
 \mathcal C_\infty=(\mathcal C_\ell\cup\mathcal C_j)^c.
\]
On $\mathcal C_\ell$, $V_j\leq CD_0^{-a}$.  Since $2<b_0$,
\begin{align}
 \|P_\ell U_j^\beta\|_{L^r(\mathcal C_\ell)}
 &\leq CD_0^{-(n-2+a)}\Phi_2(L/2)
 \leq CD_0^{-e_\alpha}\theta_*^{\,2-b_0}.
                                                        \label{eq:strong-case3b-PlUj-Cell}
\end{align}
On $\mathcal C_j$, one has
$P_\ell\leq CD_0^{-(n-2)}L^{-2}$.  The substitution
$W=\theta_*^2(X-c)$ gives
\begin{align}
 \|P_\ell U_j^\beta\|_{L^r(\mathcal C_j)}
 &\leq CD_0^{-(n-2)}L^{-2}\theta_*^{a-2b_0}
 \Phi_a(\theta_*^2L).                         \label{eq:strong-case3b-PlUj-Cj-raw}
\end{align}
If $a<b_0$, the right-hand side of
\eqref{eq:strong-case3b-PlUj-Cj-raw} equals, up to a fixed factor,
\begin{equation*}
 D_0^{-(n-2)}(D_0/\theta_*)^{-2}
 \theta_*^{a-2b_0}(D_0\theta_*)^{b_0-a}
 =D_0^{-e_\alpha}\theta_*^{\,2-b_0}.
\end{equation*}
If $a=b_0$, then $(n,\alpha)=(6,4)$ and the bound is
\begin{equation*}
 CD_0^{-n}\theta_*^{\,2-b_0}
 (1+\log(D_0\theta_*))^{1/r}.
\end{equation*}

On $\mathcal C_\infty\cap\{|X|\leq2L\}$ the rescaling $X=LZ$ gives
$$
 CD_0^{-(n-2)}\theta_*^{-a}
 \left(\int_{\substack{L/2<|X|\leq2L\\|X-c|>L/2}}
 |X|^{-2r}|X-c|^{-ar}\,dX\right)^{1/r}\leq CD_0^{-(n-2)}\theta_*^{-a}L^{b_0-2-a}
 =CD_0^{-e_\alpha}\theta_*^{\,2-b_0}.
$$
On $\mathcal C_\infty\cap\{|X|>2L\}$, the same change of variables,
written out with all powers, gives
\begin{align}
 CD_0^{-(n-2)}L^{\alpha-2}\theta_*^{-a}
 \left(\int_{\substack{|X|>2L\\|X-c|>L/2}}
 |X|^{-\alpha r}|X-c|^{-ar}\,dX\right)^{1/r}\leq CD_0^{-e_\alpha}\theta_*^{\,2-b_0}.
                                       \label{eq:strong-case3b-PlUj-Cinf-far}
\end{align}
Because $b_0>2$, equations
\eqref{eq:strong-case3b-PlUj-Cell}--
\eqref{eq:strong-case3b-PlUj-Cinf-far} imply
\begin{equation*}
 \|P_\ell U_j^\beta\|_{L^r}
 \leq C\left(D_0^{-e_\alpha}
 +D_0^{-n}(1+\log D_0)^{1/r}
  \mathbf1_{\{(n,\alpha)=(6,4)\}}\right).
\end{equation*}

We now calculate $P_jU_\ell^\beta$.  In the $x_j$ coordinates put
\[
 L=L_j,\qquad c=\lambda_j(z_\ell-z_j),\qquad |c|=2L,
\qquad
 V_\ell(X)=\theta_*^{-a}
 \langle\theta_*^{-2}(X-c)\rangle^{-a}.
\]
With the exact partition
\[
 \mathcal D_j=\{|X|\leq L/2\},\qquad
 \mathcal D_\ell=\{|X-c|\leq L/2\},\qquad
 \mathcal D_\infty=(\mathcal D_j\cup\mathcal D_\ell)^c,
\]
the region $\mathcal D_j$ gives
\begin{align}
 \|P_jU_\ell^\beta\|_{L^r(\mathcal D_j)}
 &\leq C\bigg[D_0^{-n-a+b_0}
 +D_0^{-2a}
 \left(\int_{cD_0}^{L/2}s^{n-1-(\alpha-2)r}\,ds\right)^{1/r}\bigg].
                                                        \label{eq:strong-case3b-PjUl-Dj-raw}
\end{align}
The first summand is $CD_0^{-e_\alpha}$.  Since
\[
 (\alpha-2)-b_0=-\kappa,
\]
the second summand in \eqref{eq:strong-case3b-PjUl-Dj-raw} is bounded by
$$
  CD_0^{-e_\alpha}\theta_*^\kappa, \quad \text{if}\,\, \kappa>0; \qquad
CD_0^{-(n-2)}(1+\log\theta_*)^{1/r}, \quad \text{if}\,\, \kappa=0;\qquad
CD_0^{-e_\alpha} \quad \text{if}\,\, \kappa<0.
$$

On $\mathcal D_\ell$ one has
\begin{align*}
 \|P_jU_\ell^\beta\|_{L^r(\mathcal D_\ell)}
 \leq CD_0^{-a}L^{-(\alpha-2)}\theta_*^{-a+2b_0}
 \Phi_a(L/\theta_*^2).
\end{align*}
Here $L/\theta_*^2=D_0/(2\theta_*)\to\infty$.  If $a<b_0$, then
$$
 \|P_jU_\ell^\beta\|_{L^r(\mathcal D_\ell)}
 \leq CD_0^{-a}(D_0\theta_*)^{-(\alpha-2)}
 \theta_*^{-a+2b_0}(D_0/\theta_*)^{b_0-a}=CD_0^{-e_\alpha}\theta_*^\kappa.                                      $$
If $a=b_0$, hence $(n,\alpha)=(6,4)$ and $\kappa=2$, then
\begin{equation*}
 \|P_jU_\ell^\beta\|_{L^r(\mathcal D_\ell)}
 \leq CD_0^{-n}\theta_*^2
 (1+\log(D_0/\theta_*))^{1/r}.
\end{equation*}

On $\mathcal D_\infty\cap\{|X|\leq4L\}$, the rescaling $X=LZ$
gives
\begin{align*}
 CD_0^{-a}\theta_*^a
 \left(\int_{\substack{L/2<|X|\leq4L\\|X-c|>L/2}}
 |X|^{-(\alpha-2)r}|X-c|^{-ar}\,dX\right)^{1/r}\leq CD_0^{-a}\theta_*^a
 L^{b_0-(\alpha-2)-a}
 =CD_0^{-e_\alpha}\theta_*^\kappa.
\end{align*}
On $\mathcal D_\infty\cap\{|X|>4L\}$ one obtains
\begin{align*}
 CD_0^{-a}L^2\theta_*^a
 \left(\int_{\substack{|X|>4L\\|X-c|>L/2}}
 |X|^{-\alpha r}|X-c|^{-ar}\,dX\right)^{1/r}\leq CD_0^{-a}L^2\theta_*^aL^{b_0-\alpha-a}
 =CD_0^{-e_\alpha}\theta_*^\kappa.
\end{align*}

We now eliminate every factor $\theta_*^\kappa$ without an unquantified
appeal to \eqref{eq:strong-case3-scales}.  If $\kappa\leq0$, then
$\theta_*^\kappa\leq1$.  If $\kappa>0$, then
\begin{equation}\label{eq:strong-case3b-kappa-absorption}
 D_0^{-e_\alpha}\theta_*^\kappa
 =D_0^{-(n-2)}\left(\frac{\theta_*}{D_0}\right)^\kappa
 \leq D_0^{-(n-2)}
\end{equation}
for all sufficiently large indices of the extracted sequence.  In the
endpoint $(n,\alpha)=(6,4)$, put $s_*=D_0/\theta_*\to\infty$.  Then
\begin{equation}\label{eq:strong-case3b-endpoint-absorption}
 D_0^{-n}\theta_*^2
 (1+\log(D_0/\theta_*))^{1/r}
 =D_0^{-(n-2)}s_*^{-2}(1+\log s_*)^{1/r}
 \leq CD_0^{-(n-2)}.
\end{equation}
Equations \eqref{eq:strong-case3b-PjUl-Dj-raw}--
\eqref{eq:strong-case3b-endpoint-absorption} prove
\begin{equation*}
 \|P_jU_\ell^\beta\|_{L^r}
 \leq C\left(D_0^{-e_\alpha}+D_0^{-(n-2)}
 +D_0^{-(n-2)}(1+\log D_0)^{1/r}
   \mathbf1_{\{\kappa=0\}}\right).
\end{equation*}

It remains to estimate $P_\infty$.  For the $U_j^\beta$ output,
\eqref{eq:strong-case3b-Pinf} and
\eqref{eq:strong-exterior-integral} give
\begin{align}
 \|P_\infty U_j^\beta\|_{L^r}
 \leq C\theta_*^aL_j^{-n}\Phi_a(2L_j)+C\theta_*^aL_j^{-(n-\alpha)}
 (2L_j)^{-b_0}.                              \label{eq:strong-case3b-PinfUj-raw}
\end{align}
If $a<b_0$, both terms in
\eqref{eq:strong-case3b-PinfUj-raw} are bounded by
\begin{equation*}
 CD_0^{-e_\alpha}\theta_*^{\,2-b_0}\leq CD_0^{-e_\alpha}.
\end{equation*}
If $a=b_0$, they are bounded by
\begin{equation*}
 CD_0^{-n}\theta_*^{\,2-b_0}
 (1+\log(D_0\theta_*))^{1/r}
 +CD_0^{-e_\alpha}\theta_*^{\,2-b_0}
 \leq CD_0^{-e_\alpha}(1+\log D_0)^{1/r},
\end{equation*}
where we used $b_0>2$ and $\theta_*<D_0$ for all sufficiently large
indices.

For the $U_\ell^\beta$ output we do not invoke symmetry.  In the $x_j$
coordinates, let
\[
 X=x_j,\qquad c=\lambda_j(z_\ell-z_j),\qquad |c|=2L_j,
\]
so that
\[
 U_\ell^\beta=C\lambda_j^{a/2}\theta_*^{-a}
 \langle\theta_*^{-2}(X-c)\rangle^{-a}.
\]
The factors $\theta_*^a$ in
\eqref{eq:strong-case3b-Pinf} and $\theta_*^{-a}$ in this formula cancel.
Put $M=D_0/\theta_*=2L_j/\theta_*^2\to\infty$.  On $|X|\leq2L_j$,
the substitution $W=\theta_*^{-2}(X-c)$ gives
\begin{align}
 \|P_\infty U_\ell^\beta\|_{L^r(|X|\leq2L_j)}
 &\leq CL_j^{-n}\theta_*^{2b_0}\Phi_a(CM).
                                                        \label{eq:strong-case3b-PinfUl-near}
\end{align}
The region $2L_j<|X|\leq4L_j$ obeys the same bound, because
$L_j^{-(n-\alpha)}|X|^{-\alpha}\leq CL_j^{-n}$ there.  If $a<b_0$,
the right-hand side of \eqref{eq:strong-case3b-PinfUl-near} is
\begin{equation*}
 CL_j^{-n}\theta_*^{2b_0}(D_0/\theta_*)^{b_0-a}
 =CD_0^{-e_\alpha}\theta_*^\kappa.
\end{equation*}
If $a=b_0$, it is
\begin{equation*}
 CD_0^{-n}\theta_*^2(1+\log(D_0/\theta_*))^{1/r},
\end{equation*}
which is bounded by $CD_0^{-(n-2)}$ by
\eqref{eq:strong-case3b-endpoint-absorption}.
Finally, if $|X|>4L_j$, then $|X-c|\geq|X|/2$ and
$\theta_*^{-2}|X-c|\geq D_0/\theta_*\to\infty$.  Therefore
\begin{align}
 L_j^{-(n-\alpha)}
 \left(\int_{|X|>4L_j}|X|^{-\alpha r}
 \langle\theta_*^{-2}(X-c)\rangle^{-ar}\,dX\right)^{1/r}\leq& CL_j^{-(n-\alpha)}\theta_*^{2a}
 \left(\int_{4L_j}^\infty s^{n-1-(\alpha+a)r}\,ds\right)^{1/r}\notag\\
 =&CL_j^{-e_\alpha}\theta_*^{2a}
 =CD_0^{-e_\alpha}\theta_*^\kappa.       \label{eq:strong-case3b-PinfUl-far}
\end{align}
Here we used $\alpha+a=n+2$ and $2a-e_\alpha=\kappa$.
Equations \eqref{eq:strong-case3b-PinfUl-near}--
\eqref{eq:strong-case3b-PinfUl-far}, followed by
\eqref{eq:strong-case3b-kappa-absorption}, give
\begin{equation}\label{eq:strong-case3b-PinfUl}
 \|P_\infty U_\ell^\beta\|_{L^r}
 \leq C\left(D_0^{-e_\alpha}+D_0^{-(n-2)}\right).
\end{equation}

Since $K_{\ell j}=P_\ell+P_j+P_\infty$, equations
\eqref{eq:strong-case3b-PlUl-raw}--
\eqref{eq:strong-case3b-PinfUl} prove
\begin{equation}\label{eq:strong-case3b-result}
 \|K_{\ell j}U_\ell^\beta\|_{L^r}
 +\|K_{\ell j}U_j^\beta\|_{L^r}
 \leq C\left(D_0^{-(n-2)}+D_0^{-e_\alpha}
 +\mathcal L_{\rm loc}(D_0)\right).
\end{equation}
We verify the final comparison explicitly.  By
\eqref{eq:strong-parameter-identities},
\[
 n-2\geq b_\alpha,
 \qquad
 e_\alpha-b_\alpha
 =(n-\alpha)+(b_0-b_\alpha)>0.
\]
If $\kappa=0$, then $\alpha=(n+6)/2<n$, so $n\geq7$ and
\[
 (n-2)-b_\alpha=(n-2)-b_0=\frac{n-6}{2}>0.
\]
Consequently
\[
 D_0^{-(n-2)}(1+\log D_0)^{1/r}=O(D_0^{-b_\alpha})
 \qquad(\kappa=0).
\]
At $(n,\alpha)=(6,4)$ one has $b_\alpha=4<n=6$, and hence
\[
 D_0^{-n}(1+\log D_0)^{1/r}=O(D_0^{-b_\alpha}).
\]
Finally, \eqref{eq:strong-case3-scales} implies
$D=\max\{\theta,\theta^{-1},D_0\}=D_0$ for all sufficiently large
indices.  Thus \eqref{eq:strong-case3a-result} and
\eqref{eq:strong-case3b-result} prove
\eqref{eq:strong-B-same-index} in Case 3.

We finally treat \(k\notin\{\ell,j\}\).  Since \(\beta\leq1\) and
\(U_\ell\geq U_j\) on \(\Omega_\ell\), Cauchy--Schwarz for the Riesz
potential and the one-bubble identity give
\[
 K_{\ell j}\leq
 \mathcal I_\alpha(U_\ell^{p_\alpha})^{1/2}
 \mathcal I_\alpha(U_j^{p_\alpha})^{1/2}
 =U_\ell^{\gamma/2}U_j^{\gamma/2}.
\]
Partitioning according to the largest of \(U_\ell,U_j,U_k\), the
arithmetic--geometric mean inequality on the \(U_k\)-region and
weighted Young inequality on the other two regions reduce
\(K_{\ell j}U_k^\beta\) to pair products with ordered exponent pairs
\((\beta,\gamma)\) and
\((\gamma/2,\beta+\gamma/2)\), respectively.  After taking the
\(r\)-th power, the scaled exponents in each pair sum to \(n\);
when \(\alpha\ne(n+2)/2\), Lemma \ref{A.4} gives the respective orders
\(R^{-\min\{\alpha,(n+2)/2\}}\) and \(R^{-(n+2)/2}\).  Since
\(b_\alpha\leq\min\{\alpha,(n+2)/2\}\), it follows that
\begin{equation}\label{eq:strong-B-three-Lr}
 \|K_{\ell j}U_k^\beta\|_{L^r}\leq CR^{-b_\alpha},
 \qquad k\notin\{\ell,j\}.
\end{equation}
At \(\alpha=(n+2)/2\), one has \(\beta=\gamma\) and
\[
 K_{\ell j}U_k^\beta
 \leq\frac12(U_\ell^\beta U_k^\beta+U_j^\beta U_k^\beta).
\]
Lemma \ref{A.6} and the Newtonian representation then yield
\begin{equation}\label{eq:strong-B-three-Hminus}
 \|K_{\ell j}U_k^\beta\|_{\dot H^{-1}}
 \leq CR^{-b_\alpha}(1+\log R)^{1/2}.
\end{equation}
Summing \eqref{eq:strong-B-same-index},
\eqref{eq:strong-B-three-Lr}, and
\eqref{eq:strong-B-three-Hminus} in \eqref{eq:strong-B-point} gives
\[
 \|B\|_{L^r}\leq CR^{-b_\alpha},\qquad
 (n,\alpha)\notin\mathcal C_{\log},
\]
and
\[
 \|B\|_{\dot H^{-1}}
 \leq CR^{-b_\alpha}(1+\log R)^{1/2},\qquad
 (n,\alpha)\in\mathcal C_{\log}.
\]
Together with \eqref{eq:strong-A-short} and
\eqref{eq:strong-A-short-log}, this proves the lemma.
\end{proof}

\section{Estimates for the Error Function}\label{Sec5}
This section establishes the estimates for $\rho$.
The linear problem used below is
\begin{equation}\label{3.1}
\left\{\begin{array}{l}
\Delta \phi+N_{1}(\phi,\sigma)
\\= g+\sum_{i=1}^\nu \sum_{a=1}^{n+1} c_a^i \bigg[\Big(|x|^{-\alpha}\ast |U_{i}|^{p_{\alpha}}\Big)(p_{\alpha}-1   )U_{i}^{ p_{\alpha}-2}Z_{i}^{a}+p_{\alpha}\Big(|x|^{-\alpha}\ast (U_{i}^{p_{\alpha}-1}
		Z_{i}^{a})\Big)U^{p_{\alpha}-1}_{i}\bigg], \\
\int_{\mathbb R^n} \bigg[\Big(|x|^{-\alpha}\ast |U_{i}|^{p_{\alpha}}\Big)(p_{\alpha}-1   )U_{i}^{ p_{\alpha}-2}Z_{i}^{a}+p_{\alpha}\Big(|x|^{-\alpha}\ast (U_{i}^{p_{\alpha}-1}
		Z_{i}^{a})\Big)U^{p_{\alpha}-1}_{i}\bigg] \phi\,\dx=0,
\end{array}\right.
\end{equation}
where $\sigma=\sum_{i=1}^\nu U[z_i,\lambda_i]$ is a sum of $\delta$-interacting bubbles, $i=1, \ldots, \nu$, $a=1, \ldots, n+1$. We assume throughout that $\delta$ is sufficiently small. Proposition \ref{pro4.1.1} establishes the required a priori estimate by finite-dimensional reduction.

We first isolate the operator-decoupling assertion needed in the
finite-dimensional reduction.  This assertion includes the case
$p_\alpha<2$, in which the exponent $p_\alpha-2$ is negative.

\begin{lemma}
\label{lem:linearized-form-decoupling}
Set
\[
 p=p_\alpha,\qquad \beta=p-1,\qquad
 t_\alpha=\frac{2n}{n+2-\alpha}.
\]
For a positive function $V$ satisfying
\begin{equation}\label{eq:linearized-form-domain}
 V^\beta\in L^{t_\alpha}(\mathbb R^n),
 \qquad
 \mathcal I_\alpha(V^p)V^{\beta-1}\in L^{n/2}(\mathbb R^n),
\end{equation}
define
\begin{align*}
 \mathscr B_V(v,w):=D_{n,\alpha}p
 \iint_{\mathbb R^n\times\mathbb R^n}
 \frac{V^\beta(y)v(y)V^\beta(x)w(x)}{|x-y|^\alpha}\,\dy\dx
 +D_{n,\alpha}\beta\int_{\mathbb R^n}
 \left(\int_{\mathbb R^n}\frac{V^p(y)}{|x-y|^\alpha}\,\dy\right)
 V^{\beta-1}(x)v(x)w(x)\,\dx .
\end{align*}
The HLS, H\"older, and Sobolev inequalities imply that
\begin{equation}\label{eq:linearized-form-bounded}
 |\mathscr B_V(v,w)|
 \leq C_V\|v\|_{\dot H^1}\|w\|_{\dot H^1}.
\end{equation}
Every bubble and every finite positive sum of bubbles satisfies
\eqref{eq:linearized-form-domain}; these are the only functions $V$
used below.
There is a modulus $\omega_4$, depending only on
$n,\alpha,\nu$, such that $\omega_4(t)\to0$ as $t\downarrow0$
and, for every $\delta$-interacting family
$\{V_1,\ldots,V_\nu\}$, $S=\sum_iV_i$, and every mode $Z_j^b$,
\begin{equation}\label{eq:uniform-mode-decoupling}
 \left|\mathscr B_S(v,Z_j^b)-\mathscr B_{V_j}(v,Z_j^b)\right|
 \leq \omega_4(\delta)\|v\|_{\dot H^1},
 \qquad v\in\dot H^1(\mathbb R^n).
\end{equation}
The estimate is uniform in $j$, $b$, and all bubble parameters.

More precisely, let $\delta_m\to0$, let $S_m=\sum_kV_{k,m}$ be
$\delta_m$-interacting, and normalize about $V_{i,m}$ by
\[
 \widehat v_m(y)=\lambda_{i,m}^{-(n-2)/2}
 v_m(\lambda_{i,m}^{-1}y+z_{i,m}),
 \qquad
 \widehat S_m(y)=\lambda_{i,m}^{-(n-2)/2}
 S_m(\lambda_{i,m}^{-1}y+z_{i,m}).
\]
If $\sup_m\|v_m\|_{\dot H^1}<\infty$, then, with $U=U[0,1]$,
\begin{equation}\label{eq:test-linear-decoupling}
 \mathscr B_{\widehat S_m}(\widehat v_m,\psi)
 -\mathscr B_U(\widehat v_m,\psi)\longrightarrow0
 \qquad\text{for every }\psi\in\dot H^1(\mathbb R^n).
\end{equation}
\end{lemma}

\begin{proof}
Let \(T=\sum_{k=1}^NW_k\), \(N\leq\nu\), and
\(\gamma=2^*-p\).  Since \(\beta t_\alpha=2^*\), we have
\[
\|T^\beta\|_{L^{t_\alpha}}\leq C.
\]
On \(\Omega_k=\{W_k=\max_\ell W_\ell\}\),
\(\mathcal I_\alpha(T^p)\leq C\sum_\ell W_\ell^\gamma\).
If \(\beta\geq1\), then
\(T^{\beta-1}W_\ell^\gamma\leq CW_k^{2^*-2}\); if
\(0<\beta<1\), then
\(T^{\beta-1}W_\ell^\gamma\leq W_\ell^{2^*-2}\).
Thus
\[
\|\mathcal I_\alpha(T^p)T^{\beta-1}\|_{L^{n/2}}\leq C.
\]
The bilinear HLS, H\"older, and Sobolev inequalities give
\[
|\mathscr B_T(v,w)|\leq
C(n,\alpha,N)\|v\|_{\dot H^1}\|w\|_{\dot H^1},
\]
which proves \eqref{eq:linearized-form-domain} and
\eqref{eq:linearized-form-bounded}.

For the decoupling, normalize about the tested bubble.  One profile is
\(U=U[0,1]\); every other profile has the form
\[
W_m(x)=A_n\left(\frac{\mu_m}
{1+\mu_m^2|x-\xi_m|^2}\right)^{(n-2)/2},
\qquad
\mu_m+\mu_m^{-1}+\mu_m|\xi_m|^2\longrightarrow\infty.
\]
After passage to a subsequence, one of the following alternatives holds:
\[
\mu_m\to0,\qquad
\mu_m\to\infty,\qquad
0<c\leq\mu_m\leq C<\infty,\quad |\xi_m|\to\infty.
\]
Set
\[
s=\frac{2n}{2n-\alpha},\qquad
t=t_\alpha=\frac{2n}{n+2-\alpha},\qquad
\mathfrak q=\frac{2n}{n+2}.
\]
Thus
\[
\frac1s=\frac1t+\frac1{2^*},\qquad
\beta t=2^*,\qquad
\frac1{\mathfrak q}+\frac1{2^*}=1.
\]
For every compact \(K\subset\mathbb R^n\), the explicit bubble formula
gives
\begin{equation}\label{eq:new26-foreign-local}
\|W_m^\beta\|_{L^s(K)}
+\|W_m^\beta\|_{L^{\mathfrak q}(K)}
+\|W_m^{2^*-2}\|_{L^{\mathfrak q}(K)}
\longrightarrow0.
\end{equation}
We verify this assertion.  If \(\mu_m\to0\), then
$
\sup_KW_m\leq C_K\mu_m^{(n-2)/2}$.
If \(c\leq\mu_m\leq C\) and \(|\xi_m|\to\infty\), then
$
\sup_KW_m\leq C_K|\xi_m|^{2-n}$.
If \(\mu_m\to\infty\), the change of variables
\(y=\mu_m(x-\xi_m)\) gives, for
\[
A_r:=\frac{n-2}{2}\beta r,\qquad r\in\{s,\mathfrak q\},
\qquad A_0:=2\mathfrak q,
\]
the bounds
\[
\int_KW_m^{\beta r}\dx
\leq C_K\left(\mu_m^{A_r-n}
+\mu_m^{-A_r}(1+\log\mu_m)\right),
\qquad
\int_KW_m^{(2^*-2)\mathfrak q}\dx
\leq C_K\left(\mu_m^{A_0-n}
+\mu_m^{-A_0}(1+\log\mu_m)\right).
\]
Since \(0<A_r<n\) and \(0<A_0<n\), all the right-hand sides tend to
zero.  This proves \eqref{eq:new26-foreign-local}.

The same calculation with the exponent \(1\) in place of \(\beta\)
gives
\begin{equation}\label{eq:new29-foreign-linear-local}
 \|W_m\|_{L^s(K)}
 +\|W_m\|_{L^{\mathfrak q}(K)}
 \longrightarrow0.
\end{equation}
After a further subsequence,
\begin{equation}\label{eq:new29-foreign-ae}
 W_m(x)\longrightarrow0
 \qquad\text{for almost every }x\in\mathbb R^n.
\end{equation}

Write
\[
 T_m=U+P_m,\qquad
 P_m:=\sum_{j=1}^{N-1}W_{j,m},\qquad
 A_m:=T_m^\beta-U^\beta,
\]
where \(W_{j,m}\) are the profiles different from \(U\).  Since
\(\beta t=2^*\),
\begin{equation}\label{eq:new29-A-global}
 \sup_m\|A_m\|_{L^t}<\infty.
\end{equation}
If \(0<\beta\leq1\), then
\[
 0\leq A_m\leq P_m^\beta
 \leq\sum_{j=1}^{N-1}W_{j,m}^\beta.
\]
If \(\beta>1\), then
\[
 0\leq A_m\leq
 C_{\beta,\nu}\bigl(U^{\beta-1}P_m+P_m^\beta\bigr).
\]
Equations \eqref{eq:new26-foreign-local} and
\eqref{eq:new29-foreign-linear-local}, together with
\(\beta s<2^*\) and \(\beta\mathfrak q<2^*\), imply that, for every
compact \(K\subset\mathbb R^n\),
\begin{equation}\label{eq:new29-A-local}
 \|A_m\|_{L^s(K)}
 +\|A_m\|_{L^{\mathfrak q}(K)}
 \longrightarrow0.
\end{equation}

Fix \(\psi\in C_c^\infty(\mathbb R^n)\), put
\(K:=\operatorname{supp}\psi\), and assume
\(\|v_m\|_{\dot H^1}\leq1\).  Define
\[
 F:=\mathcal I_\alpha(U^\beta\psi),\qquad
 s'=\frac{2n}{\alpha}.
\]
The HLS inequality and the smooth compact support of \(\psi\) give
\[
 F\in L^{s'}(\mathbb R^n)\cap L^\infty_{\rm loc}(\mathbb R^n).
\]
Since \(1/\mathfrak q=1/t+1/s'\), for \(L>1\),
\begin{align*}
 \|A_mF\|_{L^{\mathfrak q}}
 &\leq
 \|F\|_{L^\infty(B_L)}
 \|A_m\|_{L^{\mathfrak q}(B_L)}
 +\|A_m\|_{L^t}\|F\|_{L^{s'}(B_L^c)}.
\end{align*}
First let \(m\to\infty\) and use
\eqref{eq:new29-A-global}--\eqref{eq:new29-A-local}; then let
\(L\to\infty\).  It follows that
\begin{equation}\label{eq:new29-AF}
 \|A_mF\|_{L^{\mathfrak q}}\longrightarrow0.
\end{equation}

The difference of the first lines of the two bilinear forms has the
exact expansion
\begin{align*}
 &\int_{\mathbb R^n}
 \left[
 \mathcal I_\alpha(T_m^\beta v_m)T_m^\beta
 -\mathcal I_\alpha(U^\beta v_m)U^\beta
 \right]\psi
 \notag\\
 =&
 \int_{\mathbb R^n}A_mv_mF
 +\int_{\mathbb R^n}A_mv_m\mathcal I_\alpha(A_m\psi)
 +\int_{\mathbb R^n}U^\beta v_m\mathcal I_\alpha(A_m\psi).
\end{align*}
By Sobolev duality and \eqref{eq:new29-AF}, the first term on the
right tends to zero.  For the other two terms, bilinear HLS and
H\"older give
\begin{align*}
 \left|\int A_mv_m\mathcal I_\alpha(A_m\psi)\right|
 \leq C\|A_m\|_{L^t}\|v_m\|_{L^{2^*}}
 \|A_m\psi\|_{L^s},\qquad
 \left|\int U^\beta v_m\mathcal I_\alpha(A_m\psi)\right|
 \leq C\|U^\beta\|_{L^t}\|v_m\|_{L^{2^*}}
 \|A_m\psi\|_{L^s}.
\end{align*}
Equation \eqref{eq:new29-A-local} implies
\(\|A_m\psi\|_{L^s}\to0\).  Hence
\begin{equation}\label{eq:new29-first-form-limit}
 \int
 \left[
 \mathcal I_\alpha(T_m^\beta v_m)T_m^\beta
 -\mathcal I_\alpha(U^\beta v_m)U^\beta
 \right]\psi\longrightarrow0.
\end{equation}

It remains to estimate directly the coefficient in the second line.
Set
\[
 C_m:=\mathcal I_\alpha(T_m^p)T_m^{\beta-1}-U^{2^*-2}.
\]
Since \(p=\beta+1>1\),
\begin{equation}\label{eq:new29-p-difference}
 0\leq T_m^p-U^p
 \leq C_{\alpha,\nu}\left(
 U^\beta P_m+\sum_{j=1}^{N-1}W_{j,m}^p\right).
\end{equation}
For every foreign profile \(W_{j,m}\), split \(\mathbb R^n\) into
\(\{U\geq W_{j,m}\}\) and \(\{U<W_{j,m}\}\), and apply
Lemma~\ref{A.4} with the two scaled exponents
\[
 a=\frac{n-2}{2}\beta s,\qquad
 b=\frac{n-2}{2}s,\qquad a+b=n.
\]
Interchanging \(a\) and \(b\) on the second region gives
\begin{equation*}
 \|U^\beta W_{j,m}\|_{L^s}\longrightarrow0.
\end{equation*}
The HLS inequality therefore implies
$
 \|\mathcal I_\alpha(U^\beta W_{j,m})\|_{L^{s'}}
 \longrightarrow0$.
After a further subsequence, these potentials converge to zero almost
everywhere.  Moreover,
\[
 \mathcal I_\alpha(W_{j,m}^p)=W_{j,m}^{\gamma}
 \longrightarrow0
 \quad\text{a.e. in }\mathbb R^n
\]
by \eqref{eq:new29-foreign-ae}.  It follows from
\eqref{eq:new29-p-difference} that
\[
 \mathcal I_\alpha(T_m^p)\longrightarrow
 \mathcal I_\alpha(U^p)=U^\gamma
 \quad\text{a.e. in }\mathbb R^n.
\]
Since \(T_m\to U>0\) almost everywhere,
\begin{equation}\label{eq:new29-second-coeff-ae}
 C_m\longrightarrow
 U^{\gamma+\beta-1}-U^{2^*-2}=0
 \quad\text{a.e. in }\mathbb R^n.
\end{equation}

For \(W_{0,m}:=U\), define the measurable partition
\[
 \Omega_{k,m}:=
 \left\{W_{k,m}=\max_{0\leq\ell\leq N-1}W_{\ell,m}\right\},
 \qquad 0\leq k\leq N-1,
\]
with ties assigned by the smallest index.  On \(\Omega_{k,m}\),
\[
 \mathcal I_\alpha(T_m^p)
 \leq C_{\alpha,\nu}\sum_{\ell=0}^{N-1}W_{\ell,m}^\gamma.
\]
If \(\beta\geq1\), then
\[
 T_m^{\beta-1}\leq C_\nu W_{k,m}^{\beta-1},
 \qquad
 W_{\ell,m}^\gamma W_{k,m}^{\beta-1}
 \leq W_{k,m}^{2^*-2}.
\]
If \(0<\beta<1\), then
\[
 T_m^{\beta-1}\leq W_{k,m}^{\beta-1},
 \qquad
 W_{\ell,m}^\gamma W_{k,m}^{\beta-1}
 \leq W_{\ell,m}^{\gamma+\beta-1}
 =W_{\ell,m}^{2^*-2}.
\]
Consequently, in both ranges of \(\beta\),
\begin{equation}\label{eq:new29-second-coeff-majorant}
 |C_m|^{\mathfrak q}
 \leq C_{\alpha,\nu}
 \sum_{\ell=0}^{N-1}
 W_{\ell,m}^{(2^*-2)\mathfrak q}.
\end{equation}
The term with \(\ell=0\) belongs to \(L^1(K)\), while every foreign
term tends to zero in \(L^1(K)\) by
\eqref{eq:new26-foreign-local}.  Hence the right-hand side of
\eqref{eq:new29-second-coeff-majorant} is uniformly integrable on
\(K\).  Equations \eqref{eq:new29-second-coeff-ae} and
\eqref{eq:new29-second-coeff-majorant}, followed by Vitali's theorem,
give
\begin{equation*}
 \|C_m\|_{L^{\mathfrak q}(K)}\longrightarrow0.
\end{equation*}
Thus Sobolev inequality yields
\begin{equation}\label{eq:new29-second-form-limit}
 \left|\int_{\mathbb R^n}C_mv_m\psi\right|
 \leq C_\psi\|C_m\|_{L^{\mathfrak q}(K)}
 \|v_m\|_{L^{2^*}}\longrightarrow0.
\end{equation}

Combining \eqref{eq:new29-first-form-limit} and
\eqref{eq:new29-second-form-limit}, with the fixed constants
\(D_{n,\alpha}p\) and \(D_{n,\alpha}\beta\), gives
\begin{equation*}
 \mathscr B_{T_m}(v_m,\psi)-\mathscr B_U(v_m,\psi)
 \longrightarrow0
 \qquad\bigl(\psi\in C_c^\infty(\mathbb R^n)\bigr).
\end{equation*}
The estimates used above are uniform for
\(\|v_m\|_{\dot H^1}\leq1\), and they apply without alteration in the
two ranges \(\beta\geq1\) and \(0<\beta<1\).

Let \(\psi\in\dot H^1\) and choose
\(\psi_\ell\in C_c^\infty(\mathbb R^n)\) with
\(\|\psi_\ell-\psi\|_{\dot H^1}\to0\).  The uniform form bound gives
\[
\begin{aligned}
 \left|
 \mathscr B_{T_m}(v_m,\psi)-\mathscr B_U(v_m,\psi)
 \right|\leq
 \left|
 \mathscr B_{T_m}(v_m,\psi_\ell)
 -\mathscr B_U(v_m,\psi_\ell)
 \right|
 +C\|v_m\|_{\dot H^1}\|\psi_\ell-\psi\|_{\dot H^1}.
\end{aligned}
\]
Taking first \(m\to\infty\) and then \(\ell\to\infty\) proves
\eqref{eq:test-linear-decoupling}.

If \eqref{eq:uniform-mode-decoupling} failed, there would be
\(\delta_m\downarrow0\), normalized functions
\(\|v_m\|_{\dot H^1}=1\), and indices \((j_m,b_m)\) for which its
left-hand side is bounded below by a positive constant.  Passing to a
subsequence makes \((j_m,b_m)\) constant.  Normalization about
\(V_{j_m,m}\) converts its mode into the corresponding fixed
one-bubble mode in \(\dot H^1\), and
\eqref{eq:test-linear-decoupling} gives a contradiction.  Hence
\eqref{eq:uniform-mode-decoupling} holds.
\end{proof}

\begin{proposition}\label{pro4.1.1}
Assume that $\nu\geq1$, $n\geq3$, and $0<\alpha<n$. There exist
positive constants $\delta_0=\delta_0(n,\nu,\alpha)$ and
$C=C(n,\nu,\alpha)$ such that, for every $0\leq\delta\leq\delta_0$
and every $\delta$-interacting bubble family
$\{U_i:1\leq i\leq\nu\}$, with $\sigma=\sum_{i=1}^\nu U_i$, if
$\phi$ solves the equation
\begin{equation}\label{7.2}
    \left\{\begin{array}{l}
\Delta \phi+D_{n,\alpha} \bigg(\int_{\mathbb{R}^{n}}\frac{p_{\alpha}\sigma^{p_{\alpha}-1}(y)\phi(y)   }{|x-y|^{\alpha}}\dy \sigma^{p_{\alpha}-1}(x) +\int_{\mathbb{R}^{n}}\frac{\sigma^{p_{\alpha}}(y)   }{|x-y|^{\alpha}}\dy (p_{\alpha}-1)\sigma^{p_{\alpha}-2}(x)\phi(x)   \bigg)
\\=g+\sum_{i=1}^\nu \sum_{a=1}^{n+1} c_a^i \bigg[\Big(|x|^{-\alpha}\ast |U_{i}|^{p_{\alpha}}\Big)(p_{\alpha}-1   )U_{i}^{ p_{\alpha}-2}Z_{i}^{a}+p_{\alpha}\Big(|x|^{-\alpha}\ast (U_{i}^{p_{\alpha}-1}
		Z_{i}^{a})\Big)U^{p_{\alpha}-1}_{i}\bigg]\quad \text{ in }\,\,\mathbb{R}^{n},  \\
\int_{\mathbb R^n} \bigg[\Big(|x|^{-\alpha}\ast |U_{i}|^{p_{\alpha}}\Big)(p_{\alpha}-1   )U_{i}^{ p_{\alpha}-2}Z_{i}^{a}+p_{\alpha}\Big(|x|^{-\alpha}\ast (U_{i}^{p_{\alpha}-1}
		Z_{i}^{a})\Big)U^{p_{\alpha}-1}_{i}\bigg] \phi\,\dx=0,
  \\\quad i=1, \ldots, \nu; \,\, a=1, \ldots, n+1
\end{array}\right.
\end{equation}
for some $g$ with $\|g\|_{\dot{H}^{-1}(\mathbb{R}^{n})}<\infty$. Then
\begin{equation}\label{eq:511}
    \|\phi\|_{\dot{H}^{1}(\mathbb{R}^{n}) } \leq C\|g\|_{\dot{H}^{-1}(\mathbb{R}^{n})}.
\end{equation}
\end{proposition}
\begin{proof}
\medskip \noindent\doublebox{\textsc{Step 1.}} We claim that there is a constant $C > 0$ depending only on $n$, $\nu$, and $\alpha$ such that
\[
 \varepsilon_{\rm int}(\delta)\longrightarrow0
 \qquad\text{as }\delta\longrightarrow0,
\]
where $\varepsilon_{\rm int}(\delta)$ denotes a quantity depending only on
$n,\nu,\alpha$ and the interaction bound $\delta$, and
\begin{equation}\label{eq:512}
\sum_{i=1}^{\nu} \sum_{a=1}^{n+1} |c_a^i| \le C\left(\|g\|_{\dot{H}^{-1}(\mathbb{R}^{n})} + \varepsilon_{\rm int}(\delta)\|\phi\|_{\dot{H}^{1}(\mathbb{R}^{n}) }\right).
\end{equation}

To show \eqref{eq:512}, test \eqref{7.2} against $Z_j^b$ in the
$\dot H^{-1}$--$\dot H^1$ duality.  This gives
\begin{equation}\label{re4.2.1}
    \begin{split}
       &\int D_{n,\alpha} \bigg(\int_{\mathbb{R}^{n}}\frac{p_{\alpha}\sigma^{p_{\alpha}-1}(y)\phi (y)   }{|x-y|^{\alpha}}\dy \sigma^{p_{\alpha}-1}(x) +\int_{\mathbb{R}^{n}}\frac{\sigma^{p_{\alpha}}(y)   }{|x-y|^{\alpha}}\dy (p_{\alpha}-1)\sigma^{p_{\alpha}-2}(x) \phi (x)   \bigg)Z_j^b
       \\=&\langle g,Z_j^b\rangle_{\dot H^{-1},\dot H^1}+\sum_{i, a} \int c_a^i \bigg[\Big(|x|^{-\alpha}\ast |U_{i}|^{p_{\alpha}}\Big)(p_{\alpha}-1   )U_{i}^{ p_{\alpha}-2}Z_{i}^{a}+p_{\alpha}\Big(|x|^{-\alpha}\ast (U_{i}^{p_{\alpha}-1}
		Z_{i}^{a})\Big)U^{p_{\alpha}-1}_{i}\bigg] Z_j^b
    \end{split}
\end{equation}
for any $1 \leq j \leq \nu, 1 \leq b \leq n+1$. Here we used the orthogonality condition in \eqref{7.2}.
By Lemma \ref{A.5} and the identity $-\Delta Z_i^a=D_{n,\alpha}\Psi_i^a$, for $a,b\leq n+1$ there exist constants $\gamma^b>0$ such that
\begin{equation*}
    \begin{split}
        \sum_{i, a} \int c_a^i \bigg[\Big(|x|^{-\alpha}\ast |U_{i}|^{p_{\alpha}}\Big)(p_{\alpha}-1   )U_{i}^{ p_{\alpha}-2}Z_{i}^{a}+p_{\alpha}\Big(|x|^{-\alpha}\ast (U_{i}^{p_{\alpha}-1}
		Z_{i}^{a})\Big)U^{p_{\alpha}-1}_{i}\bigg] Z_j^b
  =c_b^j \gamma^b+\sum_{i \neq j} \sum_{a=1}^{n+1} c_a^i O\left(q_{i j}\right).
    \end{split}
\end{equation*}
Plugging the above estimates into \eqref{re4.2.1}, we see that $\left\{c_b^j\right\}$ satisfies the linear system
\begin{equation*}
    \begin{split}
    &c_b^j \gamma^b+\sum_{i \neq j} \sum_{a=1}^{n+1} c_a^i O\left(q_{i j}\right)\\=
    &-\langle g,Z_j^b\rangle_{\dot H^{-1},\dot H^1}
    +  \int D_{n,\alpha} \bigg(\int_{\mathbb{R}^{n}}\frac{p_{\alpha}\sigma^{p_{\alpha}-1}(y)\phi (y)   }{|x-y|^{\alpha}}\dy \sigma^{p_{\alpha}-1}(x) +\int_{\mathbb{R}^{n}}\frac{\sigma^{p_{\alpha}}(y)   }{|x-y|^{\alpha}}\dy (p_{\alpha}-1)\sigma^{p_{\alpha}-2}(x) \phi (x)   \bigg)Z_j^b.
    \end{split}
\end{equation*}
Denote $\vec{c}^j=(c_1^j,\ldots,c_{n+1}^j)$ and concatenate these
vectors into $\vec c\in\mathbb R^{\nu(n+1)}$.  The coefficient matrix has
the form
\[
 \mathsf G=\mathsf G_0+\mathsf E,\qquad
 \mathsf G_0=\operatorname{diag}
 \bigl(\gamma^1,\ldots,\gamma^{n+1},
       \ldots,\gamma^1,\ldots,\gamma^{n+1}\bigr),
\]
where
\[
 \|\mathsf E\|_{\ell^2\to\ell^2}
 \leq C_{n,\alpha,\nu}\max_{i\ne j}q_{ij}
 \leq C_{n,\alpha,\nu}\delta .
\]
Let $\gamma_*=\min_{1\leq b\leq n+1}\gamma^b>0$ and choose
$\delta_0$ so that
$C_{n,\alpha,\nu}\delta_0\leq\gamma_*/2$.  Then
\[
 \|\mathsf G_0^{-1}\mathsf E\|_{\ell^2\to\ell^2}\leq\frac12,
 \qquad
 \mathsf G^{-1}
 =\sum_{k=0}^{\infty}(-\mathsf G_0^{-1}\mathsf E)^k
 \mathsf G_0^{-1},
\]
 and therefore
 $\|\mathsf G^{-1}\|_{\ell^2\to\ell^2}\leq2/\gamma_*$.
The finite-dimensional system \eqref{re4.2.1} consequently gives the
projection-sensitive estimate
\begin{equation}\label{eq:multiplier-projection-bound}
 \sum_{i=1}^{\nu}\sum_{a=1}^{n+1}|c_a^i|
 \leq C\sum_{j=1}^{\nu}\sum_{b=1}^{n+1}
 \left(
 \left|\langle g,Z_j^b\rangle_{\dot H^{-1},\dot H^1}\right|
 +\left|\int_{\mathbb R^n}
 N_1(\phi,\sigma)Z_j^b\,\dx\right|
 \right).
\end{equation}

 For each $j$ and $b$, the left-hand side produced by testing
 \eqref{7.2} with $Z_j^b$ is the bilinear form
$\mathscr B_\sigma(\phi,Z_j^b)$.  By definition,
\[
 \mathscr B_\sigma(\phi,Z_j^b)-\mathscr B_{U_j}(\phi,Z_j^b)
 =\int_{\mathbb R^n}\bigl(N_1(\phi,\sigma)
       -N_1(\phi,U_j)\bigr)Z_j^b\,\dx .
\]
The orthogonality condition and $-\Delta Z_j^b=N_1(Z_j^b,U_j)$ imply
\[
 \mathscr B_{U_j}(\phi,Z_j^b)
 =\int_{\mathbb R^n}\phi N_1(Z_j^b,U_j)\,\dx
 =\int_{\mathbb R^n}\nabla\phi\cdot\nabla Z_j^b\,\dx=0.
\]
Consequently Lemma \ref{lem:linearized-form-decoupling} gives, for the
entire range $0<\alpha<n$,
\[
 \left|\mathscr B_\sigma(\phi,Z_j^b)
 -\mathscr B_{U_j}(\phi,Z_j^b)\right|
 \leq \omega_4(\delta)\|\phi\|_{\dot H^1(\mathbb R^n)}.
\]
Furthermore,
\[
 \left|\langle g,Z_j^b\rangle_{\dot H^{-1},\dot H^1}\right|
 \leq\|g\|_{\dot H^{-1}}\|Z_j^b\|_{\dot H^1}
 \leq C\|g\|_{\dot H^{-1}}.
\]
Substitution in the finite system proves \eqref{eq:512}.

\medskip \noindent\doublebox{\textsc{Step 2.}} We assert that \eqref{eq:511} holds.
Suppose not. There exist sequences of small positive numbers $\{\delta'_m\}_{m \in \N}$,
$\delta'_m$-interacting families $\{\{U_{i,m} = U[z_{i,m},\lambda_{i,m}]\}_{i = 1,\ldots,\nu}\}_{m \in \N}$, functions $\{g_m\}_{m \in \N} \subset \dot{H}^{-1}(\mathbb{R}^{n})$ and $\{\phi_{m}\}_{m \in \N} \subset \dot{H}^1(\R^n)$,
 and numbers $\{c_{a,m}^i\}_{i = 1,\ldots,\nu,\, a = 1,\ldots, n+1,\, m \in \N}$ satisfying $\|\phi_{m}\|_{\dot{H}^{1}(\mathbb{R}^{n})}=1 $, $\|g_m\|_{  \dot{H}^{-1}(\mathbb{R}^{n})}=o(1) $ and \eqref{7.2}. By \eqref{eq:512},
\begin{equation}\label{eq:513}
\sum_{i=1}^{\nu} \sum_{a=1}^{n+1} \left|c_{a,m}^i\right| \to 0, \qquad \text{as } \,\, m \to \infty.
\end{equation}
Testing \eqref{7.2} with $\phi_{m}$ and using H\"older inequality, Sobolev inequality, Hardy--Littlewood--Sobolev inequality, and \eqref{eq:513}, we obtain
\begin{equation}\label{reeq:513}
    \begin{split}
        \int_{\R^n} N_{1}( \phi_{m},\sigma_m )\phi_{m} = \|\phi_{m}\|_{\dot{H}^1(\R^n)}^2 + O\left(\|g_m\|_{\dot{H}^{-1}(\mathbb{R}^{n})} + \sum_{i=1}^{\nu} \sum_{a=1}^{n+1} \left|c_{a,m}^i\right|\right) \to 1, \qquad \text{as } \,\, m \to \infty.
    \end{split}
\end{equation}
Let
\[\hphi_{i,m}(y) =  \lambda_{i,m}^{-\frac{n-2}{2}} \phi_m\big(\lambda_{i,m}^{-1}y+z_{i,m}\big) \quad \text{for } y \in \R^n .\]
We verify $\hphi_{i,m}(y) \to 0 \,\, \text{ weakly in } \dot{H}^{1}(\R^n) $.

For $\varphi\in C_c^\infty(\mathbb R^n)$ define the rescaled
distribution $\widehat g_{i,m}\in\dot H^{-1}(\mathbb R^n)$ by
\[
 \langle\widehat g_{i,m},\varphi\rangle_{\dot H^{-1},\dot H^1}
 :=\left\langle g_m,
 \lambda_{i,m}^{\frac{n-2}{2}}
 \varphi\bigl(\lambda_{i,m}(\,\cdot-z_{i,m})\bigr)
 \right\rangle_{\dot H^{-1},\dot H^1}.
\]
The critical dilation is an isometry of $\dot H^1$, so
$\|\widehat g_{i,m}\|_{\dot H^{-1}}=\|g_m\|_{\dot H^{-1}}$.
By \eqref{7.2},
\begin{equation}\label{eq:3498}
\begin{split}
    (-\Delta) \hphi_{i,m} -   N_{1}(\hphi_{i,m},\hat{\sigma})
=  -\left[\widehat g_{i,m} + \sum_{j=1}^{\nu} \sum_{a=1}^{n+1}
c_{a,m}^jN_{2}(\hat{U}_{j,i,m},\hat{Z}_{j,m}^a ) \right] \quad \text{in } \R^n,
\end{split}
\end{equation}
where $\hat{U}_{j,i,m}(y) = \lambda_{i,m}^{-\frac{n-2}{2}} U_{j,m}(\lambda_{i,m}^{-1}y+z_{i,m}\big) $, $\hat{Z}_{j,m}^a(y) = \lambda_{i,m}^{-\frac{n-2}{2}}Z_{j,m}^a(\lambda_{i,m}^{-1}y+z_{i,m}\big)$, $\hat{\sigma}(y) =\sum_{j=1}^{\nu}\hat{U}_{j,i,m}(y) $ and
\begin{align*}
    N_{2}(\eta,\xi)=\bigg[\Big(|x|^{-\alpha}\ast |\eta|^{p_{\alpha}}\Big)(p_{\alpha}-1   )|\eta|^{ p_{\alpha}-2}\xi+p_{\alpha}\Big(|x|^{-\alpha}\ast (|\eta|^{p_{\alpha}-2}\eta
		\xi)\Big)|\eta|^{p_{\alpha}-1}\bigg].
\end{align*}
In addition,
\begin{equation*}
\begin{aligned}
\|\hphi_{i,m}\|_{\dot{H}^{1}(\R^n)}
= \|\phi_m\|_{\dot{H}^{1}(\R^n)} =1.\\
\end{aligned}
\end{equation*}
 Therefore, we may assume that
\[\hphi_{i,m} \rightharpoonup \hphi_{i,\infty} \quad \text{weakly in } \dot{H}^{1}(\R^n)
\quad \text{ and } \quad \hphi_{i,m} \to \hphi_{i,\infty} \quad \text{a.e.} \quad \text{as } m \to \infty\]
for some $\hphi_{i,\infty} \in \dot{H}^{1}(\R^n)$.
Let $\psi\in C_c^\infty(\mathbb R^n)$ be a test function. Put
$U:=U[0,1]$, $v_m:=\hphi_{i,m}$, and $S_m:=\hat\sigma$.  The exact
difference of the two linearized operators is
\begin{align*}
&N_1(v_m,S_m)-N_1(v_m,U)\notag\\
=&D_{n,\alpha}p_\alpha
\mathcal I_\alpha^0\big((S_m^{p_\alpha-1}-U^{p_\alpha-1})v_m\big)
S_m^{p_\alpha-1}+D_{n,\alpha}p_\alpha
\mathcal I_\alpha^0(U^{p_\alpha-1}v_m)
(S_m^{p_\alpha-1}-U^{p_\alpha-1})\notag\\
&+D_{n,\alpha}(p_\alpha-1)
\mathcal I_\alpha^0(S_m^{p_\alpha}-U^{p_\alpha})
S_m^{p_\alpha-2}v_m\notag+D_{n,\alpha}(p_\alpha-1)
\mathcal I_\alpha^0(U^{p_\alpha})
(S_m^{p_\alpha-2}-U^{p_\alpha-2})v_m,
\end{align*}
where
$
 \mathcal I_\alpha^0F(y):=\int_{\mathbb R^n}
 \frac{F(x)}{|x-y|^\alpha}\dx$ .
By the definition of $\mathscr B$,
\begin{equation*}
 \int_{\mathbb R^n}
 \bigl(N_1(v_m,S_m)-N_1(v_m,U)\bigr)\psi\,dy
 =\mathscr B_{S_m}(v_m,\psi)-\mathscr B_U(v_m,\psi)
 \longrightarrow0.
\end{equation*}
Here the convergence is precisely
\eqref{eq:test-linear-decoupling}.  In particular, for $p_\alpha<2$ it
follows from the dominance-region coefficient estimate in the proof of
Lemma \ref{lem:linearized-form-decoupling}, in which
\(\mathcal I_\alpha(S_m^{p_\alpha})S_m^{p_\alpha-2}\) is kept as one
factor.  The single-bubble linearized operator is symmetric,
so weak convergence of $v_m$ gives
\[
 \int N_1(v_m,U)\psi\,dy
 =\int v_mN_1(\psi,U)\,dy
 \longrightarrow
 \int N_1(\hphi_{i,\infty},U)\psi\,dy.
\]
Moreover,
\[
 \left|\langle\widehat g_{i,m},\psi\rangle_{\dot H^{-1},\dot H^1}\right|
 \leq\|\widehat g_{i,m}\|_{\dot H^{-1}}\|\psi\|_{\dot H^1}
 =\|g_m\|_{\dot H^{-1}}\|\psi\|_{\dot H^1}=o(1).
\]

For every $j,a,m$, differentiation of the rescaled one-bubble equation
and conformal scale invariance give the exact identity
\[
 D_{n,\alpha}N_2(\hat U_{j,i,m},\hat Z_{j,m}^a)
 =-\Delta\hat Z_{j,m}^a
 \quad\hbox{in }\dot H^{-1}(\mathbb R^n).
\]
Since the $\dot H^1$ norm of every translation or scale mode is
invariant under translations and critical dilations,
\[
 \sup_{i,j,a,m}
 \|N_2(\hat U_{j,i,m},\hat Z_{j,m}^a)\|_{\dot H^{-1}}
 \leq D_{n,\alpha}^{-1}
 \sup_{j,a,m}\|\hat Z_{j,m}^a\|_{\dot H^1}\leq C.
\]
Consequently, for every $\psi\in\dot H^1(\mathbb R^n)$,
\[
 \left|\sum_{j=1}^{\nu}\sum_{a=1}^{n+1}c_{a,m}^j
 \left\langle N_2(\hat U_{j,i,m},\hat Z_{j,m}^a),\psi
 \right\rangle_{\dot H^{-1},\dot H^1}\right|
 \leq C\sum_{j,a}|c_{a,m}^j|\,\|\psi\|_{\dot H^1}=o(1)\|\psi\|_{\dot H^1}
\]
by \eqref{eq:513}.  Thus the multiplier sources vanish in the weak
limit.

Accordingly, by taking the limit $m \to \infty$ in \eqref{eq:3498}, we find that
\begin{equation}\label{eq:349b}
(-\Delta) \hphi_{i,\infty} - N_{1}(\hphi_{i,\infty},U[0,1]) = 0 \quad \text{in } \R^n.
\end{equation}
Write $\Psi_{i,m}^a$ for the bracketed coefficient in the
orthogonality condition in \eqref{7.2}, with $U_i=U_{i,m}$.
The corresponding condition is invariant under the preceding rescaling:
\[
 0=\int_{\mathbb R^n}\Psi_{i,m}^a\phi_m\,\dx
 =\int_{\mathbb R^n}\Psi_U^a\hphi_{i,m}\,\dy .
\]
Passing to the weak limit gives
\[
 \int_{\mathbb R^n}\Psi_U^a\hphi_{i,\infty}\,dy=0,
 \qquad a=1,\ldots,n+1,
\]
where $D_{n,\alpha}\Psi_U^a=-\Delta Z^a$.  Thus
$\hphi_{i,\infty}$ is $\dot H^1$-orthogonal to every single-bubble kernel
direction.  By \eqref{eq:349b} and Lemma \ref{lenondegeneracy} (see also
\cite[Theorem 1.5]{Li-Liu-Tang-Xu}), it belongs to the span of those
directions.  Their Gram matrix is positive definite, and hence
$\hphi_{i,\infty}=0$ in $\mathbb R^n$.

We now derive a contradiction to \eqref{reeq:513}.  For
$x\in\mathbb R^n$, define
\[
 y_{i,m}(x):=\lambda_{i,m}(x-z_{i,m}),
 \qquad i=1,\ldots,\nu.
\]
For $L>0$, let
$\Omega_{(m),L}:=\bigcup_{i=1}^{\nu}\{x:|y_{i,m}(x)|\leq L\}$ and let
$\Omega_{(m),L}^c$ denote its complement. Then
\begin{equation*}
    \begin{split}
&\int_{\mathbb{R}^{n}}\int_{\mathbb{R}^{n}}\frac{p_{\alpha}\sigma_m^{p_{\alpha}-1}(y)\phi_{m}(y)   }{|x-y|^{\alpha}}\dy\sigma_m^{p_{\alpha}-1}(x)\phi_{m}(x) \dx
        \\= & \int_{\Omega_{(m),L}}\int_{\mathbb{R}^{n}}\frac{p_{\alpha}\sigma_m^{p_{\alpha}-1}(y)\phi_{m}(y)   }{|x-y|^{\alpha}}\dy\sigma_m^{p_{\alpha}-1}(x)\phi_{m}(x) \dx+\int_{\Omega_{(m),L}^{c}}\int_{\mathbb{R}^{n}}\frac{p_{\alpha}\sigma_m^{p_{\alpha}-1}(y)\phi_{m}(y)   }{|x-y|^{\alpha}}\dy\sigma_m^{p_{\alpha}-1}(x)\phi_{m}(x) \dx
        \\\lesssim &\|\sigma_m^{p_{\alpha}-1}(y)\phi_{m}(y)  \|_{L^{\frac{2^{*}}{p_{\alpha}}}(\mathbb{R}^{n}) }\|\sigma_m^{p_{\alpha}-1}(x)\phi_{m}(x)  \|_{L^{\frac{2^{*}}{p_{\alpha}}}(\Omega_{(m),L}) }
        +\|\sigma_m^{p_{\alpha}-1}(y)\phi_{m}(y)  \|_{L^{\frac{2^{*}}{p_{\alpha}}}(\mathbb{R}^{n}) }\|\sigma_m^{p_{\alpha}-1}(x)\phi_{m}(x)  \|_{L^{\frac{2^{*}}{p_{\alpha}}}(\Omega_{(m),L}^{c}) }
        \\\lesssim& o(1)+\frac{1}{L^{\frac{n(p_{\alpha}-1  )}{2^*}}},
    \end{split}
\end{equation*}
For the second part of $N_1$, the single-bubble identity and H\"older's
inequality give
\begin{align*}
&\int_{\mathbb R^n}\int_{\mathbb R^n}
 \frac{\sigma_m^{p_\alpha}(y)}{|x-y|^\alpha}\dy\,
 (p_\alpha-1)\sigma_m^{p_\alpha-2}(x)\phi_m^2(x)\dx\notag \\
 \lesssim&\int_{\mathbb R^n}\sigma_m^{2^*-2}\phi_m^2\dx
 =\int_{\Omega_{(m),L}}\sigma_m^{2^*-2}\phi_m^2\dx
 +\int_{\Omega_{(m),L}^c}\sigma_m^{2^*-2}\phi_m^2\dx
 \notag\\
\leq &o_m(1)
 +C\left(\int_{\Omega_{(m),L}^c}\sigma_m^{2^*}\dx\right)^{2/n}
 \|\phi_m\|_{L^{2^*}}^2
 \leq o_m(1)+CL^{-2}.
\end{align*}
This yields
\begin{equation*}
    \begin{split}
        \lim_{m \to \infty} \int_{\R^n} N_{1}( \phi_{m},\sigma_m )\phi_{m} = 0.
    \end{split}
\end{equation*} Because $L > 0$ can be taken arbitrarily large.
This contradicts \eqref{reeq:513} and proves Proposition \ref{pro4.1.1}.
\end{proof}

\begin{proposition}\label{pro4.2.1}
Suppose that $\nu\geq1$, $n\geq3$, and $0<\alpha<n$. There exist
positive constants $\delta_0=\delta_0(n,\nu,\alpha)$ and
$C=C(n,\nu,\alpha)$ such that, for every
$0\leq\delta\leq\delta_0$, every $\delta$-interacting bubble family
$\{U_i:1\leq i\leq\nu\}$, and every
$g\in\dot H^{-1}(\mathbb R^n)$, problem \eqref{3.1}, with
$\sigma=\sum_{i=1}^\nu U_i$, has a unique solution
$\phi=\mathcal{L}(g)$ and uniquely determined multipliers
$c_a^i$. Moreover, for $1\leq i\leq\nu$ and $1\leq a\leq n+1$,
    \begin{align*}
     \|\mathcal{L}(g)\|_{\dot{H}^{1}(\mathbb{R}^{n})} \leqslant C\|g\|_{\dot{H}^{-1}(\mathbb{R}^{n})}, \quad\left|c_a^i\right| \leqslant C \|g\|_{\dot{H}^{-1}(\mathbb{R}^{n})}
    \end{align*}
\end{proposition}
\begin{proof}
For $1\leq i\leq\nu$ and $1\leq a\leq n+1$, set
\[
\Psi_i^a:=
\Big(|x|^{-\alpha}\ast U_i^{p_\alpha}\Big)(p_\alpha-1)
U_i^{p_\alpha-2}Z_i^a
+p_\alpha\Big(|x|^{-\alpha}\ast
(U_i^{p_\alpha-1}Z_i^a)\Big)U_i^{p_\alpha-1},
\]
and define the closed subspace
\[
\mathcal H_\sigma:=
\left\{\phi\in\dot H^1(\mathbb R^n):
\int_{\mathbb R^n}\Psi_i^a\phi\,\dx=0
\text{ for every }i,a\right\}.
\]
For $\phi,\psi\in\mathcal H_\sigma$, let
\[
\mathcal B_\sigma(\phi,\psi):=
\int_{\mathbb R^n}\nabla\phi\cdot\nabla\psi\,\dx
-\int_{\mathbb R^n}N_1(\phi,\sigma)\psi\,\dx .
\]
The Hardy--Littlewood--Sobolev and Sobolev inequalities give
\begin{equation*}
 |\mathcal B_\sigma(\phi,\psi)|
 \leq C\|\phi\|_{\dot H^1}\|\psi\|_{\dot H^1}.
\end{equation*}
Let $\mathcal K_\sigma:\mathcal H_\sigma\to\mathcal H_\sigma$ be defined by
\[
\int\nabla(\mathcal K_\sigma\phi)\cdot\nabla\psi\,\dx
=\int N_1(\phi,\sigma)\psi\,\dx
\qquad(\psi\in\mathcal H_\sigma).
\]
This operator is compact.  Put
\[
 a_\sigma:=\sigma^{p_\alpha-1},\qquad
 b_\sigma:=(p_\alpha-1)\mathcal I_\alpha(\sigma^{p_\alpha})
 \sigma^{p_\alpha-2},\qquad
 t_\alpha:=\frac{2n}{n+2-\alpha}.
\]
Then $a_\sigma\in L^{t_\alpha}$ and $b_\sigma\in L^{n/2}$.  Indeed,
on each $\Omega_i$, the dominance inequalities and the one-bubble
potential identity give
\[
 a_\sigma\leq C\sum_{i=1}^{\nu}U_i^{p_\alpha-1},\qquad
 0\leq b_\sigma\leq C\sum_{i=1}^{\nu}U_i^{2^*-2}.
\]
Choose $\chi\in C_c^\infty(\mathbb R^n)$ with
$0\leq\chi\leq1$, $\chi=1$ on $B_1$, and $\chi=0$ outside $B_2$, and
set
\[
 \chi_M(x):=1-\prod_{i=1}^{\nu}
 \left(1-\chi\left(\frac{\lambda_i(x-z_i)}{M}\right)\right).
\]
Define $\mathcal K_{\sigma,M}$ by
\[
\begin{aligned}
 \int\nabla(\mathcal K_{\sigma,M}\phi)\cdot\nabla\psi\,\dx
 :={}&p_\alpha\int
 \mathcal I_\alpha(\chi_Ma_\sigma\phi)\,
 \chi_Ma_\sigma\psi\,\dx
 +\int\chi_Mb_\sigma\phi\psi\,\dx .
\end{aligned}
\]
For fixed $M$, choose bounded, compactly supported functions
$a_{M,j}$ and $b_{M,j}$ such that
\[
 a_{M,j}\longrightarrow\chi_Ma_\sigma\quad\text{in }L^{t_\alpha},
 \qquad
 b_{M,j}\longrightarrow\chi_Mb_\sigma\quad\text{in }L^{n/2}.
\]
Let $\mathcal K_{M,j}$ be defined by
\[
 \int\nabla(\mathcal K_{M,j}\phi)\cdot\nabla\psi\,\dx
 :=p_\alpha\int\mathcal I_\alpha(a_{M,j}\phi)a_{M,j}\psi\,\dx
 +\int b_{M,j}\phi\psi\,\dx .
\]
We verify the compactness.  Suppose that
$\phi_m\rightharpoonup\phi$ in $\mathcal H_\sigma$, and put
\[
 K_{M,j}:=\operatorname{supp}a_{M,j}
 \cup\operatorname{supp}b_{M,j},
 \qquad
 s_\alpha:=\frac{2n}{2n-\alpha}<2.
\]
Choose a ball $B_{M,j}$ containing $K_{M,j}$.  The restrictions of
$(\phi_m)$ are bounded in $H^1(B_{M,j})$.  Indeed, the Sobolev
inequality and the finite measure of \(B_{M,j}\) give
\[
 \|\phi_m\|_{L^2(B_{M,j})}
 \leq |B_{M,j}|^{1/n}\|\phi_m\|_{L^{2^*}(\mathbb R^n)}
 \leq C_{M,j}\|\phi_m\|_{\dot H^1}.
\]
Together with
\(\|\nabla\phi_m\|_{L^2(B_{M,j})}\leq\|\phi_m\|_{\dot H^1}\),
this proves the asserted local \(H^1\)-boundedness.
After passage to a subsequence, the local Rellich theorem gives
\[
 \phi_m\longrightarrow\phi
 \qquad\text{in }L^2(B_{M,j}).
\]
Since $K_{M,j}$ has finite measure and $a_{M,j},b_{M,j}$ are bounded,
\[
 a_{M,j}(\phi_m-\phi)\longrightarrow0
 \qquad\text{in }L^{s_\alpha}(\mathbb R^n).
\]
For $\psi\in\mathcal H_\sigma$ with
$\|\psi\|_{\dot H^1}\leq1$, the bilinear HLS inequality gives
\begin{align*}
 \left|\int
 \mathcal I_\alpha\bigl(a_{M,j}(\phi_m-\phi)\bigr)
 a_{M,j}\psi\,\dx\right|\leq
 C\|a_{M,j}(\phi_m-\phi)\|_{L^{s_\alpha}}
 \|a_{M,j}\psi\|_{L^{s_\alpha}}
 \longrightarrow0,
\end{align*}
uniformly in such $\psi$.  Moreover, local Sobolev embedding yields
\[
 \left|\int b_{M,j}(\phi_m-\phi)\psi\,\dx\right|
 \leq C\|\phi_m-\phi\|_{L^2(B_{M,j})}
 \|\psi\|_{L^2(B_{M,j})}
 \longrightarrow0.
\]
Taking the supremum over the unit ball of $\mathcal H_\sigma$ proves
\[
 \|\mathcal K_{M,j}(\phi_m-\phi)\|_{\dot H^1}\longrightarrow0.
\]
Thus every $\mathcal K_{M,j}$ is compact.  In addition,
the HLS, H\"older, and Sobolev inequalities give
\[
\begin{aligned}
 \|\mathcal K_{M,j}-\mathcal K_{\sigma,M}\|_
 {\mathcal L(\mathcal H_\sigma)}
 \leq C\Big[&
 \bigl(\|a_{M,j}\|_{L^{t_\alpha}}
 +\|\chi_Ma_\sigma\|_{L^{t_\alpha}}\bigr)
 \|a_{M,j}-\chi_Ma_\sigma\|_{L^{t_\alpha}}\\
 &+\|b_{M,j}-\chi_Mb_\sigma\|_{L^{n/2}}\Big]\longrightarrow0.
\end{aligned}
\]
Thus $\mathcal K_{\sigma,M}$ is compact.  Applying HLS, H\"older, and
Sobolev to the two tail coefficients yields
\begin{align}
 \|\mathcal K_\sigma-\mathcal K_{\sigma,M}\|_
 {\mathcal L(\mathcal H_\sigma)}
 &\leq C\Bigl(
 \|a_\sigma\|_{L^{t_\alpha}}
 \|(1-\chi_M)a_\sigma\|_{L^{t_\alpha}}
 +\|(1-\chi_M)b_\sigma\|_{L^{n/2}}
 \Bigr)\longrightarrow0.                         \label{eq:fredholm-tail}
\end{align}
Thus $\mathcal K_\sigma$ is the operator-norm limit of compact
operators.

The weak form of the projected equation is
\begin{equation}\label{eq:fredholm-equation}
 (I-\mathcal K_\sigma)\phi=\widetilde g
 \quad\text{in }\mathcal H_\sigma,\qquad
 \int\nabla\widetilde g\cdot\nabla\psi\,\dx=-\langle g,\psi\rangle ,
\end{equation}
and $\|\widetilde g\|_{\dot H^1}\leq\|g\|_{\dot H^{-1}}$.
The operator $I-\mathcal K_\sigma$ is Fredholm of index zero by
\eqref{eq:fredholm-tail}.  If
$(I-\mathcal K_\sigma)\phi=0$, then $\phi$ solves \eqref{7.2} with
$g=0$ for suitable multipliers.  Proposition \ref{pro4.1.1} gives
\[
 \|\phi\|_{\dot H^1}\leq C\|g\|_{\dot H^{-1}}=0,
\]
so the kernel is trivial.  Fredholm's alternative therefore gives a
unique solution of \eqref{eq:fredholm-equation}.  Applying
Proposition \ref{pro4.1.1} to this solution yields
\[
 \|\phi\|_{\dot H^1}\leq C\|g\|_{\dot H^{-1}}.
\]

The residual of the unprojected equation annihilates
$\mathcal H_\sigma$, hence it belongs to
$\operatorname{span}\{\Psi_i^a\}_{i,a}$ and determines the multipliers
$c_a^i$.  Estimate \eqref{eq:512} and the preceding bound give, after
choosing $\delta_0$ so that
$\varepsilon_{\rm int}(\delta_0)\leq1$,
\[
 \sum_{i=1}^{\nu}\sum_{a=1}^{n+1}|c_a^i|
 \leq C\left(\|g\|_{\dot H^{-1}}
 +\varepsilon_{\rm int}(\delta)
 \|\phi\|_{\dot H^1}\right)
 \leq C\|g\|_{\dot H^{-1}}.
\]
This proves existence, uniqueness, and both asserted estimates.
\end{proof}

\begin{lemma}\label{le:strong-nonlinear}
Assume that $\nu\geq1$, $4<\alpha<n$, and
$\sigma=\sum_{i=1}^\nu U_i$ is a sum of Aubin--Talenti bubbles. Set
\[
 \eta:=\min\{p_\alpha-1,2^*-2\}>0.
\]
There exists $C=C(n,\nu,\alpha)>0$ such that, whenever
$\|\rho\|_{\dot H^1}$ is
sufficiently small,
\[
    \|N(\rho,\sigma)\|_{\dot H^{-1}}\leq C\|\rho\|_{\dot H^1}^{1+\eta}.
\]
\end{lemma}

\begin{rmk}
Assume that \(4<\alpha<n\).  Then \(1<p_\alpha<2\), and the nonlinear
remainder satisfies only a H\"older-type continuity estimate in the
\(\dot H^1\)-topology.  More precisely, if
\[
 \varepsilon_H:=\|\rho_1\|_{\dot H^1}+\|\rho_2\|_{\dot H^1}\leq1,
 \qquad d_H:=\|\rho_1-\rho_2\|_{\dot H^1},
\]
then
\begin{equation}\label{eq:strong-N-mixed-continuity}
\|N(\rho_1,\sigma)-N(\rho_2,\sigma)\|_{\dot H^{-1}}
\leq C\left(
 \varepsilon_H^\eta d_H
 +\varepsilon_H^{2-p_\alpha+\eta}d_H^{p_\alpha-1}
\right).
\end{equation}
Since \(0<p_\alpha-1<1\), the second term in
\eqref{eq:strong-N-mixed-continuity} does not yield a local Lipschitz
estimate for \(N(\,\cdot\,,\sigma)\).  Consequently, for the fixed-point
map
\[
 \mathcal A_R(\rho)
 :=-\mathcal L_R\bigl(h+N(\rho,\sigma)\bigr)
\]
used in Proposition~\ref{prop:sharp-correction}, this estimate does not
imply a contraction bound on \(\mathcal E_R\), and hence the Banach
fixed-point theorem is not directly applicable.  The proof of
Proposition~\ref{prop:sharp-correction} instead establishes the weak
continuity of \(\mathcal A_R\) on the weakly compact convex set
\(\mathcal E_R\) and applies the Schauder--Tychonoff fixed-point theorem.
\end{rmk}

\begin{proof}[Proof of Lemma~\ref{le:strong-nonlinear}]
Put
\[
 p:=p_\alpha=\frac{2n-\alpha}{n-2},\qquad
 \gamma:=2^*-p=\frac{\alpha}{n-2}.
\]
The assumptions imply
$
 1<p<2$, $\gamma>0$.
For
\[
 r:=\frac{2n}{n+2},\qquad
 s:=\frac{2n}{2n-\alpha},\qquad
 t:=\frac{2n}{n+2-\alpha},\qquad
 q:=\frac{2n}{\alpha},
\]
one has
\[
 1<s<\frac{n}{n-\alpha},\qquad 1<q<\infty,
\]
and
\[
 \frac1q=\frac1s-\frac{n-\alpha}{n},\qquad
 \frac1r=\frac1q+\frac1t,\qquad
 ps=2^*,\qquad (p-1)t=2^*.
\]
Define
\[
 \mathcal I_\alpha F(x):=
 D_{n,\alpha}\int_{\mathbb R^n}\frac{F(y)}{|x-y|^\alpha}\,dy.
\]
The Hardy--Littlewood--Sobolev and H\"older inequalities give
\begin{equation}\label{eq:strong-N-HLS-complete}
 \|\mathcal I_\alpha(F)G\|_{L^r}
 \leq \|\mathcal I_\alpha(F)\|_{L^q}\|G\|_{L^t}
 \leq C\|F\|_{L^s}\|G\|_{L^t}.
\end{equation}

Set
\[
 E_1:=\{z\in\mathbb R^n:\sigma(z)>2|\rho(z)|\},\qquad
 E_2:=\{z\in\mathbb R^n:\sigma(z)\leq2|\rho(z)|\},\qquad
 \chi_i:=\mathbf1_{E_i}\quad(i=1,2),
\]
and
\begin{align*}
 \delta A:=|\sigma+\rho|^p-\sigma^p,\quad
 R_A:=\delta A-p\sigma^{p-1}\rho,\quad
 \delta B:=|\sigma+\rho|^{p-2}(\sigma+\rho)-\sigma^{p-1},\quad
 R_B:=\delta B-(p-1)\sigma^{p-2}\rho.
\end{align*}
The identity
\begin{align*}
 &\mathcal I_\alpha(|\sigma+\rho|^p)
 |\sigma+\rho|^{p-2}(\sigma+\rho)=
 \mathcal I_\alpha(\sigma^p+\delta A)
 (\sigma^{p-1}+\delta B)\\=&
 \mathcal I_\alpha(\sigma^p)\sigma^{p-1}
 +\mathcal I_\alpha(\delta A)\sigma^{p-1}
 +\mathcal I_\alpha(\sigma^p)\delta B
 +\mathcal I_\alpha(\delta A)\delta B
\end{align*}
and \eqref{N} imply
\begin{equation}\label{eq:strong-N-exact-identity}
 N(\rho,\sigma)
 =\mathcal I_\alpha(R_A)\sigma^{p-1}
 +\mathcal I_\alpha(\delta A)\delta B
 +\mathcal I_\alpha(\sigma^p)R_B.
\end{equation}
Since
$
 \mathbb R_y^n\times\mathbb R_x^n
 =\bigcup_{i,j=1}^2(E_i\times E_j)$,
\eqref{eq:strong-N-exact-identity} is equivalent to
\begin{equation}\label{eq:strong-N-four-regions}
\begin{aligned}
 N(\rho,\sigma)(x)
 =\sum_{i,j=1}^2\chi_j(x)\Big[
 \mathcal I_\alpha(R_A\chi_i)(x)\sigma^{p-1}(x)
 +\mathcal I_\alpha(\delta A\chi_i)(x)\delta B(x)+\mathcal I_\alpha(\sigma^p\chi_i)(x)R_B(x)\Big].
\end{aligned}
\end{equation}

Using Taylor's formula, we have
\begin{equation}\label{eq:strong-N-region-E1}
\begin{aligned}
 |\delta A|\chi_1&\leq C\sigma^{p-1}|\rho|\chi_1,
 &|R_A|\chi_1&\leq C\sigma^{p-2}\rho^2\chi_1
 \leq C|\rho|^p\chi_1,\\
 |\delta B|\chi_1&\leq C\sigma^{p-2}|\rho|\chi_1
 \leq C|\rho|^{p-1}\chi_1,
 &|R_B|\chi_1&\leq C\sigma^{p-3}\rho^2\chi_1.
\end{aligned}
\end{equation}
On $E_2$,
$
 \sigma\leq2|\rho|$, $|\sigma+\rho|\leq3|\rho|$,
and therefore
\begin{equation}\label{eq:strong-N-region-E2}
 |\delta A|\chi_2+|R_A|\chi_2\leq C|\rho|^p\chi_2,\,\,\,
 |\delta B|\chi_2\leq C|\rho|^{p-1}\chi_2,\,\,\,
 |R_B|\chi_2\leq
 C\bigl(|\rho|^{p-1}+\sigma^{p-2}|\rho|\bigr)\chi_2.
\end{equation}
In particular, the convolution variable satisfies
\begin{align*}
 |\mathcal I_\alpha(R_A)(x)|
 \leq \mathcal I_\alpha(|R_A|\chi_1)(x)
       +\mathcal I_\alpha(|R_A|\chi_2)(x)\leq C\mathcal I_\alpha
       (\sigma^{p-2}\rho^2\chi_1)(x)
       +C\mathcal I_\alpha(|\rho|^p\chi_2)(x)\leq C\mathcal I_\alpha(|\rho|^p)(x).
\end{align*}

Let
$
 \varepsilon:=\|\rho\|_{L^{2^*}}$.
Since $\nu<\infty$ and every $U_k$ has the same $L^{2^*}$-norm, so
$
 \|\sigma\|_{L^{2^*}}
 \leq\sum_{k=1}^{\nu}\|U_k\|_{L^{2^*}}\leq C.
$
Equations \eqref{eq:strong-N-region-E1}--\eqref{eq:strong-N-region-E2}
and $ps=2^*$, $(p-1)t=2^*$ give
\begin{align*}
 \|R_A\chi_i\|_{L^s}&\leq C\varepsilon^p,
 &&i=1,2, \\
 \|\delta A\chi_1\|_{L^s}
 &\leq C\|\sigma\|_{L^{2^*}}^{p-1}
          \|\rho\|_{L^{2^*}}
 \leq C\varepsilon,
 &\|\delta A\chi_2\|_{L^s}&\leq C\varepsilon^p,
 \\
 \|\delta B\chi_j\|_{L^t}&\leq C\varepsilon^{p-1},
 &\|\sigma^{p-1}\chi_j\|_{L^t}&\leq C,
 \qquad j=1,2.
\end{align*}
Applying \eqref{eq:strong-N-HLS-complete} separately to
$y\in E_i$ and $x\in E_j$ gives
\begin{align}
 \sum_{i,j=1}^2
 \|\chi_j\mathcal I_\alpha(R_A\chi_i)\sigma^{p-1}\|_{L^r}
 \leq C\varepsilon^p,\quad \sum_{i,j=1}^2
 \|\chi_j\mathcal I_\alpha(\delta A\chi_i)\delta B\|_{L^r}
 \leq C\bigl(\varepsilon^p+\varepsilon^{2p-1}\bigr).
 \label{eq:strong-N-four-dA}
\end{align}

The single-bubble identity and $\sigma^p\leq C\sum_{k=1}^{\nu}U_k^p$
give
\begin{align*}
 0\leq\mathcal I_\alpha(\sigma^p\chi_i)
 \leq\mathcal I_\alpha(\sigma^p)
 \leq C\sum_{k=1}^{\nu}\mathcal I_\alpha(U_k^p)
 =C\sum_{k=1}^{\nu}U_k^\gamma\leq C\left(\sum_{k=1}^{\nu}U_k\right)^\gamma
 =C\sigma^\gamma,
 \qquad i=1,2.
\end{align*}
For $x\in E_1$, \eqref{eq:strong-N-region-E1} and
$\gamma+p=2^*$ yield
\begin{equation*}
 \mathcal I_\alpha(\sigma^p\chi_i)|R_B|\chi_1
 \leq C\sigma^{2^*-3}\rho^2\chi_1.
\end{equation*}
If $5\leq n\leq6$, then $2^*-3\geq0$ and
\begin{equation*}
 \|\sigma^{2^*-3}\rho^2\|_{L^r(E_1)}
 \leq\|\sigma\|_{L^{2^*}}^{2^*-3}
       \|\rho\|_{L^{2^*}}^2
 \leq C\varepsilon^2.
\end{equation*}
If $n\geq7$, then $2^*-3<0$ and $\sigma>2|\rho|$ on $E_1$; hence
\begin{equation*}
 \sigma^{2^*-3}\rho^2\chi_1
 \leq C|\rho|^{2^*-1}\chi_1,
 \qquad
 \|\sigma^{2^*-3}\rho^2\|_{L^r(E_1)}
 \leq C\varepsilon^{2^*-1}.
\end{equation*}
For $x\in E_2$, \eqref{eq:strong-N-region-E2} and
$\sigma\leq2|\rho|$ give
\begin{align*}
 \mathcal I_\alpha(\sigma^p\chi_i)|R_B|\chi_2
 \leq C\bigl(
 \sigma^\gamma|\rho|^{p-1}
 +\sigma^{\gamma+p-2}|\rho|\bigr)\chi_2\leq C|\rho|^{2^*-1}\chi_2.
\end{align*}
Since $(2^*-1)r=2^*$,
\begin{equation}\label{eq:strong-N-four-RB}
 \sum_{i,j=1}^2
 \|\chi_j\mathcal I_\alpha(\sigma^p\chi_i)R_B\|_{L^r}
 \leq C\bigl(\varepsilon^2+\varepsilon^{2^*-1}\bigr).
\end{equation}

Combining \eqref{eq:strong-N-four-regions}, \eqref{eq:strong-N-four-dA}, and
\eqref{eq:strong-N-four-RB}, we have
\[
 \|N(\rho,\sigma)\|_{L^r}
 \leq C\bigl(
 \varepsilon^p+\varepsilon^{2p-1}
 +\varepsilon^2+\varepsilon^{2^*-1}\bigr).
\]
 Since $1<p<2$, one has
$
  \min\{p,2p-1,2,2^*-1\}
  =1+\min\{p-1,2^*-2\}=1+\eta$.
 For $\|\rho\|_{\dot H^1}\leq1$, Sobolev inequality and
 $L^{2n/(n+2)}(\mathbb R^n)\hookrightarrow\dot H^{-1}(\mathbb R^n)$ imply
 \[
  \|N(\rho,\sigma)\|_{\dot H^{-1}}
  \leq C\|\rho\|_{\dot H^1}^{1+\eta}.
 \]

\end{proof}

Now we establish the estimates for $\rho$.
\begin{proposition}\label{prop:error-estimate}
Suppose that $n\geq3$, $0<\alpha<n$, and $\nu\geq2$.  There exist
numbers $\delta_{{\rm int},0},\delta_{{\rm err},0}>0$ and a constant
$C<\infty$, depending only on $(n,\nu,\alpha)$, with the following
property.  Let
\[
 0<\delta_{\rm int}\leq\delta_{{\rm int},0},\qquad
 0<\delta_{\rm err}\leq\delta_{{\rm err},0},
\]
let $\sigma=\sum_{i=1}^{\nu}U[z_i,\lambda_i]$ be a
$\delta_{\rm int}$-interacting family, and set
\[
 R:=\frac12\min_{i\ne j}R_{ij}.
\]
Let $f\in\dot H^{-1}(\mathbb R^n)$ and
$\rho\in\dot H^1(\mathbb R^n)$ satisfy
$\|\rho\|_{\dot H^1}\leq\delta_{\rm err}$ and the equation together
with all the modulation orthogonality conditions
\begin{equation}\label{eq:error-fixed-point}
\left\{\begin{array}{l}
\Delta \rho+ D_{n,\alpha}\bigg(\frac{1}{|x|^{{\alpha}}}* \left|\sigma+\rho\right|^{p_{\alpha}} \bigg)g_{p_\alpha}(\sigma+\rho)-\sum_{i=1}^\nu U_i^{2^*-1}+f
= 0\quad \text { in } \mathbb{R}^n, \\
\int \bigg[\Big(|x|^{-\alpha}\ast |U_{i}|^{p_{\alpha}}\Big)(p_{\alpha}-1   )U_{i}^{ p_{\alpha}-2}Z_{i}^{a}+p_{\alpha}\Big(|x|^{-\alpha}\ast (U_{i}^{p_{\alpha}-1}
		Z_{i}^{a})\Big)U^{p_{\alpha}-1}_{i}\bigg] \rho=0,
\end{array}\right.
\end{equation}
in the distributional sense, where $i=1, \ldots, \nu; a=1, \ldots, n+1$.  Then
\[
\|\rho\|_{\dot{H}^{1}(\mathbb{R}^{n})}
\leq C\bigl(\mathfrak{H}_{\nu,n,\alpha}(R)
+\|f\|_{\dot{H}^{-1}(\mathbb{R}^{n})}\bigr).
\]
\end{proposition}
\begin{proof}
By hypothesis, $\rho$ satisfies \eqref{eq:error-fixed-point} and its
orthogonality conditions.
Recall that $N(\cdot,\cdot)$ and the full nonlocal interaction term $h$ are defined in
\eqref{N} and \eqref{h}, respectively.  Then \eqref{eq:error-fixed-point} is equivalent to
\[
\rho=A(\rho)=:-\mathcal{L}\left(N(\rho,\sigma)\right)-\mathcal{L}(h)-\mathcal{L}(f),
\]
where $\mathcal{L}$ is defined in Proposition \ref{pro4.2.1}.

We first estimate the nonlinear remainder.  Set
\[
    r=\frac{2n}{n+2},\qquad s=\frac{2n}{2n-\alpha},\qquad
    t=\frac{2n}{n+2-\alpha}.
\]
Then $2/s+\alpha/n=2$, $p_\alpha s=2^*$ and
$(p_\alpha-1)t=2^*$.  If $0<\alpha\leq4$, then $p_\alpha\geq2$. Using Taylor inequalities, we have
\[
\begin{split}
|N(\rho,\sigma)|
&\leq C\big(|x|^{-\alpha}*(\sigma^{p_\alpha-2}\rho^2+|\rho|^{p_\alpha})\big)
      \left(\sigma^{p_\alpha-1}+|\rho|^{p_\alpha-1}\right)\\
&\quad+C\big(|x|^{-\alpha}*(\sigma^{p_\alpha-1}|\rho|)\big)
      \left(\sigma^{p_\alpha-2}|\rho|+|\rho|^{p_\alpha-1}\right)+C\left(\sigma^{2^*-3}\rho^2\mathbf 1_{\{\sigma>2|\rho|\}}
      +|\rho|^{2^*-1}\mathbf 1_{\{\sigma\leq2|\rho|\}}\right).
\end{split}
\]
By Hardy--Littlewood--Sobolev inequality, H\"older inequality and the
Sobolev embedding,
\[
    \|N(\rho,\sigma)\|_{L^r}
    \leq C\left(\|\rho\|_{\dot H^1}^{2}
    +\|\rho\|_{\dot H^1}^{p_\alpha}
    +\|\rho\|_{\dot H^1}^{2p_\alpha-1}
    +\|\rho\|_{\dot H^1}^{2^*-1}\right).
\]
For the last term, if $2^*-3<0$, then on $\{\sigma>2|\rho|\}$ one uses
$\sigma^{2^*-3}\rho^2\leq C|\rho|^{2^*-1}$; if $2^*-3\geq0$, it follows
directly from H\"older inequality.
Thus, for $0<\alpha\leq4$,
\[
    \|N(\rho,\sigma)\|_{\dot H^{-1}}
    \leq C\|\rho\|_{\dot H^1}^{1+\eta_0}
\]
for some $\eta_0>0$, provided $\|\rho\|_{\dot H^1}$ is small.

If $4<\alpha<n$, Lemma \ref{le:strong-nonlinear} gives
\[
    \|N(\rho,\sigma)\|_{\dot H^{-1}}
    \leq C\|\rho\|_{\dot H^1}^{1+\eta_1}
\]
for some $\eta_1>0$.  Taking $\eta=\eta_0$ in the first case and
$\eta=\eta_1$ in the second case, in all cases
\[
    \|N(\rho,\sigma)\|_{\dot H^{-1}}
    \leq C\|\rho\|_{\dot H^1}^{1+\eta}.
\]
By Lemma \ref{le7.1}, Lemma \ref{le3.1}, Proposition \ref{pro4.1.1} and Proposition \ref{pro4.2.1}, we get that
\[
\begin{aligned}
\|\rho\|_{\dot{H}^{1}(\mathbb{R}^{n})}=\|A(\rho)\|_{\dot{H}^{1}(\mathbb{R}^{n})} & \lesssim \left\|N(\rho,\sigma)\right\|_{\dot H^{-1}(\mathbb{R}^{n})}+\|h\|_{\dot{H}^{-1}(\mathbb{R}^{n})}+\|f\|_{\dot{H}^{-1}(\mathbb{R}^{n})} \\
&  \lesssim \|\rho\|_{\dot{H}^{1}(\mathbb{R}^{n}) }^{1+\eta}+ \mathfrak{H}_{\nu,n,\alpha}(R)+\|f\|_{\dot{H}^{-1}(\mathbb{R}^{n})}.
\end{aligned}
\]
Choose $\delta_{{\rm int},0}$ small enough for Propositions
\ref{pro4.1.1}--\ref{pro4.2.1} and the interaction estimates, and then
choose $\delta_{{\rm err},0}$ so that
$C\delta_{{\rm err},0}^\eta\leq1/2$.  Since
$\|\rho\|_{\dot H^1}\leq\delta_{\rm err}\leq
\delta_{{\rm err},0}$, the first term on the right can be absorbed. Thus
\begin{equation*}
    \begin{split}
        \|\rho\|_{\dot{H}^{1}(\mathbb{R}^{n})}\lesssim \mathfrak{H}_{\nu,n,\alpha}(R)+\|f\|_{\dot{H}^{-1}(\mathbb{R}^{n})}.
    \end{split}
\end{equation*}
\end{proof}

We conclude with the estimates used in the proof of Theorem~\ref{mainthm1}.
\begin{proposition}\label{pro5.5}
Assume all the hypotheses of Proposition \ref{prop:error-estimate}; in
particular, $\sigma=\sum_{i=1}^{\nu}U[z_i,\lambda_i]$ is
$\delta_{\rm int}$-interacting, $R=\frac12\min_{i\ne j}R_{ij}$,
$f\in\dot H^{-1}(\mathbb R^n)$, and
$\rho\in\dot H^1(\mathbb R^n)$ satisfies
\eqref{eq:error-fixed-point}, all its modulation orthogonality conditions,
and $\|\rho\|_{\dot H^1}\leq\delta_{\rm err}$.  After decreasing
$\delta_{{\rm int},0}$ and $\delta_{{\rm err},0}$ if necessary, there are
a constant $C=C(n,\alpha,\nu)$ and a function
$\omega_5:[1,\infty)\to[0,\infty)$, depending only on
$(n,\alpha,\nu)$ and satisfying $\omega_5(R)\to0$ as $R\to\infty$,
such that, for every $1\leq r\leq\nu$ and $1\leq b\leq n+1$,
\[
\left|\int_{\mathbb R^n} N_{1}(\rho,\sigma)Z_r^b\,\dx\right|
+\left|\int_{\mathbb R^n} N(\rho,\sigma) Z_r^b\,\dx \right|
\leq \omega_5(R)Q^{\frac{a_\alpha}{n-2}}
+C\|f\|_{\dot{H}^{-1}}.
\]
The estimate holds uniformly in $r,b$ and in all bubble parameters
satisfying the stated hypotheses.
\end{proposition}

\begin{proof}
We first consider $4<\alpha<n$.  Fix
\[
    1\leq r\leq\nu,\qquad 1\leq b\leq n+1,\qquad Z_r:=Z_r^b,
\]
and write
\[
    \beta=p_\alpha-1,\qquad
    \gamma=2^*-p_\alpha .
\]
Then $0<\beta<1$ and $\gamma+\beta=2^*-1$.  Formula \eqref{eq2.1}
gives
\[
 |Z_r^b(x)|
 \leq\frac{n-2}{2}U_r(x),\qquad 1\leq b\leq n,
\qquad
 |Z_r^{n+1}(x)|
 \leq\frac{n-2}{2}U_r(x).
\]
Consequently,
\begin{equation}\label{eq:all-mode-kernel-bound}
 |Z_r^b|\leq C_nU_r,\qquad
 \|Z_r^b\|_{\dot H^1}\leq C_n,
 \qquad 1\leq b\leq n+1.
\end{equation}
Put
\[
 \mathcal F:=\|f\|_{\dot H^{-1}}.
\]
In the application to Theorem \ref{mainthm1}, the reduction
\(\mathcal F\leq1\) is exactly \eqref{eq:residual-at-most-one} in Lemma
\ref{lem:residual-smallness}.  For the standalone statement of the
present proposition there is no loss in making the same reduction.
Indeed, if \(\mathcal F>1\), the uniform form bound
\eqref{eq:linearized-form-bounded}, \eqref{eq:all-mode-kernel-bound}, and
the nonlinear estimates established in the proof of Proposition
\ref{prop:error-estimate} (in particular, Lemma
\ref{le:strong-nonlinear} when \(4<\alpha<n\)) give
\[
\begin{aligned}
 \left|\int N_1(\rho,\sigma)Z_r^b\,\dx\right|
 \leq C\|\rho\|_{\dot H^1}\leq C\mathcal F,\qquad
 \left|\int N(\rho,\sigma)Z_r^b\,\dx\right|
 \leq C\|\rho\|_{\dot H^1}^{\,1+\eta_*}
 \leq C\mathcal F,
\end{aligned}
\]
for some \(\eta_*>0\), because
\(\|\rho\|_{\dot H^1}\leq\delta_{{\rm err},0}\leq1\).
Thus the asserted estimate is immediate when \(\mathcal F>1\), and
henceforth we assume
\begin{equation}\label{eq:pro55-small-residual}
 \mathcal F\leq1.
\end{equation}
Moreover,
$
 -\Delta Z_r^b=N_1(Z_r^b,U_r).
$
The symmetry of the one-bubble linearized operator and
\eqref{introeq2} give
\begin{equation}\label{eq:one-bubble-N1-orthogonality}
\begin{aligned}
 \int_{\mathbb R^n}N_1(\rho,U_r)Z_r^b\dx
 =\int_{\mathbb R^n}\rho\,N_1(Z_r^b,U_r)\dx=\int_{\mathbb R^n}\nabla\rho\cdot\nabla Z_r^b\dx=0.
\end{aligned}
\end{equation}
Thus
\[
 \int N_1(\rho,\sigma)Z_r^b\dx
 =\int\bigl(N_1(\rho,\sigma)-N_1(\rho,U_r)\bigr)Z_r^b\dx.
\]
In the following estimates we write $Z_r=Z_r^b$.  The estimates are
uniform for $1\leq b\leq n+1$, because $Z_r$ is used only through
\eqref{eq:all-mode-kernel-bound}.  Define
\[
 \iota(x):=\min\left\{i:U_i(x)=\max_{1\leq j\leq\nu}U_j(x)\right\},
 \qquad
 \Omega_\ell:=\{x:\iota(x)=\ell\}.
\]
Thus $\{\Omega_\ell\}_{\ell=1}^{\nu}$ is a measurable partition of
$\mathbb R^n$ and $U_\ell=\max_iU_i$ on $\Omega_\ell$.
 Put
\[
 \mathcal I_\alpha f(x):=
 D_{n,\alpha}\int_{\mathbb R^n}\frac{f(y)}{|x-y|^\alpha}\,dy,
 \qquad
 F:=\sigma^{p_\alpha}-\sum_{i=1}^{\nu}U_i^{p_\alpha}.
\]
The first part of $N_1$ contains two factors $\sigma^\beta$.
Adding and subtracting one $U_r^\beta$ at a time, and decomposing
\[
 \mathcal I_\alpha(\sigma^{p_\alpha})
 =\mathcal I_\alpha(F)+\sum_{j=1}^{\nu}U_j^\gamma,
\]
gives
\begin{equation}\label{eq:new6-five-term-decomposition}
 \left|\int_{\mathbb R^n}
 \big(N_1(\rho,\sigma)-N_1(\rho,U_r)\big)Z_r\dx\right|
 \leq C\left(M_1+M_2+I_1+I_2+I_3\right),
\end{equation}
where
\begin{align*}
 M_1:=
 \iint_{\mathbb R^n\times\mathbb R^n}
 \frac{|\sigma^\beta(y)-U_r^\beta(y)|\,|\rho(y)|
 \sigma^\beta(x)|Z_r(x)|}{|x-y|^\alpha}\,dy\dx,
 \end{align*}
 \begin{align*}
 M_2:=
 \iint_{\mathbb R^n\times\mathbb R^n}
 \frac{U_r^\beta(y)|\rho(y)|
 |\sigma^\beta(x)-U_r^\beta(x)|\,|Z_r(x)|}
 {|x-y|^\alpha}\,dy\dx,
\end{align*}

 \begin{align*}
 I_1:=
 \int_{\mathbb R^n}U_r^{\gamma+1}
 |\sigma^{\beta-1}-U_r^{\beta-1}|\,|\rho|\dx,\quad
 I_2:=
 \sum_{j\ne r}\int_{\mathbb R^n}
 U_j^\gamma\sigma^{\beta-1}U_r|\rho|\dx,\quad
 I_3:=
 \int_{\mathbb R^n}
 \mathcal I_\alpha(F)\sigma^{\beta-1}U_r|\rho|\dx .
\end{align*}
For every measurable set $E\subset\mathbb R^n$, define the restriction in
the outer $x$-variable by
\begin{equation*}
\begin{aligned}
 M_1(E)
 &:=
 \int_E\int_{\mathbb R^n}
 \frac{|\sigma^\beta(y)-U_r^\beta(y)|\,|\rho(y)|
 \sigma^\beta(x)|Z_r(x)|}{|x-y|^\alpha}\,dy\dx .
\end{aligned}
\end{equation*}
In particular,
\begin{equation*}
 M_1=M_1(\mathbb R^n)
 =\sum_{\ell=1}^{\nu}M_1(\Omega_\ell),\qquad \mathcal I_\alpha(U_i^{p_\alpha})=U_i^\gamma,\qquad \mathcal I_\alpha(\sigma^{p_\alpha})
=\mathcal I_\alpha(F)+\sum_{i=1}^{\nu}U_i^\gamma,
\end{equation*}
which verifies \eqref{eq:new6-five-term-decomposition}.

Put
\[
 \mathfrak r:=\frac{2n}{n+2},\qquad
 s:=\frac{2n}{2n-\alpha},\qquad
 t:=\frac{2n}{n+2-\alpha},\qquad
 \eta_\alpha:=\frac{2n-\alpha}{2}.
\]
To keep track of the precise logarithmic cases, set
\[
\begin{aligned}
 \chi_0&:=\mathbf 1_{\{\nu\geq2,\ n=6\ \text{or}\
 \ \alpha=(n+2)/2\}},\qquad \chi_H:=\mathbf 1_{\{\nu\geq2,\
 (n,\alpha)\in\mathcal C_{\log}\}}\\
 \chi_1&:=\mathbf 1_{\{\nu\geq2,\ n=6\}}
 +\mathbf 1_{\{\nu\geq2,\ n\geq7,\ \alpha=(n+2)/2\}}
 +\mathbf 1_{\{n\geq7,\ \nu\geq3,\ \alpha>(n+2)/2\}}.
\end{aligned}
\]
For $4<\alpha<n$,
$
\chi_0,\chi_1,\chi_H\in\{0,1\}
$.
Moreover,
\[
 \frac2s+\frac{\alpha}{n}=2,\qquad
 \frac1s=\frac1t+\frac1{2^*},\qquad
 \frac1{\mathfrak r}+\frac1{2^*}=1,\qquad
 \beta t=2^*.
\]
We also use
\begin{equation}\label{eq:C2-restricted-potential}
 \mathcal I_\alpha(U_i^{p_\alpha}\mathbf1_E)
 \leq\mathcal I_\alpha(U_i^{p_\alpha})=U_i^\gamma.
\end{equation}

For $k\ne m$, Lemma \ref{A.4} gives
\begin{equation}\label{eq:C2-tail-s}
 \|U_k^\beta U_m\|_{L^s(\Omega_k)}
 \leq C R_{km}^{-\eta_\alpha}.
\end{equation}
Since $U_k\geq U_m$ on $\Omega_k$, Lemma \ref{A.4} gives
\[
\int_{\Omega_k}(U_k^\beta U_m)^s\dx\leq CR_{km}^{-n},
\qquad
\|U_k^\beta U_m\|_{L^s(\Omega_k)}
\leq CR_{km}^{-n/s}=CR_{km}^{-\eta_\alpha}.
\]

In addition,
\begin{align}
 \|U_k^{2^*-2}U_m\|_{L^{\mathfrak r}(\Omega_k)}
 &\leq C R_{km}^{-c_n}
 (1+\log R_{km})^{\frac1{\mathfrak r}\mathbf 1_{\{n=6\}}},
                                                               \label{eq:new6-yamabe}\\
 \|U_k^\beta U_m^\gamma\|_{L^{\mathfrak r}(\Omega_k)}
 &\leq C R_{km}^{-b_\alpha}
 (1+\log R_{km})^{\frac1{\mathfrak r}
 \mathbf 1_{\{\alpha=(n+2)/2\}}},                    \label{eq:new6-singular}\\
 \|U_m^{2^*-1}\|_{L^{\mathfrak r}(\Omega_k)}
 &\leq C R_{km}^{-(n+2)/2}.                            \label{eq:new6-pure-tail}
\end{align}
Lemmas \ref{A.4} and \ref{A.3} give
\[
 \|U_k^\beta U_m^\gamma\|_{L^{\mathfrak r}(\Omega_k)}
 \leq C
 \begin{cases}
 R_{km}^{-\alpha},&4<\alpha<(n+2)/2,\\
 R_{km}^{-(n+2)/2}(1+\log R_{km})^{1/\mathfrak r},
     &\alpha=(n+2)/2,\\
 R_{km}^{-(n+2)/2},&(n+2)/2<\alpha<n,
 \end{cases}
\]
which implies \eqref{eq:new6-singular}.  Moreover,
\[
 U_m\leq U_k\quad\text{on }\Omega_k,
\]
and Lemma \ref{A.4}, applied with exponent pair $(0,n)$, gives
\[
 \int_{\Omega_k}U_m^{2^*}\dx
 \leq\int_{\{U_k\geq U_m\}}U_m^{2^*}\dx
 \leq C R_{km}^{-n};
\]
Taking the $\mathfrak r$-th root proves
\eqref{eq:new6-pure-tail}.

We now estimate $M_1$ by splitting its outer $x$-integral.  On
$\Omega_r$,
\[
 \sigma^\beta|Z_r|\leq CU_r^{\beta+1}=CU_r^{p_\alpha}.
\]
Self-adjointness and \eqref{eq:C2-restricted-potential} give
\begin{align*}
 M_1(\Omega_r)
 \leq C\int_{\mathbb R^n}
 U_r^\gamma|\sigma^\beta-U_r^\beta|\,|\rho|\dx\notag=C\sum_{\ell=1}^{\nu}\int_{\Omega_\ell}
 U_r^{2^*-p_\alpha}
 |\sigma^\beta-U_r^\beta|\,|\rho|\dx .
\end{align*}
On $\Omega_r$,
$
 |\sigma^\beta-U_r^\beta|
 \leq C U_r^{\beta-1}\sum_{m\ne r}U_m,
$
while on $\Omega_\ell$, $\ell\ne r$,
$
 |\sigma^\beta-U_r^\beta|\leq C U_\ell^\beta.
$
Therefore the H\"older and Sobolev inequalities,
\eqref{eq:new6-yamabe}, and \eqref{eq:new6-singular} imply
\begin{equation}\label{eq:new6-M1-r}
 M_1(\Omega_r)
 \leq C\|\rho\|_{\dot H^1}\left[
 \sum_{m\ne r}
 \|U_r^{2^*-2}U_m\|_{L^{\mathfrak r}(\Omega_r)}
 +\sum_{\ell\ne r}
 \|U_\ell^\beta U_r^\gamma\|_{L^{\mathfrak r}(\Omega_\ell)}
 \right]\leq C R^{-b_\alpha}(1+\log R)^{\chi_0/\mathfrak r}
 \|\rho\|_{\dot H^1}.
\end{equation}

For the outer region $\Omega_k$, $k\ne r$,
$
 \sigma^\beta|Z_r|\leq C U_k^\beta U_r$.
Moreover, on every $\Omega_\ell$,
$|\sigma^\beta-U_r^\beta|\leq C U_\ell^\beta$, including
$\ell=r$ because all other bubbles are bounded by $U_r$ there.
Consequently
\[
 \sum_\ell\|(\sigma^\beta-U_r^\beta)
 \mathbf1_{\Omega_\ell}\|_{L^t}\leq C.
\]
The HLS, H\"older, and Sobolev inequalities, together with
\eqref{eq:C2-tail-s}, yield
\begin{align}
 \sum_{k\ne r}M_1(\Omega_k)
 &\leq C\|\rho\|_{\dot H^1}
 \sum_{k\ne r}\|U_k^\beta U_r\|_{L^s(\Omega_k)}
 \leq C R^{-\eta_\alpha}\|\rho\|_{\dot H^1}.           \label{eq:new6-M1-tail}
\end{align}
Combining \eqref{eq:new6-M1-r} and
\eqref{eq:new6-M1-tail},
\begin{equation}\label{eq:new6-M1}
 M_1\leq C\left[
 R^{-b_\alpha}(1+\log R)^{\chi_0/\mathfrak r}
 +R^{-\eta_\alpha}\right]\|\rho\|_{\dot H^1}.
\end{equation}

For $M_2$,
$
 \|U_r^\beta\rho\|_{L^s}\leq C\|\rho\|_{\dot H^1}$.
On the outer region $\Omega_r$,
\[
 |\sigma^\beta-U_r^\beta||Z_r|
 \leq C\sum_{m\ne r}U_r^\beta U_m,
\]
and, on $\Omega_k$ with $k\ne r$,
\[
 |\sigma^\beta-U_r^\beta||Z_r|
 \leq C U_k^\beta U_r.
\]
Therefore HLS and \eqref{eq:C2-tail-s} give
\begin{equation}\label{eq:new6-M2}
 M_2\leq C R^{-\eta_\alpha}\|\rho\|_{\dot H^1}.
\end{equation}

We next estimate the original three terms $I_1,I_2,I_3$.
Since $\beta-1<0$, on $\Omega_r$ the mean value theorem gives
\[
 |\sigma^{\beta-1}-U_r^{\beta-1}|
 \leq C U_r^{\beta-2}\sum_{m\ne r}U_m.
\]
On $\Omega_\ell$, $\ell\ne r$, one has
$\sigma\geq U_r$, and hence
\[
 |\sigma^{\beta-1}-U_r^{\beta-1}|
 \leq U_r^{\beta-1}.
\]
It follows from the H\"older and Sobolev inequalities,
\eqref{eq:new6-yamabe}, and \eqref{eq:new6-pure-tail} that
\begin{align}
 I_1
 \leq C\|\rho\|_{\dot H^1}\left[
 \sum_{m\ne r}
 \|U_r^{2^*-2}U_m\|_{L^{\mathfrak r}(\Omega_r)}
 +\sum_{\ell\ne r}
 \|U_r^{2^*-1}\|_{L^{\mathfrak r}(\Omega_\ell)}
 \right]\leq C R^{-b_\alpha}
 (1+\log R)^{\frac1{\mathfrak r}\mathbf 1_{\{n=6\}}}
 \|\rho\|_{\dot H^1}.                                  \label{eq:new6-I1}
\end{align}

For $I_2$, split the outer variable into
$\Omega_r$, $\Omega_j$, and $\Omega_k$ with
$k\notin\{j,r\}$.  The first two pieces satisfy
\begin{align*}
 U_j^\gamma\sigma^{\beta-1}U_r
 \leq U_j^\gamma U_r^\beta
 \quad\hbox{on }\Omega_r,\quad
 U_j^\gamma\sigma^{\beta-1}U_r
 \leq U_j^{2^*-2}U_r
 \quad\hbox{on }\Omega_j.
\end{align*}
For the remaining regions define
\[
 \tau_\alpha:=
 \begin{cases}
 1-\beta,&n=5,\\
 (\gamma-\beta)/2,&n\geq6,\ \alpha>(n+2)/2,\\
 0,&n\geq6,\ \alpha\leq(n+2)/2.
 \end{cases}
\]
Then $0\leq\tau_\alpha\leq1-\beta$.  Since
$U_j/U_k\leq1$ and $U_r/U_k\leq1$ on $\Omega_k$,
\begin{align*}
 U_j^\gamma U_k^{\beta-1}U_r
 =U_k^{\beta+\gamma}
 \left(\frac{U_j}{U_k}\right)^\gamma
 \left(\frac{U_r}{U_k}\right)\leq U_k^{\beta+\gamma}
 \left(\frac{U_j}{U_k}\right)^{\gamma-\tau_\alpha}
 \left(\frac{U_r}{U_k}\right)^{\beta+\tau_\alpha}
 =U_j^{\gamma-\tau_\alpha}U_r^{\beta+\tau_\alpha}.
\end{align*}
Put
\[
q_j:=\gamma-\tau_\alpha,\qquad
q_r:=\beta+\tau_\alpha,
\]
so that
\[
q_j+q_r=2^*-1,\qquad
\frac{n-2}{2}q_j\mathfrak r
+\frac{n-2}{2}q_r\mathfrak r=n.
\]
Moreover,
\begin{align*}
 \|U_j^{q_j}U_r^{q_r}\|_{L^{\mathfrak r}(\Omega_k)}^{\mathfrak r}
\leq
 \int_{\{U_j\geq U_r\}}
 U_j^{q_j\mathfrak r}U_r^{q_r\mathfrak r}\dx
 +\int_{\{U_r>U_j\}}
 U_j^{q_j\mathfrak r}U_r^{q_r\mathfrak r}\dx .       \label{eq:new6-three-split}
\end{align*}
If $n=5$, then $(q_j,q_r)=(4/3,1)$, and Lemma \ref{A.4}
gives
\[
 \int_{\{U_j\geq U_r\}}
 U_j^{q_j\mathfrak r}U_r^{q_r\mathfrak r}\dx
 \leq C R_{jr}^{-3\mathfrak r},\qquad
 \int_{\{U_r>U_j\}}
 U_j^{q_j\mathfrak r}U_r^{q_r\mathfrak r}\dx
 \leq C R_{jr}^{-5}.
\]
If $n\geq6$ and $4<\alpha<(n+2)/2$, then
$(q_j,q_r)=(\gamma,\beta)$, and
\[
 \int_{\{U_j\geq U_r\}}
 U_j^{q_j\mathfrak r}U_r^{q_r\mathfrak r}\dx
 \leq C R_{jr}^{-n},\qquad
 \int_{\{U_r>U_j\}}
 U_j^{q_j\mathfrak r}U_r^{q_r\mathfrak r}\dx
 \leq C R_{jr}^{-\alpha\mathfrak r}.
\]
Finally, if $n\geq6$ and $\alpha\geq(n+2)/2$, then
$q_j=q_r=(2^*-1)/2$, and Lemma \ref{A.3} gives
\[
 \int_{\mathbb R^n}
 U_j^{q_j\mathfrak r}U_r^{q_r\mathfrak r}\dx
 \leq C R_{jr}^{-n}(1+\log R_{jr}).
\]
Taking the $\mathfrak r$-th root in these three alternatives gives
\begin{equation}\label{eq:new6-three-rate}
 \|U_j^{\gamma-\tau_\alpha}U_r^{\beta+\tau_\alpha}
 \|_{L^{\mathfrak r}(\Omega_k)}
 \leq C R^{-b_\alpha}(1+\log R)^{\chi_1/\mathfrak r}.
\end{equation}
Using \eqref{eq:new6-singular}, \eqref{eq:new6-yamabe},
and \eqref{eq:new6-three-rate}, we obtain
\begin{equation*}
 I_2\leq C R^{-b_\alpha}(1+\log R)^{\chi_1/\mathfrak r}
 \|\rho\|_{\dot H^1}.
\end{equation*}

Finally, on $\Omega_\ell$ the mean value theorem gives
\[
 0\leq F\mathbf1_{\Omega_\ell}
 \leq C\sum_{m\ne\ell}
 U_\ell^\beta U_m\mathbf1_{\Omega_\ell}.
\]
Consequently \eqref{eq:C2-tail-s} yields
\[
 \|F\|_{L^s}\leq C R^{-\eta_\alpha}.
\]
On $\Omega_k$, because $\beta-1<0$ and $U_r\leq U_k\leq\sigma$,
\[
 \sigma^{\beta-1}U_r
 \leq U_k^{\beta-1}U_r\leq U_k^\beta,
\]
so the H\"older and Sobolev inequalities give
$
 \|\sigma^{\beta-1}U_r\rho\|_{L^s}\leq
 C\|\rho\|_{\dot H^1}$.
The HLS inequality therefore gives
\begin{equation}\label{eq:new6-I3}
 I_3\leq C R^{-\eta_\alpha}\|\rho\|_{\dot H^1}.
\end{equation}

Combining \eqref{eq:new6-five-term-decomposition},
\eqref{eq:new6-M1}, \eqref{eq:new6-M2}, and
\eqref{eq:new6-I1}--\eqref{eq:new6-I3}, and recalling that
$b\in\{1,\ldots,n+1\}$ was arbitrary, we obtain
\begin{equation}\label{eq:new6-N1-bound}
 \max_{1\leq b\leq n+1}
 \left|\int N_1(\rho,\sigma)Z_r^b\dx\right|
 \leq C R^{-d_\alpha}(1+\log R)^{\chi_1/\mathfrak r}
 \|\rho\|_{\dot H^1},
 \qquad
 d_\alpha:=\min\left\{b_\alpha,\eta_\alpha\right\}.
\end{equation}
Proposition \ref{prop:error-estimate} gives, for $4<\alpha<n$,
\[
\|\rho\|_{\dot H^1}
 \leq C\left[
 R^{-b_\alpha}(1+\log R)^{\chi_H/2}
 +\|f\|_{\dot H^{-1}}\right].
\]
Combining this estimate with \eqref{eq:new6-N1-bound} yields
\begin{align*}
\left|\int N_1(\rho,\sigma)Z_r\dx\right|
 \leq C R^{-(b_\alpha+d_\alpha)}
 (1+\log R)^{\frac{\chi_H}{2}+\frac{\chi_1}{\mathfrak r}}+C R^{-d_\alpha}(1+\log R)^{\chi_1/\mathfrak r}
 \|f\|_{\dot H^{-1}}.
\end{align*}
Moreover,
\[
2b_\alpha-a_\alpha
=
\begin{cases}
a_\alpha,&b_\alpha=a_\alpha,\\
n+2-a_\alpha,&b_\alpha=(n+2)/2<a_\alpha,
\end{cases}
>0,
\]
and
\[
b_\alpha+\eta_\alpha-a_\alpha
=
\begin{cases}
\dfrac{2n-\alpha}{2},&b_\alpha=a_\alpha,\\[2mm]
\dfrac{3n+2-3\alpha}{2}\geq4,
&b_\alpha=(n+2)/2,\ a_\alpha=\alpha\leq n-2,\\[2mm]
\dfrac{n+6-\alpha}{2}>3,
&b_\alpha=(n+2)/2,\ a_\alpha=n-2.
\end{cases}
\]
Consequently,
\[
b_\alpha+d_\alpha-a_\alpha
=\min\{2b_\alpha-a_\alpha,\,
b_\alpha+\eta_\alpha-a_\alpha\}>0,
\]
\[
\begin{aligned}
R^{-(b_\alpha+d_\alpha)}
(1+\log R)^{\frac{\chi_H}{2}+\frac{\chi_1}{\mathfrak r}}
=o(R^{-a_\alpha}),\qquad
R^{-d_\alpha}(1+\log R)^{\chi_1/\mathfrak r}
=o(1).
\end{aligned}
\]
Hence
\begin{equation}\label{eq:C2-N1-final}
    \left|\int N_1(\rho,\sigma)Z_r\dx\right|
    \leq o\left(Q^{\frac{a_\alpha}{n-2}}\right)
    +C\|f\|_{\dot H^{-1}}.
\end{equation}

Put
\[
 m_\alpha:=\min\{p_\alpha,2,2^*-1,2p_\alpha-1\},
\]
Lemma \ref{le:strong-nonlinear}, Proposition \ref{prop:error-estimate}, and
$\|Z_r\|_{\dot H^1}\leq C$ give, for
$\|f\|_{\dot H^{-1}}\leq1$,
\[
 \left|\int N(\rho,\sigma)Z_r\dx\right|
 \leq C\|\rho\|_{\dot H^1}^{m_\alpha}
 \leq C R^{-m_\alpha b_\alpha}
 (1+\log R)^{m_\alpha\chi_H/2}
 +C\|f\|_{\dot H^{-1}}
\]
The implication
\[
R^{-m_\alpha b_\alpha}
(1+\log R)^{m_\alpha\chi_H/2}=o(R^{-a_\alpha})
\]
requires $m_\alpha b_\alpha>a_\alpha$, whereas
\[
(n,\alpha)=(12,10)\Longrightarrow
(a_\alpha,b_\alpha,m_\alpha)=\left(10,7,\frac75\right),\qquad
m_\alpha b_\alpha=\frac{49}{5}<10=a_\alpha.
\]
Thus Lemma \ref{le:strong-nonlinear} alone is insufficient for the full
range $4<\alpha<n$.  We therefore keep the inner convolution and outer
projection regions separate.  Put $p=p_\alpha$,
\[
 E_1=\{\sigma>2|\rho|\},\qquad E_2=\{\sigma\leq2|\rho|\},
 \qquad \chi_i=\mathbf1_{E_i},
\]
\[
 \delta A=|\sigma+\rho|^p-\sigma^p,\quad
 R_A=\delta A-p\sigma^{p-1}\rho,\quad
 \delta B=|\sigma+\rho|^{p-2}(\sigma+\rho)-\sigma^{p-1},\quad
 R_B=\delta B-(p-1)\sigma^{p-2}\rho .
\]
Expanding before estimating gives the exact four-region identity
\begin{equation}\label{eq:C2-four-region-identity}
\begin{aligned}
 N(\rho,\sigma)(x)=\sum_{i,j=1}^2\chi_j(x)\big[
 \mathcal I_\alpha(R_A\chi_i)(x)\sigma^{p-1}(x)
 +\mathcal I_\alpha(\delta A\chi_i)(x)\delta B(x)+\mathcal I_\alpha(\sigma^p\chi_i)(x)R_B(x)\big].
\end{aligned}
\end{equation}
The index $i$ always refers to the convolution variable and $j$ to the
projection variable.  Let
\[
 \varepsilon=\|\rho\|_{\dot H^1},\qquad
 \|\rho\|_{L^{2^*}}\leq C\varepsilon,\qquad
 s=\frac{2n}{2n-\alpha},\qquad \eta_\alpha=\frac{2n-\alpha}{2}.
\]
Taylor's formula on $E_1$ and the elementary bounds on $E_2$ give
\begin{equation*}
\begin{aligned}
 |R_A|\chi_1&\leq C\sigma^{p-2}\rho^2\chi_1
 \leq C|\rho|^p\chi_1,
 &|R_A|\chi_2&\leq C|\rho|^p\chi_2,\\
 |\delta A|\chi_1&\leq C\sigma^{p-1}|\rho|\chi_1,
 &|\delta A|\chi_2&\leq C|\rho|^p\chi_2,\\
 |\delta B|\chi_1&\leq C\sigma^{p-2}|\rho|\chi_1,
 &|\delta B|\chi_2&\leq C|\rho|^{p-1}\chi_2,\\
 |R_B|\chi_1&\leq C\sigma^{p-3}\rho^2\chi_1,
 &|R_B|\chi_2&\leq C\bigl(|\rho|^{p-1}
 +\sigma^{p-2}|\rho|\bigr)\chi_2.
\end{aligned}
\end{equation*}
In particular,
\[
 \|R_A\chi_i\|_{L^s}\leq C\varepsilon^p,\qquad
 \|\delta A\chi_1\|_{L^s}\leq C\varepsilon,\qquad
 \|\delta A\chi_2\|_{L^s}\leq C\varepsilon^p.
\]

We write $J_{ij}^A,J_{ij}^{AB},J_{ij}^B$ for the integrals against
$Z_r$ of the three terms in \eqref{eq:C2-four-region-identity}.  We
show the only estimate that is not a direct symmetric-HLS product.
On $E_1\cap\Omega_r$, self-adjointness and the restricted one-bubble
potential give
\[
\begin{aligned}
 |J_{11}^{A,r}|+|J_{21}^{A,r}|
 &\leq C\sum_\ell\int_{E_1\cap\Omega_\ell}
 U_\ell^{p-2}U_r^\gamma\rho^2
 +C\sum_\ell\int_{E_2\cap\Omega_\ell}U_r^\gamma|\rho|^p\\
 &\leq C\sum_\ell\int U_\ell^{2^*-2}\rho^2
 +C\int|\rho|^{2^*}
 \leq C(\varepsilon^2+\varepsilon^{2^*}).
\end{aligned}
\]
Here $U_r\leq U_\ell$ on $\Omega_\ell$ and
$U_r\leq\sigma\leq2|\rho|$ on $E_2$.  On the remaining output regions,
HLS and Lemma \ref{A.4} give
\[
 \sum_{k\ne r}|J_{i1}^{A,k}|
 \leq C\|R_A\chi_i\|_{L^s}
 \sum_{k\ne r}\|U_k^\beta U_r\|_{L^s(\Omega_k)}
 \leq CR^{-\eta_\alpha}\varepsilon^p .
\]
All other cells use the same two inequalities.  For clarity, the full
list is
\begin{equation}\label{eq:C2-twelve-cells}
\begin{aligned}
 |J_{11}^{A}|+|J_{21}^{A}|
 \leq C\bigl(\varepsilon^2+\varepsilon^{2^*}
 +R^{-\eta_\alpha}\varepsilon^p\bigr),\qquad
 |J_{12}^{A}|+|J_{22}^{A}| \leq C\varepsilon^{2p},\qquad
 |J_{11}^{AB}|\leq C\varepsilon^2,\\
 |J_{12}^{AB}|+|J_{21}^{AB}|\leq C\varepsilon^{p+1},\qquad
 |J_{22}^{AB}|\leq C\varepsilon^{2p},\qquad
 |J_{11}^{B}|+|J_{21}^{B}|\leq C\varepsilon^2, \quad
 |J_{12}^{B}|+|J_{22}^{B}|\leq C\varepsilon^{2^*}.
\end{aligned}
\end{equation}
For the middle row apply symmetric HLS to
$\|\delta A\chi_i\|_{L^s}\|\delta BZ_r\chi_j\|_{L^s}$; on
$E_1\cap\Omega_k$ the negative power is paired first:
\[
 \sigma^{p-2}|\rho|\,|Z_r|
 \leq U_k^{p-2}U_r|\rho|\leq U_k^\beta|\rho|.
\]
For the last row use
$0\leq\mathcal I_\alpha(\sigma^p\chi_i)\leq C\sigma^\gamma$.
Then $E_1\cap\Omega_\ell$ is bounded by
$CU_\ell^{2^*-2}\rho^2$, while on $E_2$ it is bounded by
$C|\rho|^{2^*}$.  These observations prove every entry of
\eqref{eq:C2-twelve-cells}.

Since $p>1$, $2^*>2$, and $\varepsilon<1$, summing the twelve cells in
\eqref{eq:C2-four-region-identity} gives, uniformly for
$1\leq b\leq n+1$,
\begin{equation}\label{eq:C2-N-projection-explicit}
 \max_{1\leq b\leq n+1}
 \left|\int N(\rho,\sigma)Z_r^b\dx\right|
 \leq C\left[\varepsilon^2+
 R^{-\eta_\alpha}\varepsilon^{p_\alpha}\right].
\end{equation}
Let $L_R=(1+\log R)^{\chi_H/2}$.  Proposition
\ref{prop:error-estimate} gives
\[
 \varepsilon\leq C(R^{-b_\alpha}L_R+\mathcal F).
\]
By \eqref{eq:pro55-small-residual} and $p_\alpha>1$,
\[
 \varepsilon^q\leq
 C R^{-qb_\alpha}L_R^q+C\mathcal F,
 \qquad q\in\{2,p_\alpha\}.
\]
Moreover,
\[
 2b_\alpha>a_\alpha,\qquad
 \eta_\alpha+b_\alpha>a_\alpha.
\]
The second inequality was verified above in the proof of
\eqref{eq:C2-N1-final}.  Consequently,
\[
 2p_\alpha b_\alpha>2b_\alpha>a_\alpha,\qquad
 \eta_\alpha+p_\alpha b_\alpha>
 \eta_\alpha+b_\alpha>a_\alpha.
\]
\[
\begin{aligned}
R^{-2b_\alpha}L_R^2=o(R^{-a_\alpha}),\qquad
R^{-(\eta_\alpha+p_\alpha b_\alpha)}L_R^{p_\alpha}
=o(R^{-a_\alpha}),\qquad
R^{-a_\alpha}\asymp Q^{a_\alpha/(n-2)}.
\end{aligned}
\]
Hence,
\begin{equation}\label{eq:C2-N-final}
    \left|\int N(\rho,\sigma)Z_r\dx\right|
    \leq o\left(Q^{\frac{a_\alpha}{n-2}}\right)
    +C\|f\|_{\dot H^{-1}}.
\end{equation}
The estimates \eqref{eq:C2-N1-final} and \eqref{eq:C2-N-final} prove the
proposition for $4<\alpha<n$.

It remains to consider \(0<\alpha\leq4\).  Set
\[
 p:=p_\alpha,\qquad \beta:=p-1\geq1,\qquad
 \gamma:=2^*-p>0,
\]
and fix \(1\leq r\leq\nu\), \(1\leq b\leq n+1\).  Write
\[
 Z_r:=Z_r^b,\qquad
 A_r:=\sigma^\beta-U_r^\beta,\qquad
 G:=\sigma^p-\sum_{i=1}^{\nu}U_i^p.
\]
By \eqref{eq:one-bubble-N1-orthogonality},
\begin{equation}\label{eq:weak-orthogonal-reduction}
 \int N_1(\rho,\sigma)Z_r\,\dx
 =\int\bigl(N_1(\rho,\sigma)-N_1(\rho,U_r)\bigr)Z_r\,\dx.
\end{equation}

The first component of \(N_1\) satisfies the exact identity
\begin{align*}
 \int\left[
 \mathcal I_\alpha(\sigma^\beta\rho)\sigma^\beta
 -\mathcal I_\alpha(U_r^\beta\rho)U_r^\beta
 \right]Z_r\,\dx                                =
 \int A_r\rho\,\mathcal I_\alpha(U_r^\beta Z_r)\,\dx
 +\int A_r\rho\,\mathcal I_\alpha(A_rZ_r)\,\dx
 +\int U_r^\beta\rho\,\mathcal I_\alpha(A_rZ_r)\,\dx .
\end{align*}
For the second component, using
\(\mathcal I_\alpha(U_i^p)=U_i^\gamma\),
\begin{align*}
 \mathcal I_\alpha(\sigma^p)\sigma^{\beta-1}
 -\mathcal I_\alpha(U_r^p)U_r^{\beta-1}   =
 \mathcal I_\alpha(G)\sigma^{\beta-1}
 +\sum_{j\ne r}U_j^\gamma\sigma^{\beta-1}
 +U_r^\gamma\bigl(\sigma^{\beta-1}-U_r^{\beta-1}\bigr).
\end{align*}
Consequently,
\begin{equation}\label{eq:weak-five-integrals}
\left|\int\bigl(N_1(\rho,\sigma)-N_1(\rho,U_r)\bigr)
Z_r\,\dx\right|
\leq C(J_1+J_2+I_1+I_2+I_3),
\end{equation}
where
\begin{align*}
J_1:=\int |A_r\rho|\,
\mathcal I_\alpha(U_r^\beta|Z_r|)\,\dx,\qquad
J_2:=\left|\int A_r\rho\,\mathcal I_\alpha(A_rZ_r)\,\dx\right|
+\left|\int U_r^\beta\rho\,\mathcal I_\alpha(A_rZ_r)\,\dx\right|,
\end{align*}
\begin{align*}
I_1:=\int U_r^\gamma
\left|\sigma^{\beta-1}-U_r^{\beta-1}\right|
|\rho Z_r|\,\dx,\qquad
I_2:=\sum_{j\ne r}\int
U_j^\gamma\sigma^{\beta-1}|\rho Z_r|\,\dx,\qquad
I_3:=\int\mathcal I_\alpha(G)
\sigma^{\beta-1}|\rho Z_r|\,\dx .
\end{align*}
The terms \(J_1,J_2\), which are absent from the local-coefficient
calculation, are required by the first component of \(N_1\).

Put
\[
 q:=\frac{2n}{n+2},\qquad
 s:=\frac{2n}{2n-\alpha},\qquad
 t_\alpha:=\frac{2n}{n+2-\alpha}.
\]
Moreover,
\[
 \|A_r\|_{L^{t_\alpha}}
 +\|U_r^\beta\|_{L^{t_\alpha}}
 +\|\sigma^{\beta-1}Z_r\|_{L^{t_\alpha}}\leq C,
\]
and \eqref{eq:all-mode-kernel-bound} gives
\[
 |Z_r|\leq CU_r,\qquad
 \mathcal I_\alpha(U_r^\beta|Z_r|)
 \leq C\mathcal I_\alpha(U_r^p)=CU_r^\gamma.
\]
H\"older inequality and the bilinear HLS inequality therefore give
\begin{align}\label{eq:weak-J1}
J_1\leq C\|A_rU_r^\gamma\|_{L^q}\|\rho\|_{\dot H^1},
\qquad
J_2\leq C\|A_rZ_r\|_{L^s}\|\rho\|_{\dot H^1}.
\end{align}
\begin{align} \label{eq:weak-I1}
I_1\leq C\left\|U_r^{\gamma+1}
\left|\sigma^{\beta-1}-U_r^{\beta-1}\right|\right\|_{L^q}
\|\rho\|_{\dot H^1},     \qquad
I_2\leq C\sum_{j\ne r}
\|U_j^\gamma\sigma^{\beta-1}U_r\|_{L^q}
\|\rho\|_{\dot H^1},                 \qquad
I_3\leq C\|G\|_{L^s}\|\rho\|_{\dot H^1}.
\end{align}

Let
\[
 \Omega_\ell:=\{U_\ell=\max_iU_i\},\qquad
 S_\ell:=\sum_{j\ne\ell}U_j,
\]
with ties assigned to the smallest index.  On
\(\Omega_\ell\), \(U_\ell\leq\sigma\leq\nu U_\ell\).  The mean-value
theorem gives
\begin{equation*}
\begin{aligned}
A_r\mathbf1_{\Omega_r}
&\leq CU_r^{\beta-1}S_r,
&A_r\mathbf1_{\Omega_\ell}
&\leq CU_\ell^\beta\mathbf1_{\Omega_\ell}\quad(\ell\ne r),
&G\mathbf1_{\Omega_\ell}
&\leq CU_\ell^\beta S_\ell\mathbf1_{\Omega_\ell}.
\end{aligned}
\end{equation*}
It follows that
\begin{align}
\|A_rU_r^\gamma\|_{L^q}
&\leq C\sum_{j\ne r}
\|U_r^{2^*-2}U_j\mathbf1_{\Omega_r}\|_{L^q}
+C\sum_{\ell\ne r}
\|U_\ell^\beta U_r^\gamma\mathbf1_{\Omega_\ell}\|_{L^q},
\label{eq:weak-J1-products}\\
\|A_rZ_r\|_{L^s}
&\leq C\sum_{j\ne r}
\|U_r^\beta U_j\mathbf1_{\Omega_r}\|_{L^s}
+C\sum_{\ell\ne r}
\|U_\ell^\beta U_r\mathbf1_{\Omega_\ell}\|_{L^s},
\label{eq:weak-J2-products}\\
\|G\|_{L^s}
&\leq C\sum_{\ell=1}^{\nu}\sum_{j\ne\ell}
\|U_\ell^\beta U_j\mathbf1_{\Omega_\ell}\|_{L^s}.
\label{eq:weak-I3-products}
\end{align}

If \(\beta=1\), then \(I_1=0\).  If \(\beta>1\), then
\begin{align}
\left\|U_r^{\gamma+1}
\left|\sigma^{\beta-1}-U_r^{\beta-1}\right|\right\|_{L^q}
\leq C\sum_{j\ne r}
\|U_r^{2^*-2}U_j\mathbf1_{\Omega_r}\|_{L^q}
+C\sum_{\ell\ne r}
\|U_\ell^{\beta-1}U_r^{\gamma+1}
\mathbf1_{\Omega_\ell}\|_{L^q}.
\label{eq:weak-I1-products}
\end{align}
For \(I_2\), on \(\Omega_\ell\) the integrand coefficient is bounded by
\(CU_j^\gamma U_\ell^{\beta-1}U_r\).  If \(\ell=j\) or \(\ell=r\),
this is respectively
\[
 U_j^{2^*-2}U_r,\qquad U_j^\gamma U_r^\beta.
\]
If \(\ell\notin\{j,r\}\) and \(\beta>1\), Young inequality gives
\[
 U_j^\gamma U_\ell^{\beta-1}U_r
 \leq C U_\ell^{\beta-1}
 \left(U_j^{\gamma+1}+U_r^{\gamma+1}\right).
\]
When \(\beta=1\), it equals \(U_j^\gamma U_r\).  Hence every term in
\eqref{eq:weak-I1} is bounded by a finite sum of ordered two-bubble
products of total power \(2^*-1\).

The products in \eqref{eq:weak-J2-products} and
\eqref{eq:weak-I3-products} have total power \(p\) and are measured in
\(L^s\).  All products in
\eqref{eq:weak-J1-products}, \eqref{eq:weak-I1-products}, and the
preceding estimate for \(I_2\) have total power \(2^*-1\) and are
measured in \(L^q\).  For \(i\ne j\), set
\[
 \mathcal P_\tau(\mu,\eta;i,j)
 :=\|U_i^\mu U_j^\eta
 \mathbf1_{\{U_i\geq U_j\}}\|_{L^\tau}.
\]
If
\[
 (\tau,\mu+\eta)=(s,p)
 \quad\text{or}\quad
 (\tau,\mu+\eta)=(q,2^*-1),
\]
then
\begin{align}
 \|U_i^\mu U_j^\eta
 \mathbf1_{\{U_i\geq U_j\}}\|_{L^\tau}^{\tau}
 &=C_{n,\mu,\eta}
 \int_{\{U_i\geq U_j\}}
 \frac{\lambda_i^{\mathfrak a}}
 {(1+\lambda_i^2|x-z_i|^2)^{\mathfrak a}}
 \frac{\lambda_j^{\mathfrak b}}
 {(1+\lambda_j^2|x-z_j|^2)^{\mathfrak b}}\,\dx,
 \label{eq:weak-A4-local-map}\\
 \mathfrak a&:=\frac{n-2}{2}\mu\tau,\qquad
 \mathfrak b:=\frac{n-2}{2}\eta\tau,\qquad
 \mathfrak a+\mathfrak b=n.\notag
\end{align}
Since \(R_{ij}\geq2R\), Lemma \ref{A.4}, followed by the uniformity
principle \eqref{eq:uniform-principle-a}--\eqref{eq:uniform-principle-b},
gives the explicit estimate
\begin{equation}\label{eq:weak-ordered-product-rate}
 \mathcal P_\tau(\mu,\eta;i,j)
 \leq C
 \begin{cases}
  R^{-(n-2)\eta},&\mu>\eta,\\
  R^{-n/\tau}(1+\log R)^{1/\tau},&\mu=\eta,\\
  R^{-n/\tau},&\mu<\eta.
 \end{cases}
\end{equation}
Indeed, these are respectively the \(a>b\), \(a=b\), and \(a<b\)
conclusions of Lemma \ref{A.4}, after taking the \(\tau\)-th root in
\eqref{eq:weak-A4-local-map}.

We now calculate every norm occurring in
\eqref{eq:weak-J1}--\eqref{eq:weak-I1}.  First,
\begin{equation}\label{eq:weak-pair-yamabe-rate}
 \mathcal P_q(2^*-2,1;i,j)
 \leq C
 \begin{cases}
  R^{-(n-2)},&3\leq n\leq5,\\
  R^{-4}(1+\log R)^{2/3},&n=6,\\
  R^{-(n+2)/2},&n\geq7.
 \end{cases}
\end{equation}
For the pair \((\beta,\gamma)\), since
\((n-2)\gamma=\alpha\) and \(n/q=(n+2)/2\),
\begin{equation}\label{eq:weak-pair-beta-gamma-rate}
 \mathcal P_q(\beta,\gamma;i,j)
 \leq C
 \begin{cases}
  R^{-\alpha},&0<\alpha<(n+2)/2,\\
  R^{-(n+2)/2}(1+\log R)^{(n+2)/(2n)},
       &\alpha=(n+2)/2,\\
  R^{-(n+2)/2},&(n+2)/2<\alpha\leq4.
 \end{cases}
\end{equation}
The last line is present only when \((n+2)/2<4\).  Since
\(\beta\geq1\) and \(n/s=(2n-\alpha)/2\),
\begin{equation}\label{eq:weak-pair-beta-one-rate}
 \mathcal P_s(\beta,1;i,j)
 \leq C
 \begin{cases}
  R^{-(n-2)},&0<\alpha<4,\\
  R^{-(n-2)}(1+\log R)^{(n-2)/n},&\alpha=4.
 \end{cases}
\end{equation}
Finally, when \(0<\alpha<4\), application of
\eqref{eq:weak-ordered-product-rate} to
\((\mu,\eta)=(\beta-1,\gamma+1)\) gives
\begin{equation}\label{eq:weak-pair-remainder-rate}
 \mathcal P_q(\beta-1,\gamma+1;i,j)
 \leq C
 \begin{cases}
  R^{-(n-2+\alpha)},&
       0<\alpha<(6-n)/2,\\
  R^{-(n+2)/2}(1+\log R)^{(n+2)/(2n)},&
       \alpha=(6-n)/2>0,\\
  R^{-(n+2)/2},&
       \alpha>(6-n)/2.
 \end{cases}
\end{equation}
The first two alternatives in \eqref{eq:weak-pair-remainder-rate} are
empty when \(n\geq6\).

We substitute these four estimates term by term.  From
\eqref{eq:weak-J1}, \eqref{eq:weak-J1-products},
\eqref{eq:weak-pair-yamabe-rate}, and
\eqref{eq:weak-pair-beta-gamma-rate},
\begin{equation}\label{eq:weak-J1-explicit}
 J_1\leq C
 \begin{cases}
  R^{-a_\alpha}\|\rho\|_{\dot H^1},
       &(n,\alpha)\ne(6,4),\\
  R^{-4}(1+\log R)^{2/3}\|\rho\|_{\dot H^1},
       &(n,\alpha)=(6,4).
 \end{cases}
\end{equation}
From \eqref{eq:weak-J1}, \eqref{eq:weak-J2-products}, and
\eqref{eq:weak-pair-beta-one-rate},
\begin{equation*}
 J_2\leq C
 \begin{cases}
  R^{-(n-2)}\|\rho\|_{\dot H^1},&0<\alpha<4,\\
  R^{-(n-2)}(1+\log R)^{(n-2)/n}
       \|\rho\|_{\dot H^1},&\alpha=4.
 \end{cases}
\end{equation*}
If \(\alpha=4\), then \(I_1=0\).  If \(0<\alpha<4\), equations
\eqref{eq:weak-I1}, \eqref{eq:weak-I1-products},
\eqref{eq:weak-pair-yamabe-rate}, and
\eqref{eq:weak-pair-remainder-rate} give
\begin{equation*}
 I_1\leq C
 \begin{cases}
  R^{-(n-2)}\|\rho\|_{\dot H^1},&3\leq n\leq5,\\
  R^{-4}(1+\log R)^{2/3}\|\rho\|_{\dot H^1},&n=6,\\
  R^{-(n+2)/2}\|\rho\|_{\dot H^1},&n\geq7.
 \end{cases}
\end{equation*}
For \(0<\alpha<4\), the three types of pairs in \(I_2\) are precisely
\((2^*-2,1)\), \((\beta,\gamma)\), and
\((\beta-1,\gamma+1)\).  If \(\alpha=4\), then
\(\beta=1\), and splitting \(U_j^\gamma U_r\) over
\(\{U_j\geq U_r\}\) and \(\{U_r\geq U_j\}\) reduces \(I_2\) to the
ordered pairs \((\gamma,1)\) and \((1,\gamma)\).  Thus
\eqref{eq:weak-I1} and
\eqref{eq:weak-pair-yamabe-rate}--\eqref{eq:weak-pair-remainder-rate}
yield
\begin{equation*}
 I_2\leq C
 \begin{cases}
  R^{-a_\alpha}\|\rho\|_{\dot H^1},
       &(n,\alpha)\ne(6,4),\\
  R^{-4}(1+\log R)^{2/3}\|\rho\|_{\dot H^1},
       &(n,\alpha)=(6,4).
 \end{cases}
\end{equation*}
Finally, \eqref{eq:weak-I1}, \eqref{eq:weak-I3-products}, and
\eqref{eq:weak-pair-beta-one-rate} give
\begin{equation}\label{eq:weak-I3-explicit}
 I_3\leq C
 \begin{cases}
  R^{-(n-2)}\|\rho\|_{\dot H^1},&0<\alpha<4,\\
  R^{-(n-2)}(1+\log R)^{(n-2)/n}
       \|\rho\|_{\dot H^1},&\alpha=4.
 \end{cases}
\end{equation}
This proves a separate explicit bound for each of the five terms in
\eqref{eq:weak-five-integrals}.

Proposition \ref{prop:error-estimate} and
\eqref{eq:pro55-small-residual} give
$
 \|\rho\|_{\dot H^1}
 \leq C\bigl(\mathfrak H_{\nu,n,\alpha}(R)+\mathcal F\bigr)$.
If \((n,\alpha)\ne(6,4)\), then
\(\mathfrak H_{\nu,n,\alpha}(R)=R^{-a_\alpha}\).  Every coefficient
of \(\|\rho\|_{\dot H^1}\) on the right-hand sides of
\eqref{eq:weak-J1-explicit}--\eqref{eq:weak-I3-explicit} tends to zero.
Consequently, \eqref{eq:weak-five-integrals} gives
\[
 \left|\int\bigl(N_1(\rho,\sigma)-N_1(\rho,U_r)\bigr)
 Z_r^b\,\dx\right|
 \leq o(R^{-a_\alpha})+C\mathcal F .
\]
At \((n,\alpha)=(6,4)\), equations
\eqref{eq:weak-J1-explicit}--\eqref{eq:weak-I3-explicit} and
\(I_1=0\) give
\[
 \left|\int\bigl(N_1(\rho,\sigma)-N_1(\rho,U_r)\bigr)
 Z_r^b\,\dx\right|
 \leq
 CR^{-4}(1+\log R)^{2/3}
 \bigl(R^{-4}(1+\log R)^{1/2}+\mathcal F\bigr).
\]
Since
$
 R^{-8}(1+\log R)^{7/6}=o(R^{-4})
$
and \(R^{-4}(1+\log R)^{2/3}\leq C\) for \(R\geq2\), the same
conclusion holds at the endpoint.  Combining this conclusion with
\eqref{eq:weak-orthogonal-reduction} yields
\[
 \left|\int N_1(\rho,\sigma)Z_r^b\,\dx\right|
 \leq o(R^{-a_\alpha})+C\mathcal F.
\]
For the nonlinear remainder, the Taylor estimate used in the proof of
Proposition \ref{prop:error-estimate} gives,
\begin{equation}\label{eq:pro55-N-Hminus1}
 \|N(\rho,\sigma)\|_{\dot H^{-1}}
 \leq C\|\rho\|_{\dot H^1}^{m},
\end{equation}
where
$
 m:=\min\{2,p_\alpha,2^*-1\}>1
$. By \eqref{eq:pro55-small-residual},
\[
 \|\rho\|_{\dot H^1}^{m}
 \leq C\mathfrak H_{\nu,n,\alpha}(R)^m+C\mathcal F.
\]
Since $\mathfrak H_{\nu,n,\alpha}(R)=R^{-a_\alpha}$ except at
$(n,\alpha)=(6,4)$, where it equals
$R^{-a_\alpha}(\log R)^{1/2}$, one has
\[
 \mathfrak H_{\nu,n,\alpha}(R)^m=o(R^{-a_\alpha}).
\]
Using Proposition \ref{prop:error-estimate} once more and
$\|Z_r^b\|_{\dot H^1}\leq C$, we conclude that
\[
 \left|\int N(\rho,\sigma)Z_r^b\,dx\right|
 \leq o(R^{-a_\alpha})+C\mathcal F.
\]
Finally, \(R^{-a_\alpha}\asymp Q^{a_\alpha/(n-2)}\).  All the
little-\(o\) terms above are uniform by
\eqref{eq:uniform-principle-a}--\eqref{eq:uniform-principle-b};
taking the maximum over the finitely many modes \(r,b\) therefore
produces the modulus \(\omega_5\) in the statement and completes the
proof.
\end{proof}

\section{A Sharp Example}\label{Sec6}
In this section, we construct the two-bubble example which gives the optimality of the estimate in Theorem \ref{mainthm1}.  The argument is written in a unified form for
$
    n\geq6, \frac{n+2}{2}\leq\alpha<n,$
\[
\operatorname{dist}_{\dot H^1}\bigl(u,\mathcal T_2^{(D)}\bigr)\gtrsim
\begin{cases}
\Gamma(u)\left(1+|\log\Gamma(u)|\right)^{1/2},
&(n,\alpha)\in\mathcal C_{\log},\\
\Gamma(u)^{\theta_\alpha},
&n\geq7,\ \alpha>(n+2)/2,
\end{cases}
\qquad
\theta_\alpha=\dfrac{b_\alpha}{a_\alpha}.
\]
Let
\[
    U_1:=U[-Re_1,1],\qquad U_2:=U[Re_1,1],\qquad \sigma:=U_1+U_2,
\]
where $R\gg1$.  Then $Q=q_{12}\approx R^{-(n-2)}$.
For $i=1,2$ and $a=1,\ldots,n+1$, define the linearized kernel densities
\[
    \Psi_i^a:=
    \Big(|x|^{-\alpha}\ast U_i^{p_\alpha}\Big)(p_\alpha-1)U_i^{p_\alpha-2}Z_i^a
    +p_\alpha\Big(|x|^{-\alpha}\ast (U_i^{p_\alpha-1}Z_i^a)\Big)U_i^{p_\alpha-1}.
\]
Equivalently, $-\Delta Z_i^a=D_{n,\alpha}\Psi_i^a$.  Hence the condition $\int \Psi_i^a\rho=0$ is the same as $\dot H^1$-orthogonality of $\rho$ to $Z_i^a$.

\begin{proposition}\label{prop:sharp-correction}
  Suppose that $n\geq6$ and $\frac{n+2}{2}\leq\alpha<n$.  If $R$ is large enough, then there exist $\rho\in\dot H^1(\mathbb R^n)$ and scalars $\left(c_a^i\right)$ such that
\begin{equation}\label{eq8.1.1}
\left\{\begin{array}{l}
\Delta \rho+N_1(\rho,\sigma)+h+N(\rho,\sigma)
=\displaystyle\sum_{i=1}^2\sum_{a=1}^{n+1}c_a^i\Psi_i^a\quad\text{in }\mathbb R^n,\\[4pt]
\displaystyle\int_{\mathbb R^n}\Psi_i^a\rho\,\dx=0,\qquad i=1,2,\quad a=1,\ldots,n+1,
\end{array}\right.
\end{equation}
where $h$, $N_1$, and $N(\cdot,\cdot)$ are the quantities defined in
\eqref{2.1} with $\nu=2$. Moreover,
\begin{equation}\label{eq8.1.2}
    \sum_{i=1}^2\sum_{a=1}^{n+1}|c_a^i|\lesssim Q^{\frac{a_\alpha}{n-2}}\approx R^{-a_\alpha},
    \qquad
    \|\nabla\rho\|_{L^2}\lesssim
    \mathfrak H_{2,n,\alpha}(R).
\end{equation}
\end{proposition}

\begin{proof}
Let
\[
\mathcal Z_R=\operatorname{span}\{Z_i^a:i=1,2,\ 1\leq a\leq n+1\},
\qquad X_R=\mathcal Z_R^\perp.
\]
Proposition \ref{pro4.2.1} gives a uniformly bounded inverse
\(\mathcal L_R:\dot H^{-1}\to X_R\) for the projected linearized
operator.  Thus \eqref{eq8.1.1} is equivalent on \(X_R\) to
\[
\rho=\mathcal A_R(\rho):=-\mathcal L_R\bigl(h+N(\rho,\sigma)\bigr).
\]
Lemma \ref{le3.1} and Lemma \ref{le:strong-nonlinear} give, when
\(4<\alpha<n\),
\[
\|\mathcal A_R(\rho)\|_{\dot H^1}
\leq C\mathfrak H_{2,n,\alpha}(R)
+C\|\rho\|_{\dot H^1}^{1+\eta},
\qquad \eta=\min\{p_\alpha-1,2^*-2\}>0.
\]
For \((n,\alpha)=(6,4)\), the identity
\[
N(\rho,\sigma)=\mathcal I_4(\rho^2)\sigma
+2\mathcal I_4(\sigma\rho)\rho+\mathcal I_4(\rho^2)\rho
\]
and the bilinear HLS inequality give the same estimate with \(\eta=1\).
Choose \(M>2C\), and then \(R\) so large that
\[
C M^{1+\eta}\mathfrak H_{2,n,\alpha}(R)^\eta\leq\frac M2.
\]
Then the weakly compact convex set
\[
\mathcal E_R=\{\rho\in X_R:
\|\rho\|_{\dot H^1}\leq M\mathfrak H_{2,n,\alpha}(R)\}
\]
satisfies \(\mathcal A_R(\mathcal E_R)\subset\mathcal E_R\).

Set
\begin{equation}\label{eq:sec6-fixed-point-exponents}
 \mathscr F(v):=\mathcal I_\alpha(|v|^{p_\alpha})
 g_{p_\alpha}(v),
 \qquad
 s:=\frac{2n}{2n-\alpha},\quad
 t:=\frac{2n}{n+2-\alpha},\quad
 q:=\frac{2n}{\alpha},\quad
 r:=\frac{2n}{n+2}.
\end{equation}
Let \(v_k\rightharpoonup v\) in \(\dot H^1(\mathbb R^n)\).  For an
arbitrary subsequence, the local Rellich theorem supplies a further
subsequence, still denoted by \(v_k\), such that
\[
 v_k\longrightarrow v\quad\text{a.e. in }\mathbb R^n.
\]
Define
\[
 A_k:=|v_k|^{p_\alpha},\qquad
 B_k:=g_{p_\alpha}(v_k),\qquad
 A:=|v|^{p_\alpha},\qquad
 B:=g_{p_\alpha}(v).
\]
Sobolev inequality and \eqref{eq:sec6-fixed-point-exponents} give
\[
 \sup_k\|A_k\|_{L^s}
 +\sup_k\|B_k\|_{L^t}<\infty.
\]
Since \(1<s,t<\infty\), boundedness and almost-everywhere convergence
imply
\begin{equation}\label{eq:sec6-powers-weak}
 A_k\rightharpoonup A\quad\text{in }L^s,\qquad
 B_k\rightharpoonup B\quad\text{in }L^t.
\end{equation}

For every compact \(K\Subset\mathbb R^n\),
\begin{equation}\label{eq:sec6-riesz-local-compactness}
 \mathcal I_\alpha(A_k)\longrightarrow\mathcal I_\alpha(A)
 \quad\text{in }L^{t'}(K).
\end{equation}
Indeed, \(\alpha<n\) and \eqref{eq:sec6-fixed-point-exponents} imply
\[
 s<q,\qquad t'<q.
\]
Fix
\[
 \max\{s,t'\}<\ell<q,\qquad
 \frac1a:=1+\frac1\ell-\frac1s
 =\frac1\ell+\frac{\alpha}{2n}.
\]
Then
\begin{equation}\label{eq:sec6-a-range}
 1<a<\frac n\alpha,\qquad \frac na-\alpha>0.
\end{equation}
Put \(G_k:=A_k-A\).  By \eqref{eq:sec6-powers-weak},
\[
 G_k\rightharpoonup0\quad\text{in }L^s,\qquad
 \sup_k\|G_k\|_{L^s}<\infty.
\]
Choose \(R_K>0\) with \(K\subset B_{R_K}\), and let
\(L\geq2R_K\).  Since \(s'=2n/\alpha\), H\"older inequality yields
\begin{align}
 &\sup_k\left\|D_{n,\alpha}\int_{|y|>L}
 \frac{G_k(y)}{|\,\cdot-y|^\alpha}\,\dy
 \right\|_{L^\ell(K)}
 \leq C_KL^{-\alpha/2}.
 \label{eq:sec6-riesz-tail}
\end{align}
For \(0<\delta<1\), Young's convolution inequality and
\eqref{eq:sec6-a-range} yield
\begin{align}
 &\sup_k\left\|D_{n,\alpha}
 \int_{\substack{|y|\leq L\\|\,\cdot-y|<\delta}}
 \frac{G_k(y)}{|\,\cdot-y|^\alpha}\,\dy
 \right\|_{L^\ell(K)}
 \leq C\delta^{\,n/a-\alpha}.
 \label{eq:sec6-riesz-diagonal}
\end{align}

For fixed \(L,\delta\), define
\[
 (T_{L,\delta}g)(x):=
 D_{n,\alpha}\int_{B_L}
 \frac{\mathbf1_{\{|x-y|\geq\delta\}}}{|x-y|^\alpha}g(y)\,\dy,
 \qquad x\in K.
\]
Its kernel \(\mathcal K_{L,\delta}\) satisfies
\begin{align*}
 \|\mathcal K_{L,\delta}\|_
 {L^\ell(K;L^{s'}(B_L))}
 &\leq D_{n,\alpha}\delta^{-\alpha}
 |K|^{1/\ell}|B_L|^{1/s'}<\infty.
\end{align*}
Hence \(T_{L,\delta}\) is a Hille--Tamarkin operator.  The
Hille--Tamarkin compactness theorem
\cite[Theorem~3.2.11, pp.~383--384]{DunfordPettis1940};
see also \cite[Theorem~2.3 and Remark~2,
pp.~927--930]{Schep1985} gives
\begin{equation}\label{eq:sec6-T-compact}
 T_{L,\delta}:L^s(B_L)\longrightarrow L^\ell(K)
 \quad\text{compact}.
\end{equation}
The restriction \(G_k|_{B_L}\) converges weakly to zero in \(L^s(B_L)\).
Therefore \eqref{eq:sec6-T-compact} implies
\begin{equation}\label{eq:sec6-T-weak-strong}
 \|T_{L,\delta}(G_k|_{B_L})\|_{L^\ell(K)}
 \longrightarrow0.
\end{equation}
Equations \eqref{eq:sec6-riesz-tail},
\eqref{eq:sec6-riesz-diagonal}, and
\eqref{eq:sec6-T-weak-strong} give
\[
 \limsup_{k\to\infty}
 \|\mathcal I_\alpha(G_k)\|_{L^\ell(K)}
 \leq C_KL^{-\alpha/2}+C\delta^{\,n/a-\alpha}.
\]
Letting \(\delta\downarrow0\) and then \(L\to\infty\) proves convergence
in \(L^\ell(K)\).  Since \(\ell>t'\) and \(|K|<\infty\),
\eqref{eq:sec6-riesz-local-compactness} follows.

Let \(\zeta\in C_c^\infty(\mathbb R^n)\) and
\(K=\operatorname{supp}\zeta\).  By
\eqref{eq:sec6-riesz-local-compactness} and the boundedness of \(B_k\)
in \(L^t\),
\[
 \int_{\mathbb R^n}
 \bigl(\mathcal I_\alpha(A_k)-\mathcal I_\alpha(A)\bigr)
 B_k\zeta\,\dx\longrightarrow0.
\]
The HLS inequality, \(t'<q\), and \(|K|<\infty\) imply
\[
 \mathcal I_\alpha(A)\zeta\in L^{t'}(\mathbb R^n).
\]
Hence the weak convergence of \(B_k\) in \(L^t\) gives
\[
 \int_{\mathbb R^n}
 \mathcal I_\alpha(A)(B_k-B)\zeta\,\dx\longrightarrow0.
\]
Consequently,
\begin{equation}\label{eq:sec6-Hartree-distributional}
 \int_{\mathbb R^n}
 \bigl(\mathscr F(v_k)-\mathscr F(v)\bigr)\zeta\,\dx
 \longrightarrow0
 \qquad(\zeta\in C_c^\infty(\mathbb R^n)).
\end{equation}
The bilinear HLS inequality and
\eqref{eq:sec6-fixed-point-exponents} give
\[
 \sup_k\|\mathscr F(v_k)\|_{L^r}
 \leq C\sup_k\|v_k\|_{\dot H^1}^{\,2p_\alpha-1}<\infty,
 \qquad
 \|\mathscr F(v)\|_{L^r}<\infty.
\]
Since \(1<r<\infty\) and
\(C_c^\infty(\mathbb R^n)\) is dense in \(L^{r'}(\mathbb R^n)\),
\eqref{eq:sec6-Hartree-distributional} implies
\begin{equation}\label{eq:sec6-Hartree-weak-continuity}
 \mathscr F(v_k)\rightharpoonup\mathscr F(v)
 \quad\text{in }L^r(\mathbb R^n)
\end{equation}
along the selected subsequence.  Since every subsequence admits a
further subsequence with this limit, the convergence in
\eqref{eq:sec6-Hartree-weak-continuity} holds for the original sequence.

The identity
$
 N(\rho,\sigma)
 =\mathscr F(\sigma+\rho)-\mathscr F(\sigma)-N_1(\rho,\sigma)
$
and the bounded linearity of \(N_1(\,\cdot\,,\sigma)\) imply
\[
 \rho_k\rightharpoonup\rho\quad\text{in }X_R
 \quad\Longrightarrow\quad
 N(\rho_k,\sigma)\rightharpoonup N(\rho,\sigma)
 \quad\text{in }\dot H^{-1}.
\]
The bounded linearity of
\(\mathcal L_R:\dot H^{-1}\to X_R\) therefore gives
\[
 \mathcal A_R(\rho_k)\rightharpoonup\mathcal A_R(\rho)
 \quad\text{in }X_R.
\]
The set \(\mathcal E_R\) is convex and weakly compact.  Since \(X_R\)
is a separable Hilbert space, the weak topology restricted to
\(\mathcal E_R\) is metrizable; weak sequential continuity is therefore
equivalent to weak continuity on \(\mathcal E_R\).  The
Schauder--Tychonoff fixed-point theorem
\cite{Tychonoff1935} yields
\[
 \rho\in\mathcal E_R,\qquad \mathcal A_R(\rho)=\rho.
\]
The Lagrange multipliers \(c_a^i\) in \eqref{eq8.1.1} are then
determined uniquely.

Applying instead the projection-sensitive estimate
\eqref{eq:multiplier-projection-bound}, obtained from the same
finite-dimensional system in the proof of Proposition
\ref{pro4.1.1}, gives
\begin{equation*}
\sum_{i=1}^2\sum_{a=1}^{n+1}|c_a^i|
\leq C\sum_{j=1}^2\sum_{b=1}^{n+1}
\left(
\left|\int hZ_j^b\,\dx\right|
+\left|\int N_1(\rho,\sigma)Z_j^b\,\dx\right|
+\left|\int N(\rho,\sigma)Z_j^b\,\dx\right|
\right).
\end{equation*}

Lemma \ref{le2.1} for
\((n,\alpha)=(6,4)\) and Lemma \ref{le2.1strong} for
\(4<\alpha<n\) give
\[
\left|\int hZ_i^{n+1}\dx\right|\leq CR^{-a_\alpha}.
\]
For \(2\leq b\leq n\), parity in \(x^b\) gives
\[
\int hZ_i^b\dx=0.
\]
For \(b=1\), write \(j=3-i\), \(\beta=p_\alpha-1\), and
\(\gamma=2^*-p_\alpha\).  Taylor's formula on
\(\{U_i\geq U_j\}\) and \(\{U_j>U_i\}\), followed by the
self-adjointness identity used in \eqref{eq:new26-strong-master}, gives
\[
\int hZ_i^1\dx
=\int U_j^\gamma U_i^\beta Z_i^1\dx
+(2^*-1)\int U_i^{2^*-2}U_jZ_i^1\dx+O(R^{-n}).
\]
With \(S=2+4R^2\), \(q_{12}=S^{-(n-2)/2}\), and
\(z_i^1-z_j^1=2(-1)^iR\), Lemma \ref{A.2} yields
\[
\int hZ_i^1\dx
=-(-1)^i\frac{4R}{S}
\left(C_{2,n,\alpha}q_{12}^{\alpha/(n-2)}
+(2^*-1)C_{1,n}q_{12}\right)
+o(R^{-a_\alpha-1}).
\]
Consequently,
\begin{equation}\label{eq:sec6-h-all-projections}
\left|\int hZ_i^{n+1}\dx\right|\leq CR^{-a_\alpha},\qquad
\max_{1\leq b\leq n}\left|\int hZ_i^b\dx\right|
\leq CR^{-a_\alpha-1}.
\end{equation}
Since \(\rho\in\mathcal E_R\), put
\(\varepsilon_R:=\|\rho\|_{\dot H^1}
\leq M\mathfrak H_{2,n,\alpha}(R)\).
For \(4<\alpha<n\), set
\[
d_\alpha:=\min\left\{b_\alpha,\eta_\alpha\right\},
\qquad \eta_\alpha:=\frac{2n-\alpha}{2}.
\]
Then \eqref{eq:new6-N1-bound} and
\eqref{eq:C2-N-projection-explicit} give
\[
\max_{j,b}\left|\int
\bigl(N_1(\rho,\sigma)+N(\rho,\sigma)\bigr)Z_j^b\dx\right|
\leq C\left[
R^{-d_\alpha}(1+\log R)^C\varepsilon_R
+\varepsilon_R^2
+R^{-\eta_\alpha}\varepsilon_R^{p_\alpha}\right]
=o(R^{-a_\alpha}),
\]
where the last step uses the exponent inequalities established in the
proof of Proposition \ref{pro5.5}.  For \((n,\alpha)=(6,4)\),
\eqref{eq:weak-five-integrals},
\eqref{eq:weak-J1-explicit}--\eqref{eq:weak-I3-explicit}, and
\eqref{eq:pro55-N-Hminus1} with \(m=2\) yield
\[
\max_{j,b}\left|\int
\bigl(N_1(\rho,\sigma)+N(\rho,\sigma)\bigr)Z_j^b\dx\right|
\leq C\left[
R^{-4}(1+\log R)^{2/3}\varepsilon_R+\varepsilon_R^2\right]
=o(R^{-4}).
\]
Indeed,
\(\varepsilon_R\leq MR^{-4}(1+\log R)^{1/2}\), so the two displayed
terms are bounded respectively by
\(CR^{-8}(1+\log R)^{7/6}\) and
\(CR^{-8}(1+\log R)\).
Thus, throughout the parameter range of this section,
\[
\max_{j,b}\left|\int
\bigl(N_1(\rho,\sigma)+N(\rho,\sigma)\bigr)Z_j^b\dx\right|
=o(R^{-a_\alpha}).
\]
Consequently,
\[
\sum_{i=1}^2\sum_{a=1}^{n+1}|c_a^i|
\leq CR^{-a_\alpha}
\asymp Q^{a_\alpha/(n-2)}.
\]
This proves \eqref{eq8.1.2} and the proposition.
\end{proof}

Now let $u:=\sigma+\rho$.  By \eqref{eq8.1.1},
\[
    \Delta u+D_{n,\alpha}\left(|x|^{-\alpha}\ast |u|^{p_\alpha}\right)
    g_{p_\alpha}(u)
    =\sum_{i=1}^2\sum_{a=1}^{n+1}c_a^i\Psi_i^a=:f .
\]
Using $|\Psi_i^a|\lesssim U_i^{2^*-1}$, Sobolev embedding, and \eqref{eq8.1.2}, we obtain
\begin{equation}\label{eq8.1.3}
    \Gamma(u)=\|f\|_{\dot H^{-1}}
    \lesssim \sum_{i,a}|c_a^i|
    \lesssim Q^{\frac{a_\alpha}{n-2}}
    \approx R^{-a_\alpha}.
\end{equation}

\begin{lemma}\label{le8.1.1}
For $R$ large enough, the solution constructed above satisfies
\[
    \|\nabla\rho\|_{L^2}\gtrsim
    \mathfrak H_{2,n,\alpha}(R).
\]
\end{lemma}
\begin{proof}
Let $P_R$ be the $\dot H^1$-orthogonal projection onto
$\mathcal Z_R^\perp$, and define
\[
 \mathcal A_Rv:=v-P_R(-\Delta)^{-1}N_1(v,\sigma),
 \qquad v\in\mathcal Z_R^\perp.
\]
Proposition~\ref{pro4.1.1} gives constants $0<c<C$, independent of $R$,
such that
\begin{equation*}
 c\|v\|_{\dot H^1}\leq\|\mathcal A_Rv\|_{\dot H^1}
 \leq C\|v\|_{\dot H^1},
 \qquad v\in\mathcal Z_R^\perp.
\end{equation*}
Projecting \eqref{eq8.1.1} onto $\mathcal Z_R^\perp$ gives
\begin{equation*}
 \mathcal A_R\rho
 =P_R(-\Delta)^{-1}\bigl(h+N(\rho,\sigma)\bigr).
\end{equation*}
Consequently,
\begin{equation}\label{eq:sec6-master-lower}
 C\|\rho\|_{\dot H^1}
 \geq\|P_R(-\Delta)^{-1}h\|_{\dot H^1}
 -C\|N(\rho,\sigma)\|_{\dot H^{-1}}.
\end{equation}
By \eqref{eq8.1.2} and the nonlinear estimate,
\begin{equation}\label{eq:sec6-nonlinear-small}
 \|N(\rho,\sigma)\|_{\dot H^{-1}}
 \leq C\|\rho\|_{\dot H^1}^{1+\eta}
 \leq C\mathfrak H_{2,n,\alpha}(R)^{1+\eta}
 =o_R\bigl(\mathfrak H_{2,n,\alpha}(R)\bigr).
\end{equation}
It remains to prove, uniformly in both parameter ranges,
\begin{equation}\label{eq:sec6-h-lower}
 \|P_R(-\Delta)^{-1}h\|_{\dot H^1}
 \geq c\mathfrak H_{2,n,\alpha}(R).
\end{equation}

Set
\begin{equation}\label{eq:sec6-beta-gamma}
 \beta:=p_\alpha-1=\frac{n+2-\alpha}{n-2}>0,
 \qquad
 \gamma:=2^*-p_\alpha=\frac{\alpha}{n-2}>0.
\end{equation}
Let $\Omega_i:=\{U_i\geq U_j\}$, where $\{i,j\}=\{1,2\}$.  Since
$
 h=A+B$, $B\geq0$, and
 $A=\sum_{\ell=1}^2
 U_\ell^\gamma\bigl(\sigma^\beta-U_\ell^\beta\bigr)$,
the two summands of $A$ satisfy on $\Omega_i$
\begin{align}
 U_i^\gamma\bigl((U_i+U_j)^\beta-U_i^\beta\bigr)
 =U_i^{\beta+\gamma}
 \left[\left(1+\frac{U_j}{U_i}\right)^\beta-1\right]\geq c_\beta U_i^{2^*-2}U_j,                    \label{eq:sec6-positive-first}\\
 U_j^\gamma\bigl((U_i+U_j)^\beta-U_j^\beta\bigr)
 =U_i^\beta U_j^\gamma
 \left[
 \left(1+\frac{U_j}{U_i}\right)^\beta
 -\left(\frac{U_j}{U_i}\right)^\beta
 \right]\geq c_\beta U_i^\beta U_j^\gamma.              \label{eq:sec6-positive-second}
\end{align}
Define the single comparison function
\begin{equation}\label{eq:sec6-unified-g}
g_R:=
 \sum_{i=1}^2
 \left(U_i^{2^*-2}U_{3-i}+U_i^\beta U_{3-i}^\gamma\right)
 \mathbf1_{\Omega_i}.
\end{equation}
Equations \eqref{eq:sec6-positive-first} and
\eqref{eq:sec6-positive-second} give
\begin{equation}\label{eq:sec6-h-dominates-g}
 h\geq c_\beta g_R\geq0.
\end{equation}

If $(n,\alpha)\in\mathcal C_{\log}$, then
\[
 g_R\geq
 \begin{cases}
  U_1U_2,&n=6,\\[1mm]
  U_1^\vartheta U_2^\vartheta,
  &n\geq7,\ \alpha=(n+2)/2,
 \end{cases}
 \qquad
 \vartheta=\frac{n+2}{2(n-2)}.
\]
Indeed, the first term in \eqref{eq:sec6-unified-g} equals $U_1U_2$
when $n=6$, while the second term equals
$U_1^\vartheta U_2^\vartheta$ when
$n\geq7$ and $\beta=\gamma=\vartheta$.
For nonnegative measurable functions
\(f_1,f_2\in L^{2n/(n+2)}(\mathbb R^n)\) satisfying
\(f_1\geq c f_2\) almost everywhere,
the Newtonian representation gives
\[
 \|f_1\|_{\dot H^{-1}}^2
 =\kappa_n\iint\frac{f_1(x)f_1(y)}{|x-y|^{n-2}}\,dx\,dy
 \geq c^2\kappa_n\iint
 \frac{f_2(x)f_2(y)}{|x-y|^{n-2}}\,dx\,dy
 =c^2\|f_2\|_{\dot H^{-1}}^2.
\]
Applying this inequality to \eqref{eq:sec6-h-dominates-g} and then using
Lemma~\ref{A.6} gives
\begin{equation}\label{eq:sec6-h-unprojected-log}
 \|h\|_{\dot H^{-1}}
 \geq cR^{-b_\alpha}(\log R)^{1/2}
 =c\mathfrak H_{2,n,\alpha}(R).
\end{equation}

If $n\geq7$ and $\alpha>(n+2)/2$, let
\[
 \mathcal B_R^{\mathrm{mid}}:=B_{R/4}(0).
\]
For $x\in\mathcal B_R^{\mathrm{mid}}$ and $i=1,2$,
\[
 cR^{-(n-2)}\leq U_i(x)\leq CR^{-(n-2)}.
\]
Using \eqref{eq:sec6-beta-gamma} in either term of
\eqref{eq:sec6-unified-g}, we obtain
\begin{equation*}
 g_R(x)\geq cR^{-(n+2)},\qquad x\in\mathcal B_R^{\mathrm{mid}}.
\end{equation*}
Therefore
\begin{align*}
 \|h\|_{\dot H^{-1}}^2
 =\kappa_n\iint_{\mathbb R^n\times\mathbb R^n}
 \frac{h(x)h(y)}{|x-y|^{n-2}}\,dx\,dy\geq cR^{-2(n+2)}
 \iint_{\mathcal B_R^{\mathrm{mid}}\times\mathcal B_R^{\mathrm{mid}}}
 |x-y|^{2-n}\,dx\,dy=cR^{-(n+2)},
\end{align*}
and consequently
\begin{equation}\label{eq:sec6-h-unprojected-power}
 \|h\|_{\dot H^{-1}}
 \geq cR^{-(n+2)/2}
 =cR^{-b_\alpha}
 =c\mathfrak H_{2,n,\alpha}(R).
\end{equation}
Thus \eqref{eq:sec6-h-unprojected-log} and
\eqref{eq:sec6-h-unprojected-power} yield in both cases
\begin{equation}\label{eq:sec6-h-unprojected-unified}
 \|(-\Delta)^{-1}h\|_{\dot H^1}
 =\|h\|_{\dot H^{-1}}
 \geq c\mathfrak H_{2,n,\alpha}(R).
\end{equation}

Let $v=(-\Delta)^{-1}h$ and write
\[
 (I-P_R)v=\sum_{i=1}^2\sum_{a=1}^{n+1}d_i^aZ_i^a.
\]
The kernel Gram matrix is uniformly invertible, and
\[
 \langle v,Z_j^b\rangle_{\dot H^1}=\int hZ_j^b\,\dx.
\]
Consequently, \eqref{eq:sec6-h-all-projections} gives
\begin{equation}\label{eq:sec6-kernel-removal-unified}
 \|(I-P_R)(-\Delta)^{-1}h\|_{\dot H^1}
 \leq C\sum_{j=1}^2\sum_{b=1}^{n+1}
 \left|\int hZ_j^b\,\dx\right|
 \leq CR^{-a_\alpha}
 =o\bigl(\mathfrak H_{2,n,\alpha}(R)\bigr).
\end{equation}
Here the last equality follows from
\[
 \begin{cases}
  a_\alpha=b_\alpha,\quad
  \mathfrak H_{2,n,\alpha}(R)
  =R^{-b_\alpha}(\log R)^{1/2},
  &(n,\alpha)\in\mathcal C_{\log},\\
  a_\alpha>b_\alpha,\quad
  \mathfrak H_{2,n,\alpha}(R)=R^{-b_\alpha},
  &n\geq7,\ \alpha>(n+2)/2.
 \end{cases}
\]
Combining \eqref{eq:sec6-h-unprojected-unified} and
\eqref{eq:sec6-kernel-removal-unified}, we have
\begin{align*}
 \|P_R(-\Delta)^{-1}h\|_{\dot H^1}
 &\geq \|(-\Delta)^{-1}h\|_{\dot H^1}
 -\|(I-P_R)(-\Delta)^{-1}h\|_{\dot H^1}\\
 &\geq c\mathfrak H_{2,n,\alpha}(R)
 -o\bigl(\mathfrak H_{2,n,\alpha}(R)\bigr)
 \geq \frac c2\mathfrak H_{2,n,\alpha}(R)
\end{align*}
for all sufficiently large $R$.  This proves
\eqref{eq:sec6-h-lower}.  Finally,
\eqref{eq:sec6-master-lower}, \eqref{eq:sec6-nonlinear-small}, and
\eqref{eq:sec6-h-lower} give
$
 \|\rho\|_{\dot H^1}\geq c\mathfrak H_{2,n,\alpha}(R).
$
\end{proof}

\begin{proof}[Proof of Theorem \ref{mainthm3}]
Let \(u=\sigma+\rho\) be furnished by Proposition
\ref{prop:sharp-correction}.  Equations \eqref{eq8.1.2},
\eqref{eq8.1.3}, and Lemma \ref{le8.1.1} give
\begin{equation}\label{eq:new26-sharp-basic}
c\mathfrak H_{2,n,\alpha}(R)\leq\|\rho\|_{\dot H^1}
\leq C\mathfrak H_{2,n,\alpha}(R),\qquad
\Gamma(u)\leq C_0R^{-a_\alpha}.
\end{equation}

Let \(\widetilde\sigma=\widetilde U_1+\widetilde U_2\) be a minimizing
two-bubble sum for \(u\).  Since
\(\|u-\sigma\|_{\dot H^1}
\leq C\mathfrak H_{2,n,\alpha}(R)\to0\), the standard
parameter-matching assertion in
\cite[Appendix A]{BahriCoron} (see also
\cite[(1.12)--(1.14)]{DSW}) shows, after relabeling, that
\[
\sum_{i=1}^2\left(
\left|\log\frac{\widetilde\lambda_i}{\lambda_i}\right|
+\lambda_i|\widetilde z_i-z_i|\right)
\longrightarrow0.
\]
Moreover, the minimizing property and the admissible competitor
\(\sigma\) give
\[
\|\widetilde\sigma-\sigma\|_{\dot H^1}
\leq\|u-\widetilde\sigma\|_{\dot H^1}
+\|u-\sigma\|_{\dot H^1}
\leq2\|u-\sigma\|_{\dot H^1}
\leq C\mathfrak H_{2,n,\alpha}(R).
\]
Taylor's formula for the \(C^2\)-map
\((z,\log\lambda)\mapsto U[z,\lambda]\in\dot H^1\) gives
\begin{equation}\label{eq:new26-modulation-Taylor}
\widetilde U_i-U_i=\sum_{a=1}^{n+1}t_i^aZ_i^a+A_i,\qquad
\sum_{i=1}^2\|A_i\|_{\dot H^1}\leq C|t|^2,
\end{equation}
where \(|t|^2=\sum_{i,a}|t_i^a|^2\).
The one-bubble Gram matrix
\((\langle Z_i^a,Z_i^b\rangle_{\dot H^1})_{a,b}\) is positive
definite.  Lemma \ref{A.5} gives
\[
\langle Z_1^a,Z_2^b\rangle_{\dot H^1}=o_R(1),
\]
so the full two-bubble Gram matrix is uniformly positive definite.
Therefore, for \(R\) large and \(|t|\) small,
\begin{equation}\label{eq:new26-modulation-equivalence}
c|t|\leq\|\widetilde\sigma-\sigma\|_{\dot H^1}\leq C|t|,
\qquad
\sum_i\|A_i\|_{\dot H^1}
\leq C\|\widetilde\sigma-\sigma\|_{\dot H^1}^2.
\end{equation}

The orthogonality in \eqref{eq8.1.1} is equivalent to
\(\langle\rho,Z_i^a\rangle_{\dot H^1}=0\).  Hence
\eqref{eq:new26-modulation-Taylor}--%
\eqref{eq:new26-modulation-equivalence} imply
\[
|\langle\rho,\widetilde\sigma-\sigma\rangle_{\dot H^1}|
\leq C\|\rho\|_{\dot H^1}
\|\widetilde\sigma-\sigma\|_{\dot H^1}^2.
\]
Since both norms tend to zero,
\[
\begin{aligned}
\|u-\widetilde\sigma\|_{\dot H^1}^2
=\|\rho\|_{\dot H^1}^2
+\|\sigma-\widetilde\sigma\|_{\dot H^1}^2
+2\langle\rho,\sigma-\widetilde\sigma\rangle_{\dot H^1}\geq\frac12\|\rho\|_{\dot H^1}^2.
\end{aligned}
\]
Together with the admissible competitor \(\sigma\) and
\eqref{eq:new26-sharp-basic}, this proves
\[
\operatorname{dist}_{\dot H^1}(u,\mathcal T_2^{(D)})
\asymp\mathfrak H_{2,n,\alpha}(R).
\]

If \((n,\alpha)\in\mathcal C_{\log}\), then \(a_\alpha=b_\alpha\).
For \(g(t)=t(1+|\log t|)^{1/2}\), one has
\[
g'(t)=\frac{1-2\log t}{2(1-\log t)^{1/2}}>0,
\qquad0<t\leq e^{-1}.
\]
Thus \eqref{eq:new26-sharp-basic} gives
\[
\Gamma(u)(1+|\log\Gamma(u)|)^{1/2}
\leq CR^{-a_\alpha}(\log R)^{1/2}
=C\mathfrak H_{2,n,\alpha}(R).
\]
If \(n\geq7\) and \(\alpha>(n+2)/2\), then
\[
\Gamma(u)^{\theta_\alpha}
\leq CR^{-a_\alpha\theta_\alpha}
=CR^{-b_\alpha}=C\mathfrak H_{2,n,\alpha}(R).
\]
This proves Theorem \ref{mainthm3}.
\end{proof}


  \appendix
\section{Auxiliary Results}\label{appA}
\subsection{Non-degeneracy Result}
\begin{lemma}\label{lenondegeneracy}
    Let $n\geq 3$, $0<\alpha<n$ and
    $p_{\alpha}=\frac{2n-\alpha}{n-2}$.  Let
    $U_i=U[z_i,\lambda_i]$ be an Aubin--Talenti bubble.  Assume that
    $Z\in\dot H^1(\mathbb R^n)$ is a weak solution of
    \begin{equation*}
    \begin{split}
             \Delta Z(x) +D_{n,\alpha} \bigg(\int_{\mathbb{R}^{n}}\frac{p_{\alpha}U_{i}^{p_{\alpha}-1}(y)Z(y)   }{|x-y|^{\alpha}}\dy U_{i}^{p_{\alpha}-1}(x) +\int_{\mathbb{R}^{n}}\frac{U_{i}^{p_{\alpha}}(y)   }{|x-y|^{\alpha}}\dy (p_{\alpha}-1)U_{i}^{p_{\alpha}-2}(x) Z(x)   \bigg)=0 \,\,\, \text{ in } \, \mathbb{R}^{n},
    \end{split}
    \end{equation*}
  Then $Z\in \operatorname{span}\{Z^{a}_{i}:a=1,\ldots,n+1\}$.
\end{lemma}

\begin{proof}
By translation and dilation invariance, it suffices to take
\(U=U[0,1]=A_n(1+|x|^2)^{-(n-2)/2}\).  Put
\(\beta=p_\alpha-1\), \(s_\alpha=2n/(2n-\alpha)\), and
\(\mathcal I_\alpha f=D_{n,\alpha}|x|^{-\alpha}*f\).  The equation becomes
\begin{equation}\label{eq:A1-normalized-linearized}
 -\Delta Z=\beta U^{2^*-2}Z
 +p_\alpha U^\beta\mathcal I_\alpha(U^\beta Z).
\end{equation}
Sobolev and Hardy--Littlewood--Sobolev give
\(Z\in L^{2^*}\), \(U^\beta Z\in L^{s_\alpha}\), and
\(\mathcal I_\alpha(U^\beta Z)\in L^{2n/\alpha}\).  We use the
standard localized Riesz-potential bootstrap: on concentric balls,
split \(\mathcal I_\alpha f=\mathcal I_\alpha(\chi f)+
\mathcal I_\alpha((1-\chi)f)\).  The second term is smooth on the
smaller ball, whereas the first is convolution with a compactly
truncated \(L^1\)-kernel.  Hence it preserves local \(H^m\)-regularity.
Interior elliptic regularity applied to \eqref{eq:A1-normalized-linearized}
therefore gives \(Z\in H^m_{\rm loc}\Rightarrow Z\in H^{m+2}_{\rm loc}\).
Starting from \(Z\in H^1_{\rm loc}\) and iterating yields
\(Z\in C^\infty(\mathbb R^n)\); see
\cite[Chapter~3]{AdamsHedberg} and
\cite[Chapter~V, Section~1]{SteinSingular}.

To control infinity, apply the Kelvin transform
\[
 \widetilde Z(x)=|x|^{2-n}Z\!\left(\frac{x}{|x|^2}\right).
\]
The Kelvin transform is an isometry of \(\dot H^1\); conformal
covariance of the bubble and the Riesz kernel shows that \(\widetilde Z\)
satisfies \eqref{eq:A1-normalized-linearized} on
\(\mathbb R^n\setminus\{0\}\).  Since a point has zero
\(\dot H^1\)-capacity for \(n\geq3\), a standard cutoff removes the
origin.  Repeating the preceding bootstrap gives
\(\widetilde Z\in C^\infty\) and \(\widetilde Z\in L^\infty(B_1)\).
Consequently
\[
 |Z(x)|\leq C|x|^{2-n}\quad(|x|\geq1),
 \qquad Z\in C^\infty\cap L^\infty,
\]
so \(Z\) is a bounded classical solution vanishing at infinity.  The
nondegeneracy theorem \cite[Theorem~1.5]{Li-Liu-Tang-Xu}, followed by
the inverse translation and dilation, yields
\(Z\in\operatorname{span}\{Z_i^a:1\leq a\leq n+1\}\).
\end{proof}

\subsection{Some Important Estimates}

\begin{lemma}\label{A.2}
Let $n\geq3$ and $0<\alpha<n$. Let
\(\{(U_{i,(k)},U_{j,(k)})\}_{k\geq1}\) be any sequence of
two-bubble configurations such that \(q_{ij,(k)}\to0\).  Suppressing the
index \((k)\), there are constants \(C_{1,n}>0\),
\(C_{2,n,\alpha}>0\), and a nondecreasing modulus
\(\omega_6:[0,1)\to[0,\infty)\), depending only on \((n,\alpha)\)
and independent of the sequence,
such that \(\omega_6(t)\to0\) as \(t\downarrow0\) and
\begin{align}
\int_{\mathbb R^n}U_i^{2^*-2}U_jZ_i^{n+1}\dx
&=-C_{1,n}\mathfrak d_{ij}q_{ij}
+r_{ij}^{(1)},\notag\\
\int_{\mathbb R^n}U_j^{2^*-p_\alpha}
U_i^{p_\alpha-1}Z_i^{n+1}\dx
&=-C_{2,n,\alpha}\mathfrak d_{ij}
q_{ij}^{\frac{\alpha}{n-2}}
+r_{ij}^{(\alpha)}.
\label{eq:A2-signed-singular}
\end{align}
Here the remainders are uniform in all bubble parameters in the precise
sense that
\begin{equation}\label{eq:A2-uniform-modulus}
 |r_{ij}^{(1)}|\leq\omega_6(q_{ij})q_{ij},
 \qquad
 |r_{ij}^{(\alpha)}|
 \leq\omega_6(q_{ij})q_{ij}^{\frac{\alpha}{n-2}}.
\end{equation}
If $\lambda_i=\lambda_j=1$,
$D_{ij}:=|z_i-z_j|\to\infty$, set
\[
 \mathfrak c_1:=C_{1,n},\qquad
 \mathfrak c_{\alpha/(n-2)}:=C_{2,n,\alpha},
 \qquad C_{2,n,n-2}:=C_{1,n}.
\]
Then, uniformly for $1\leq b\leq n$ and
$\beta\in\{1,\alpha/(n-2)\}$,
\begin{equation}\label{eq:A2-centre-direct}
 \int_{\mathbb R^n}
 U_i^{2^*-\beta-1}U_j^\beta Z_i^b\dx
 =-2\mathfrak c_\beta
 \frac{z_i^b-z_j^b}{S_{ij}}q_{ij}^\beta
 +r_{ij}^{(b,\beta)},
\end{equation}
where, after enlarging the same modulus if necessary,
\begin{equation}\label{eq:A2-centre-uniform-modulus}
 \max_{1\leq b\leq n,\ \beta\in\{1,\alpha/(n-2)\}}
 |r_{ij}^{(b,\beta)}|
 \leq\omega_6(D_{ij}^{-1})
 D_{ij}^{-1}q_{ij}^\beta.
\end{equation}
\end{lemma}
\begin{proof}
Put
\[
 m:=\frac{n-2}{2},\qquad
 \mathcal B:=\left\{1,\frac{\alpha}{n-2}\right\}.
\]
For $\beta\in\mathcal B$, set
\[
 a_\beta:=m\beta,\qquad
 p_\beta:=2^*-\beta=\frac{n-a_\beta}{m},
\]
\[
 F_\beta(y):=(1+|y|^2)^{-(n-a_\beta)},\qquad
 W_\beta(y):=
 \frac{1-|y|^2}{(1+|y|^2)^{n-a_\beta+1}}.
\]
For
\[
 \mathscr I_{\beta,ij}:=
 \int_{\mathbb R^n}U_i^{p_\beta}U_j^\beta\,\dx,
 \qquad
 P_{\beta,ij}:=
 \int_{\mathbb R^n}
 U_i^{p_\beta-1}U_j^\beta Z_i^{n+1}\,\dx,
\]
one has the exact scale-projection identity
\begin{equation*}
 \lambda_i\partial_{\lambda_i}\mathscr I_{\beta,ij}
 =p_\beta P_{\beta,ij}.
\end{equation*}
For $\beta=1$, integration by parts and
$-\Delta U_k=U_k^{2^*-1}$ give
\begin{equation}\label{eq:A2-DSW-identification}
 (2^*-1)P_{1,ij}
 =\int_{\mathbb R^n}U_j^{2^*-1}Z_i^{n+1}\,\dx.
\end{equation}
Formula (F16) in \cite{Bahri1989}, in the form recorded in
\cite[Lemma~A.2]{DSW} and with $i,j$ interchanged, gives a constant
$C_{1,n}>0$ such that
\begin{align}
 \int_{\mathbb R^n}U_j^{2^*-1}Z_i^{n+1}\,\dx
 =
 -(2^*-1)C_{1,n}
 \left(q_{ij}^{-\frac2{n-2}}
 -2\frac{\lambda_j}{\lambda_i}\right)
 q_{ij}^{\frac n{n-2}}
 +O\left(q_{ij}^{\frac n{n-2}}
 |\log q_{ij}|\right).
 \label{eq:A2-Bahri-F16}
\end{align}
Since
\[
 q_{ij}^{-\frac2{n-2}}=S_{ij},\qquad
 S_{ij}-2\frac{\lambda_j}{\lambda_i}
 =\mathfrak d_{ij}S_{ij},\qquad
 S_{ij}q_{ij}^{\frac n{n-2}}=q_{ij},
\]
\eqref{eq:A2-DSW-identification} and \eqref{eq:A2-Bahri-F16} yield
\begin{align}
 P_{1,ij}
 =-C_{1,n}\mathfrak d_{ij}q_{ij}
 +O\left(q_{ij}^{\frac n{n-2}}|\log q_{ij}|\right)=-C_{1,n}\mathfrak d_{ij}q_{ij}+o(q_{ij}).
 \label{eq:A2-F16-yamabe-conclusion}
\end{align}
Put
\[
 \beta_\alpha:=\frac{\alpha}{n-2},\qquad
 a:=a_{\beta_\alpha}=\frac{\alpha}{2},\qquad
 p:=p_{\beta_\alpha}=p_\alpha.
\]
For $\beta\in\mathcal B$, define
\begin{equation*}
M_\beta:=\int_{\mathbb R^n}F_\beta(y)\,dy
 =\pi^{n/2}
 \frac{\Gamma(\frac n2-a_\beta)}{\Gamma(n-a_\beta)},
 \qquad
 \mathbf C_\beta:=A_n^{2^*}M_\beta, \qquad
 \mathfrak c_\beta:=
 \frac{a_\beta\mathbf C_\beta}{p_\beta}>0.
\end{equation*}
Since
\[
 W_\beta
 =-\frac{a_\beta}{n-a_\beta}F_\beta
 +\frac1{n-a_\beta}\operatorname{div}(yF_\beta),
\]
one has
\begin{equation}\label{eq:A2-direct-W-mass}
 \int_{\mathbb R^n}W_\beta(y)\,dy
 =-\frac{a_\beta}{n-a_\beta}M_\beta.
\end{equation}
For $r,s>0$,
\[
 A^{-r}B^{-s}
 =\frac{\Gamma(r+s)}{\Gamma(r)\Gamma(s)}
 \int_0^1
 \frac{t^{r-1}(1-t)^{s-1}}
 {[tA+(1-t)B]^{r+s}}\,dt.
\]
Consequently, for $\xi\in\mathbb R^n$,
\begin{align*}
 H_\beta(\xi)
 &:=
 \int_{\mathbb R^n}
 F_\beta(y)|y-\xi|^{-2a_\beta}\,dy=\frac{\pi^{n/2}\Gamma(\frac n2)}
 {\Gamma(n-a_\beta)\Gamma(a_\beta)}
 \int_0^1
 \frac{v^{a_\beta-1}
 (1-v)^{\frac n2-a_\beta-1}}
 {(1+v|\xi|^2)^{n/2}}\,dv=M_\beta(1+|\xi|^2)^{-a_\beta}.
\end{align*}
For $\varepsilon>0$, integration by parts gives
\begin{align*}
 \int_{\mathbb R^n}\frac{W_\beta(y)}
 {(\varepsilon^2+|y-\xi|^2)^{a_\beta}}\,dy=\frac1{n-a_\beta}
 \left[
 a_\beta H_{\beta,\varepsilon}(\xi)
 +\xi\cdot\nabla_\xi H_{\beta,\varepsilon}(\xi)
 -2a_\beta\varepsilon^2
 \int_{\mathbb R^n}\frac{F_\beta(y)}
 {(\varepsilon^2+|y-\xi|^2)^{a_\beta+1}}\,dy
 \right],
\end{align*}
where
\[
 H_{\beta,\varepsilon}(\xi):=
 \int_{\mathbb R^n}F_\beta(y)
 (\varepsilon^2+|y-\xi|^2)^{-a_\beta}\,dy.
\]
Put $d_\beta:=n-2a_\beta>0$.  For $0<\varepsilon\leq1/2$, polar
coordinates about $\xi$ give
\begin{align*}
 \varepsilon^2\int_{B_1(\xi)}
 (\varepsilon^2+|y-\xi|^2)^{-a_\beta-1}\,dy=C_n\varepsilon^{d_\beta}
 \int_0^{1/\varepsilon}
 t^{n-1}(1+t^2)^{-a_\beta-1}\,dt\leq C_{n,\beta}
 \begin{cases}
  \varepsilon^{d_\beta},&0<d_\beta<2,\\
  \varepsilon^2|\log\varepsilon|,&d_\beta=2,\\
  \varepsilon^2,&d_\beta>2.
 \end{cases}
\end{align*}
Each quantity on the right tends to zero as $\varepsilon\downarrow0$.
Moreover, $2a_\beta<n$ and the explicit formulas for $F_\beta$ and
$\nabla F_\beta$ give the integrable majorants
\[
 F_\beta(y)|y-\xi|^{-2a_\beta},
 \qquad
 |\nabla F_\beta(y)|\,|y-\xi|^{-2a_\beta}.
\]
Integration by parts in $y$ yields
\[
 \nabla_\xi H_{\beta,\varepsilon}(\xi)
 =\int_{\mathbb R^n}\nabla F_\beta(y)
 (\varepsilon^2+|y-\xi|^2)^{-a_\beta}\,dy.
\]
Dominated convergence therefore gives
\[
 H_{\beta,\varepsilon}(\xi)\longrightarrow H_\beta(\xi),
 \qquad
 \nabla_\xi H_{\beta,\varepsilon}(\xi)\longrightarrow
 \nabla_\xi H_\beta(\xi).
\]
The preceding bounds, together with
\[
 \varepsilon^2
 \int_{\mathbb R^n\setminus B_1(\xi)}
 F_\beta(y)(\varepsilon^2+|y-\xi|^2)^{-a_\beta-1}\,dy
 \leq C\varepsilon^2
\]
imply
\begin{align}
 L_\beta(\xi)
 :=
 \int_{\mathbb R^n}
 W_\beta(y)|y-\xi|^{-2a_\beta}\,dy=\frac{a_\beta M_\beta}{n-a_\beta}
 (1-|\xi|^2)(1+|\xi|^2)^{-a_\beta-1}.
 \label{eq:A2-direct-W-riesz}
\end{align}

Set
\[
 \tau:=\frac{\lambda_j}{\lambda_i},\qquad
 \xi:=\lambda_i(z_j-z_i),\qquad
 \eta:=\lambda_j(z_i-z_j)=-\tau\xi,
 \qquad y:=\lambda_i(x-z_i).
\]
Recalling that
\(a=a_{\beta_\alpha}=m\beta_\alpha=\alpha/2\), direct substitution yields
\begin{align}
 P_{\beta_\alpha,ij}
 :=
 \int_{\mathbb R^n}
 U_i^{p-1}U_j^{\beta_\alpha}Z_i^{n+1}\,\dx=mA_n^{2^*}\tau^a
 \int_{\mathbb R^n}
 W_{\beta_\alpha}(y)
 (1+|\tau y+\eta|^2)^{-a}\,dy,
 \label{eq:A2-A3-exact-projection}
\end{align}
and
\begin{align*}
 S_{ij}
 =\tau^{-1}(1+\tau^2+|\eta|^2)
 =\tau(\tau^{-2}+1+|\xi|^2),\qquad
 \mathfrak d_{ij}
 =\frac{1-\tau^2+|\eta|^2}{1+\tau^2+|\eta|^2}
 =\frac{\tau^{-2}-1+|\xi|^2}
 {\tau^{-2}+1+|\xi|^2},\qquad
 q_{ij}^{\beta_\alpha}=S_{ij}^{-a}.
\end{align*}

Let $S_{ij}\to\infty$.  By the bubble-tree trichotomy
\eqref{confi7}--\eqref{confi8}, after extraction of a subsequence,
exactly one of the following alternatives occurs:
\[
 \mathrm{(I)}\quad j\prec i,\qquad
 \mathrm{(II)}\quad i\prec j,\qquad
 \mathrm{(III)}\quad i\text{ and }j\text{ are incomparable}.
\]

\noindent\emph{Case I.}
Put
\[
 R_i:=2(1+|\eta|)\tau^{-1}
 =2\bigl(1+\lambda_j|z_i-z_j|\bigr)
 \frac{\lambda_i}{\lambda_j},
 \qquad E_i:=B_{R_i}(0),\qquad
 E_i^c:=\mathbb R^n\setminus E_i.
\]
Under $y=\lambda_i(x-z_i)$, these are the images of
\[
 \widehat\Omega_i:=
 B_{\frac{2(1+\lambda_j|z_i-z_j|)}{\lambda_j}}(z_i),
 \qquad
 \widehat\Omega_i^c.
\]
For
$
 K_a(w):=(1+|w|^2)^{-a}$,
one has $K_a\in C_b(\mathbb R^n)$, $K_a$ is uniformly continuous,
$R_i\to\infty$, and $W_{\beta_\alpha}\in L^1(\mathbb R^n)$.  Hence
\begin{align*}
 \int_{E_i^c}
 |W_{\beta_\alpha}(y)|K_a(\tau y+\eta)\,dy
 \leq\int_{|y|>R_i}|W_{\beta_\alpha}(y)|\,dy=o(1),
\qquad
 \int_{E_i}
 W_{\beta_\alpha}(y)
 \bigl[K_a(\tau y+\eta)-K_a(\eta)\bigr]\,dy
 =o(1),
\end{align*}
\begin{align*}
 \int_{E_i}W_{\beta_\alpha}(y)\,dy
 =\int_{\mathbb R^n}W_{\beta_\alpha}(y)\,dy+o(1).
\end{align*}
Thus \eqref{eq:A2-A3-exact-projection} and
\eqref{eq:A2-direct-W-mass} give
\begin{align}
 P_{\beta_\alpha,ij}
 =-\mathfrak c_{\beta_\alpha}
 \tau^a(1+|\eta|^2)^{-a}+o(\tau^a)
 =-\mathfrak c_{\beta_\alpha}
 \mathfrak d_{ij}q_{ij}^{\beta_\alpha}
 +o(q_{ij}^{\beta_\alpha}).
 \label{eq:A2-A3-I-final}
\end{align}

\noindent\emph{Case II.}
Put
\[
 u:=\lambda_j(x-z_j)=\tau(y-\xi),
\]
\[
 R_j:=2(1+|\eta|)\tau
 =2\bigl(1+\lambda_j|z_i-z_j|\bigr)
 \frac{\lambda_j}{\lambda_i},
 \qquad F_j:=B_{R_j}(0),\qquad F_j^c:=\mathbb R^n\setminus F_j.
\]
In the $y$-variables,
\[
 \widetilde F_j:=B_{2(1+|\eta|)}(\xi),\qquad
 \widetilde F_j^c:=\mathbb R^n\setminus\widetilde F_j.
\]
Equation \eqref{eq:A2-A3-exact-projection} becomes
\begin{align*}
 P_{\beta_\alpha,ij}
 &=mA_n^{2^*}\tau^{-a}
 \sum_{\ell\in\{j,j^c\}}
 \int_{\widetilde F_\ell}
 W_{\beta_\alpha}(y)
 (\tau^{-2}+|y-\xi|^2)^{-a}\,dy.
\end{align*}
After extraction,
\[
 \xi\to\bar\xi,\qquad |\eta|\to\ell\in[0,\infty].
\]
If $\ell<\infty$, put
$
 F_\infty:=B_{2(1+\ell)}(\bar\xi)$;
if $\ell=\infty$, put $F_\infty:=\mathbb R^n$.  For $0<r<1$ and
$R>2(1+|\bar\xi|)$,
\begin{align*}
 \sup_k\int_{B_r(\xi)}
 (\tau^{-2}+|y-\xi|^2)^{-a}\,dy
 \leq Cr^{n-2a},
\qquad
 \sup_k\int_{|y|>R}|W_{\beta_\alpha}(y)|
 (\tau^{-2}+|y-\xi|^2)^{-a}\,dy
 \leq CR^{-n}.
\end{align*}
It follows that
$$
 \int_{\widetilde F_j}W_{\beta_\alpha}(y)
 (\tau^{-2}+|y-\xi|^2)^{-a}\,dy
 \longrightarrow
 \int_{F_\infty}W_{\beta_\alpha}(y)
 |y-\bar\xi|^{-2a}\,dy,$$
 $$
 \int_{\widetilde F_j^c}W_{\beta_\alpha}(y)
 (\tau^{-2}+|y-\xi|^2)^{-a}\,dy
 \longrightarrow
 \int_{\mathbb R^n\setminus F_\infty}
 W_{\beta_\alpha}(y)|y-\bar\xi|^{-2a}\,dy.
$$
Hence, by \eqref{eq:A2-direct-W-riesz},
\begin{align}
 P_{\beta_\alpha,ij}
 &=mA_n^{2^*}\tau^{-a}
 L_{\beta_\alpha}(\xi)+o(\tau^{-a})=\mathfrak c_{\beta_\alpha}\tau^{-a}
 (1-|\xi|^2)(1+|\xi|^2)^{-a-1}
 +o(\tau^{-a})\notag\\
 &=-\mathfrak c_{\beta_\alpha}
 \mathfrak d_{ij}q_{ij}^{\beta_\alpha}
 +o(q_{ij}^{\beta_\alpha}).
 \label{eq:A2-A3-II-final}
\end{align}

\noindent\emph{Case III.}
Put
\[
 d_{ij}:=\frac{|z_i-z_j|}{2},\qquad
 \Omega_i:=B_{d_{ij}}(z_i),\qquad
 \Omega_j:=B_{d_{ij}}(z_j),
\]
\[
 \Omega_0:=\mathbb R^n\setminus(\Omega_i\cup\Omega_j),\qquad
 L_i:=\lambda_i|z_i-z_j|=|\xi|,\qquad
 L_j:=\lambda_j|z_i-z_j|=|\eta|=\tau L_i.
\]
In the $y$-variables,
\[
 E_i=B_{L_i/2}(0),\qquad
 E_j=B_{L_i/2}(\xi),\qquad
 E_0=\mathbb R^n\setminus(E_i\cup E_j).
\]
Write
\[
 P_{\beta_\alpha,ij}
 =P_{\beta_\alpha,ij}(\Omega_i)
 +P_{\beta_\alpha,ij}(\Omega_j)
 +P_{\beta_\alpha,ij}(\Omega_0).
\]
On $E_i$, dominated convergence and
\eqref{eq:A2-direct-W-mass} yield
\begin{align}
 P_{\beta_\alpha,ij}(\Omega_i)
 =mA_n^{2^*}\tau^aL_j^{-2a}
 \left(\int_{\mathbb R^n}W_{\beta_\alpha}(y)\,dy+o(1)\right)=-\mathfrak c_{\beta_\alpha}
 \tau^{-a}L_i^{-2a}(1+o(1)).
 \label{eq:A2-A3-III-i}
\end{align}
With $u=\tau(y-\xi)$,
$$
 |P_{\beta_\alpha,ij}(\Omega_j)|
 \leq C\tau^aL_i^{-2(n-a)}\tau^{-n}
 \int_{|u|<L_j/2}(1+|u|^2)^{-a}\,du\leq C\tau^{-a}L_i^{-n}.
 \label{eq:A2-A3-III-j}
$$
With $\widehat\xi:=\xi/L_i$ and $y=L_iv$,
\begin{align}
 |P_{\beta_\alpha,ij}(\Omega_0)|
 \leq C\tau^{-a}L_i^{-n}
 \int_{\substack{|v|\geq1/2
 |v-\widehat\xi|\geq1/2}}
 |v|^{-2(n-a)}|v-\widehat\xi|^{-2a}\,dv\leq C\tau^{-a}L_i^{-n}.
 \label{eq:A2-A3-III-zero}
\end{align}
Since $2a<n$,
$
 \tau^{-a}L_i^{-n}
 =o(\tau^{-a}L_i^{-2a})$.
Moreover,
\[
 S_{ij}=\tau L_i^2(1+o(1)),\qquad
 q_{ij}^{\beta_\alpha}
 =\tau^{-a}L_i^{-2a}(1+o(1)),\qquad
 \mathfrak d_{ij}=1+o(1).
\]
Therefore
\begin{equation}\label{eq:A2-A3-III-final}
 P_{\beta_\alpha,ij}
 =-\mathfrak c_{\beta_\alpha}
 \mathfrak d_{ij}q_{ij}^{\beta_\alpha}
 +o(q_{ij}^{\beta_\alpha}).
\end{equation}
Equations \eqref{eq:A2-A3-I-final}, \eqref{eq:A2-A3-II-final}, and
\eqref{eq:A2-A3-III-final}, together with the subsequence principle, give
\begin{equation}\label{eq:A2-A3-conclusion}
 P_{\beta_\alpha,ij}
 =-\mathfrak c_{\beta_\alpha}
 \mathfrak d_{ij}q_{ij}^{\beta_\alpha}
 +o(q_{ij}^{\beta_\alpha})
\end{equation}
uniformly in the bubble parameters.  Thus
$
 C_{2,n,\alpha}:=\mathfrak c_{\beta_\alpha}
 =\frac{\alpha\mathbf C_{\beta_\alpha}}{2p_\alpha}>0,
$
and \eqref{eq:A2-signed-singular} follows.

For $\beta=1$, the computation
\eqref{eq:A2-A3-III-i}--\eqref{eq:A2-A3-III-zero}, with
$\tau=1$ and $|z_i-z_j|\to\infty$, gives
$
 P_{1,ij}=-\mathfrak c_1q_{ij}+o(q_{ij}).
$
Comparison with \eqref{eq:A2-F16-yamabe-conclusion} gives
\[
 \mathfrak c_1=C_{1,n},\qquad
 \alpha=n-2\Longrightarrow C_{2,n,n-2}=C_{1,n}.
\]

It remains to prove \eqref{eq:A2-centre-direct}.  Let
$\lambda_i=\lambda_j=1$, set
\[
 \xi:=z_j-z_i,\qquad D:=|\xi|=D_{ij},\qquad
 S:=S_{ij}=2+D^2,
\]
and, for $1\leq b\leq n$ and $\beta\in\mathcal B$, set
\[
 T_{\beta,ij}^b:=
 \int_{\mathbb R^n}
 U_i^{2^*-\beta-1}U_j^\beta Z_i^b\,\dx.
\]
With $y=x-z_i$,
\begin{equation*}
 T_{\beta,ij}^b
 =2mA_n^{2^*}
 \int_{\mathbb R^n}
 \frac{y^b\,dy}
 {(1+|y|^2)^{n-a_\beta+1}
 (1+|y-\xi|^2)^{a_\beta}}.
\end{equation*}
Put
\[
 \Omega_i:=B_{D/2}(0),\qquad
 \Omega_j:=B_{D/2}(\xi),\qquad
 \Omega_0:=\mathbb R^n\setminus(\Omega_i\cup\Omega_j).
\]
For $y\in\Omega_i$, Taylor's formula gives
\begin{align}
 (1+|y-\xi|^2)^{-a_\beta}
 =(1+D^2)^{-a_\beta}
 +2a_\beta(1+D^2)^{-a_\beta-1}\xi\cdot y
 +\mathcal R_\beta(y),\qquad
 |\mathcal R_\beta(y)|
 \leq CD^{-2a_\beta-2}|y|^2.
 \label{eq:A2-centre-equal-Taylor}
\end{align}
Furthermore,
\begin{equation*}
 \int_{\mathbb R^n}
 \frac{y^by^c}{(1+|y|^2)^{n-a_\beta+1}}\,dy
 =\frac{\delta_{bc}}{2(n-a_\beta)}M_\beta,
\end{equation*}
and
\[
 \int_{|y|<D/2}
 \frac{|y|^3}{(1+|y|^2)^{n-a_\beta+1}}\,dy=o(D),
\qquad
 \int_{|y|>D/2}
 \frac{|y|^2}{(1+|y|^2)^{n-a_\beta+1}}\,dy=o(1).
\]
The constant term in \eqref{eq:A2-centre-equal-Taylor} is odd on
$\Omega_i$.  Hence
\begin{align}
 T_{\beta,ij}^b(\Omega_i)
 &=2\mathfrak c_\beta\xi^b
 (1+D^2)^{-a_\beta-1}
 +o(D^{-2a_\beta-1}).
 \label{eq:A2-centre-equal-i}
\end{align}
On $\Omega_j$, the substitution $u=y-\xi$ gives
$$
 |T_{\beta,ij}^b(\Omega_j)|
 \leq CD^{-2(n-a_\beta+1)+1}
 \int_{|u|<D/2}(1+|u|^2)^{-a_\beta}\,du
 \leq CD^{-n-1}.
$$
On $\Omega_0$, put $\widehat\xi:=\xi/D$ and $y=Dv$.  Then
\begin{align}
 |T_{\beta,ij}^b(\Omega_0)|
 \leq CD^{-n-1}
 \int_{\substack{|v|\geq1/2\\|v-\widehat\xi|\geq1/2}}
 |v|^{-2(n-a_\beta+1)+1}
 |v-\widehat\xi|^{-2a_\beta}\,dv\leq CD^{-n-1}.
 \label{eq:A2-centre-equal-zero}
\end{align}
Since $2a_\beta<n$,
\[
 D^{-n-1}=o(D^{-2a_\beta-1}),\qquad
 (1+D^2)^{-a_\beta-1}
 =S^{-a_\beta-1}(1+o(1)).
\]
Equations \eqref{eq:A2-centre-equal-i}--
\eqref{eq:A2-centre-equal-zero} imply
\begin{align*}
 T_{\beta,ij}^b
 =2\mathfrak c_\beta\xi^bS^{-a_\beta-1}
 +o(D^{-1}S^{-a_\beta})=-2\mathfrak c_\beta
 \frac{z_i^b-z_j^b}{S_{ij}}q_{ij}^\beta
 +o(D_{ij}^{-1}q_{ij}^\beta),
\end{align*}
which is \eqref{eq:A2-centre-direct}.  For \(0<t<1\), define
\begin{align*}
 e_\alpha(t)
 &:=
 \sup\left\{
 q_{ij}^{-\beta_\alpha}
 \left|P_{\beta_\alpha,ij}
 +C_{2,n,\alpha}\mathfrak d_{ij}q_{ij}^{\beta_\alpha}\right|
 \,\middle|\,
 \begin{gathered}
 (z_i,\lambda_i),(z_j,\lambda_j)
 \in\mathbb R^n\times(0,\infty),\\
 i\ne j,\quad 0<q_{ij}\leq t
 \end{gathered}
 \right\},
 \\
 e_c(t)
 &:=
 \sup\left\{
 D_{ij}q_{ij}^{-\beta}
 \left|T_{\beta,ij}^b
 +2\mathfrak c_\beta
 \frac{z_i^b-z_j^b}{S_{ij}}q_{ij}^{\beta}\right|
 \,\middle|\,
 \begin{gathered}
 z_i,z_j\in\mathbb R^n,\quad z_i\ne z_j,\\
 \lambda_i=\lambda_j=1,\quad D_{ij}^{-1}\leq t,\\
 1\leq b\leq n,\quad\beta\in\mathcal B
 \end{gathered}
 \right\}.
\end{align*}
Then
\[
 \eqref{eq:A2-A3-conclusion}
 \Longrightarrow\lim_{t\downarrow0}e_\alpha(t)=0,
 \qquad
 \eqref{eq:A2-centre-equal-i}\text{--}
 \eqref{eq:A2-centre-equal-zero}
 \Longrightarrow\lim_{t\downarrow0}e_c(t)=0.
\]
Consequently, after choosing \(C_n>0\) sufficiently large,
\[
 \omega_6(t):=
 \begin{cases}
  \max\{e_\alpha(t),e_c(t),C_nt^{1/(n-2)}\},&0<t<1,\\
  0,&t=0,
 \end{cases}
 \qquad
 \lim_{t\downarrow0}\omega_6(t)=0.
\]
The functions \(e_\alpha\), \(e_c\), and \(\omega_6\) are
nondecreasing.  Since
\(t^{2/(n-2)}|\log t|\leq C_nt^{1/(n-2)}\) for \(0<t<1\),
\eqref{eq:A2-F16-yamabe-conclusion} and the definitions above imply
\eqref{eq:A2-uniform-modulus} and
\eqref{eq:A2-centre-uniform-modulus}.
\end{proof}

\begin{lemma}\label{A.3}
    Let $n\geq3$, let $U_1,U_2$ be two Aubin--Talenti bubbles, fix
    $\varepsilon>0$, and let $\mu,\eta\geq0$ satisfy
    $\mu+\eta=2^*$.  Assume in addition that either
    $\mu=\eta$ or $|\mu-\eta|\geq\varepsilon$.  Then
    \[
\int_{\mathbb{R}^n} U_1^\mu U_2^\eta \approx_{n,\varepsilon}
\begin{cases}
q_{12}^{\min\{\mu,\eta\}},&|\mu-\eta|\geq\varepsilon,\\
q_{12}^{\frac{n}{n-2}}|\log q_{12}|,&\mu=\eta.
\end{cases}
    \]
\end{lemma}

\begin{proof}
In the notation of \cite[Proposition~B.2]{Figalli-Glaudo2020}, the
interaction parameter is
\[
 \varepsilon_{12}
 =
 \left(
 \frac{\lambda_1}{\lambda_2}
 +\frac{\lambda_2}{\lambda_1}
 +\lambda_1\lambda_2|z_1-z_2|^2
 \right)^{-\frac{n-2}{2}}
 =q_{12}.
\]
Thus the two alternatives in Lemma \ref{A.3} are exactly the
two-bubble interaction estimates in that proposition, with the same
uniform dependence on $n$ and $\epsilon$.
\end{proof}

\begin{lemma}\label{A.4}
Let $n\geq3$, $a,b\geq0$, and $a+b=n$.  Suppose
$
 \left\{(\lambda_{1,(k)},\lambda_{2,(k)},
 z_{1,(k)},z_{2,(k)})\right\}_{k=1}^{\infty}
$
satisfies $R_{12,(k)}\to\infty$.  Then
\begin{align*}
 I_{12,(k)}&:=
 \int_{\{U_{1,(k)}\geq U_{2,(k)}\}}
 \frac{\lambda_{1,(k)}^a}
 {(1+\lambda_{1,(k)}^2|y-z_{1,(k)}|^2)^a}
 \frac{\lambda_{2,(k)}^b}
 {(1+\lambda_{2,(k)}^2|y-z_{2,(k)}|^2)^b}\,dy.
\end{align*}
The two-bubble construction of Section~\ref{Sec2}, specifically
\eqref{confi2}--\eqref{confi3} and \eqref{confi7}--\eqref{confi8}, yields,
after restoring the original labels (with $i\succ j$ meaning $j\prec i$),
a subsequence on which exactly one of the following mutually exclusive
alternatives holds:
\begin{align}\label{eq:A4-order-greater}
 1\succ 2
 &\quad\Longleftrightarrow\quad
 \frac{\lambda_{1,(k)}}{\lambda_{2,(k)}}\longrightarrow\infty
 \quad\text{and}\quad
 z_{21,(k)}:=\lambda_{2,(k)}(z_{1,(k)}-z_{2,(k)})
 \longrightarrow\bar z_{21}\in\mathbb R^n,
\end{align}
$$
 1\prec 2
 \quad\Longleftrightarrow\quad
 \frac{\lambda_{2,(k)}}{\lambda_{1,(k)}}\longrightarrow\infty
 \quad\text{and}\quad
 z_{12,(k)}:=\lambda_{1,(k)}(z_{2,(k)}-z_{1,(k)})
 \longrightarrow\bar z_{12}\in\mathbb R^n,
$$
\begin{align}\label{eq:A4-incomparable}
 1\not\prec 2,\quad 1\not\succ 2
 \quad\Longleftrightarrow\quad
 \sqrt{\frac{\lambda_{1,(k)}}{\lambda_{2,(k)}}}
 +\sqrt{\frac{\lambda_{2,(k)}}{\lambda_{1,(k)}}}
 =o(R_{12,(k)}).
\end{align}
Moreover,
$$
 a>b\Longrightarrow
 I_{12,(k)}\approx_{n,a,b}R_{12,(k)}^{-2b},\qquad a<b\Longrightarrow
 I_{12,(k)}=O\!\left(R_{12,(k)}^{-n}\right)
$$
$$
 a=b=\frac n2\Longrightarrow
 I_{12,(k)}=O\!\left(R_{12,(k)}^{-n}
 (1+\log R_{12,(k)})\right).
$$
If $a>b$, put
\[
 \mathcal M(\xi):=\max\{1,|\xi|^2\},\qquad
 \mathcal J_a:=\int_{\mathbb R^n}\frac{dy}{(1+|y|^2)^a},
\qquad
 \mathcal J_{a,b}(\xi):=
 \int_{\mathbb R^n}\frac{dy}{(1+|y|^2)^a|y-\xi|^{2b}}.
\]
Along this subsequence,
\begin{equation}\label{eq:A4-original-limits}
 \lim_{k\to\infty}R_{12,(k)}^{2b}I_{12,(k)}
 =\begin{cases}
 \displaystyle
 \mathcal J_a\left(\frac{\mathcal M(\bar z_{21})}
 {1+|\bar z_{21}|^2}\right)^b,
 &1\succ 2,\quad z_{21,(k)}\to\bar z_{21},\\[3mm]
 \displaystyle
 \mathcal M(\bar z_{12})^b\mathcal J_{a,b}(\bar z_{12}),
 &1\prec 2,\quad z_{12,(k)}\to\bar z_{12},\\[2mm]
 \mathcal J_a,
 &1\not\prec 2\ \text{and}\ 1\not\succ 2.
 \end{cases}
\end{equation}

\end{lemma}
\begin{proof}
If \(a=b=n/2\), Lemma~\ref{A.3} and
\(q_{12}\asymp R_{12}^{-(n-2)}\) give
\begin{align*}
 I_{12}
 \leq C_n\int_{\mathbb R^n}
 U_1^{\frac n{n-2}}U_2^{\frac n{n-2}}\,dy\leq C_nq_{12}^{\frac n{n-2}}(1+|\log q_{12}|)
 \leq C_nR_{12}^{-n}(1+\log R_{12}).
\end{align*}
It remains to consider \(a\ne b\).  We suppress the index \(k\) and put
\(d_{12}:=|z_1-z_2|\).

\medskip
\noindent\textbf{Case A: \(1\succ2\).}
For all sufficiently large \(k\),
\[
 |y-z_1|\leq\frac1{2\sqrt{\lambda_1\lambda_2}}
 \quad\Longrightarrow\quad
 \lambda_1(1+\lambda_2^2|y-z_2|^2)
 \geq\lambda_1
 \geq\lambda_2+\frac{\lambda_1}{4}
 \geq\lambda_2(1+\lambda_1^2|y-z_1|^2),
\]
so \(U_1(y)\geq U_2(y)\).  If
$
 |y-z_1|\geq
 \frac{2(1+|z_{21}|)}{\sqrt{\lambda_1\lambda_2}}
$
and \(y\in B_{d_{12}/2}(z_1)\cup B_{d_{12}/2}(z_2)\), then
\[
 \lambda_1(1+\lambda_2^2|y-z_2|^2)
 \leq\lambda_1(1+4\lambda_2^2d_{12}^2)
 \leq\lambda_2(1+\lambda_1^2|y-z_1|^2).
\]
If the same lower bound for \(|y-z_1|\) holds and
\(y\notin B_{d_{12}/2}(z_1)\cup B_{d_{12}/2}(z_2)\), then
\(|y-z_2|\leq3|y-z_1|\), hence
\[
 \lambda_1(1+\lambda_2^2|y-z_2|^2)
 \leq\lambda_1+9\lambda_1\lambda_2^2|y-z_1|^2
 \leq\lambda_2(1+\lambda_1^2|y-z_1|^2).
\]
Therefore, with \(x=\lambda_1(y-z_1)\),
\begin{align}
 I_{12}
 &=\left(\frac{\lambda_2}{\lambda_1}\right)^b
 \int_{\{|x|\leq\frac12\sqrt{\lambda_1/\lambda_2}\}}
 \frac{dx}{(1+|x|^2)^a
 \bigl(1+|(\lambda_2/\lambda_1)(x-z_{12})|^2\bigr)^b}
 \notag\\
 &\quad+
 O\left(
 \left(\frac{\lambda_2}{\lambda_1}\right)^b
 \int_{\substack{\frac12\sqrt{\lambda_1/\lambda_2}\leq|x|\\
 |x|\leq2(1+|z_{21}|)\sqrt{\lambda_1/\lambda_2}}}
 \frac{dx}{(1+|x|^2)^a
 \bigl(1+|(\lambda_2/\lambda_1)(x-z_{12})|^2\bigr)^b}
 \right).
 \label{eq:A4-case-A-direct}
\end{align}
Since
$
 R_{12}^2
 =\max\left\{\frac{\lambda_1}{\lambda_2},
 \frac{\lambda_2}{\lambda_1},
 \frac{\lambda_1}{\lambda_2}|z_{21}|^2\right\}
 =\frac{\lambda_1}{\lambda_2}\mathcal M(z_{21})$,
dominated convergence in the first integral of
\eqref{eq:A4-case-A-direct}, together with the radial estimate on its
annular error, yields
\begin{equation}\label{eq:A4-case-A-limit}
 I_{12}=
 \begin{cases}
 \displaystyle
 \mathcal J_a
 \left(\frac{\mathcal M(\bar z_{21})}
 {1+|\bar z_{21}|^2}\right)^bR_{12}^{-2b}
 +o(R_{12}^{-2b}),&a>b,\\[2mm]
 O(R_{12}^{-n}),&a<b.
 \end{cases}
\end{equation}
For \(a<b\), the required radial bound is
\[
 \left(\frac{\lambda_2}{\lambda_1}\right)^b
 \int_0^{C\sqrt{\lambda_1/\lambda_2}}
 \frac{r^{n-1}}{(1+r^2)^a}\,dr
 \leq C\left(\frac{\lambda_2}{\lambda_1}\right)^{n/2}
 \leq CR_{12}^{-n}.
\]

\medskip
\noindent\textbf{Case B: \(1\prec2\).}
For all sufficiently large \(k\), if
$
 |y-z_2|\geq
 \frac{2(1+|z_{12}|)}{\sqrt{\lambda_1\lambda_2}}
$
and \(y\in B_{d_{12}/2}(z_1)\cup B_{d_{12}/2}(z_2)\), then
\[
 \lambda_1(1+\lambda_2^2|y-z_2|^2)
 \geq\lambda_1+4\lambda_2+4\lambda_1^2\lambda_2d_{12}^2
 \geq\lambda_2(1+\lambda_1^2|y-z_1|^2).
\]
If the same lower bound holds outside these two balls, then
\(|y-z_2|\geq|y-z_1|/3\), and
\[
 \lambda_1(1+\lambda_2^2|y-z_2|^2)
 \geq\lambda_1+2\lambda_2+
 \frac1{18}\lambda_1\lambda_2^2|y-z_1|^2
 \geq\lambda_2(1+\lambda_1^2|y-z_1|^2).
\]
Moreover,
\[
 |y-z_2|\leq\frac1{4\sqrt{\lambda_1\lambda_2}}
 \quad\Longrightarrow\quad
 \lambda_1(1+\lambda_2^2|y-z_2|^2)
 \leq\lambda_1+\frac{\lambda_2}{16}
 \leq\lambda_2(1+\lambda_1^2|y-z_1|^2).
\]
Consequently, after \(x=\lambda_1(y-z_1)\),
\begin{align}
 I_{12}
 &=\left(\frac{\lambda_2}{\lambda_1}\right)^b
 \int_{\{|x-z_{12}|\geq
 2(1+|z_{12}|)\sqrt{\lambda_1/\lambda_2}\}}
 \frac{dx}{(1+|x|^2)^a
 \bigl(1+|(\lambda_2/\lambda_1)(x-z_{12})|^2\bigr)^b}
 \notag\\
 &\quad+
 O\left(
 \left(\frac{\lambda_2}{\lambda_1}\right)^b
 \int_{\substack{\frac14\sqrt{\lambda_1/\lambda_2}
 \leq|x-z_{12}|\\
 |x-z_{12}|\leq
 2(1+|z_{12}|)\sqrt{\lambda_1/\lambda_2}}}
 \frac{dx}{(1+|x|^2)^a
 \bigl(1+|(\lambda_2/\lambda_1)(x-z_{12})|^2\bigr)^b}
 \right).
 \label{eq:A4-case-B-direct}
\end{align}
Here
\[
 R_{12}^2
 =\max\left\{\frac{\lambda_1}{\lambda_2},
 \frac{\lambda_2}{\lambda_1},
 \frac{\lambda_2}{\lambda_1}|z_{12}|^2\right\}
 =\frac{\lambda_2}{\lambda_1}\mathcal M(z_{12}).
\]
Multiplying \eqref{eq:A4-case-B-direct} by \(R_{12}^{2b}\), followed
by dominated convergence on
\(\mathbb R^n\setminus\{\bar z_{12}\}\) and the radial estimate for
the error annulus, gives
\begin{equation}\label{eq:A4-case-B-limit}
 I_{12}=
 \begin{cases}
 \displaystyle
 \mathcal M(\bar z_{12})^b
 \mathcal J_{a,b}(\bar z_{12})R_{12}^{-2b}
 +o(R_{12}^{-2b}),&a>b,\\[2mm]
 O(R_{12}^{-n}),&a<b.
 \end{cases}
\end{equation}

\medskip
\noindent\textbf{Case C: \(1\not\prec2\) and \(1\not\succ2\).}
Set
$
 L:=\sqrt{\lambda_1\lambda_2}\,d_{12}.
$
Then, for all sufficiently large \(k\), we have
\begin{equation}\label{eq:A4-case-C-L}
 \sqrt{\frac{\lambda_1}{\lambda_2}}
 +\sqrt{\frac{\lambda_2}{\lambda_1}}=o(L),
 \qquad R_{12}=L.
\end{equation}

\smallskip
\noindent\emph{Subcase C.1:
\(\lambda_2\leq\lambda_1\leq16\lambda_2\).}
If
$
 |y-z_1|\leq
 \frac{\sqrt{\lambda_2}\,d_{12}}{8\sqrt{\lambda_1}},
$
then \(|y-z_2|\geq d_{12}/2\), and
\[
 \lambda_1(1+\lambda_2^2|y-z_2|^2)
 \geq\lambda_1+\frac14\lambda_1\lambda_2^2d_{12}^2
 \geq\lambda_2(1+\lambda_1^2|y-z_1|^2).
\]
If \(|y-z_2|\leq d_{12}/4\), then
\(|y-z_1|\geq3d_{12}/4\), and
\[
 \lambda_1(1+\lambda_2^2|y-z_2|^2)
 \leq\lambda_1+\frac1{16}\lambda_1\lambda_2^2d_{12}^2
 \leq\lambda_2(1+\lambda_1^2|y-z_1|^2).
\]
Thus, after \(x=\lambda_1(y-z_1)\),
\begin{align}
 I_{12}
 &=\left(\frac{\lambda_2}{\lambda_1}\right)^b
 \int_{\{|x|\leq L/8\}}
 \frac{dx}{(1+|x|^2)^a
 \bigl(1+|(\lambda_2/\lambda_1)(x-z_{12})|^2\bigr)^b}
 \notag\\
 &\quad+
 O\left(
 \left(\frac{\lambda_2}{\lambda_1}\right)^b
 \int_{\substack{|x|\geq L/8\\
 |x-z_{12}|\geq\lambda_1d_{12}/4}}
 \frac{dx}{(1+|x|^2)^a
 \bigl(1+|(\lambda_2/\lambda_1)(x-z_{12})|^2\bigr)^b}
 \right).
 \label{eq:A4-case-C1-direct}
\end{align}

\smallskip
\noindent\emph{Subcase C.2: \(\lambda_1\geq16\lambda_2\).}
The first dominance inclusion in Subcase C.1 remains valid.  If
$
 |y-z_1|\geq
 \frac{2\sqrt{\lambda_2}\,d_{12}}{\sqrt{\lambda_1}}
$
and \(y\in B_{d_{12}}(z_1)\cup B_{d_{12}}(z_2)\), then
\[
 \lambda_1(1+\lambda_2^2|y-z_2|^2)
 \leq\lambda_1+\frac94\lambda_1\lambda_2^2d_{12}^2
 \leq\lambda_2(1+\lambda_1^2|y-z_1|^2).
\]
For \(y\notin B_{d_{12}}(z_1)\cup B_{d_{12}}(z_2)\),
\(|y-z_2|\leq3|y-z_1|\); by \eqref{eq:A4-case-C-L},
for all sufficiently large \(k\),
\[
 \lambda_1(1+\lambda_2^2|y-z_2|^2)
 \leq\lambda_1+9\lambda_1\lambda_2^2|y-z_1|^2
 \leq\lambda_2(1+\lambda_1^2|y-z_1|^2).
\]
Hence
\begin{align*}
 I_{12}
 &=\left(\frac{\lambda_2}{\lambda_1}\right)^b
 \int_{\{|x|\leq L/8\}}
 \frac{dx}{(1+|x|^2)^a
 \bigl(1+|(\lambda_2/\lambda_1)(x-z_{12})|^2\bigr)^b}
 \notag\\
 &\quad+
 O\left(
 \left(\frac{\lambda_2}{\lambda_1}\right)^b
 \int_{\{L/8\leq|x|\leq2L\}}
 \frac{dx}{(1+|x|^2)^a
 \bigl(1+|(\lambda_2/\lambda_1)(x-z_{12})|^2\bigr)^b}
 \right).
 \label{eq:A4-case-C2-direct}
\end{align*}

\smallskip
\noindent\emph{Subcase C.3:
\(\lambda_1\leq\lambda_2\leq16\lambda_1\).}
If \(|y-z_1|\leq d_{12}/4\), then
\(|y-z_2|\geq3d_{12}/4\), and
\[
 \lambda_1(1+\lambda_2^2|y-z_2|^2)
 \geq\lambda_1+\frac9{16}\lambda_1\lambda_2^2d_{12}^2
 \geq\lambda_2(1+\lambda_1^2|y-z_1|^2).
\]
If \(|y-z_2|\leq d_{12}/64\), then
\(|y-z_1|\geq d_{12}/2\), and
\[
 \lambda_1(1+\lambda_2^2|y-z_2|^2)
 \leq\lambda_1+\frac1{64^2}\lambda_1\lambda_2^2d_{12}^2
 \leq\lambda_2(1+\lambda_1^2|y-z_1|^2).
\]
It follows that
\begin{align*}
 I_{12}
 &=\left(\frac{\lambda_2}{\lambda_1}\right)^b
 \int_{\{|x|\leq\lambda_1d_{12}/4\}}
 \frac{dx}{(1+|x|^2)^a
 \bigl(1+|(\lambda_2/\lambda_1)(x-z_{12})|^2\bigr)^b}
 \notag\\
 &\quad+
 O\left(
 \left(\frac{\lambda_2}{\lambda_1}\right)^b
 \int_{\substack{|x|\geq\lambda_1d_{12}/4\\
 |x-z_{12}|\geq\lambda_1d_{12}/64}}
 \frac{dx}{(1+|x|^2)^a
 \bigl(1+|(\lambda_2/\lambda_1)(x-z_{12})|^2\bigr)^b}
 \right).
 \label{eq:A4-case-C3-direct}
\end{align*}

\smallskip
\noindent\emph{Subcase C.4: \(\lambda_2\geq16\lambda_1\).}
If
$
 |y-z_2|\geq
 \frac{2\sqrt{\lambda_1}\,d_{12}}{\sqrt{\lambda_2}}
$
and \(y\in B_{d_{12}}(z_1)\cup B_{d_{12}}(z_2)\), then
\[
 \lambda_1(1+\lambda_2^2|y-z_2|^2)
 \geq\lambda_1+4\lambda_1\lambda_2^2d_{12}^2
 \geq\lambda_2(1+\lambda_1^2|y-z_1|^2).
\]
Outside these two balls, \(|y-z_1|\leq3|y-z_2|\); hence, for all
sufficiently large \(k\),
\[
 \lambda_1(1+\lambda_2^2|y-z_2|^2)
 \geq\lambda_2+9\lambda_2\lambda_1^2|y-z_2|^2
 \geq\lambda_2(1+\lambda_1^2|y-z_1|^2).
\]
If
$
 |y-z_2|\leq
 \frac{\sqrt{\lambda_1}\,d_{12}}{8\sqrt{\lambda_2}},
$
then \(|y-z_1|\geq d_{12}/2\), and
\[
 \lambda_1(1+\lambda_2^2|y-z_2|^2)
 \leq\lambda_2+\frac14\lambda_2\lambda_1^2d_{12}^2
 \leq\lambda_2(1+\lambda_1^2|y-z_1|^2).
\]
Therefore
\begin{align}
 I_{12}
 &=\left(\frac{\lambda_2}{\lambda_1}\right)^b
 \int_{\{|x-z_{12}|\geq
 2\lambda_1^{3/2}\lambda_2^{-1/2}d_{12}\}}
 \frac{dx}{(1+|x|^2)^a
 \bigl(1+|(\lambda_2/\lambda_1)(x-z_{12})|^2\bigr)^b}
 \notag\\
 &\quad+
 O\left(
 \left(\frac{\lambda_2}{\lambda_1}\right)^b
 \int_{\substack{\lambda_1^{3/2}\lambda_2^{-1/2}d_{12}/8
 \leq|x-z_{12}|\\
 |x-z_{12}|\leq
 2\lambda_1^{3/2}\lambda_2^{-1/2}d_{12}}}
 \frac{dx}{(1+|x|^2)^a
 \bigl(1+|(\lambda_2/\lambda_1)(x-z_{12})|^2\bigr)^b}
 \right).
 \label{eq:A4-case-C4-direct}
\end{align}

In each of \eqref{eq:A4-case-C1-direct}--\eqref{eq:A4-case-C4-direct},
\eqref{eq:A4-case-C-L}, \(a+b=n\), and the radial decomposition about
\(z_1\) and \(z_2\) give
\begin{equation}\label{eq:A4-case-C-limit}
 I_{12}=
 \begin{cases}
 \displaystyle
 \mathcal J_aR_{12}^{-2b}+o(R_{12}^{-2b}),&a>b,\\[2mm]
 O(R_{12}^{-n}),&a<b.
 \end{cases}
\end{equation}
Indeed, \eqref{eq:A4-case-C-L} implies
\(\lambda_1d_{12}\to\infty\) and \(\lambda_2d_{12}\to\infty\), and
for every fixed \(x\in\mathbb R^n\),
\[
 R_{12}^{2b}\left(\frac{\lambda_2}{\lambda_1}\right)^b
 \left(1+
 \left|\frac{\lambda_2}{\lambda_1}(x-z_{12})\right|^2
 \right)^{-b}\longrightarrow1
\]
locally uniformly in \(x\).  If \(a>b\), the principal region in each
subcase exhausts every fixed ball, while the corresponding error
region contributes \(o(R_{12}^{-2b})\).  If \(a<b\), splitting each
error integral at distance \(R_{12}\) gives
\[
 C R_{12}^{-2b}\int_0^{R_{12}}
 \frac{r^{n-1}}{(1+r^2)^a}\,dr
 +C\int_{R_{12}}^\infty r^{-n-1}\,dr
 \leq CR_{12}^{-n}.
\]

Equations \eqref{eq:A4-case-A-limit}, \eqref{eq:A4-case-B-limit}, and
\eqref{eq:A4-case-C-limit} prove \eqref{eq:A4-original-limits} and all
three estimates in the statement.  The limiting coefficients in the
case \(a>b\) are positive.  The three subsequence alternatives are
exhaustive; hence a contradiction argument applied to any sequence on
which either side of the asserted two-sided estimate fails yields
\(I_{12}\asymp_{n,a,b}R_{12}^{-2b}\).
\end{proof}

\begin{lemma}\label{A.5}
     Let $n\geq3$ and let $U_i$ be an Aubin--Talenti bubble. For the
     functions $Z_i^a$ defined in \eqref{eq2.1}, there exist constants
     $\gamma^a=\gamma^a(n)>0$ such that
\[
\int_{\mathbb R^n} U_i^{\frac{4}{n-2}} Z_i^a Z_i^b\,\dx
= \begin{cases}0, & \text { if } a \neq b, \\ \gamma^a, & \text { if } 1 \leq a=b \leq n+1,\end{cases}
\]
If $U_j$ is another Aubin--Talenti bubble, $i\ne j$, and
$1\leq a,b\leq n+1$, then
\[
\left|\int_{\mathbb R^n} U_i^{\frac{4}{n-2}} Z_i^a Z_j^b\,\dx\right| \lesssim q_{i j}.
\]
\end{lemma}

\begin{proof}
To compare normalizations with
\cite[Formulas~(F1)--(F6)]{Bahri1989}, set
\[
 \delta_i=\lambda_i^{-1},\qquad a_i=z_i.
\]
Under this identification, the translation and scale tangent fields
used there agree with
\[
 Z_i^a=\lambda_i^{-1}\partial_{z_i^a}U[z_i,\lambda_i],
 \qquad
 Z_i^{n+1}=\lambda_i\partial_{\lambda_i}U[z_i,\lambda_i],
\]
up to nonzero dimensional constants and signs.  These constants are
absorbed into $\gamma^a(n)>0$ and the implicit constant in the
cross-interaction estimate.  The stated conclusions therefore follow
from those formulas with the present normalization and with interaction
parameter $q_{ij}$.
\end{proof}

\begin{lemma}\label{A.6}
Assume that $n\geq3$ and that
$\{U_{1,(k)},U_{2,(k)}\}_{k=1}^\infty$ is a two-bubble sequence satisfying
$R_{12,(k)}\to\infty$.  If
\[
 \vartheta:=\frac{n+2}{2(n-2)},
\]
then
\begin{equation}\label{eq:A6-two-sided}
 \bigl\|U_{1,(k)}^\vartheta U_{2,(k)}^\vartheta\bigr\|_{\dot H^{-1}}
 \asymp_n R_{12,(k)}^{-\frac{n+2}{2}}
 \bigl(\log R_{12,(k)}\bigr)^{1/2}.
\end{equation}
\end{lemma}
\begin{proof}
We use the two-case, four-region decomposition of the Newtonian energy.
For readability we suppress the subscript $(k)$.  Put
    \[
    m:=\frac{n+2}{2},\qquad
    \mathfrak b_n:=\kappa_n[n(n-2)]^m.
    \]
Set \(s:=2n/(n+2)\) and
\(f:=U_1^\vartheta U_2^\vartheta\).  Since
\(2s\vartheta=2^*\), H\"older inequality gives \(f\in L^s\).
By the standard Newtonian representation of
\((-\Delta)^{-1}:L^s(\mathbb R^n)\to\dot H^1(\mathbb R^n)\)
\cite[Section~1.2.2]{AdamsHedberg} and the
Hardy--Littlewood--Sobolev inequality \cite{Lieb},
\begin{equation*}
\begin{split}
 \|U_1^\vartheta U_2^\vartheta\|_{\dot H^{-1}}^2
 =\int_{\mathbb R^n}f(-\Delta)^{-1}f\,\dx=\kappa_n\iint_{\mathbb R^n\times\mathbb R^n}
 \frac{U_1^\vartheta(x)U_2^\vartheta(x)
 U_1^\vartheta(y)U_2^\vartheta(y)}
 {|x-y|^{n-2}}\,\dx\dy .
\end{split}
\end{equation*}
Put
\[
 G(x):=(1+|x|^2)^{-m/2},\qquad
 \mathscr B(f,g):=\iint_{\mathbb R^n\times\mathbb R^n}
 \frac{f(x)g(y)}{|x-y|^{n-2}}\,dx\,dy.
\]
Newton's shell formula gives the standard estimate
\begin{equation}\label{eq:A6-annulus-energy}
 \mathscr B(G\mathbf1_{B_H},G\mathbf1_{B_H})
 \asymp_n\log H,\qquad H\geq8.
\end{equation}

    By Lemma~\ref{A.4}, every subsequence has a further subsequence satisfying
    one of the three alternatives \eqref{eq:A4-order-greater}--
    \eqref{eq:A4-incomparable}.  Since the Newtonian energy is invariant under
    interchanging the two bubbles, it is enough to treat one ordered
    configuration and the incomparable configuration; the reverse ordered
    configuration is then obtained by exchanging the indices.

    \noindent\textbf{Case 1. $1\succ2$.}
    Put
    \[
     \rho:=\frac{\lambda_2}{\lambda_1},\qquad
     \zeta:=z_{21},\qquad
     H:=\frac{4(1+|\zeta|)}{\rho},\qquad P:=B_H(0).
    \]
    Then \(\rho\to0\), \(\zeta\to\bar z_{21}\), and
    \(R_{12}^2=\rho^{-1}\max\{1,|\zeta|^2\}\).  After the change
    \(x_1=\lambda_1(x-z_1)\), set
    \[
     F_{\rho,\zeta}(x):=
     \frac1{(1+|x|^2)^{m/2}(1+|\rho x+\zeta|^2)^{m/2}}.
    \]
    Decompose
    \begin{align*}
     \mathfrak b_n\rho^m
     \left(\iint_{P\times P}+\iint_{P\times P^c}
     +\iint_{P^c\times P}+\iint_{P^c\times P^c}\right)
     \frac{F_{\rho,\zeta}(x)F_{\rho,\zeta}(y)}
     {|x-y|^{n-2}}\,dx\,dy=:I_1+I_2+I_3+I_4.
    \end{align*}
    On \(P\), \(F_{\rho,\zeta}\asymp_{\mathscr U}G\); on \(P^c\),
    \(F_{\rho,\zeta}(x)\leq C\rho^{-m}|x|^{-2m}\).  Hence
    \begin{align*}
     I_1\asymp_{n,\mathscr U}\rho^m\log H
          \asymp_{n,\mathscr U}R_{12}^{-2m}\log R_{12},\end{align*}
          \begin{align*}
     I_2+I_3\leq
          CH^{n-m}\int_H^\infty r^{1-2m}\,dr
          \leq C_{n,\mathscr U}R_{12}^{-2m},\qquad
     I_4\leq C\rho^{-m}H^{-2m}
          \leq C_nR_{12}^{-2m}.
    \end{align*}
    Thus
    \[
     \|U_1^\vartheta U_2^\vartheta\|_{\dot H^{-1}}^2
     \asymp_{n,\mathscr U}R_{12}^{-2m}\log R_{12}.
    \]

    \noindent\textbf{Case 2. $1\not\succ2$ and $1\not\prec2$.}
    Put
    \[
     \rho:=\frac{\lambda_2}{\lambda_1},\qquad
     w:=z_{12},\qquad d:=|w|,\qquad L:=\sqrt{\rho}\,d .
    \]
    Then \(d\to\infty\), \(\rho d\to\infty\), and \(R_{12}=L\).  Set
    \[
     F_{\rho,w}(x):=
     \frac1{(1+|x|^2)^{m/2}(1+\rho^2|x-w|^2)^{m/2}},
     \qquad P:=B_{4d}(0),
    \]
    and decompose
    \begin{align*}
     \mathfrak b_n\rho^m
     \left(\iint_{P\times P}+\iint_{P\times P^c}
     +\iint_{P^c\times P}+\iint_{P^c\times P^c}\right)
     \frac{F_{\rho,w}(x)F_{\rho,w}(y)}
     {|x-y|^{n-2}}\,dx\,dy=:I_1+I_2+I_3+I_4.
    \end{align*}
    For
    \[
     P_0:=B_{d/4}(0),\qquad
     P_1:=B_{d/4}(w),\qquad
     P_m:=P\setminus(P_0\cup P_1),
    \]
    one has
    \[
     F_{\rho,w}\asymp(\rho d)^{-m}G\quad\hbox{on }P_0,
     \qquad
     F_{\rho,w}(x)\asymp d^{-m}G(\rho(x-w))\quad\hbox{on }P_1,
     \qquad
     F_{\rho,w}\leq CL^{-2m}\quad\hbox{on }P_m.
    \]
    Splitting the \(y\)-integral in \(I_1\) over \(P_0,P_1,P_m\) and
    using \eqref{eq:A6-annulus-energy} gives
    \[
     I_1\asymp_n
     L^{-2m}\bigl(\log d+\log(\rho d)\bigr)
     =2L^{-2m}\log L.
    \]
    On \(P^c\),
    \[
     F_{\rho,w}(y)\leq C\rho^{-m}|y|^{-2m},\qquad
     \int_PF_{\rho,w}\leq C\rho^{-m}d^{-2}.
    \]
    Hence Newton's shell formula and scaling give
    \begin{align*}
     I_2+I_3
     \leq Cd^{-n}\int_PF_{\rho,w}(x)\,dx
       \leq C\rho^{-m}d^{-n-2}=CL^{-2m},\qquad
     I_4
     \leq C\rho^{-m}d^{\,n+2-4m}
       =CL^{-2m}.
    \end{align*}
    Therefore
    \[
     \|U_1^\vartheta U_2^\vartheta\|_{\dot H^{-1}}^2
     \asymp_nR_{12}^{-2m}\log R_{12}.
    \]
    The reverse ordered case follows by exchanging the indices.
    Lemma~\ref{A.4} and the subsequence principle make the constants
    uniform.  Taking square roots proves \eqref{eq:A6-two-sided}.
\end{proof}

\section{Projection Expansions}\label{appB}
\begin{lemma}\label{le:weak-projection-expansion}
Let $n\geq3$, $0<\alpha<n$, $\alpha\leq4$, $\nu\geq2$, and let
$\{U_i:1\leq i\leq\nu\}$ be a $\delta$-interacting family. Set
$\sigma=\sum_{i=1}^\nu U_i$, let $h$ be defined by \eqref{h}, and set
$R=\frac12\min_{i\ne j}R_{ij}$. Then, for
$1\leq r\leq\nu$,
\begin{equation}\label{eq:weak-range-master-expansion}
\begin{split}
 \int_{\mathbb R^n}hZ_r^{n+1}\dx
 =\sum_{j\ne r}\int_{\mathbb R^n}
  U_j^{2^*-p_\alpha}U_r^{p_\alpha-1}Z_r^{n+1}\dx
 +(2^*-1)\sum_{j\ne r}\int_{\mathbb R^n}
  U_r^{2^*-2}U_jZ_r^{n+1}\dx
 +\mathcal R_{{\rm weak},r}.
\end{split}
\end{equation}
There is a function
$\omega_7:[2,\infty)\to[0,\infty)$, depending only on
$(n,\alpha,\nu)$, such that
\[
 \omega_7(R)\longrightarrow0,
 \qquad
 |\mathcal R_{{\rm weak},r}|
 \leq\omega_7(R)R^{-a_\alpha}
\]
as $R\to\infty$, uniformly for all $r$ and all admissible bubble
families.
\end{lemma}

\begin{proof}
Put
\[
\beta=p_\alpha-1,\qquad\gamma=2^*-p_\alpha,\qquad
F=\sigma^{p_\alpha}-\sum_iU_i^{p_\alpha},
\qquad
A=\sum_iU_i^\gamma(\sigma^\beta-U_i^\beta),\qquad
B=\mathcal I_\alpha(F)\sigma^\beta.
\]
Then \(h=A+B\).  Let
\(\Omega_\ell=\{U_\ell=\max_iU_i\}\) and \(Z_r=Z_r^{n+1}\).
All constants below depend only on \((n,\alpha,\nu)\).  Every estimate
obtained from Lemma \ref{A.4} is uniform in the bubble parameters:
otherwise a sequence with \(R_{ij}\to\infty\) would contradict one of
the three alternatives in that lemma.
On \(\Omega_r\), Taylor's formula gives
\[
\begin{aligned}
U_r^\gamma(\sigma^\beta-U_r^\beta)
&=\beta U_r^{2^*-2}\sum_{j\ne r}U_j
+O\!\left(U_r^{2^*-3}\Big(\sum_{j\ne r}U_j\Big)^2\right),\\
U_j^\gamma(\sigma^\beta-U_j^\beta)
&=U_j^\gamma U_r^\beta
+O(U_j^{2^*-1})
+O(U_j^{\gamma+1}U_r^{\beta-1})+O\!\left(\sum_{m\ne r,j}U_j^\gamma U_mU_r^{\beta-1}\right).
\end{aligned}
\]
Indeed,
\[
\sigma^\beta
=U_r^\beta+
O\!\left(U_r^{\beta-1}\sum_{m\ne r}U_m\right)
+O\!\left(\sum_{m\ne r}U_m^\beta\right)
\quad\hbox{on }\Omega_r.
\]
Since \(|Z_r|\leq CU_r\), the integral of the Taylor remainder on
\(\Omega_r\) is bounded by
\[
\begin{aligned}
\mathcal E_{A,r}^{\rm in}:={}&
C\sum_{j,m\ne r}\int_{\Omega_r}U_r^{2^*-2}U_jU_m\dx
+C\sum_{j\ne r}\int_{\Omega_r}U_j^{2^*-1}U_r\dx\\
&+C\sum_{j\ne r}\int_{\Omega_r}U_j^{\gamma+1}U_r^\beta\dx
+C\sum_{j\ne r}\sum_{\substack{m\ne r\\m\ne j}}
\int_{\Omega_r}U_j^\gamma U_mU_r^\beta\dx .
\end{aligned}
\]
Put
\[
\mathcal Y_n(R)=
\begin{cases}
R^{-2},&n=3,\\
R^{-4}\log R,&n=4,\\
R^{-n},&n\geq5,
\end{cases}
\qquad
\mathcal S_{n,\alpha}(R)=
\begin{cases}
R^{-(n-2+\alpha)},&0<\alpha<2,\\
R^{-n}\log R,&\alpha=2,\\
R^{-n},&2<\alpha\leq4.
\end{cases}
\]
For the first sum, \(U_jU_m\leq(U_j^2+U_m^2)/2\), and Lemma
\ref{A.4}, with the scaled exponent pair \((2,n-2)\), gives
\[
\int_{\Omega_r}U_r^{2^*-2}U_j^2\dx
\leq C
\begin{cases}
R_{rj}^{-2},&n=3,\\
R_{rj}^{-4}\log R_{rj},&n=4,\\
R_{rj}^{-n},&n\geq5.
\end{cases}
\]
The second and third sums correspond respectively to
\[
\left(\frac{n-2}{2},\frac{n+2}{2}\right),
\qquad
\left(\frac{n+2-\alpha}{2},
\frac{n+\alpha-2}{2}\right).
\]
Consequently,
\[
\int_{\Omega_r}U_j^{2^*-1}U_r\dx\leq CR_{rj}^{-n},
\qquad
\int_{\Omega_r}U_j^{\gamma+1}U_r^\beta\dx
\leq C
\begin{cases}
R_{rj}^{-(n-2+\alpha)},&0<\alpha<2,\\
R_{rj}^{-n}\log R_{rj},&\alpha=2,\\
R_{rj}^{-n},&2<\alpha\leq4.
\end{cases}
\]
For the fourth sum,
\[
U_j^\gamma U_m
\leq C\bigl(U_j^{\gamma+1}+U_m^{\gamma+1}\bigr),
\]
so it satisfies the same \(\mathcal S_{n,\alpha}\)-bound.  Hence
\[
\mathcal E_{A,r}^{\rm in}
\leq C\bigl(\mathcal Y_n(R)+\mathcal S_{n,\alpha}(R)+R^{-n}\bigr).
\]
On \(\Omega_\ell\), \(\ell\ne r\),
\[
0\leq A\leq C\left(
U_\ell^{2^*-2}\sum_{j\ne\ell}U_j+
\sum_{j\ne\ell}U_j^\gamma U_\ell^\beta\right).
\]
Since \(|Z_r|\leq CU_r\),
\[
\begin{aligned}
\left|\int_{\Omega_r^c}AZ_r\dx\right|
&\leq C\sum_{\ell\ne r}\int_{\Omega_\ell}
U_\ell^{2^*-2}U_r^2\dx+C\sum_{\ell\ne r}\sum_{\substack{j\ne\ell\\j\ne r}}
\int_{\Omega_\ell}U_\ell^{2^*-2}U_jU_r\dx\\
&\quad+C\sum_{\ell\ne r}\int_{\Omega_\ell}
U_r^{\gamma+1}U_\ell^\beta\dx+C\sum_{\ell\ne r}\sum_{\substack{j\ne\ell\\j\ne r}}
\int_{\Omega_\ell}U_j^\gamma U_\ell^\beta U_r\dx .
\end{aligned}
\]
The two pointwise inequalities
\[
U_jU_r\leq\frac12(U_j^2+U_r^2),\qquad
U_j^\gamma U_r\leq
C(U_j^{\gamma+1}+U_r^{\gamma+1})
\]
reduce these integrals to the scaled exponent pairs
\[
(2,n-2),\qquad
\left(\frac{n+2-\alpha}{2},
\frac{n+\alpha-2}{2}\right).
\]
Lemma \ref{A.4} therefore gives
\[
\left|\int_{\Omega_r^c}AZ_r\dx\right|
\leq C\bigl(\mathcal Y_n(R)+\mathcal S_{n,\alpha}(R)\bigr).
\]
The parts of the two leading terms outside \(\Omega_r\) satisfy
\[
\begin{aligned}
&\sum_{j\ne r}\left|\int_{\Omega_r^c}
U_j^\gamma U_r^\beta Z_r\dx\right|
+\sum_{j\ne r}\left|\int_{\Omega_r^c}
U_r^{2^*-2}U_jZ_r\dx\right|\\
\leq& C\sum_{\ell\ne r}\int_{\Omega_\ell}
\bigl(U_\ell^\gamma U_r^{\beta+1}
+U_\ell U_r^{2^*-1}\bigr)\dx+C\sum_{\ell\ne r}\sum_{\substack{j\ne\ell\\j\ne r}}
\int_{\Omega_\ell}\bigl(
U_j^\gamma U_\ell^\beta U_r+
U_\ell^{2^*-2}U_jU_r\bigr)\dx .
\end{aligned}
\]
In addition to the preceding two pairs, the last display uses
\[
\left(\frac{n-2}{2},\frac{n+2}{2}\right),\qquad
\left(\frac{\alpha}{2},n-\frac{\alpha}{2}\right).
\]
Thus Lemma \ref{A.4} bounds it by
\(C(\mathcal Y_n+\mathcal S_{n,\alpha}+R^{-n})\).  Combining the
estimates on \(\Omega_r\) and \(\Omega_r^c\) gives
\begin{equation}\label{eq:new26-weak-A}
\begin{aligned}
\int AZ_r\dx
={}&\sum_{j\ne r}\int U_j^\gamma U_r^\beta Z_r\dx
+\beta\sum_{j\ne r}\int U_r^{2^*-2}U_jZ_r\dx
+\mathcal R_{A,r}^{\rm weak},\\
|\mathcal R_{A,r}^{\rm weak}|
\leq{}&C\bigl(\mathcal Y_n(R)+
\mathcal S_{n,\alpha}(R)+R^{-n}\bigr).
\end{aligned}
\end{equation}
Since \(a_\alpha=\min\{\alpha,n-2\}\),
$
R^{a_\alpha}\bigl(\mathcal Y_n(R)+
\mathcal S_{n,\alpha}(R)+R^{-n}\bigr)\longrightarrow0$.
For \(B_r=\int BZ_r\dx\), put
\[
S_r=\sum_{i\ne r}U_i,\qquad
G_r=F-p_\alpha U_r^{p_\alpha-1}S_r.
\]
Self-adjointness of \(\mathcal I_\alpha\) and differentiation of the
single-bubble identity give
\[
B_r=\gamma\int U_r^{2^*-2}S_rZ_r\dx
+\frac{\gamma}{p_\alpha}\int G_rU_r^{\gamma-1}Z_r\dx
+\mathcal R_B,
\quad
\mathcal R_B=\int\mathcal I_\alpha(F)
(\sigma^\beta-U_r^\beta)Z_r\dx.
\]
Although \(\gamma-1\) may be negative,
\[
U_r^{\gamma-1}|Z_r|\leq CU_r^\gamma.
\]
On \(\Omega_r\), Taylor's formula and \(S_r\leq(\nu-1)U_r\) give
\[
|G_r|\leq C\left(
U_r^{p_\alpha-2}S_r^2+\sum_{i\ne r}U_i^{p_\alpha}\right),
\]
and therefore
\[
\begin{aligned}
\int_{\Omega_r}|G_r|U_r^\gamma\dx
\leq C\sum_{i\ne r}\int_{\Omega_r}
U_r^{2^*-2}U_i^2\dx\leq C\mathcal Y_n(R).
\end{aligned}
\]
Fix \(\ell\ne r\).  On \(\Omega_\ell\),
\[
0\leq F\leq CU_\ell^\beta S_\ell,
\qquad S_\ell=\sum_{i\ne\ell}U_i.
\]
Since \(S_r\leq\nu U_\ell\) there, the pointwise estimate
\[
\begin{aligned}
|G_r|U_r^\gamma
\leq CU_\ell^\beta S_\ell U_r^\gamma
+CU_r^{\beta+\gamma}S_r\leq C\sum_{k\ne\ell}U_\ell^\beta U_k^{\gamma+1}
+CU_\ell U_r^{2^*-1}
\end{aligned}
\]
follows from
\(U_kU_r^\gamma\leq
C(U_k^{\gamma+1}+U_r^{\gamma+1})\).
The scaled exponent pairs are
\[
\left(\frac{n+2-\alpha}{2},
\frac{n+\alpha-2}{2}\right),
\qquad
\left(\frac{n-2}{2},\frac{n+2}{2}\right).
\]
Lemma \ref{A.4} yields
$
\sum_{\ell\ne r}\int_{\Omega_\ell}|G_r|U_r^\gamma\dx
\leq C\bigl(\mathcal S_{n,\alpha}(R)+R^{-n}\bigr)$.
Thus
\[
\left|\frac{\gamma}{p_\alpha}
\int G_rU_r^{\gamma-1}Z_r\dx\right|
\leq C\bigl(\mathcal Y_n(R)+
\mathcal S_{n,\alpha}(R)+R^{-n}\bigr).
\]

It remains to estimate \(\mathcal R_B\).  On \(\Omega_\ell\),
$
|F|\leq CU_\ell^\beta\sum_{i\ne\ell}U_i$.
Moreover,
\[
|\sigma^\beta-U_r^\beta|\,|Z_r|
\leq
\begin{cases}
CU_r^\beta\sum_{m\ne r}U_m,&\Omega_r,\\
CU_k^\beta U_r,&\Omega_k,\quad k\ne r.
\end{cases}
\]
Let \(s=2n/(2n-\alpha)\).  Applying the bilinear HLS inequality on
each \(\Omega_\ell\times\Omega_k\) block gives
\[
\begin{aligned}
|\mathcal R_B|
\leq C\sum_{\ell}\sum_{i\ne\ell}\sum_{m\ne r}
\|U_\ell^\beta U_i\|_{L^s(\Omega_\ell)}
\|U_r^\beta U_m\|_{L^s(\Omega_r)}+C\sum_{\ell}\sum_{i\ne\ell}\sum_{k\ne r}
\|U_\ell^\beta U_i\|_{L^s(\Omega_\ell)}
\|U_k^\beta U_r\|_{L^s(\Omega_k)}.
\end{aligned}
\]
Here
\[
(\beta+1)s=2^*,\qquad
\frac{n-2}{2}(\beta+1)s=n,\qquad
\beta=1\Longleftrightarrow\alpha=4.
\]
Lemma \ref{A.4} therefore gives
\[
\|U_\ell^\beta U_i\|_{L^s(\Omega_\ell)}
\leq CR_{\ell i}^{-(n-2)}
(1+\log R_{\ell i})^{s^{-1}\mathbf1_{\{\alpha=4\}}}.
\]
and hence
\[
|\mathcal R_B|
\leq C\mathcal H_{n,\alpha}(R),\qquad
\mathcal H_{n,\alpha}(R):=
R^{-2(n-2)}
(1+\log R)^{2s^{-1}\mathbf1_{\{\alpha=4\}}}.
\]
Consequently,
\begin{equation}\label{eq:new26-weak-B}
B_r=\gamma\sum_{j\ne r}\int
U_r^{2^*-2}U_jZ_r\dx+\mathcal R_{B,r}^{\rm weak},\qquad
|\mathcal R_{B,r}^{\rm weak}|
\leq C\bigl(\mathcal Y_n(R)+\mathcal S_{n,\alpha}(R)
+R^{-n}+\mathcal H_{n,\alpha}(R)\bigr).
\end{equation}
Since \(a_\alpha<n\) and \(2(n-2)>a_\alpha\),
$
R^{a_\alpha}\bigl(\mathcal Y_n+\mathcal S_{n,\alpha}
+R^{-n}+\mathcal H_{n,\alpha}\bigr)\longrightarrow0.
$
Define
\[
\omega_7(R)
:=C R^{a_\alpha}\bigl(\mathcal Y_n(R)+
\mathcal S_{n,\alpha}(R)+R^{-n}+
\mathcal H_{n,\alpha}(R)\bigr).
\]
Then \(\omega_7(R)\to0\), uniformly in \(r\) and in all
bubble parameters.  Since \(\beta+\gamma=2^*-1\), adding
\eqref{eq:new26-weak-A} and \eqref{eq:new26-weak-B} proves
\eqref{eq:weak-range-master-expansion}.
\end{proof}

In the strongly singular range $4<\alpha<n$, the localized HLS estimates needed
for the scale projection are combined with the self-adjointness identity for
the Riesz potential in the proof of Lemma
\ref{le:strong-projection-proof}.

Assume $4<\alpha<n$ and use the decomposition $h=A+B$.  Put
\[
    \beta=p_\alpha-1=\frac{n+2-\alpha}{n-2}\in(0,1),\qquad
    \gamma=2^*-p_\alpha=\frac{\alpha}{n-2}.
\]
Thus $\gamma+\beta=2^*-1$ and $\gamma+\beta-1=2^*-2$.

\begin{lemma}\label{le:strong-projection-proof}
Assume that $4<\alpha<n$, $\nu\geq2$, and
$\{U_i:1\leq i\leq\nu\}$ is a $\delta$-interacting family. Set
$\sigma=\sum_{i=1}^\nu U_i$, let $h$ be defined by \eqref{h}, and set
$R=\frac12\min_{i\ne j}R_{ij}$. Then
\eqref{eq:strong-projection-expansion} holds for every
$r\in\{1,\ldots,\nu\}$.  More precisely, there is a function
$\omega_8:[1,\infty)\to[0,\infty)$, depending only on
$n,\alpha,\nu$, such that
\[
 \omega_8(R)\longrightarrow0
 \quad\text{as }R\longrightarrow\infty
\]
and the absolute value of the remainder in
\eqref{eq:strong-projection-expansion} is bounded by
$\omega_8(R)R^{-a_\alpha}$, uniformly for all $r$ and all admissible
bubble families.
\end{lemma}

\begin{proof}
Throughout the proof \(R\geq2\), and all constants depend only on
\((n,\alpha,\nu)\).  As in Lemma \ref{le:weak-projection-expansion},
every use of Lemma \ref{A.4} is uniform in the ordered pair of bubbles.
Write \(Z_r=Z_r^{n+1}\),
\[
F=\sigma^{p_\alpha}-\sum_iU_i^{p_\alpha},\qquad
A=\sum_iU_i^\gamma(\sigma^\beta-U_i^\beta),\qquad
B=\mathcal I_\alpha(F)\sigma^\beta.
\]
On \(\Omega_r=\{U_r=\max_iU_i\}\),
\[
\begin{aligned}
U_r^\gamma(\sigma^\beta-U_r^\beta)
&=\beta U_r^{2^*-2}\sum_{j\ne r}U_j
+O\!\left(U_r^{2^*-3}\Big(\sum_{j\ne r}U_j\Big)^2\right),\\
U_j^\gamma(\sigma^\beta-U_j^\beta)
&=U_j^\gamma U_r^\beta+O(U_j^{2^*-1})
+O(U_j^{\gamma+1}U_r^{\beta-1})+O\!\left(\sum_{m\ne r,j}
U_j^\gamma U_mU_r^{\beta-1}\right).
\end{aligned}
\]
Since \(|Z_r|\leq CU_r\), the contribution of these Taylor remainders
on \(\Omega_r\) is bounded by
\[
\begin{aligned}
\mathcal E_{A,r}^{\rm in}:={}&
C\sum_{j,m\ne r}\int_{\Omega_r}U_r^{2^*-2}U_jU_m\dx
+C\sum_{j\ne r}\int_{\Omega_r}U_j^{2^*-1}U_r\dx\\
&+C\sum_{j\ne r}\int_{\Omega_r}U_j^{\gamma+1}U_r^\beta\dx
+C\sum_{j\ne r}\sum_{\substack{m\ne r\\m\ne j}}
\int_{\Omega_r}U_j^\gamma U_mU_r^\beta\dx .
\end{aligned}
\]
For the first sum,
\[
U_jU_m\leq\frac12(U_j^2+U_m^2),\qquad
\int_{\Omega_r}U_r^{2^*-2}U_j^2\dx\leq CR_{rj}^{-n},
\]
where Lemma \ref{A.4} is used with \((a,b)=(2,n-2)\).
For the second and third sums the scaled exponent pairs are
\[
\left(\frac{n-2}{2},\frac{n+2}{2}\right),
\qquad
\left(\frac{n+2-\alpha}{2},
\frac{n+\alpha-2}{2}\right),
\]
respectively.  Their entries add to \(n\), and the dominant bubble
corresponds to the smaller entry; hence
\[
\int_{\Omega_r}U_j^{2^*-1}U_r\dx
+\int_{\Omega_r}U_j^{\gamma+1}U_r^\beta\dx
\leq CR_{rj}^{-n}.
\]
Finally,
$
U_j^\gamma U_m
\leq C(U_j^{\gamma+1}+U_m^{\gamma+1})$,
and the preceding estimate applies to the fourth sum.  Thus
$
\mathcal E_{A,r}^{\rm in}\leq CR^{-n}.
$
On \(\Omega_\ell\), \(\ell\ne r\), use
\[
0\leq A\leq C\left(
U_\ell^{2^*-2}\sum_{j\ne\ell}U_j+
\sum_{j\ne\ell}U_j^\gamma U_\ell^\beta\right).
\]
It follows that
\[
\begin{aligned}
\left|\int_{\Omega_r^c}AZ_r\dx\right|
&\leq C\sum_{\ell\ne r}\int_{\Omega_\ell}
U_\ell^{2^*-2}U_r^2\dx+C\sum_{\ell\ne r}\sum_{\substack{j\ne\ell\\j\ne r}}
\int_{\Omega_\ell}U_\ell^{2^*-2}U_jU_r\dx\\
&\quad+C\sum_{\ell\ne r}\int_{\Omega_\ell}
U_r^{\gamma+1}U_\ell^\beta\dx+C\sum_{\ell\ne r}\sum_{\substack{j\ne\ell\\j\ne r}}
\int_{\Omega_\ell}U_j^\gamma U_\ell^\beta U_r\dx .
\end{aligned}
\]
The parts of the two leading terms over \(\Omega_r^c\) satisfy
\[
\begin{aligned}
&\sum_{j\ne r}\left|\int_{\Omega_r^c}
U_j^\gamma U_r^\beta Z_r\dx\right|
+\sum_{j\ne r}\left|\int_{\Omega_r^c}
U_r^{2^*-2}U_jZ_r\dx\right|\\
\leq &C\sum_{\ell\ne r}\int_{\Omega_\ell}
\bigl(U_\ell^\gamma U_r^{\beta+1}
+U_\ell U_r^{2^*-1}\bigr)\dx+C\sum_{\ell\ne r}\sum_{\substack{j\ne\ell\\j\ne r}}
\int_{\Omega_\ell}\bigl(
U_j^\gamma U_\ell^\beta U_r+
U_\ell^{2^*-2}U_jU_r\bigr)\dx .
\end{aligned}
\]
For the double sums, use
\[
U_jU_r\leq\frac12(U_j^2+U_r^2),\qquad
U_j^\gamma U_r\leq C(U_j^{\gamma+1}+U_r^{\gamma+1}).
\]
All resulting integrals are instances of Lemma \ref{A.4} with
\[
(2,n-2),\quad
\left(\frac{n+2-\alpha}{2},\frac{n+\alpha-2}{2}\right),\quad
\left(\frac{n-2}{2},\frac{n+2}{2}\right),\quad
\left(\frac{\alpha}{2},n-\frac{\alpha}{2}\right).
\]
Each pair has sum \(n\), and the dominance set fixes the correct
ordering.  Consequently the two preceding displays are \(O(R^{-n})\),
and
\begin{equation}\label{eq:new26-strong-A}
\int AZ_r\dx
=\sum_{j\ne r}\int U_j^\gamma U_r^\beta Z_r\dx
+\beta\sum_{j\ne r}\int U_r^{2^*-2}U_jZ_r\dx
+\mathcal R_{A,r},
\qquad |\mathcal R_{A,r}|\leq CR^{-n}.
\end{equation}

For the nonlocal part, put
\[
S_\ell=\sum_{m\ne\ell}U_m,\qquad
L_\ell=p_\alpha U_\ell^\beta S_\ell,\qquad
E_\ell=F-L_\ell\quad\hbox{on }\Omega_\ell,
\]
and \(H=\sigma^\beta-U_r^\beta\).  On \(\Omega_\ell\),
\[
|E_\ell|\leq C\left(U_\ell^{\beta-1}S_\ell^2+
\sum_{m\ne\ell}U_m^{p_\alpha}\right),\qquad
|F|\leq CU_\ell^\beta S_\ell.
\]
Moreover,
\[
|H|\leq
\begin{cases}
CU_r^{\beta-1}\sum_{m\ne r}U_m,&\Omega_r,\\
CU_\ell^\beta,&\Omega_\ell,\ \ell\ne r.
\end{cases}
\]
Decompose the double integral according to
\(\Omega_k\) in the \(x\)-variable and \(\Omega_\ell\) in the
\(y\)-variable.  The first-order contribution of the
\(\Omega_r\times\Omega_r\) block is
\[
T_B:=D_{n,\alpha}p_\alpha\sum_{j\ne r}
\int_{\Omega_r}\int_{\Omega_r}
\frac{U_r^\beta(y)U_j(y)U_r^\beta(x)Z_r(x)}
{|x-y|^\alpha}\dy\dx .
\]
Writing
$
\int BZ_r\dx=T_B+\mathcal E_B,
$
the Taylor identities imply
\[
|\mathcal E_B|
\leq \mathcal E_{B,1}+\mathcal E_{B,2}
+|\mathcal E_{B,3}^{\rm lin}|
+\mathcal E_{B,3}^{\rm rem}+\mathcal E_{B,4},
\]
where
\[
\begin{aligned}
\mathcal E_{B,1}
&:=C\int_{\Omega_r}\int_{\Omega_r}
\frac{\bigl(U_r^{p_\alpha-2}S_r^2+
\sum_{j\ne r}U_j^{p_\alpha}\bigr)(y)
U_r^\beta(x)|Z_r(x)|}{|x-y|^\alpha}\dy\dx,\\
\mathcal E_{B,2}
&:=C\int_{\Omega_r}\int_{\Omega_r}
\frac{U_r^\beta(y)S_r(y)
U_r^{\beta-1}(x)S_r(x)|Z_r(x)|}
{|x-y|^\alpha}\dy\dx,
\end{aligned}
\]
\[
\mathcal E_{B,3}^{\rm lin}
:=D_{n,\alpha}p_\alpha\sum_{\ell\ne r}
\int_{\Omega_r}\int_{\Omega_\ell}
\frac{U_\ell^\beta(y)S_\ell(y)U_r^\beta(x)Z_r(x)}
{|x-y|^\alpha}\dy\dx,
\]
\[
\begin{aligned}
\mathcal E_{B,3}^{\rm rem}
:=&C\sum_{\ell\ne r}\int_{\Omega_r}\int_{\Omega_\ell}
\frac{|E_\ell(y)|U_r^\beta(x)|Z_r(x)|}{|x-y|^\alpha}\dy\dx+C\sum_{\ell\ne r}\int_{\Omega_r}\int_{\Omega_\ell}
\frac{L_\ell(y)|H(x)Z_r(x)|}{|x-y|^\alpha}\dy\dx\\
&+C\sum_{\ell\ne r}\int_{\Omega_r}\int_{\Omega_\ell}
\frac{|E_\ell(y)H(x)Z_r(x)|}{|x-y|^\alpha}\dy\dx ,
\end{aligned}
\]
and \(\mathcal E_{B,4}\) is the contribution of \(x\in\Omega_k\),
\(k\ne r\).  On such a block,
\[
\mathcal E_{B,4}
\leq C\sum_{k\ne r}\sum_{\ell}\sum_{m\ne\ell}
\int_{\Omega_k}\int_{\Omega_\ell}
\frac{U_\ell^\beta(y)U_m(y)U_k^\beta(x)U_r(x)}
{|x-y|^\alpha}\dy\dx .
\]
The following three estimates control every resulting block:
\[
\iint\frac{G(y)U_r^\beta(x)|Z_r(x)|}{|x-y|^\alpha}
\dy\dx\leq C\int G\,U_r^\gamma\dx,
\]
\[
\iint\frac{\Phi(y)\Psi(x)}{|x-y|^\alpha}\dy\dx
\leq C\|\Phi\|_{L^s}\|\Psi\|_{L^s},
\qquad s=\frac{2n}{2n-\alpha},
\]
\[
\|U_a^\beta U_b\|_{L^s(\Omega_a)}
+\|U_b^{p_\alpha}\|_{L^s(\Omega_a)}
\leq CR_{ab}^{-\eta_\alpha},
\qquad \eta_\alpha=\frac{2n-\alpha}{2}.
\]
For completeness, the first norm follows from
\[
A_s:=\frac{n-2}{2}\beta s
=\frac{n(n+2-\alpha)}{2n-\alpha},
\qquad
B_s:=\frac{n-2}{2}s
=\frac{n(n-2)}{2n-\alpha}.
\]
Since \(A_s+B_s=n\) and \(\alpha>4\) implies \(A_s<B_s\),
Lemma \ref{A.4} gives
\[
\int_{\Omega_a}U_a^{\beta s}U_b^s\dx\leq CR_{ab}^{-n}.
\]
Taking the \(s\)-th root yields
\(\|U_a^\beta U_b\|_{L^s(\Omega_a)}
\leq CR_{ab}^{-\eta_\alpha}\).  Similarly,
\[
\|U_b^{p_\alpha}\|_{L^s(\Omega_a)}^s
=\int_{\Omega_a}U_b^{2^*}\dx\leq CR_{ab}^{-n}.
\]
We now estimate the five error classes above.  By self-adjointness,
\[
\begin{aligned}
\mathcal E_{B,1}
\leq C\int_{\Omega_r}
\left(U_r^{p_\alpha-2}S_r^2+
\sum_{j\ne r}U_j^{p_\alpha}\right)U_r^\gamma\dx\leq C\sum_{j\ne r}\int_{\Omega_r}
\left(U_r^{2^*-2}U_j^2+U_j^{p_\alpha}U_r^\gamma\right)\dx .
\end{aligned}
\]
The two scaled exponent pairs are
\[
(2,n-2),\qquad
\left(\frac{\alpha}{2},\frac{2n-\alpha}{2}\right),
\]
and both have sum \(n\).  Lemma \ref{A.4} gives
\[
\mathcal E_{B,1}\leq CR^{-n}.
\]
For \(\mathcal E_{B,2}\), bilinear HLS and the tail estimate give
$$
\mathcal E_{B,2}
\leq C\sum_{j,m\ne r}
\|U_r^\beta U_j\|_{L^s(\Omega_r)}
\|U_r^\beta U_m\|_{L^s(\Omega_r)}\leq CR^{-2\eta_\alpha}.
$$

For the first and third integrals in
\(\mathcal E_{B,3}^{\rm rem}\), self-adjointness reduces the estimate
to
\[
C\sum_{\ell\ne r}\int_{\Omega_\ell}|E_\ell|U_r^\gamma\dx.
\]
On \(\Omega_\ell\),
\[
\begin{aligned}
|E_\ell|U_r^\gamma
\leq C\left(U_\ell^{\beta-1}S_\ell^2+
\sum_{m\ne\ell}U_m^{\beta+1}\right)U_r^\gamma\leq C\sum_{m\ne\ell}
U_\ell^\beta\bigl(U_m^{\gamma+1}+U_r^{\gamma+1}\bigr).
\end{aligned}
\]
The relevant pair is
\[
\left(\frac{n+2-\alpha}{2},
\frac{n+\alpha-2}{2}\right),
\]
whose entries add to \(n\), with the first strictly smaller than the
second.  Hence Lemma \ref{A.4} gives
\[
\sum_{\ell\ne r}\int_{\Omega_\ell}|E_\ell|U_r^\gamma\dx
\leq CR^{-n}.
\]
For the middle integral in \(\mathcal E_{B,3}^{\rm rem}\), HLS gives
\[
\begin{aligned}
\sum_{\ell\ne r}\int_{\Omega_r}\int_{\Omega_\ell}
\frac{L_\ell(y)|H(x)Z_r(x)|}{|x-y|^\alpha}\dy\dx\leq C\sum_{\ell\ne r}\sum_{m\ne\ell}\sum_{q\ne r}
\|U_\ell^\beta U_m\|_{L^s(\Omega_\ell)}
\|U_r^\beta U_q\|_{L^s(\Omega_r)}
\leq CR^{-2\eta_\alpha}.
\end{aligned}
\]
Therefore
\[
\mathcal E_{B,3}^{\rm rem}
\leq C(R^{-n}+R^{-2\eta_\alpha}).
\]

To preserve the sign structure of
\(\mathcal E_{B,3}^{\rm lin}\), set
\[
K_r(y):=D_{n,\alpha}p_\alpha
\int_{\mathbb R^n}
\frac{U_r^\beta(x)Z_r(x)}{|x-y|^\alpha}\dx.
\]
Differentiating
\(\mathcal I_\alpha(U_r^{p_\alpha})=U_r^\gamma\) yields
\[
K_r=\gamma U_r^{\gamma-1}Z_r,\qquad |K_r|\leq CU_r^\gamma.
\]
Extending the \(x\)-integration from \(\Omega_r\) to
\(\mathbb R^n\) gives
\[
\mathcal E_{B,3}^{\rm lin}
=\sum_{\ell\ne r}\int_{\Omega_\ell}
U_\ell^\beta S_\ell K_r\dx-\mathcal Q_{B,3}.
\]
For the full-space term,
\[
U_mU_r^\gamma\leq
C(U_m^{\gamma+1}+U_r^{\gamma+1})
\]
and the preceding \((\beta,\gamma+1)\)-pair give
\[
\left|\sum_{\ell\ne r}\int_{\Omega_\ell}
U_\ell^\beta S_\ell K_r\dx\right|\leq CR^{-n}.
\]
On \(\Omega_k\), \(k\ne r\),
\(|U_r^\beta Z_r|\leq CU_r^{p_\alpha}\).  Hence
\[
\begin{aligned}
|\mathcal Q_{B,3}|
\leq C\sum_{\ell\ne r}\sum_{m\ne\ell}\sum_{k\ne r}
\|U_\ell^\beta U_m\|_{L^s(\Omega_\ell)}
\|U_r^{p_\alpha}\|_{L^s(\Omega_k)}\leq CR^{-2\eta_\alpha}.
\end{aligned}
\]
Thus
\[
|\mathcal E_{B,3}^{\rm lin}|
\leq C(R^{-n}+R^{-2\eta_\alpha}).
\]
Finally,
\[
\begin{aligned}
\mathcal E_{B,4}
\leq C\sum_{k\ne r}\sum_{\ell}\sum_{m\ne\ell}
\|U_\ell^\beta U_m\|_{L^s(\Omega_\ell)}
\|U_k^\beta U_r\|_{L^s(\Omega_k)}\leq CR^{-2\eta_\alpha}.
\end{aligned}
\]
Consequently,
\[
|\mathcal E_B|\leq C(R^{-n}+R^{-2\eta_\alpha}).
\]

It remains to replace the restricted main term \(T_B\) by its
full-space analogue
\[
\widetilde T_B:=D_{n,\alpha}p_\alpha\sum_{j\ne r}
\iint_{\mathbb R^n\times\mathbb R^n}
\frac{U_r^\beta(y)U_j(y)U_r^\beta(x)Z_r(x)}
{|x-y|^\alpha}\dy\dx .
\]
If \(x\in\Omega_k\), \(k\ne r\), and \(y\in\Omega_r\), bilinear HLS
gives
\[
\sum_{j\ne r}\sum_{k\ne r}
\|U_r^\beta U_j\|_{L^s(\Omega_r)}
\|U_r^{p_\alpha}\|_{L^s(\Omega_k)}
\leq CR^{-2\eta_\alpha}.
\]
If \(y\in\Omega_\ell\), \(\ell\ne r\), self-adjointness gives
\[
\int_{\Omega_\ell}U_r^\beta U_jU_r^\gamma\dx
\leq C\int_{\Omega_\ell}U_\ell U_r^{2^*-1}\dx
\leq CR^{-n},
\]
where the last pair is
\(((n-2)/2,(n+2)/2)\).  The same estimate controls the part
\(x\in\Omega_r\), \(y\notin\Omega_r\).  Hence
\[
|T_B-\widetilde T_B|
\leq C(R^{-n}+R^{-2\eta_\alpha}).
\]

Differentiating the one-bubble identity in the scale direction gives
\[
D_{n,\alpha}p_\alpha
\left(\int\frac{U_r^\beta(y)Z_r(y)}
{|x-y|^\alpha}\dy\right)U_r^\beta(x)
=\gamma U_r^{2^*-2}(x)Z_r(x).
\]
By self-adjointness,
\[
\widetilde T_B
=\gamma\sum_{j\ne r}\int U_r^{2^*-2}U_jZ_r\dx .
\]
Combining the last three estimates yields
\[
\int BZ_r\dx
=\gamma\sum_{j\ne r}\int U_r^{2^*-2}U_jZ_r\dx
+\mathcal R_{B,r},\qquad
|\mathcal R_{B,r}|\leq C(R^{-n}+R^{-2\eta_\alpha}).
\]
Combining this identity with \eqref{eq:new26-strong-A} yields
\begin{equation}\label{eq:new26-strong-master}
\begin{aligned}
\int hZ_r\dx
=\sum_{j\ne r}\int U_j^\gamma U_r^\beta Z_r\dx
+(2^*-1)\sum_{j\ne r}\int U_r^{2^*-2}U_jZ_r\dx
+O(R^{-n}+R^{-2\eta_\alpha}).
\end{aligned}
\end{equation}

Lemma \ref{A.2}, uniformly for \(R_{rj}\geq2R\), gives
\[
\int U_j^\gamma U_r^\beta Z_r\dx
=-c_1\mathfrak d_{rj}q_{rj}^{\alpha/(n-2)}
+o\!\left(q_{rj}^{\alpha/(n-2)}\right),
\]
\[
\int U_r^{2^*-2}U_jZ_r\dx
=-\widetilde c_2\mathfrak d_{rj}q_{rj}
+o(q_{rj}),
\qquad c_1,\widetilde c_2>0.
\]
Equivalently, there is a modulus
\(\overline\omega_{A.2}(T)\to0\) such that, whenever
\(R_{rj}\geq T\),
\[
\left|\int U_j^\gamma U_r^\beta Z_r\dx
+c_1\mathfrak d_{rj}q_{rj}^{\alpha/(n-2)}\right|
\leq\overline\omega_{A.2}(T)
q_{rj}^{\alpha/(n-2)},
\qquad
\left|\int U_r^{2^*-2}U_jZ_r\dx
+\widetilde c_2\mathfrak d_{rj}q_{rj}\right|
\leq\overline\omega_{A.2}(T)q_{rj}.
\]
Since
\[
\sum_{j\ne r}\bigl(q_{rj}^{\alpha/(n-2)}+q_{rj}\bigr)
\leq C R^{-a_\alpha},\qquad
a_\alpha<n,\qquad a_\alpha<2\eta_\alpha,
\]
the choice
\[
\omega_8(R)
=C\left(R^{a_\alpha-n}
+R^{a_\alpha-2\eta_\alpha}
+\overline\omega_{A.2}(2R)\right)
\]
satisfies \(\omega_8(R)\to0\).  Equation
\eqref{eq:new26-strong-master} and the two quantified pair expansions
prove \eqref{eq:strong-projection-expansion}, uniformly in \(r\) and in
all bubble parameters.
\end{proof}

\section*{Data Availability Statement}
No datasets were generated or analysed during the current study.

\section*{Conflict of Interest}
The authors declare no conflict of interest. All authors approved the final manuscript.

\section*{Acknowledgements}
The authors are deeply grateful to the anonymous referees for their careful reading and very valuable comments and suggestions that improved the presentation of the paper.

\raggedbottom

\end{document}